\documentclass[11pt]{article}

\usepackage{amsmath}
\usepackage{latexsym}
\usepackage{amssymb}
\usepackage[all]{xy}
\usepackage{tikz}
\usetikzlibrary{calc,intersections,through,backgrounds}
\usetikzlibrary{arrows.meta,positioning}
\usetikzlibrary{decorations.pathmorphing}
\usetikzlibrary{plothandlers}
\usetikzlibrary{plotmarks}
\usetikzlibrary{shadings}
\usetikzlibrary{shapes.geometric}
\usetikzlibrary{datavisualization.formats.functions}

\newfont{\gothique}{eufm10 scaled 1100}  

\newcommand{\PP}{{\mathbb{P}}}
\newcommand{\CC}{{\bf{C}}}
\newcommand{\RR}{{\mathbb{R}}}
\newcommand{\ZZ}{{\mathbb{Z}}}
\newcommand{\NN}{{\mathbb{N}}}
\newcommand{\QQ}{{\mathbb{Q}}}   

\newcommand{\AMX}{Amp(X)}
\newcommand{\NSX}{NS(X)_{\RR}}
\newcommand{\TET}{\Theta_{X}}
\newcommand{\TT}{{\cal{T}}}
\newcommand{\FF}{{\cal{F}}}
\newcommand{\OM}{\Omega}

\newcommand{\OO}{{\cal{O}}}
\newcommand{\HH}{{\cal{H}}}
\newcommand{\DD}{{\cal{D}}}
\newcommand{\AC}{{\cal{A}}}

\newcommand{\HA}{\frac{1}{2}}

\newcommand{\EXT}{{{\cal E}xt}}

\newcommand{\om}{\omega}

\newtheorem{thm}{Theorem}[section]
\newtheorem{lem}[thm]{Lemma}
\newtheorem{pro}[thm]{Proposition}
\newtheorem{cor}[thm]{Corollary}
\newtheorem{rem}[thm]{Remark}
\newtheorem{cl} [thm]{Claim}
\newtheorem{example}[thm]{Example}
\newtheorem{defi}[thm]{Definition}  

\newlength{\myskip}
\newenvironment{pf}{
     \addvspace{\myskip}  

     \noindent {\it Proof.$\, $}}
     {$\Box$

     \addvspace{\myskip}
     }

\makeatletter
\renewcommand{\@seccntformat}[1]{\S \/ {\csname the#1\endcsname}\hspace{0.5em}}
\makeatother

\title{\bf BRIDGELAND STABILITY CONDITIONS AND THE TANGENT BUNDLE OF SURFACES OF GENERAL TYPE}
\author{Igor Reider}

\begin{document}

\bibliographystyle{amsplain}
\maketitle

\setcounter{section}{0}
\numberwithin{equation}{section}
\begin{abstract}
 Let $X$ be a smooth compact complex surface with the canonical divisor
 $K_X$ ample and let $\Theta_X$ be its holomorphic tangent bundle.
 Bridgeland stability conditions are used to study the space
 $H^1 (\TET)$ of infinitesimal deformations of complex structures on $X$ and
 its relation to the geometry/topology of $X$. The main observation is that 
 for $X$ with $H^1 (\TET)\neq 0$ and the Chern numbers $(c_2 (X),K^2_X)$ subject to
 $$
 \tau_X =2ch_2 (\TET)=K^2_X -2c_2 (X) >0 \,(\text{$X$ is of positive index}),
 $$
 the object $\TET[1] $ of the derived category $\DD$ of bounded complexes of
 coherent sheaves on $X$ is Bridgeland {\it unstable} in a certain part of the space of Bridgeland stability conditions. The Harder-Narasimhan
 filtrations of $\TET[1] $ for those stability conditions are expected to provide new insights into geometry of surfaces of general type and the study of their moduli. The paper provides a certain body of evidence that this is indeed the case.
\end{abstract}
 \tableofcontents 
\section{Main problem and basic idea}
The main concern of the paper is the study of the cohomology group $H^1 (\TET)$ for a compact complex surface $X$ with the canonical divisor $K_X$ ample, that is, $X$ is a surface of general type for which a sufficiently large multiple of $K_X$ gives an embedding into a projective space, see \cite{Bo} for precise results. Throughout the paper $\TET$ and its dual $\Omega_X$ stand for the holomorphic tangent and the cotangent
bundles of $X$.
Since the works of Kodaira and Spencer on the deformation of complex structures, \cite{KS}, it is known that this group parametrizes the infinitesimal deformations of complex structure of $X$. Hence the interest to understand this group. Throughout the paper we assume its nonvanishing 
\begin{equation}\label{H1non0}
H^1 (\TET) \neq 0.
\end{equation}

\vspace{0.2cm}
\noindent
{\bf 1.1. Two identifications.}  Our approach is based on the following identifications
\begin{equation}\label{ident}
H^1 (\TET)\cong Ext^1 (\OO_X, \TET) \cong Hom_{\cal D} (\OO_X,\TET[1]),
\end{equation}
where ${\cal D} := {\cal D} (X)$ is the derived category of bounded complexes of coherent sheaves on $X$ and $[1]$ is the autoequivalence of $\DD$ defined by the shift of complexes to the left.

The first identification above has been used in several previous works, see
\cite{R}. The approach there consists of viewing a nonzero cohomology class $\xi$ in $H^1 (\TET )$, via the first identification above, as an extension sequence, that is, as a short exact sequence
\begin{equation}\label{ext}
\xymatrix@R=12pt@C=12pt{
0\ar[r]& \TET \ar[r]^(0.57){i_{\xi}}&{\cal G}_{\xi} \ar[r]^{j_{\xi}}& \OO_X \ar[r]&0}
\end{equation}
in the abelian category ${\cal A}=Coh(X)$ of coherent sheaves on $X$ and then extracting topological/geometric information by studying the sheaf ${\cal G}_{\xi}$ in the middle of that sequence from the point of view of the classical stability conditions - Mumford-Takemoto or Bogomolov (in)stability.

In this paper we will use the second identification in \eqref{ident}. With this identification the cohomology classes in $H^1 (\TET)$ are thought of as morphisms
$$
\OO_X \longrightarrow \TET[1]
$$
in the derived category ${\cal D}$.  This idea, in different contexts, has appeared in the works \cite{AB}, \cite{Ba}. This in turn is largely inspired by Bridgeland's theory of stability conditions. Namely, the category of coherent sheaves ${\cal A}$ is replaced by a new full abelian subcategory ${\cal A'}$ of ${\cal D}$  so that $\TET[1]$ and $\OO_X $ are objects of this new subcategory. So  our identification becomes
$$
H^1 (\TET)\cong  Hom_{\cal D} (\OO_X,\TET[1])\cong Hom_{\cal A'} (\OO_X,\TET[1]) $$
Thus the cohomology classes in $H^1 (\TET)$ are viewed now as morphisms
$$
\OO_X \longrightarrow \TET[1]
$$
in the abelian category  ${\cal A'}$. The essential new feature is that  ${\cal A'}$ admits a stability function
$$
Z: K_0 ({\cal A'}) \longrightarrow \CC
$$
in a sense of Bridgeland. This stability function, according to Bridgeland's theory, lifts to a stability condition $\sigma=(Z, {\cal P})$ on ${\cal D}$. For each $\sigma$ the objects $\TET[1]$ and $\OO_X $ will come with Harder-Narasimhan filtrations and nonzero morphisms in $Hom_{\cal A'} (\OO_X,\TET[1])$ are expected to produce restrictions on the semistable factors of those filtrations.
 One would expect to obtain interesting topological/geometric information from those restrictions. In any case, there seem to be a source for some interesting new invariants of surfaces. The purpose of the paper is to give some evidence that these ideas are fruitful for surfaces of general type.

\vspace{0.2cm}
\noindent
{\bf 1.2. Extensions versus morphisms.} Before going further, let us spell out the relation between two approaches. To go from the extension \eqref{ext} to the corresponding morphism
in $Hom_{\cal D} (\OO_X,\TET[1])$, one replaces \eqref{ext} by the distinguished triangle
\begin{equation}\label{mapcone}
\xymatrix@R=12pt@C=12pt{
\TET \ar[r]^(0.57){i_{\xi}}& {\cal G}_{\xi}\ar[r]& M(i_{\xi}) \ar[r]& \TET [1]
}
\end{equation}
in $\DD$, where $M(i_{\xi})$ is the mapping cone of the morphism $i_{\xi}:\TET \longrightarrow {\cal G}_{\xi}$. In the derived category ${\cal D}$ the mapping cone $M(i_{\xi})$ is isomorphic to $\OO_X$, the cokernel of $i_{\xi}$ in $\AC$. Hence the resulting morphism $f(\xi): \OO_X \longrightarrow \TET[1]$.

 Conversely, given a morphism $\eta: \OO_X  \longrightarrow \TET[1]$ in ${\cal D}$, we complete it to a distinguished triangle
$$
\xymatrix@R=12pt@C=12pt{
\OO_X \ar[r]^(0.4){\eta}& {\TET[1]}\ar[r]& C \ar[r] &\OO_X [1]
}
$$
in ${\cal D}$ and apply the homological functor ${\cal H}^0:{\cal D}\longrightarrow \AC $ to obtain the exact complex of sheaves
$$
\xymatrix@R=12pt@C=12pt{
0\ar[r]&\TET\ar[r]&{\cal H}^{-1} (C)\ar[r]&\OO_X \ar[r]&0,
}
$$
a representative of an element of $Ext^1 (\OO_X , \TET)$. Obviously, different ways of completing $\eta: \OO_X  \longrightarrow \TET[1]$ to a distinguished triangle in ${\cal D}$  give rise to isomorphic extensions in $\AC $. 

The two constructions are inverse to each other. If we start with an extension
sequence \eqref{ext}, the morphism $f(\xi): \OO_X  \longrightarrow \TET[1]$
is a part of the distinguished triangle
$$
\xymatrix@R=12pt@C=12pt{
\OO_X \ar[r]^(0.45){f(\xi)}& {\TET[1]}\ar[r]^(.55){-i_{\xi}}& {\cal G}_{\xi} [1] \ar[r]^(0.45){-j_{\xi}} &\OO_X [1]
}
$$
obtained from the triangle \eqref{mapcone} by applying the axiom of `turning' triangles two times (observe the change of signs occurring according to that axiom ). Applying the homological functor ${\cal H}^0$ we obtain the exact sequence
$$
\xymatrix@R=12pt@C=12pt{
0\ar[r]&\TET \ar[r]^(0.47){-i_{\xi}}& {\cal G}_{\xi}\ar[r]^(0.45){-j_{\xi}} & \OO_X \ar[r]& 0
}
$$
Observe, the sequence obtained is not the same, but only isomorphic to the one we started with via the following isomorphism of exact sequences
$$
\xymatrix@R=19pt@C=12pt{
0\ar[r]&\TET \ar[r]^(0.47){i_{\xi}} \ar[d]_{id}& {\cal G}_{\xi}\ar[r]^(0.45){j_{\xi}} \ar[d]_{-id}& \OO_X \ar[r] \ar[d]^{id}& 0\\
0\ar[r]&\TET \ar[r]^(0.47){-i_{\xi}}& {\cal G}_{\xi}\ar[r]^(0.45){-j_{\xi}} & \OO_X \ar[r]& 0
}
$$

The discussion above shows that the two points of view on studying $H^1 (\Theta_X)$ are equivalent and one can freely go from one to another. But as we will see shortly, going to the derived category $\DD$ and using Bridgeland stability conditions is much more powerful - it brings in the continuity arguments which under suitable topological assumptions on $X$ {\it guarantee the instability of $\TET[1]$ in a certain part of the space $Stab_X$ of Bridgeland stability conditions.} This observation seems to be of independent interest and should have many interesting applications in the theory of surfaces of general type.

\vspace{0.2cm}
\noindent
{\bf 1.3. Main results of the paper.}
 The topological assumption on $X$ we make in the paper is
\begin{equation}\label{tau-int}
\tau_X :=K^2_X - 2c_2 (X)>0,
\end{equation}
where $K^2_X$ and $c_2 (X)$ are the Chern numbers of $X$; this amounts to saying that we will be concerned with surfaces of positive topological index  (see a rapid summary on surfaces of general type in \S\ref{S-surf}).
Inside the cone of ample divisor classes $Amp(X)$ we consider the subcone
$$
Amp^s (X):=\{h \in Amp(X) | \text{$\Theta_X$ is $h$-stable}\}.
$$
Here $h$-stability is the classical Mumford-Takemoto slope stability and the cone is nonempty because the ray generated by the canonical class $K_X$ is in $Amp^s (X)$. The starting point of our results is the following.
\begin{thm}\label{th-Bregion}
	For every integral divisor class $H\in Amp^s (X)$ there is a region
	 $B^{un}_H$ in the space of Bridgeland stability conditions for which
	 the object $\TET[1]$ is Bridgeland unstable.
	 \end{thm} 
For the reader who is not familiar with Bridgeland stability we tried to give a minimum background to make the paper self contained, see \S\ref{S-BrS}. For now it will be enough to say that  $B^{un}_H$ is the set of pairs $(cH,bH)$ where $(c,b)$ form a part of the semicircular region in the upper half plane $\Pi_{H}$:
$$
D_H=\Big\{
(c,b)\in \Pi_{H}\,\, \Big{|} \,\, \Big(b+\frac{\tau_X}{2H\cdot K_X}\Big)^2 +c^2 <\Big(\frac{\tau_X}{2H\cdot K_X}\Big)^2 
\Big \},
$$
where $\tau_X$ is as in \eqref{tau-int}. Let us spell out what the above result means:

\vspace{0.2cm}
$\bullet$ for every $(cH,bH)$ in $B^{un}_H$ there is a full abelian subcategory
$\AC_{H,b H}$ of $\DD$ together with the stability function $Z_{cH,bH}$;

$\bullet$
 $\TET[1]$ is an object of $\AC_{H,b H}$ and it is unstable with respect to the stability function $Z_{cH,bH}$, that is, there is the Harder-Narasimhan filtration (abbreviated to HN filtration in the sequel)
 \begin{equation}\label{filtTh1-intro}
 \TET[1]=E^{c,b}_n \hookleftarrow E^{c,b}_{n-1} \hookleftarrow \cdots \hookleftarrow E^{c,b}_1 \hookleftarrow E^{c,b}_0 =0. 
\end{equation}
in the abelian category $\AC_{H,b H}$.

The above filtrations are intrinsically associated with $(cH,bH)$, for every $H\in Amp^s (X)$, and should contain interesting information about the surface $X$. For example, if we return to the cohomology, we obtain the `shadow' of the above filtration
$$
\begin{gathered}
 H^1(\TET)=Hom_{\AC_{H,b H}} (\OO_X,\TET[1])=Hom_{\AC_{H,b H}} (\OO_X,E^{c,b}_n) \supset
 \\
  Hom_{\AC_{H,b H}} (\OO_X,E^{c,b}_{n-1}) \supset \cdots \supset Hom_{\AC_{H,b H}} (\OO_X,E^{c,b}_1) \supset Hom_{\AC_{H,b H}} (\OO_X,E^{c,b}_0) =0.
 \end{gathered} 
 $$ 
 on the cohomology level. This can be stated as follows.
 \begin{thm}\label{th-H-filt}
 	Let $X$ be a smooth compact complex surface with the canonical divisor $K_X$ ample. Assume that $H^1 (\Theta_X)$ is nonzero and $\tau_X=K^2_X -2c_2 (X)$ is positive. Then $H^1 (\TET)$ admits
 	a filtration
 	$$
  H^1(\TET) =H^1(\TET)^n _{cH,bH} \supset H^1(\TET)^{n-1}_{cH,bH} \supset \cdots \supset H^1(\TET)^1_{cH,bH} \supset H^1(\TET)^0_{cH,bH} =0
  $$
  intrinsically attached to every integral divisor class $H \in Amp^s (X)$ and every $(cH,bH) \in B^{un}_H$.
 \end{thm}

This result could be viewed as a part of the duality
$$
\text{K\"ahler moduli}\leftrightarrow \text{Complex structure moduli}
$$
inspired by the Mirror symmetry. T. Bridgeland discovered his space of stability conditions in response to string theorists' quest to understand
superconformal field theories (SCFT) on Calabi-Yau manifolds. We refer the reader to the original paper of T. Bridgeland, \cite{Br}, as well as to works of Douglas, \cite{Do}, \cite{Do1}, for physicists' perspective.

It is believed  that the Mirror Symmetry principle relating symplectic and complex structures on a smooth manifold and controlled respectively by  K\"ahler moduli and moduli of complex structures goes beyond the Calabi-Yau manifolds and, in particular, should manifest itself in manifolds of general type, see for example \cite{KKOY}. Our Theorem \ref{th-H-filt} could be regarded as such a manifestation.

  Namely, the cohomology space $H^{\ast}(X,\CC)$ is a functorial linearization of a smooth manifold structure of $X$. A K\"ahler and  complex structures of $X$ give Hodge structures on $H^{\ast}(X,\CC)$; viewed over $\{\text{Complex moduli}\}$ we have the
  variation of Hodge structures on $H^{\ast}(X,\CC)$, see \cite{G} for an overview.
  
  The cohomology space $H^1 (\TET)$ is a linearization of the moduli space of complex structures on $X$ and
  the above result tells us that  K\"ahler and  complex structures of $X$ give a sort of `Hodge structure' on $H^1 (\TET)$ - the filtration of Theorem \ref{th-H-filt}; viewed over
    (a part of)  
 $\{\text{ K\"ahler moduli}\}$ we also expect a theory of variation of these
 `Hodge structures'. 
 
 We do not know how far the above analogy goes. But it should be clear
 that there are additional interesting structures on the space $H^1 (\TET)$
 of the infinitesimal
 deformations of complex structures of $X$. Those structures in turn
 arise from additional structures of the derived category $\DD$ of bounded complexes of coherent sheaves on $X$. In this our results inscribe in the quest for understanding the derived categories of coherent sheaves on projective varieties. 
  
  Besides proving the above results, the paper focuses on the geometric aspects of Bridgeland instability of $\TET[1]$ and its cohomological counterpart, the filtration on $H^1 (\TET)$. To give an idea for a kind of geometry contained in the filtration \eqref{filtTh1-intro} and its cohomological counterpart in Theorem \ref{th-H-filt}, let us take
  an intermediate
  subobject $E^{c,b}_k$, for $k\in [1,n-1]$, in the filtration \eqref{filtTh1-intro}. The inclusion $E^{c,b}_k \hookrightarrow \TET[1]$ in the abelian category $\AC_{H,bH}$ can be completed to the exact sequence
   $$
   \xymatrix@R=12pt@C=12pt{
  0\ar[r]&E^{c,b}_k \ar[r]& \TET[1]\ar[r]& F_k [1] \ar[r]&0 
}
  $$ 	
in $\AC_{H,bH}$. It turns out that $F_k$ is a torsion free sheaf of rank at least $2$.  We suggest to view the epimorphism
\begin{equation}\label{cat-diff}
\TET[1]\longrightarrow F_k [1]
\end{equation}
 in the above sequence as a categorical version of the differential of a morphism of $X$ into some other variety so that $F_k$ is a substitute for the pull back of the tangent sheaf of that variety. In the sequel we refer to the morphism \eqref{cat-diff} as {\it categorical differential}.
 
 To see the analogy more clearly, apply the homological functor $\HH^0$ to the above exact sequence
 to obtain the exact complex of coherent sheaves on $X$
  $$
 \xymatrix@R=12pt@C=12pt{
 	0\ar[r]&\HH^{-1}(E^{c,b}_k )\ar[r]& \TET\ar[r]& F_k  \ar[r]&\HH^0 (E^{c,b}_k ) \ar[r]&0. 
 }
 $$ 	
This is now a sheaf version of the differential of a morphism from $X$ to some variety with $F_k$ playing the role of the tangent sheaf of the target variety. From the complex the geometric role of the complex $E^{c,b}_k$ becomes transparent:

the cohomology sheaves $\HH^{-1}(E^{c,b}_k )$ and $\HH^{0}(E^{c,b}_k )$ of the complex $E^{c,b}_k$ are respectively the kernel and the cokernel of the categorical differential; if $\HH^{-1}(E^{c,b}_k )$ is nonzero, then $X$ and the target variety are foliated; if $\HH^{-1}(E^{c,b}_k )$ is zero, then the categorical differential is generically an immersion, the subobject $E^{c,b}_k$ is
a sheaf and it could be envisaged as the categorical version of the normal sheaf of that immersion.

Turning to the filtration on $H^1 (\TET)$ in Theorem \ref{th-H-filt}, the subspace $Hom (\OO_X, E^{c,b}_k )$ could be interpreted as the infinitesimal deformations of $X$ for which
the categorical differential \eqref{cat-diff} deforms along. 

The subspaces of the filtration $H^1 (\TET)$ could also be thought as distributions on the moduli stack of surfaces of general type. If integrable
this would tell us that the moduli stack is foliated.

The above discussion indicates the interest to study the filtrations in \eqref{filtTh1-intro} and Theorem \ref{th-H-filt}. Part of our effort
goes into uncovering the relations between the Chern numbers $(c_2(X),K^2_X)$
and the invariants of the sheaves $F_k$ in \eqref{cat-diff}. The most probing results in this direction are obtained when the divisor class 
$H$ in Theorem \ref{th-Bregion} and all subsequent considerations is taken to be the canonical divisor $K_X$. To state our result we let
$$
\alpha_X :=\frac{c_2 (X)}{K^2_X}
$$
the ratio of the Chern numbers of $X$. We recall that $\alpha_X \geq \displaystyle{\frac{1}{3}}$. This is the celebrated Bogomolov-Miyaoka-Yau inequality for surfaces of general type, see \cite{Mi} and a brief introduction in \S2. Our assumption $\tau_X>0$ in \eqref{tau-int} and the nonvanishing of $H^1 (\TET)$ in \eqref{H1non0} imply
$$
\frac{1}{3} <\alpha_X <\HA.
$$
\begin{thm}\label{th:filtKX-intro}
	Let $\alpha_X < \displaystyle{\frac{3}{8}}$ and let $B^{un}_{K_X}$ be the region in Theorem \ref{th-Bregion}. There is a constant $\beta_{K_X}$ such that for every stability condition $(cK_X,bK_X)$ with $(c,b)\in B^{un}_{K_X}$ and $b>\beta_{K_X}$ the HN filtration
	$$
	\TET[1]=E^{c,b}_n \supset E^{c,b}_{n-1}\supset \cdots \supset E^{c,b}_{1}  \supset E^{c,b}_{0}=0
	$$
	has the following properties.
	
	1) All objects $E^{c,b}_{k}$, for $k\in [1,n-1]$, are sheaves.
	
	2) In the quotient objects
	$$
	F_k [1] =\TET[1] / E^{c,b}_k,
	$$
	for every $k \in [1,n-1]$, the sheaves $F_k$ are locally free and either
	Bogomolov unstable or have rank $2$ or $3$.
	
	3) If $F_k$ is Bogomolov unstable, then its HN filtration
	$$
	F_k =F^{l_k} _k \supset F^{l_k-1} _k \supset \cdots \supset F^{1} _k \supset F^{0} _k =0
	$$
with respect to $K_X$ has the last semistable quotient
	$$
	Q^{l_k}_k= F^{l_k} _k / F^{l_k-1}_k 
	$$ 
	of rank $2$ or $3$ and the divisor class $(-c_1(Q^{l_k}_k))$ lies in the positive cone of $NS(X)_{\RR}$, the real N\'eron-Severi group of $X$.
\end{thm}

This somewhat technical statement can be informally explained as follows.
The assumption on the ratio $\alpha_X$ of the Chern numbers of $X$ allows us to control the properties
of the categorical differentials
$$
\TET[1] \longrightarrow F_k [1]:
$$

the property 1) says that it is always generically an immersion and 
$ E^{c,b}_k$ is a categorical version of the normal sheaf of that immersion;

the property 2) controls the rank of $F_k$, or the dimension of the target variety of the immersion: it is either a surface or a 3-fold, unless
$F_k$ is Bogomolov unstable; in this case the property 3) intervenes and it tells us that we always have categorical immersions into a surface or a 3-fold by replacing $F_k$ by the last semistable quotient $Q^{l_k}_k$ in
Theorem \ref{th:filtKX-intro}, 3); namely, we compose the morphisms
$$
 \xymatrix@R=12pt@C=12pt{
\TET[1]\ar[r]&F_k [1] \ar[d]\\
&Q^{l_k}_k [1]
}
$$  
to obtain the categorical differential
\begin{equation}\label{Th1-Qk}
\xymatrix@R=12pt@C=12pt{
	\TET[1]\ar[r]&Q^{l_k}_k [1]
}	
\end{equation}
  whose target is of dimension $2$ or $3$. In addition, the property
  3) of Theorem \ref{th:filtKX-intro} says that the divisor class
  $(-c_1 (Q^{l_k}_k))$ is in the positive cone of $X$. Remember, that the geometric analogue of $Q^{l_k}_k$ is the pull back of the tangent sheaf of the target variety. So the positivity of  $(-c_1 (Q^{l_k}_k))$ could be viewed as the statement 
  that the target variety is of general type.
  
  Thus with the dimensional reduction technique outlined above, Theorem \ref{th:filtKX-intro} tells us that $X$ admits categorical differentials
  to surfaces or $3$-folds of general type. The first case could be viewed as
   {\it categorical ramified covering} and the second as the composition
   of the categorical ramified covering with the {\it categorical immersion}
  into a $3$-fold.
  
  On the level the cohomological filtration of $H^1 (\TET)$ we have the following `reduced' version of Theorem \ref{th-H-filt}.
  \begin{thm}\label{th:H1-filt-intro}
  	Let $X$ be a smooth compact complex surface with the canonical divisor
  	$K_X$ ample and the ratio $\alpha_X$ of the Chern numbers
  	$$
  	\alpha_X =\frac{c_2(X)}{K^2_X} <\frac{3}{8}.
  	$$
  	Assume $h^1 (\TET)\geq 2$. Then $H^1 (\TET) $ admits a one-step filtration
  	$$
  	H^1 (\TET) \supset H^1 (\TET)^0
  	$$
  	where the subspace $ H^1 (\TET)^0$ controls the categorical ramified coverings of $X$ and the quotient
  	$	H^1 (\TET)/ H^1 (\TET)^0$ controls the categorical immersions into
  	a $3$-fold.
  \end{thm} 

In the main body of the paper the reader will find an expanded version of this theorem, see Theorem \ref{th:H1Theta}. Let us unravel the above
version just enough to see how the above links to the geometry of $X$.

The reductional mechanism allows us to construct the following exact sequence
\begin{equation}\label{TET-Q'X-intro}
\xymatrix@R=12pt@C=12pt{
	0\ar[r]&\TET \ar[r]& Q'_X \ar[r]& S_X \ar[r]&0
}
\end{equation}
where $Q'_X$ is a $K_X$-semistable torsion free sheaf of rank $2$ or $3$.
 
In the rank $2$ case $S_X$ is a torsion sheaf and we have
\begin{equation}\label{H1=H10-intro}
	\begin{gathered}
\bullet \,\,H^1 (\TET)= H^1 (\TET)^0 \cong H^0 (S_X),
\\
\bullet\,\,\text{$S_X$ is globally generated,}
\end{gathered}
\end{equation}
and everything is controlled by the geometry of the divisor
$$
D^0_X =c_1 (S_X).
$$
 This divisor could be viewed as the categorical version of the ramification divisor and the identification in the first item of \eqref{H1=H10-intro} means that all infinitesimal deformations of $X$ are obtained from the infinitesimal deformations of $D^0_X$. The exact sequence \eqref{TET-Q'X-intro} gives the decomposition of the canonical divisor $K_X$ of $X$
 $$
 K_X =(-c_1 (Q'_X)) +D^0_X
 $$
 which is reminiscent of the formula for the canonical divisor of $X$ under
 a ramified covering from $X$ to some other surface, call it $X'$. The sheaf
 $Q'_X$ looks like the pull back of  the tangent sheaf of $X'$. It can be shown that $(-c_1 (Q'_X))$ is again in the positive cone of $X$ and this, as we already mentioned, could be viewed as $X'$ being of general type.
 Furthermore, the divisor $D^0_X$ is subject to 
 $$
 D^0_X \cdot K_X < 4(3c_2 (X)-K^2_X),
 $$  
 while the Chern numbers of $Q'_X$ satisfy
 $$
 3c_2 (X)-K^2_X \geq 3c_2 (Q'_X)-c^2_1 (Q'_X).
 $$
 This last inequality means that if our hypothetical $X'$ exists, then its Chern numbers
 would lie closer to the Bogomolov-Miyaoka-Yau line
 $$
 3c_2 =K^2
 $$
 than the Chern numbers of $X$.
 
 The analogy with ramified covering can be pushed further by observing that
$D^0_X$ can be decomposed into two parts
$$
D^0_X=\big(D^0_X\big)' + \big(D^0_X\big)''.
$$ 
Each part is composed of the reduced, irreducible components of $D^0_X$ which sheaf theoretically look like components on which a ramified covering is respectively finite or a contraction to a point.

In addition, it turns out that $D^0_X$ is nef unless there are smooth rational curves $C$ on $X$ with 
$$
C\cdot D^0_X =-1 \,\, \text{or} \,\, -2.
$$
Such curves must be components of $\big(D^0_X\big)''$ and if they exist one can construct special filtrations of 
$H^1 (\TET) $ labeled by rational curves. Thus what comes out of the above discussion is the relation
\begin{equation}\label{flag-X}
	\text{\it Flag varieties of $H^1(\TET)^0$}\leftrightarrow \text{\it Geometry of $X$.}
\end{equation}

\vspace{0.2cm}
Let us now turn to the case when $rk(Q'_X)=3$. The sheaf $S_X$ in \eqref{TET-Q'X-intro} is of rank $1$. It may have the torsion part 
denoted $Tor(S_X)$ and we write ${\cal J}_{A_X} (L_X)$ for its torsion free part. So we have the exact sequence
$$
\xymatrix@R=12pt@C=12pt{
	0\ar[r]& Tor(S_X)\ar[r]& S_X \ar[r]&{\cal J}_{A_X} (L_X)\ar[r]&0.
}
$$
As before we have
\begin{equation}\label{H1=H0SX-3}
	\begin{gathered}
	\bullet \,\,H^1 (\TET)\cong H^0 (S_X),
	\\
	\bullet\,\,\text{$S_X$ is globally generated.}
\end{gathered}	 
\end{equation} 
In addition, the above exact sequence for $S_X$ brakes the exact sequence
\eqref{TET-Q'X-intro} as follows.
\begin{equation}\label{TET-Q'X-3-diag}
\xymatrix@R=12pt@C=12pt{
	&&0\ar[d]&0\ar[d]&\\
0\ar[r]&\TET\ar[r]\ar@{=}[d]&Q''_X \ar[r]\ar[d]&Tor( S_X) \ar[d]\ar[r]&0\\ 	
0\ar[r]&\TET\ar[r]&Q'_X \ar[r]\ar[d]& S_X \ar[r]\ar[d]&0\\
&&{\cal J}_{A_X} (L_X)\ar@{=}[r]\ar[d]& {\cal J}_{A_X} (L_X)\ar[d]&\\
&&0&0&
}
\end{equation}
The subspace $H^0 (Tor(S_X)) \subset H^0(S_X)$,
via the isomorphism in the first item in \eqref{H1=H0SX-3}, gives rise to the subspace of $H^1 (\TET)$ which is denoted $H^1(\TET)^0$. If this subspace
is nonzero, then the top row of the diagram \eqref{TET-Q'X-3-diag} puts us back into the rank $2$ case with the divisor
$$
D^0_X =c_1 (Tor(S_X))
$$
controlling the infinitesimal deformations parametrized by the subspace $H^1(\TET)^0$. 

The quotient space
$H^1(\TET)/H^1(\TET)^0$ is a new aspect in our considerations. It is identified with a subspace of 
global sections of ${\cal J}_{A_X} (L_X)$ and that subspace globally generates the sheaf. Hence the geometric meaning of $H^1(\TET)/H^1(\TET)^0$ is in its relation to the properties of the linear system $|L_X|$.

In addition, we have the numerical estimates similar to the rank $2$ case:
$$
\begin{gathered}
(D^0_X +L_X)\cdot K_X \leq 2(3c_2(X)-K^2_X),
\\
3c_2(X)-K^2_X \geq 3c_2 (Q''_X)-c^2_1 (Q''_X).
\end{gathered}
$$
Heuristically, the diagram \eqref{TET-Q'X-3-diag} tells us that $X$ 

$\bullet$ admits a categorical ramified covering onto a surface, call it $X''$; the exact
sequence on the top of the diagram playing the role of the differential of that ramified covering, the divisor $D^0_X$ playing the role of the ramification divisor and its geometric properties controlling the subspace
$H^1 (\TET)^0$ according to the principle stated in \eqref{flag-X}; furthermore, the above estimates tell us that $X''$ is of general type with the Chern numbers placed closer to the Bogomolov-Miyaoka-Yau line (BMY line in the sequel), 

$\bullet$ the surface $X''$ admits a categorical immersion into a $3$-fold;
the middle column of the diagram playing the role of the normal sequence of that immersion and the sheaf ${\cal J}_{A_X} (L_X)$ playing the role of the normal sheaf, and geometric properties of the linear system $|L_X|$ controlling the quotient space $H^1(\TET)/H^1(\TET)^0$.

\vspace{0.2cm} 
We have seen that the geometric output of our results, the divisors
$D^0_X$ and $L_X$, are subject to the inequalities depending on
$(3c_2 (X)-K^2_X)$, the quantity which measures how far the Chern numbers $(c_2 (X),K^2_X)$ are away from the BMY line. So by imposing the additional assumptions on $(3c_2 (X)-K^2_X)$, we can control the geometric properties of these divisors. In \S9 we impose the condition
\begin{equation}\label{qbound-intro}
3c_2 (X)-K^2_X <\HA (K^2_X)^{\HA}
\end{equation}
 to deduce the following:
 \begin{equation}
 	\begin{gathered}
 	\text{\it the sublattice of the N\'eron-Severi group $NS(X)$ of $X$
 		spanned}
 	\\
 	\text{\it by the irreducible components of the divisor
 		$(D^0_X +L_X)$ is negative semidefinite.}
 	\end{gathered}
 \end{equation}
This imposes serious restrictions on the linear system $L_X$.
\begin{thm}\label{th:fX-fibration-intro}
	Let $X$ be a smooth compact complex surface with the canonical divisor
	$K_X$ ample and the Chern numbers $(c_2 (X),K^2_X)$ of $X$ subject to
	$$
3c_2 (X)-K^2_X <\HA (K^2_X)^{\HA}.
$$
Assume furthermore, $dim(H^1(\TET)/H^1(\TET)^0)\geq 2$.	Then $X$ admits a
fibration
$$
f_X :X \longrightarrow B
$$
onto a smooth projective curve $B$ with connected fibres of genus $g_{{}_X}$ subject to
$$
2(g_{{}_X}-1)(dim(H^1(\TET)/H^1(\TET)^0)-1)\leq 3c_2 (X)-K^2_X.
$$
This fibration is essentially unique.
\end{thm}

We refer the reader to \S9 for details and corollaries of this result. Here we wish to make several informal remarks. The assumptions of the theorem tell us that $X$ has the structure of a fibration and the genus of the fibration is an integer $g$ subject to
$$
2(g-1)\leq 3c_2 (X)-K^2_X
$$
{\it intrinsically} associated to $X$. 
This intrinsic fibration gives rise to a unique morphism
$$
\overline{m}_{{}_X}: B \longrightarrow \overline{\mathfrak M}_g
$$
into the Deligne-Mumford stack $\overline{\mathfrak M}_g$ of stable curves of genus $g$. This should be functorial, in particular, a deformation of $X$
should give rise to a deformation of $\overline{m}_{{}_X}$. Thus the moduli stack of surfaces subject to Theorem \ref{th:fX-fibration-intro} give rise
to the substack of curves on $\overline{\mathfrak M}_g$. This could be viewed as a manifestation of Gromov-Witten type invariants. In particular, if $X$ a regular surface, that is, $q(X)=h^1 (\OO_X)=0$, the curve $B$ in Theorem \ref{th:fX-fibration-intro} is $\PP^1$ and we obtain the substack of ${\bf Maps}(\PP^1,\overline{\mathfrak M}_g)$ of maps from $\PP^1$ to $\overline{\mathfrak M}_g$, that is the Gromov-Witten theory of genus zero.
  
  Of course the condition \eqref{qbound-intro} may look completely artificial and the reader may wonder about the utility of Theorem \ref{th:fX-fibration-intro}. In fact, it becomes quite natural if
  one looks at the problem of geography/moduli of surfaces of general type
  from the point of view of the quantity
  $$
  (3c_2 -K^2).
  $$
  In \S2, for convenience of the reader we recall the `world' of numerical invariants of surfaces of general type, see the picture \eqref{Pworld}.
  We suggest to `foliate' the region on that picture by lines parallel to BMY line and then to study the surfaces with the Chern invariants lying on a {\it fixed} line. In other words, we fix a positive integer $d$ and look at the line
  $$
  BMY_d:=\{3c_2 -K^2 =d\}.
  $$
  The assumption \eqref{qbound-intro} now can be restated as `going sufficiently far' along that line (precisely, all admissible integral points of the line with $K^2 >4d^2$). So Theorem \ref{th:fX-fibration-intro}
  tells us that {\it all} surfaces with invariants lying sufficiently far on the line $BMY_d$ are either fibrations as described in the theorem or are subject to
  $$
  dim(H^1 (\TET)/H^1 (\TET)^0) \leq 1.
  $$
 As we explained above, having this intrinsic fibration should be useful for understanding of surfaces and their moduli stacks. So the admissible points
 on $BMY_d$ subject to Theorem \ref{th:fX-fibration-intro} could be viewed
 as amenable for a better understanding.
 
  The other possibility, $dim(H^1 (\TET)/H^1 (\TET)^0) \leq 1$, always under the assumption $h^1 (\TET) \geq 2$, shifts the focus on the subspace $H^1 (\TET)^0$ and, according to our heuristic principle, this controls the categorical ramified coverings of some surfaces $X'$ lying on the lines $BMY_{d'}$ with $d'< d$. Thus we on the lines closer to the BMY line and this provides an inductive procedure for studying surfaces of general type.
 
 We end the discussion of our results with a word of warning: there are no examples in this paper. We have
 attempted to present here several general phenomena in the theory of surfaces of general type stemming from the Brindgeland's instability of 
 $\TET[1]$. A separate paper will be devoted to examples.
 
 \vspace{0.2cm}
 \noindent
 {\bf 1.4. Organization of the paper.} We tried to make the paper relatively
 self contained. This partly explains its length. There are two main themes involved:
 
 1) surfaces of general type with the problems of geography and the moduli as
 the main subjects of concern,
 
 2) Bridgeland stability conditions.
 
 In a nut shell the paper is about the relation of these two themes and it is written with the hope that the reader with interests in one of those themes
 will find the paper accessible:
 on the one hand we included
 standard facts known to experts of surfaces of general type and on the other
 hand the basics on Bridgeland stability.
 
 Below we give a detailed guide to the contents of each section, so the interested reader could go directly to the ones corresponding to her interests.
 
 \vspace{0.2cm}
 
 \S2 is a brief account of basics on surfaces of general type. In particular, we recall the region where the Chern numbers of smooth minimal surfaces of general type reside, see \eqref{Pworld}, the problem of `populating' the region is the geography problem. For more background, see \cite{P1} and references there in.
 
 Another aspect of this background material is the particular features of the
 cotangent bundle of surfaces of general type which are used throughout the paper. One feature says that the cotangent bundle of a surface of general type can not have rank one subsheaves of the Iitaka dimension $2$. This is known as the Bogomolov lemma, see \cite{Mi}. The other is generic positivity of the cotangent bundle of surfaces of general type, a result of Y. Miyaoka, see
 \cite{Mi2}.
 
 In the subsection {\bf 2.3} of \S2 we recall classical Mumford-Takemoto 
 slope stability in the abelian category of coherent sheaves on $X$. Then, in {\bf 2.4}, we
 go on to describe the structure induced by the slope (in)stability of $\TET$ on the ample cone of $X$: 
 
 there is the subcone of polarizations with respect to which $\TET$ is stable; it is surrounded by codimension $1$ faces of strictly semistable polarizations; the faces are labeled by destabilizing sequences for $\TET$;
 these are foliations on $X$.

 All this could be considered as a classical version of the results in Theorem \ref{th-H-filt} and is intended to make contrast with a much finer
 structures coming from Bridgeland stability conditions.
 
 The last subsection, {\bf 2.5}, reviews Bogomolov-Gieseker inequality and
 the Bogomolov \linebreak(in)stability.
 
 \S3 is a fast paced review of the elements of Bridgeland stability conditions. We introduce the bare minimum needed to proceed with the study of $\TET$ and $H^1(\TET)$. The reader is invited to consult Bridgeland's original paper \cite{Br} as well as the overview \cite{Ma-S} and the references therein.
 
 In \S4 we apply Bridgeland stability conditions to study the objects
 $\OO_X$ and $\TET[1]$. Here we define, for every integral polarization
 $H$ a suitable tilted abelian subcategory ${\AC_{H,bH}}$ of $\DD$ which contains both objects. The main result is Corollary \ref{cor:ThXBrunst}
 which is formulated as Theorem \ref{th-Bregion} in the subsection {\bf 1.3}.
 The key intermediate results here are Proposition \ref{pro:Phi_V-inj} and
 Proposition \ref{pro:semicircH}. The first gives a criterion for morphisms
 $$
 \OO_X \longrightarrow \TET[1]
 $$
 to be injective in ${\AC_{H,bH}}$ and the second identifies the semicircular wall in the upper half plane defined by $\OO_X$: inside the wall the phase of $\OO_X$ is bigger than the phase of $\TET[1]$ and hence Bridgeland instability of one of those objects.
 
 \S5 studies the HN filtration of $\TET[1]$, the one appearing in \eqref{filtTh1-intro}. Every nonzero intermediate object of the filtration
 gives rise to a semicircular wall in the upper half plane $\Pi_H$.
 The subsection {\bf 5.2} collects basic facts about the numerical and actual walls for $\TET[1]$, see Proposition \ref{pro:Th-walls}, Corollary \ref{cor:semicirwalls}, Lemma \ref{lem:awall-finite}. To keep the paper self
 contained we supplied the proofs for most of the facts. We refer the reader
 to \cite{Ma-S} and references therein for general facts.
 
 The subsection {\bf 5.3} shows that the walls for $\TET[1]$ defined by the
 intermediate subobjects of the filtration \eqref{filtTh1-intro} are at least as high as the wall defined by $\OO_X$ in \S4. This technical result is used to deduce the relations of the numerical invariants of the surface
 $X$ and the Chern invariants of the sheaves $F_k$ involved in the categorical differentials \eqref{cat-diff}, see Corollary \ref{cor:tau-bound}, Proposition \ref{pro:no1rankbound}, Corollary \ref{cor:tau-bound1}, Corollary \ref{cor:-c1Fi}. The main application of those results is for the case of canonical polarization: Proposition \ref{pro:HNalpha} contains all the statements of Theorem \ref{th:filtKX-intro}.
 
 In \S6 we take up the question of Bridgeland stability of $\OO_X$.
 The main result of the section is Proposition \ref{pro:OXBrsemist}:
 it identifies a part of the region found in \S4 where $\OO_X$ is Bridgeland
 semistable. 
 
 \S7 combines the results of previous two sections to develop the reduction of the categorical differential to the rank $2$ or $3$.
 The subsections {\bf 7.2} and {\bf 7.3} are devoted to each of those possibilities. 
 
 Proposition \ref{pro:QVrk2} summarizes all the properties of rank $2$ case.
 In the subsections {\bf 7.2.1} and {\bf 7.2.2} are treated the geometric properties of the categorical ramification divisor: Proposition \ref{pro:VDfilt} contains the results encapsulated in the `dictionary' \eqref{flag-X}, while Proposition \ref{pro:nefpair} treats the case of nef ramification divisor.
 
 The properties of the rank $3$ are collected in Proposition \ref{pro:QVrk3}. 
   
   All the result discussed until now are formulated for an arbitrary nonzero subspace $V\subset H^1 (\TET) $. In \S8 we specialize to 
   $V=H^1 (\TET)$. This is Theorem \ref{th:H1Theta} which is an expanded version of Theorem \ref{th:H1-filt-intro} formulated above.
   
  \S9 treats the surfaces with $(3c_2 (X)-K^2_X) $ small, where `small' means the condition \eqref{qbound-intro}. The main results here are 
  Corollary \ref{cor:f_X-fibration} and Corollary \ref{cor:fX-unique}
  which, taken together, correspond to Theorem \ref{th:fX-fibration-intro} stated above.
  The last result of the section, Proposition \ref{pro:imdp}, tells us that for regular surfaces the images of the maps $\PP^1 \longrightarrow \overline{\mathfrak{M}}_g$ resulting from Theorem \ref{th:fX-fibration-intro} sweep out at most a two dimensional subscheme of $\overline{\mathfrak{M}}_g$. It is given for future references and will be used in the sequel to this paper.
     
\section{Basics on surfaces of general type}\label{S-surf} 

In this section we collect some basic facts about surfaces of general type which are used throughout the paper. We also recall Mumford-Takemoto slope stability of coherent sheaves as well as Bogomolov (in)stability.

\vspace{0.2cm}
\noindent
{\bf 2.1.} Let $X$ be a smooth, compact, complex surface with the canonical divisor
$K_X$ ample. This means that $X$ is a minimal surface of general type. Its holomorphic tangent (resp. cotangent) bundle will be denoted by $\TET$ (resp. $\Omega_X$).
The Chern classes of $\TET$ (resp. $\Omega_X$) are
\begin{equation}\label{Chernclasses}
\mbox{$c_1 (\TET)=-K_X$ (resp. $c_1 (\Omega_X)=K_X$) and $c_2 (\TET)=c_2 (\Omega_X):=c_2 (X)$.}
\end{equation}
It is standard to identify the singular cohomology group $H^4 (X,\ZZ)$ with $\ZZ$ via the orientation class corresponding to the complex structure of $X$. With this identification in mind the cohomology classes
$$
K^2_X =c^2_1 (\TET), \,\, c_2 (X)
$$
in $H^4 (X,\ZZ)$ are thought of as integers and called the Chern numbers of $X$.
Their ratio is denoted
\begin{equation}\label{ratio-Chernnum}
\alpha_X := \frac{c_2 (X)}{K^2_X}
\end{equation}

We are concerned with minimal surfaces of general type. This puts restrictions on the Chern numbers:
\begin{equation}\label{restChnum}
	\begin{gathered}
	K^2_X, c_2  (X) \geq 1;\hspace{0.2cm} K^2_X \geq \frac{1}{5}(c_2(X)-36);
	\\
	 K^2_X \leq 3c_2 (X).
	 \end{gathered}
 \end{equation} 
 The restrictions in the first line are classical and go back to the classical Italian school of algebraic geometry and Max Noether, the inequality in the second line is the celebrated Bogomolov-Miyaoka-Yau inequality. All of the above cuts out the following region in the $(c_2,K^2)$-plane:
  
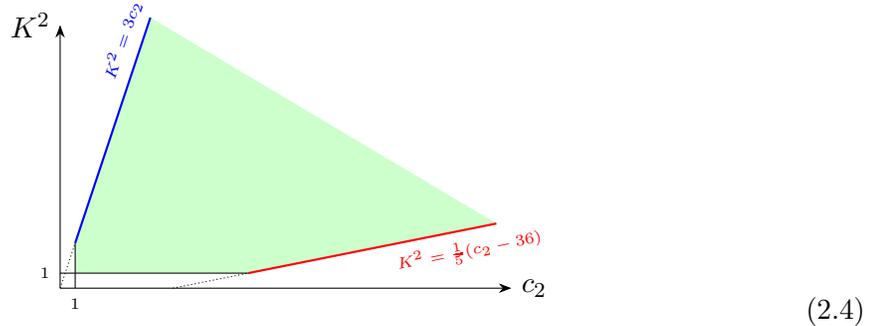
\begin{equation}\label{Pworld}
\begin{tikzpicture}[>=Stealth]
	\draw[->](0,0)--(6,0)
	node[right]{$c_2$};
	\path 
	(0.2,0) node[below]{\tiny{$1$}}
	(0,0.2) node[left]{\tiny{$1$}};
	\draw[->](0,0)--(0,3.5) node[left]{$K^2$};
	\draw[densely dotted](0,0)--(0.2,0.6);
	\draw[blue,thick](0.2,0.6)--(1.2,3.6)
	node[very near end,sloped,above]{\tiny{$K^2=3c_2$}};
	\draw[red,thick](2.5,0.2)--(5.8,0.86)
	node[very near end,sloped,below]{\tiny{$K^2=\frac{1}{5}(c_2 -36)$}};
	\draw[densely dotted](1.5,0)--(2.5,0.2);
	\draw[thin](0.2,0)--(0.2,0.6);
	\draw[thin](0,0.2)--(2.5,0.2);
	\colorlet{Pworld}{green!20!white}
	\begin{pgfonlayer}{background}
		\fill[Pworld](0.2,0.2)--(0.2,0.6)--(1.2,3.6)--(5.8,0.86)--(2.5,0.2)--(0.2,0.2);
	\end{pgfonlayer}
\end{tikzpicture}
\end{equation}

A point $(m,n)$ with integer coordinates lying in the above region is called {\it admissible} if there is a minimal surface $X$ of general type
with Chern numbers $c_2 (X)=m$ and $K^2_X =n$. 

An admissible point $(m,n)$ must in addition satisfy the divisibility condition
$$
m+n\equiv 0\!\!\!\!\mod(12).
$$
This comes from the Noether's formula
$$
\chi(\OO_X)=\frac{K^2_X +c_2(X)}{12},
$$
where the left hand side is the Euler characteristic of the structure sheaf $\OO_X$ of $X$; the formula is the Riemann-Roch for $\OO_X$, see \cite{GH}.

It is generally believed that {\it all} integer points in the area depicted in \eqref{Pworld} subject to the above mentioned divisibility condition are admissible. But the question does not seem to be completely settled.  
 
\vspace{0.2cm}

Another numerical invariant which will be used in our considerations is 
\begin{equation}\label{ch2}
\tau_X :=2ch_2 (X)=K^2_X -2c_2 (X).
\end{equation}
It is often called the index of $X$ for the reason we briefly recall. 

We have the intersection pairing 
$$
Q_X: H^2 (X,\RR) \times H^2 (X,\RR) \longrightarrow H^4 (X,\RR)\cong \RR.
$$
This is symmetric and non-degenerate. The corresponding quadratic form on $H^2 (X,\RR)$ will be also denoted by $Q_X$. This form is indefinite and its signature and index can be computed as follows.

It is well-known that in the Hodge decomposition
\begin{equation}\label{HDX}
H^2(X,\CC)=H^2 (X,\RR)\otimes \CC=H^{2,0} (X) \oplus H^{1,1} (X)\oplus H^{0,2} (X)
\end{equation}
the subspaces $(H^{2,0} (X)\oplus H^{0,2} (X))$ and $H^{1,1} (X)$ are invariant under the complex conjugation in $H^2(X,\CC)$. Their respective subspaces  $(H^{2,0} (X)\oplus H^{0,2} (X))^{fixed}$ and $(H^{1,1} (X))^{fixed}$ {\it fixed} under the complex conjugation are vector subspaces of
$H^2 (X,\RR)$. These two subspaces are orthogonal with respect to the intersection pairing $Q_X$ and furthermore $Q_X$ is positive definite on $(H^{2,0} (X)\oplus H^{0,2} (X))^{fixed}$ and has signature $(1,h^{1,1} (X)-1)$ on $(H^{1,1} (X))^{fixed}$, see \cite{GH} and \cite{Vo}. Thus the signature of $Q_X$ is 
$$
sign(Q_X)=(2p_g(X)+1, h^{1,1} (X)-1),
$$ 
where $p_g(X):=h^{2,0} (X)=h^0(\OO_X (K_X))$ is the geometric genus of $X$.
The index of $Q_X$, the difference of positive and negative definite parts, is
\begin{equation}\label{ind-calc}
\begin{gathered}
2p_g(X)+2 -h^{1,1} (X)=2p_g(X)+2 -(h^2(X)-2p_g(X))=4p_g(X)+2-h^2(X)\\
=4\chi(\OO_X)+4q(X) -2-h^2(X)=4\chi(\OO_X)-c_2(X)=4\frac{K^2_X +c_2(X)}{12}-c_2(X)\\=\frac{1}{3}(K^2_X -2c_2 (X))=\frac{1}{3}\tau_X,
\end{gathered}
\end{equation}
where 
$$
\chi(\OO_X):=1-q(X)+p_g(X)
$$
 is the Euler characteristic of $\OO_X$, $q(X):=h^1(\OO_X)$ is the irregularity of $X$ and where we used (the fourth equality from the left) that $c_2 (X)$ is the topological Euler characteristic of $X$
$$
c_2(X)=\sum^4_{i=0}(-1)^i h^i (X)=2-2h^1 (X) +h^2 (X)=2-4q(X)+h^2(X),
$$
and the Noether formula
$$
\chi(\OO_X)=\frac{K^2_X +c_2(X)}{12}
$$
in the fifth equality. From \eqref{ind-calc} it follows that $\tau_X$ essentially controls the index of $Q_X$. Abusing the terminology we will often refer to $\tau_X$ as the index of $X$. Surfaces for which $\tau_X >0$ are called surfaces of positive index. This is equivalent to the inequality
\begin{equation}\label{p-index-in}
K^2_X >2c_2 (X).
\end{equation}
In terms of the ratio $\alpha_X$ in \eqref{ratio-Chernnum} this translates into
\begin{equation}\label{p-index-in}
\alpha_X < \HA.
\end{equation}

\vspace{0.2cm}
\noindent
{\bf 2.2.} $NS(X)$ denotes the N\'eron-Severi group of $X$. By definition this is the image of the first Chern class map
$$
c_1: Pic(X) \longrightarrow H^2(X,\ZZ),
$$
where $Pic(X)$ is the Picard group of $X$, the group of the isomorphism classes of line bundles on $X$, and $c_1$ is the map which sends an isomorphism class of a line bundle ${\cal L}$ to its first Chern class $c_1({\cal L})$.

Tensoring $NS(X)$ with $\QQ$ (resp. $\RR$) we obtain the vector space
$$
\mbox{$NS(X)_{\QQ} := NS(X)\otimes_{\ZZ} \QQ$ \hspace{0.2cm}  (resp. $NS(X)_{\RR} := NS(X)\otimes_{\ZZ} \RR$).}
$$

Inside $NS(X)_{\RR}$ we will distinguish several cones. The first one is the cone of ample divisor classes
$$
\mbox{$Amp(X):=$ the convex cone in $NS(X)_{\RR}$ generated by ample divisor classes in $NS(X)$.}
$$
This cone is contained in the cone
$$
\mbox{$Eff(X):=$ the convex cone in $NS(X)_{\RR}$ generated by effective divisor classes.}
$$ 

The vector space $NS(X)_{\RR}$ comes equipped with the intersection form, the restriction to $NS(X)_{\RR} \subset H^2(X, \RR)$ of the pairing $Q_X$ in {\bf 2.1}. We continue to denote it by $Q_X$. By Hodge Index theorem, the signature of $Q_X$ on $NS(X)_{\RR}$ is
$$
(1, \rho_X -1),
$$
where $\rho_X =dim_{\RR} (NS(X)_{\RR})$ is the Picard number of $X$, see \cite{GH}. 

 The intersection pairing $Q_X$ gives the identification of  $NS(X)_{\RR}$ with its dual vector space $NS(X)^{\ast}_{\RR}$:
\begin{equation}\label{dual}
NS(X)_{\RR} \ni x\mapsto Q_X (x,\bullet) \in NS(X)^{\ast}_{\RR} .
\end{equation}
In the sequel we write $x\cdot y$ instead of $Q_X (x,y)$. The linear form
$Q_X (x,\bullet)$ will be denoted $\hat{x}$. In this notation we have
\begin{equation}\label{dualf}
\hat{x} (y)=x\cdot y,
\end{equation}
for $x,y \in NS(X)_{\RR}$. 

The cone of the effective divisor classes defines the cone $NE(X)$ of numerically effective (nef) divisor classes on $X$:
\begin{equation}\label{NE}
NE(X):=\{ x\in NS(X)_{\RR} | \,\,x\cdot C \geq 0, \,\forall C\in Eff(X)\}.
\end{equation}
This cone clearly contains the ample cone $Amp(X)$. In fact, it is known that $NE(X)$ is the closure of $Amp(X)$. More precisely, fix a norm $\parallel \bullet \parallel$ in $NS(X)_{\RR}$. All such norms are equivalent. Hence the resulting topology of $NS(X)_{\RR}$ is the same. We call this topology the norm topology of $NS(X)_{\RR}$ and whenever $NS(X)_{\RR}$ is viewed as a topological space, it is with respect to the norm topology.  

The closure $\overline{Amp(X)}$ of $Amp(X)$ in $NS(X)_{\RR}$, equipped with the norm topology, is the cone $NE(X)$ of the numerically effective divisor classes:
$$
\overline{Amp(X)}=NE(X).
$$

Another cone in $NS(X)_{\RR}$ is the cone of {\it big} or {\it positive} divisor classes
\begin{equation}\label{C+}
C^{+}(X):=\{ x\in NS(X)_{\RR} |\,\, \text{$x^2=x\cdot x > 0$ and $x\cdot h>0$, for some $h\in Amp(X)$} \}.
\end{equation}
We call it the positive cone of $X$. Observe, by the Hodge Index the cone $C^{+}(X)$ is independent of the choice of an ample divisor $h$ in the above definition. The divisor classes in $C^{+}(X)$ are called positive or big. It is obvious that ample divisors are positive. Thus the inclusion
$$
 C^{+}(X) \supset Amp(X).
$$

\vspace{0.2cm}
\noindent
{\bf 2.3.} Denote by ${\cal A}_X:=Coh(X)$, the abelian category of coherent $\OO_X$-modules on $X$. Fix an ample divisor class $h\in Amp(X)$. Following Mumford and Takemoto one defines the slope function on ${\cal A}_X$ by the rule
\begin{equation}\label{slope}
\mu_h ({\cal S})=\begin{cases} \frac{c_1({\cal S}) \cdot h}{rk({\cal S})},&\text{if $rk({\cal S})\neq 0$},\\
+\infty, &\text{if $rk({\cal S})= 0$},
\end{cases}
\end{equation}
where $rk({\cal S})$ (resp. $c_1({\cal S})$) is the rank (resp. the first Chern class) of a sheaf ${\cal S} \in {\cal A}_X$. 

The following is the standard notion of slope (semi)stability of sheaves in
${\cal A}_X$.

\begin{defi}\label{MT-stability}
A nonzero sheaf ${\cal S}$ in ${\cal A}_X$ is called $h$-semistable if for any nonzero proper subsheaf 
${\cal S'}$ of ${\cal S}$ the inequality 
$$
\mu_h ({\cal S'}) \leq \mu_h ({\cal S}/{\cal S'})
$$
holds; ${\cal S}$ is called $h$-stable, if the inequality above is strict.
\end{defi}

\begin{example}\label{ex:semist}
1) Any torsion sheaf on $X$ is $h$-semistable, for any $h\in Amp(X)$.

2) Any torsion free sheaf of rank $1$ is $h$-stable, for any $h\in Amp(X)$.

3) For $X$ with the canoncal divisor $K_X$ ample, the holomorphic tangent bundle $\TET$ (resp. $\Omega_X$) is $K_X$-semistable, a result of H. Tsuji, \cite{T}.
\end{example}

It is well known that every nonzero sheaf admits a filtration with $h$-semistable quotients having strictly decreasing slopes. Namely, for a nonzero ${\cal S}$ in ${\cal A}_X$ there is a finite step filtration
$$
{\cal S}={\cal S}_n \supset {\cal S}_{n-1} \supset \cdots \supset {\cal S}_{1} 
\supset {\cal S}_0 =0
$$
such that the quotient sheaves ${\cal Q}_i :={\cal S}_i /{\cal S}_{i-1}$, for $i=1,\ldots, n$, are all $h$-semistable and their slopes $\mu_h({\cal Q}_i)$ form a strictly decreasing function of $i$. In other words we have the inequalities
$$
\mu_h({\cal Q}_1)> \mu_h({\cal Q}_2) > \cdots > \mu_h({\cal Q}_n).
$$
The filtration with above properties is unique and is called the Harder-Narasimhan (HN, for short) filtration of ${\cal S}$ with respect to $h$.
Observe: the quotient sheaf ${\cal Q}_1$ is equal to ${\cal S}_1$, the first sheaf from the right of the HN filtration of ${\cal S}$. That sheaf is called {\it the maximal $h$-destabilizing subsheaf} of ${\cal S}$.

In Example \ref{ex:semist}, 3), we recalled that $\TET$ (resp. $\Omega_X$) is $K_X$-semistable. If $X$ is a surface of positive index the following holds.
\begin{lem}\label{l:TX-Kst}
For a surface $X$ with the index $\tau_X >0$ the (co)tangent bundle $\TET$ (resp. $\Omega_X$) is $K_X$-stable.
\end{lem}
\begin{pf} The sheaves $\TET$ and $\Omega_X$ are dual to each other, so it is enough to give a prove for one of them. We do it for $\Omega_X$. Let ${\cal L}$ be a nonzero proper subsheaf of $\Omega_X$. Then it is torsion free of rank $1$. By replacing ${\cal L}$ by its reflexive hall, that is, by taking the double dual
of ${\cal L}$, we may assume ${\cal L}$ to be a line bundle. Denote it by $\OO_X (L)$. Furthermore, we may assume the inclusion $\OO_X (L) \subset \Omega_X$ to be saturated, that is, the quotient sheaf is torsion free. This gives the following exact sequence
\begin{equation}\label{dest-L}
\xymatrix@R=12pt@C=12pt{
0\ar[r]& \OO_X ( L)\ar[r]& \Omega_X\ar[r]& {\cal I}_Z (M)\ar[r]& 0
}
\end{equation}
where ${\cal I}_Z $ is the ideal sheaf of a subscheme $Z$ of dimension at most $0$ and $M$ is a divisor. We know that $\Omega_X $ is $K_X$-semistable meaning that $L\cdot K_X \leq \HA K^2_X $. We need to check that the inequality must be strict. Assume the equality:
\begin{equation}\label{half}
L\cdot K_X = \HA K^2_X.
\end{equation}
Combining this with $K^2_X >0$ we obtain $L\cdot K_X >0$. Since line subsheaves of $\Omega_X $ can not be positive, Bogomolov lemma, see \cite{Mi}, we deduce $L^2 \leq 0$.
The equality \eqref{half} also implies
$$
M\cdot K_X =(K_X -L) \cdot K_X =K^2_X - L \cdot K_X =\HA K^2_X.
$$
From this it follows
$$
0=M\cdot K_X -L\cdot K_X=(M-L)\cdot K_X=(M-L)\cdot (M+L)=M^2 - L^2.
$$
But the exact sequence \eqref{dest-L} tells us
$$
\tau_X =2ch_2(\Omega_X) =2ch_2(\OO_X ( L))+2ch_2({\cal I}_Z (M))=L^2 +M^2 -2deg(Z) =2L^2-2deg(Z) \leq 0,
$$
contradicting the positivity of $\tau_X$.
\end{pf}

\vspace{0.2cm}
\noindent
{\bf 2.4.} Define
\begin{equation}\label{Unst}
Un_X := \{h \in Amp(X) | \, \text{$\Omega_X$ is $h$-unstable} \}.
\end{equation}
From Example \ref{ex:semist}, 3), this is a proper subset of the ample cone $Amp(X)$. We assume it to be nonempty and wish to  understand this subset.

Let $h\in Un_X $. The maximal $h$-destabilizing subsheaf of $\Omega_X$ is a locally free subsheaf of rank $1$. We denote it by $\OO_X (L_h)$. It gives rise to the $h$-destabilizing sequence
\begin{equation}\label{h-destseq}
\xymatrix@R=12pt@C=12pt{
0\ar[r]& \OO_X ( L_h)\ar[r]& \Omega_X\ar[r]& {\cal I}_{Z_h} (M_h)\ar[r]& 0.
}
\end{equation}
Observe: the quotient sheaf is torsion free by $h$-maximality of $\OO_X ( L_h)$.
The numerical condition for $\OO_X ( L_h)$ being $h$-destabilizing is
\begin{equation}\label{unst-cond}
L_h \cdot h > \HA K_X \cdot h.
\end{equation}
The properties of the destabilizing sequence above are summarized in the following statement.
\begin{lem}\label{l:h-unstseq}
1) $K_X = L_h+M_h$.

\vspace{0.2cm}
2) $2c_2(X)=K^2_X  - L^2_h -M^2_h +2deg(Z_h).$

\vspace{0.2cm}
3) $L_h \cdot h > \HA K_X \cdot h$ and $L^2_h \leq 0$.

\vspace{0.2cm}
4) $M_h\cdot h' \geq 0$, for any $h'\in Amp(X)$. In addition, if the index $\tau_X >0$, then
$M^2_h >0$ and $M_h$ lies in the positive cone $C^{+}(X)$ of $X$.
\end{lem}
\begin{pf}
The items $1) - 2)$ follows immediately from the destabilizing sequence \eqref{h-destseq}. The first assertion in $3)$ is \eqref{unst-cond} and the second assertion follows from the first and the Bogomolov lemma: $L_h$ can not lie in the positive cone of $X$.

The first assertion in $4)$ is the result of Miyaoka on the generic semipositivity of $\Omega_X$, see \cite{Mi2}. For the second assertion we write
$$
\tau_X=2ch_2(\Omega_X)= 2ch_2(\OO_X ( L_h))+2ch_2({\cal I}_{Z_h} (M_h))=L^2_h +M^2_h -2deg(Z_h) \leq M^2_h -2deg(Z_h),
$$
where the last inequality uses $L^2_h \leq 0$ in $3)$ of the lemma. From the above and $\tau_X >0$ we deduce
$$
M^2_h >2deg(Z_h)\geq 0.
$$
This together with the first assertion of 4) and the Hodge index imply $M_h\cdot h' > 0$, for any ample divisor class $h'$. Thus $M_h \in C^{+}(X)$ as asserted. 
\end{pf}

From now on we assume that $X$ is a surface of positive index. That is
\begin{equation}\label{cond-pi}
\tau_X >0.
\end{equation}
For $h\in Un_X$ the divisor $M_h$ in the quotient of the destabilizing sequence \eqref{h-destseq} gives rise to the linear map
$$
\widehat{(M_h -\HA K_X)} :NS(X)_{\RR} \longrightarrow \RR.
$$
The hyperplane
\begin{equation}\label{H_h}
H_{M_h}:=ker(\widehat{(M_h -\HA K_X)})=\{x\in NS(X)_{\RR} |\,(M_h -\HA K_X)\cdot x =0\}
\end{equation}
cuts the ample cone $Amp(X)$ into two disjoint nonempty open subcones
$$
\begin{gathered}
Amp^{+}_{M_h}:=\{x\in Amp(X)|\,(M_h -\HA K_X)\cdot x >0\},\\ 
 Amp^{-}_{M_h}:=\{x\in Amp(X)|\,(M_h -\HA K_X)\cdot x <0\};
\end{gathered} 
$$
that these are disjoint open convex cones is obvious and the fact that they are nonempty comes from $K_X \in Amp^{+}_{M_h}$, see Lemma \ref{l:TX-Kst}, and $h\in Amp^{-}_{M_h}$, see \eqref{unst-cond}.

Motivated by the above discussion we call an extension sequence
\begin{equation}\label{dest-seq}
\xymatrix@R=12pt@C=12pt{
0\ar[r]& \OO_X ( L)\ar[r]& \Omega_X\ar[r]& {\cal I}_{Z} (M)\ar[r]& 0
}
\end{equation}
{\it destabilizing} for $\Omega_X$ if the divisor $M$ has the following properties:

$(i)$ $ M \in C^{+}(X)$,

$(ii)$ the cones  $Amp^{\pm}_{M} $ are nonempty, where we let
\begin{equation}\label{M-cones}
\begin{gathered}
Amp^{+}_{M}:=\{x\in Amp(X)|\,(M -\HA K_X)\cdot x >0\},\\ 
 Amp^{-}_{M}:=\{x\in Amp(X)|\,(M -\HA K_X)\cdot x <0\};
\end{gathered} 
\end{equation}
\begin{rem}\label{rem:L2} From the nonemptiness of $ Amp^{-}_{M}$ it follows that
$L\cdot h >\HA K_X \cdot h$, for all $h\in Amp^{-}_{M}$. This and the fact that 
$L$ can not be in the positive cone of $X$ implies
$$
L^2 \leq 0,
$$ 
for every destabilizing sequence for $\Omega_X$.
\end{rem}   
We denote
\begin{equation}\label{DestX}
\mbox{$\mathfrak{Dest}_X:=$ the set of the destabilizing for $\Omega_X$ sequences.}
\end{equation}
Taking the first Chern class of the quotient term of the destabilizing sequences gives the map
\begin{equation}\label{c1-map}
c^{quot}_1: \mathfrak{Dest}_X \longrightarrow NS^{+} (X):=NS (X)\cap C^{+}(X).
\end{equation}
 The image of this map will be denoted by $Dest_X $. In particular, a divisor class $M\in Dest_X$ has the properties $(i)-(ii)$ listed after the extension sequence \eqref{dest-seq}. Furthermore, Lemma \ref{l:h-unstseq} gives a relation between $Dest_X $ and the set $Un_X$ defined in \eqref{Unst}. Namely, define the map
\begin{equation}\label{mu-map}
\mu: Un_X \longrightarrow Dest_X
\end{equation}
by the formula
$$
\mu (h)=M_h, \,\forall h\in Un_X
$$
where $M_h$ is as in the $h$-destabilizing sequence \eqref{h-destseq}. The map 
$\mu$ partitions $Un_X$ into disjoint subsets
$$
 Un_X =\bigcup_{M\in  Dest_X} \mu^{-1} (M).
 $$
 \begin{lem}\label{l:parts}
   Each part $\mu^{-1} (M)$ equals the negative cone $Amp^{-}_{M}$ defined in \eqref{M-cones}.
   \end{lem}
   \begin{pf}
     We have the destabilizing sequence
     \begin{equation}\label{M-seq}
   \xymatrix@R=12pt@C=12pt{
0\ar[r]& \OO_X ( L)\ar[r]& \Omega_X\ar[r]& {\cal I}_{Z} (M)\ar[r]& 0
}
\end{equation}
for $\Omega_X$ corresponding to $M$. By definition  $h\in \mu^{-1} (M)$ if and only if
$M_h$ in the $h$-destabilizing sequence \eqref{h-destseq} equals $M$. Hence $h\cdot(M -\HA K_X)=h\cdot(M_h -\HA K_X)<0$. Thus $h\in Amp^{-}_{M}$ and we have the inclusion
$$
\mu^{-1} (M) \subset Amp^{-}_{M}.
$$
Let $h\in Amp^{-}_{M}$. This means  $\OO_X ( L)$ in \eqref{M-seq} is an $h$-destabilizing subsheaf of $\Omega_X$. So $\Omega_X$ is $h$-unstable and we need to check that $\OO_X ( L)$ is the maximal $h$-destabilizing subsheaf. In other words the $h$-destabilizing sequence
$$
\xymatrix@R=12pt@C=12pt{
  0\ar[r]& \OO_X ( L_h)\ar[r]& \Omega_X\ar[r]& {\cal I}_{Z_h} (M_h)\ar[r]& 0
}
$$
is equal to the one in \eqref{M-seq}. For this we put two sequences together to obtain the diagram
$$
\xymatrix@R=12pt@C=12pt{
  &&0\ar[d]&&\\
  &&\OO_X ( L)\ar[d]\ar[dr]&&\\
  0\ar[r]& \OO_X ( L_h)\ar[r]& \Omega_X\ar[r]\ar[d]& {\cal I}_{Z_h} (M_h)\ar[r]& 0\\
  &&{\cal I}_{Z} (M)\ar[d]&&\\
  &&0&&
}
$$
where the slanted arrow is the composition of the relevant arrows of the horizontal and vertical exact sequences. To prove the equality of two sequences is equivalent to showing that the slanted arrow is zero. Suppose it is not. Then the divisor
$$
E=M_h-L
$$
is effective. Hence
$$
0\leq E\cdot h=(M_h-L)\cdot h
$$
But $M_h\cdot h < \HA K_X \cdot h$ and $L\cdot h >\HA K_X \cdot h$ implying
$$
(M_h-L)\cdot h < 0.
$$
\end{pf}

Thus we obtain a partition
 \begin{equation}\label{partUn}
   Un_X =\bigcup_{M\in  Dest_X} Amp^{-}_{M}
 \end{equation}
 of $ Un_X$ into the disjoint union of convex open subcones $ Amp^{-}_{M}$ of $Amp(X)$. The subcones are labeled by the divisor classes $M$ in $Dest_X$. The latter, as we have seen above, is a discrete subset in $NS^{+} (X)=NS (X)\cap C^{+}(X) \subset C^{+}(X)$.

 If we fix $x_0 \in  C^{+}(X)$, we can stratify $Dest_X$ according to the `height' of the divisors with respect to $x_0$. Namely, the intersection with $x_0$ gives a positive valued function
 $$
 \widehat{x_0}:Dest_X \longrightarrow \RR_{>0}:=(0,+\infty),
 $$
measuring the  height (= degree) of the divisors in $Dest_X$ with respect to $x_0$. For every positive value $a$ define
 \begin{equation}\label{a-hight}
   Dest^{\leq a}_X (x_0):=\widehat{x_0}^{-1}((0,a])=\{ M\in Dest_X |\, M\cdot x_0 \leq a \}.
 \end{equation}
 \begin{lem} 1) For every positive value $a$ the set $Dest^{\leq a}_X (x_0)$ is at most finite.

   2) For $a'>a$ in $\RR_{>0}$ we have the inclusion
   $$
   Dest^{\leq a}_X (x_0) \subset  Dest^{\leq a'}_X (x_0).
   $$

   3) $\displaystyle{Dest_X =\bigcup_{a\in \RR_{>0}} Dest^{\leq a}_X (x_0)}$.
 \end{lem}
 \begin{pf}
   The items 2) and 3) are obvious from the definitions of the function $\widehat{x_0}$ and the sets $Dest^{\leq a}_X (x_0)$.

   The first statement of the lemma is an easy consequence of the Hodge index theorem.  Assume $Dest^{\leq a}_X (x_0)$ to be nonempty, otherwise there is nothing to prove. The signature of the intersection form $Q_X$ on $NS(X)_{\RR}$ is $(1,\rho_X -1)$. The vector $x_0$ spans the subspace of $NS(X)_{\RR}$ on which $Q_X$ is positive definite. On the orthogonal complement
   $x^{\perp}_0$ the form $Q_X$ is negative definite. Choose an orthonormal basis $x_1, \ldots, x_{\rho_X -1}$ on $(x^{\perp}_0, -Q_X)$. Then any $M\in Dest^{\leq a}_X (x_0)$ is expressed uniquely as follows
   $$
   M=\frac{M\cdot x_0}{x^2_0} x_0 +\sum^{\rho_X -1}_{i=1} a^M_i x_i,
   $$
   for some $a^M_i \in \RR$. Taking the self-intersection of $M$ gives
   $$
   M^2 = \frac{(M\cdot x_0)^2}{x^2_0} -\sum^{\rho_X -1}_{i=1} (a^M_i )^2.
   $$
   From $M^2 >0$ (because $M\in C^{+}(X)$) we obtain
   $$
   \frac{(M\cdot x_0)^2}{x^2_0} -\sum^{\rho_X -1}_{i=1} (a^M_i )^2>0,
   $$
   or, equivalently,
   $$
   \sum^{\rho_X -1}_{i=1} (a^M_i )^2 <\frac{(M\cdot x_0)^2}{x^2_0} \leq \frac{a^2}{x^2_0},
   $$
   where the last inequality uses the assumption $M \in Dest^{\leq a}$, that is,
   $M\cdot x_0 \leq a$.

   Now choose the norm $\parallel \bullet \parallel$ defined by the positive definite form $Q=Q_X \Big|_{\RR x_0} \oplus (-Q_X)\Big|_{x^{\perp}_0}$ on $NS(X)_{\RR}$. The above inequality tells us that the norm of $M \in Dest^{\leq a}$ is bounded
   $$
   \parallel M\parallel^2 = \frac{(M\cdot x_0)^2}{x^2_0}+\sum^{\rho_X -1}_{i=1} (a^M_i )^2 <2\frac{(M\cdot x_0)^2}{x^2_0}\leq 2\frac{a^2}{x^2_0}.
   $$
   Thus $Dest^{\leq a}$ is contained in the compact ball
   $$
   B_a =\Big \{b\in NS(X)_{\RR} | \parallel b \parallel \leq  a\Big(\frac{2}{x^2_0}\Big)^{\HA} \Big \}.
   $$
 Since $Dest^{\leq a}$ is discrete it must be a finite subset of $B_a$.  
   \end{pf}

From the definition of the partition \eqref{partUn} we know that the cones $\{Amp^{-}_M\}_{M\in Dest_X}$ are pairwise disjoint. It will be instructive to see this property from a more geometric perspective.
\begin{lem}\label{l:EMM'}
Let $M$ and $M'$ be two distinct divisor classes in $Dest_X$. Then the divisor
$$
E_{M,M'}:=M+M'-K_X
$$
is a nonzero effective divisor.
\end{lem}
\begin{pf}
Each $M$ and  $M'$ comes from its own destabilizing sequence
$$
\xymatrix@R=12pt@C=12pt{
0\ar[r]& \OO_X ( L)\ar[r]& \Omega_X\ar[r]& {\cal I}_{Z} (M)\ar[r]& 0,
\\
0\ar[r]& \OO_X ( L')\ar[r]& \Omega_X\ar[r]& {\cal I}_{Z'} (M')\ar[r]& 0.
}
$$
Putting them together gives the following diagram
$$
\xymatrix@R=12pt@C=12pt{
&&0\ar[d]&&\\
& &\OO_X ( L')\ar[d]\ar[dr]&&\\
0\ar[r]& \OO_X ( L)\ar[r]\ar[dr]& \Omega_X\ar[r]\ar[d]& {\cal I}_{Z} (M)\ar[r]& 0\\
& &{\cal I}_{Z'} (M')\ar[d]&&\\
&&0&&
}
$$
where the slanted arrows are the compositions of the relevant arrows in the vertical and the horizontal sequences. The slanted arrow in the upper right corner is nonzero: otherwise the monomorphism in the vertical sequence factors through
$\OO_X ( L)$ giving the isomorphism $\OO_X ( L)\cong \OO_X ( L')$ and hence
the equality $M=M'$ and this is contrary to the assumption that $M$ and $M'$ are distinct.

Once the slanted arrow in the upper right corner is nonzero it produces the effective divisor
$$
E_{M,M'}=M-L'=M-(K_X-M')=M+M'-K_X.
$$
It is nonzero: otherwise $M=L'$; but this is impossible, since $M$ is in $C^{+}(X)$ and $(L')^2 \leq 0$ by Remark \ref{rem:L2}. 
\end{pf}

Denote 
\begin{equation}\label{AmpM0}
Amp^0_M =\{x\in Amp(X) | \, x\cdot (M-\HA K_X)=0\}
\end{equation}
the intersection of the hyperplane 
$$
H_M=ker(\widehat{(M-\HA K_X)})
$$
with the ample cone $Amp(X)$. We also set
$$
\begin{gathered}
Amp^{\leq 0}_M=Amp^0_M \cup Amp^{-}_M,\\
Amp^{\geq 0}_M=Amp^0_M \cup Amp^{+}_M,
\end{gathered}
$$
the closure of the cones $Amp^{\pm}_M$ in $Amp(X)$.
\begin{lem}\label{l:MM'}
Let $M$ and $M'$ be two distinct divisor classes in $Dest_X$. Then the subcones
$Amp^{\leq 0}_M$ and $Amp^{\leq 0}_{M'}$ are disjoint.
\end{lem}
\begin{pf}
From Lemma \ref{l:EMM'} we know
$$
E_{M,M'}=M+M'-K_X
$$
is an effective nonzero divisor. We can write it
$$
E_{M,M'}=(M-\HA K_X)+(M'-\HA K_X).
$$
From this it follows
$$
0<E_{M,M'} \cdot h= (M-\HA K_X)\cdot h+(M'-\HA K_X)\cdot h,
$$
for any $h\in Amp(X)$. In particular, taking $h\in Amp^{\leq 0}_M $, the above inequality gives 
$$
(M'-\HA K_X)\cdot h >0,
$$
that is, $Amp^{\leq 0}_M $ is contained in $Amp^{+}_{M'}$. The latter by definition is the complement of $Amp^{\leq 0}_{M'}$ in $Amp(X)$. Hence 
$$
Amp^{\leq 0}_M \cap Amp^{\leq 0}_{M'} =\emptyset.
$$
\end{pf}
\begin{cor}\label{c:AM-intM'}
  For every $M\in Dest_X$ we have the inclusion
  $$
  Amp^{\leq 0}_M  \subset \bigcap_{M'\neq M} Amp^{+}_{M'}.
  $$
\end{cor}
\begin{pf}
  From Lemma \ref{l:MM'} we know that $Amp^{\leq 0}_M \subset  Amp^{+}_{M'}$ for every $M'\neq M$. Hence the asserted inclusion.
\end{pf}

Set
\begin{equation}\label{SemStcone}
  \begin{gathered}
    Amp^{\geq 0}_X:= \bigcap_{ M\in Dest_X} Amp^{\geq 0}_{M}\\
    Amp^{+}_X:= \bigcap_{ M\in Dest_X} Amp^{+}_{M}.
  \end{gathered}
\end{equation}
\begin{pro}\label{pro:stcone}
 1) $Amp^{\geq 0}_X$ is the nonempty convex cone of polarizations of $X$ for which $\Omega_X$ is semistable.

2) The complement  of $Amp^{ +}_X$ in $Amp^{\geq 0}_X $ is the disjoint union of codimension $1$ subcones $Amp^{0}_M$, where $M\in Dest_X$. That is, we have the equality
$$
Amp^{\geq 0}_X \setminus Amp^{+}_X = \coprod_{M\in Dest_X} Amp^{0}_M.
$$ 
This disjoint union is a subset of polarizations of $X$ for which $\Omega_X$ is strictly semistable.

3) $Amp^{+}_X$ is the nonempty convex cone of polarizations of $X$ for which $\Omega_X$ is stable.
\end{pro}
 \begin{pf}
By definition $Amp^{\geq 0}_X$ (resp. $Amp^{+}_X$ ) is the intersection of convex cones. So both, $Amp^{\geq 0}_X$ and $Amp^{+}_X$, convex cones. The fact that $Amp^{+}_X$ is nonempty is Lemma \ref{l:TX-Kst} asserting that $\Omega_X$ is $K_X$-stable. 
 Hence $K_X \in Amp^{+}_X$. The inclusion $Amp^+_X \subset Amp^{\geq 0}_X$ implies that both cones are nonempty.

By definition of $Un_X$, see \eqref{Unst}, the complement of this set in $Amp(X)$ is the set of polarizations for which
$\Omega_X$ is semistable. From the partition of $Un_X$ in \eqref{partUn} we have
$$
(Un_X)^c=(\bigcup_{M\in Dest_X} Amp^{-}_{M})^c= \bigcap_{M\in Dest_X} (Amp^{-}_{M})^c = \bigcap_{M\in Dest_X} Amp^{\geq 0}_{M}=Amp^{\geq 0}_X.
  $$
The inclusion 
$$
Amp^{\geq 0}_X \setminus Amp^{+}_X \subset \bigcup_{M\in Dest_X} Amp^{0}_M
$$
follows from the definition of  $Amp^{\geq 0}_X$. The opposite inclusion comes from Corollary \ref{c:AM-intM'}. Thus the equality
$$
Amp^{\geq 0}_X \setminus Amp^{+}_X = \bigcup_{M\in Dest_X} Amp^{0}_M
$$
 By Lemma \ref{l:MM'} the union above is disjoint. Hence the formula
$$
Amp^{\geq 0}_X \setminus Amp^{+}_X = \coprod_{M\in Dest_X} Amp^{0}_M
$$ 
in the second item of the proposition.

From the inclusion  $ Amp^{0}_M \subset Amp^{\geq 0}_X$ and the part 1) of the proposition, we know that $ Amp^{0}_M$ is a subset of polarizations for which $\Omega_X$ is semistable. Since we have an extension sequence
$$
\xymatrix@R=12pt@C=12pt{
0\ar[r]& \OO_X ( L)\ar[r]& \Omega_X\ar[r]& {\cal I}_{Z} (M)\ar[r]& 0
}
$$
subject to $(M-\HA K_X)\cdot h=0$, for all $h\in Amp^{0}_M$, it follows that
$\Omega_X$ is strictly semistable for all $h\in Amp^{0}_M$.

It remains to show that  $ Amp^{+}_X$ is the set of polarizations for which $\Omega_X$ is stable. Clearly, all stable polarizations are in
$ Amp^{+}_X$. Conversely, let $h\in  Amp^{+}_X$, we need to check that
$\Omega_X$ is $h$-stable. We already know that it is $h$-semistable.
Suppose it is strictly semistable with respect to $h$. This means that there is an exact sequence
\begin{equation}\label{semst-h}
\xymatrix@R=12pt@C=12pt{
	0\ar[r]& \OO_X ( L_0)\ar[r]& \Omega_X\ar[r]& {\cal I}_{Z} (M_0)\ar[r]& 0
}
\end{equation}
for which $L_0\cdot h=M_0\cdot h=\HA K_X \cdot h$. This tells us that
$h$ lies in the face
$$
A^0_{M_0} =Amp(X)\bigcap ker(\widehat{(M_0-\HA K_X)}).
$$
This is a closed proper subset of codimension $1$ in $Amp(X)$, so its complement consists of two nonempty disjoint subcones denoted
$A^+_{M_0}$ and $A^-_{M_0}$, where the linear function
$\widehat{(M_0-\HA K_X)}$ is respectively positive and negative.
Furthermore, from the strict semistability of the exact sequence \eqref{semst-h} it follows that the divisor $M_0$ lies in the positive cone of $X$. Hence that exact sequence lies in the set $\mathfrak{Dest}_X$
of destabilizing sequences for $\Omega_X$ and the divisor $M_0$ lies in the set $Dest_X$. But then, by the part 2) of the proposition,
$h$ lies $Amp^{\geq 0}_X \setminus Amp^{+}_X$ and this is contrary to our assumption.  
\end{pf}

From now on we set
\begin{equation}\label{st-cone}
Amp^s_X:=\{h\in Amp(X) | \mbox{\,$\Omega_X$ is $h$-stable}\}
\end{equation}
and what we have learned from Proposition \ref{pro:stcone} is the equality
\begin{equation}\label{st-cone=+}
	Amp^s_X =Amp^+_X
\end{equation}
The discussion above gives us the `classical' picture of the structure of the ample cone $Amp(X)$ of $X$ with respect to the slope stability
of $\Omega_X$ (resp. $\TET$). Let us summarize:

\vspace{0.2cm}
$\bullet$ inside $Amp(X)$ there is the subcone
$Amp^s_X$, see \eqref{st-cone} and \eqref{st-cone=+}, of polarizations with respect to which $\Omega_X$ is slope-stable;

\vspace{0.2cm}
$\bullet$ the cone $Amp^s_X$ is surrounded by faces, $Amp^0_M$ labeled
by the divisor classes $M \in C^+(X)$ appearing in a destabilizing sequence
of the form
\begin{equation}\label{L-M-seq}
\xymatrix@R=12pt@C=12pt{
	0\ar[r]& \OO_X ( L)\ar[r]& \Omega_X\ar[r]& {\cal I}_{Z} (M)\ar[r]& 0,
}
\end{equation}
where $L$ and $M$ are subject to the properties of Lemma \ref{l:h-unstseq};

the face $Amp^0_M$ is defined as the intersection of $Amp(X)$ with the hyperplane
$$
H_M=ker(\widehat{(M-\HA K_X)})=\{x\in NS(X)_{\RR} |\,\, (M-\HA K_X)\cdot x=0\};
$$

\vspace{0.2cm}
$\bullet$ each $M$ in the previous item decomposes the ample cone into the disjoint union of subcones
$$
Amp(X)=Amp^0_M \bigcup Amp^{+}_M \bigcup Amp^{-}_M,
$$
where $Amp^{\pm}_M$ are the connected components of the intersection of
$
Amp(X)
$
with the complement of $H_M$; 

the faces $Amp^0_M$ are codimension $1$ subcones of
$Amp(X)$ and parameterize polarizations for which the exact sequence 
\eqref{L-M-seq} is strictly semistable; 

crossing a face $Amp^0_M$ from  $ Amp^{+}_M $ to $ Amp^{-}_M$ turns
\eqref{L-M-seq} into a destabilizing sequence of $\Omega_X$: for a continuous path $\gamma:[0,1] \longrightarrow Amp(X) $ starting at
$h_0=\gamma(0)\in A^+_M$ and ending at $h_1=\gamma(1)\in A^-_M$, the subsheaf
$\OO_X (L)$ becomes maximal $h_1$-destabilizing subsheaf of $\Omega_X$.

\vspace{0.2cm}
What kind of geometry can be extracted from the above? Part of the geometric data is contained in the decomposition
$$
K_X =L+M
$$
of the canonical divisor and its properties collected in Lemma \ref{l:h-unstseq}. Another part comes from dualizing the exact sequence
\eqref{L-M-seq}

\begin{equation}\label{L-M-seq-foli}
	\xymatrix@R=12pt@C=12pt{
		0\ar[r]& \OO_X ( -M)\ar[r]& \Theta_X\ar[r]& {\cal I}_{Z} (-L)\ar[r]& 0.
	}
\end{equation}
This is a (singular along $Z$) foliation of $X$. So we have the identification
$$
\mathfrak{Dest}_X  \cong \mathfrak{Fol}_X
$$
where on the left is the set of the destabilizing sequences introduced in \eqref{DestX} and on the right
$$
\mathfrak{Fol}_X:= \text{the set of exact sequences dual to the ones in
$\mathfrak{Dest}_X$}
$$	
The foliation in \eqref{L-M-seq-foli} is what the experts in the subject call foliations of Kodaira dimension $2$, see \cite{McQ}, that is, the Kodaira dimension
of the {\it dual} of the subsheaf in \eqref{L-M-seq-foli} is $2$: this is because $M$ according to Lemma \ref{l:h-unstseq} is in the positive
cone of $X$.

From the foliation \eqref{L-M-seq-foli} we also deduce the filtration:
\begin{equation}\label{KS-filt-1}
	H^1(\TET)\supset H^1 (\OO_X (-M)).
\end{equation}
This is a sort of classical counterpart of the filtration in Theorem \ref{th-H-filt}. This is nonvacuous if $H^1 (\OO_X (-M))$ is nonzero and
this nonvanishing has a precise geometric meaning.
\begin{pro}\label{pro:-M}
If $H^1 (\OO_X (-M))$ is nonzero, then there are smooth rational curves
$R$ on $X$ subject to the following:

a) $R\cdot M=-1$,

b) each such curve 
$R$ is an integral curve of the foliation and it passes through a single point $z_R$ of the singularity subscheme $Z$,

c) all such curves are mutually disjoint.
\end{pro} 
\begin{pf}
We know that $M$ is in the positive cone of $X$. So $M$ is big and the nonvanishing of $H^1 (\OO_X (-M))$ implies that $M$ is not nef. Hence there are curves on $X$ intersecting $M$ negatively. Furthermore, all such curves
are in the negative part of the Zariski decomposition of $M$. Hence they form a negative definite sublattice of the N\'eron-Severi group $NS(X)$ of $X$. 

	Let $C$ be a reduced irreducible curve with $M\cdot C <0$. From the above
	it follows
	$C^2 < 0$. The restriction to $C$ of the monomorphism in \eqref{L-M-seq-foli} put together with the normal sequence of $C$ gives
	the diagram
	$$
	\xymatrix@R=12pt@C=12pt{
		&&0\ar[d]\\
		&&\Theta_C \ar[d]\\
	0\ar[r]& \OO_C ( -M)\ar[r]& \Theta_X \otimes \OO_C \ar[d]\\
	&&\OO_C (C)
}
	$$
	where $\Theta_C$ is the tangent sheaf of $C$. Furthermore, the horizontal arrow factors through $\Theta_C$: $Hom(\OO_C (-M),\OO_C (C))=H^0 (\OO_C (M+C))=0$, because $C\cdot (M+C) <0$. The resulting morphism
	$$
	\OO_C ( -M) \longrightarrow \Theta_C
	$$
	is nonzero and vanishes precisely on the $0$-dimensional subscheme
	$$
	Z^C=Z\cap C
	$$
	where the singular scheme of the foliation $Z$ intersects $C$. This means that $C$ is an integral curve of the foliation. The morphism, being an isomorphism outside $Z^C$, also tells us that $\Theta_C$ is locally free outside $Z^C$. By Lipman's criterion of smoothness, \cite{Li}, the curve $C$ is smooth outside $Z^C$.
	
	Taking the normalization
	$$
	\eta_C: \widetilde{C} \longrightarrow C
	$$
	the above morphism gives the morphism
	$$
	\eta^{\ast}(\OO_C ( -M) )\longrightarrow \Theta_{\widetilde{C}}
	$$
	vanishing along the divisor $\eta^{\ast}(Z^C)$. Hence an isomorphism
	$$
	\Theta_{\widetilde{C}} \cong 	\eta^{\ast}(\OO_C ( -M) )\otimes \OO_{\widetilde{C}} (\eta^{\ast}(Z^C)).
	$$
	Comparing the degrees on the both sides gives the equality
	$$
	2-2g_{\widetilde{C}}=-M\cdot C +deg(\eta^{\ast}(Z^C)),
	$$
	where $g_{\widetilde{C}}$ is the genus of $\widetilde{C}$. This and the assumption that $M\cdot C$ is negative imply:
	$$
	g_{\widetilde{C}}=0,\,\,M\cdot C \geq -2, \,\, deg(\eta^{\ast}(Z^C)) \leq 1.
	$$
	The last inequality tells us that $Z^C$ is either empty or a smooth point of $C$. Since we already know that $C$ is smooth everywhere else, we deduce that $C$ is a smooth rational curve and
	$M\cdot C =-2$ or $-1$. The first possibility means that $Z^C=\emptyset$, while the second says: $Z^C$ is a single point.
	This together with the result of M. McQuillan which gives an estimate of such curves $\{C_i\}$:
	$$
	\text{the number of curves in $\{C_i\}$} \leq deg (Z \cap \big(\cup C_i\big)),
	$$
	see \cite{McQ}, Fact I.3.2, imply that there are no curves $C$
	with $Z^C=\emptyset$ and through every point of $Z$ passes at most
	one curve $C$ with $M\cdot C =-1$.	
\end{pf}
\begin{rem}\label{rem:moreintcurves}
	It should be noted that the nonvanishing of $H^1(\OO_X (-M))$ in Proposition \ref{pro:-M} contains more geometry than stated in the proposition: the marked curves $(R,z_R)$ in Proposition \ref{pro:-M}, b), are `connected' by some other integral curves of the foliation. This is to say that there are integral curves $\Gamma$ with
	$M\cdot \Gamma \geq 0$ and intersecting certain number of curves $R$; more precisely, one has
	$$
	\Gamma \cdot \Big(\sum \frac{R}{(-R^2)}\Big) > M\cdot \Gamma,
	$$
	where the sum is taken over all rational curves $R$ in Proposition \ref{pro:-M}. Furthermore, the intersection
	$$
	\Gamma \cdot  {R}=  m_{{}_{\Gamma R}} z_R
	$$
	is supported at $z_R$ for every rational curve $R$. All of this comes from the study of the divisor
	$$
M':=	M-\sum \frac{R}{(-R^2)}
		$$
		and the observation that $H^1 (\OO_X (-M))$ vanishes unless the above divisor fails to be nef. The curves $\Gamma$ above are precisely the ones intersecting $M'$ negatively.
\end{rem}

More informally, the above considerations say that the part
$H^1 (\OO_X (-M))$ of the filtration \eqref{KS-filt-1}, if nonzero, detects a collection of smooth rational curves. Those are integral curves of the foliation  \eqref{L-M-seq-foli} and each comes with a marked point, the unique point of the singularity of the foliation lying on that curve.

All of the above shows that on the classical level of slope-stability
of $\Omega_X$ (or, equivalently, of $\TET$), the cone of ample divisors $Amp(X)$, the K\"ahler moduli attached to $X$, acquires a `wall' structure, this is the collection of faces $Amp^0_M$ in Proposition \ref{pro:stcone}; those come from the collection of objects in the category of complexes of coherent sheaves on $X$, the exact sequences in \eqref{L-M-seq} or their duals \eqref{L-M-seq-foli}, and they carry a certain amount of geometry of $X$.

But the classical level is completely mute about the subcone 
$Amp^s_X$ of stable polarizations, see \eqref{st-cone} and \eqref{st-cone=+}. Nor says it anything about the quotient part of the filtration of $H^1 (\TET)$ in \eqref{KS-filt-1}. All of this changes dramatically when we use Bridgeland stability conditions: we will see
that the `black box' of stable polarizations opens up to produce much finer wall structure on the auxiliary upper half plane bundle over $Amp^s_X$. This in turn comes from changing the heart in the derived category $\DD$ of coherent sheaves and producing many more complexes intrinsically attached to
$\TET$ and the cohomology space $H^1 (\TET)$. Before we go to those topics, see \S\ref{S-BrS} and further on, we need to recall one more classical aspect of
stability of coherent sheaves.

\vspace{0.5cm}
\noindent
{\bf 2.5. Bogomolov (in)stability.} The semistability with respect to a polarization $h$ of a vector bundle or, more generally a torsion free sheaf, imposes a constraint on its Chern invariants. Namely, let $E$ be a torsion free sheaf on $X$ which is $h$-semistable, then the following inequality holds
\begin{equation}\label{B-Gineq}
2rk(E)c_2 (E) \geq (rk(E)-1)c^2_1 (E),
\end{equation}
where $c_1(E)$, $c_2 (E)$ and $rk(E)$ are respectively the Chern classes and the rank of $E$. This is the famous Bogomolov-Gieseker inequality. Inverting the logic one obtains the following.
\begin{lem}\label{lem:Boguns}
If $E$ is a torsion free sheaf on $X$ for which the Bogomolov-Gieseker inequality fails, that is, one has
$$
 (rk(E)-1)c^2_1 (E)> 2rk(E)c_2 (E),
$$
then $E$ is unstable with respect to any $h\in Amp(X)$.
\end{lem}

In view of the above it is convenient to introduce a purely numerical notion of stability.
\begin{defi}\label{def:Bst}
A torsion free sheaf $E$ on $X$ is Bogomolov semistable if the inequality
$$
2rk(E)c_2 (E) \geq (rk(E)-1)c^2_1 (E)
$$
holds. Otherwise $E$ is called Bogomolov unstable.
\end{defi}

For the subsequent discussion we recast the Bogomolov-Gieseker inequality in terms of the Chern character. Recall that for a sheaf $E$ in $\AC$ one defines the Chern character 
$$
ch(E)=ch_0 (E)+ch_1 (E) +ch_2 (E),
$$
where $ch_0 (E)=rk(E)$ is the rank of $E$, $ch_1(E)=c_1 (E)$ is the first Chern class of $E$ and $ch_2 (E)=\HA (c^2_1 (E)-2c_2(E))$. From this it follows that the Bogomolov-Gieseker inequality, see \eqref{B-Gineq}, can be written as follows:
\begin{equation}\label{BG-ch}
ch^2_1(E) -2ch_0 (E)ch_2 (E) \geq 0.
\end{equation}
The expression on the left hand side 
\begin{equation}\label{deltE}
	\Delta(E):=ch^2_1(E) -2ch_0 (E)ch_2 (E)
\end{equation} 
is invariant with respect to tensoring $E$ with a line bundle. That is, one has 
$$
\Delta(E)=\Delta (E\otimes \OO_X(D))
$$
for any line bundle $\OO_X(D)$. This can be checked formally using the multiplicative property of the Chern character
$$
ch(E\otimes F)=ch(E)ch(F)
$$
for a pair of coherent sheaves $E$ et $F$ on $X$. Perhaps, a more conceptual way to see that invariance is to observe that for a locally free sheaf $E$ we have
$$
-\Delta(E)=2ch_0(E)ch_2(E)-ch^2_1 (E)=ch_2 (E\otimes E^{\ast}),
$$
where $E^{\ast} $ is the dual of $E$; the invariance of $ch_2 (E\otimes E^{\ast})$ with respect to tensoring $E$ with a line bundle is now obvious.

Recall that the Chern character extends to a homomorphism 
$$
ch: K_0(\DD) \longrightarrow H^{\ast} (X,\QQ),
$$
where $K_0(\DD)$ is the Grothendieck group of the derived category $\DD$ of coherent sheaves on $X$. The inequality \eqref{BG-ch} now makes sense for objects
of the derived category $\DD$. Thus Definition \ref{def:Bst} can be generalized  to objects of $\DD$, that is, a nonzero object $E$ of $\DD$ is called {\it Bogomolov semistable} if the inequality
\eqref{BG-ch} holds. Otherwise it is called {\it Bogomolov unstable}.

For a polarization $h\in Amp(X)$ the inequality \eqref{BG-ch} gives
\begin{equation}\label{h-BGineq}
(ch_1(E)\cdot h)^2 -2h^2 ch_0 (E) ch_2 (E) \geq 0,
\end{equation}
where we used the Hodge Index inequality 
$$
h^2 ch^2_1 (E) \leq (ch_1(E)\cdot h)^2.
$$
 We will meet quite often the expression on the left hand side of the  inequality \eqref{h-BGineq}. So we set
\begin{equation}\label{Deltah}
	\Delta^h (E):=(ch_1(E)\cdot h)^2 -2h^2 ch_0 (E) ch_2 (E),
\end{equation}
for any object $E$ in $\DD$. This could be viewed as a polarized version of $\Delta(E)$ in \eqref{deltE} and one often refers to it as $h$-{\it discriminant} of $E$, see \cite{Ma-S}. The terminology can be explained as follows.

 We tensor $E$ with a line bundle $\OO_X (-D)$ to obtain
$$
ch(E(-D))=ch(E)ch(\OO_X (-D))=ch(E)(1-D+\HA D^2).
$$
Hence the formulas
\begin{equation}\label{chE(-D)}
	\begin{gathered}
		ch_0 (E(-D))=ch_0 (E), \hspace{0.2cm} ch_1 (E(-D))=ch_1 (E)-ch_0(E)D,\\
		ch_2(E(-D))=ch_2(E)-ch_1(E)\cdot D +\HA ch_0 (E)D^2.
	\end{gathered}
\end{equation}
Observe that the expressions on the right hand side in the above formulas make sense for divisor classes $D\in NS(X)_{\RR}$. With this
in mind
 we replace $D$ in the above formulas by $tD$, where $t$ is a formal variable, to obtain the parameter version of the Chern character. In particular, $ch_2(E(-tD))$ is a quadratic polynomial in $t$
$$
ch_2 (E(-tD))=ch_2(E)-(ch_1(E)\cdot D ) t +\HA ch_0 (E)D^2 t^2.
$$
Its discriminant is $\Delta^D (E)$. This explains the terminology.

In addition, the change of variable $t\mapsto (t+a)$ does not effect the discriminant. Hence the equality
\begin{equation}\label{DelataEa}
	\Delta^D (E(-aD))=\Delta^D(E),
\end{equation}
for $a$, a formal variable.

 The above discussion can be now summarized as follows:
\begin{equation}\label{Bsobj}
	\begin{gathered}
	\text{\it a Bogomolov semistable object $E$ in $\DD$ is subject to the inequality} 
	\\
	\Delta^h (E)\geq 0, \,\forall h\in Amp(X).
\end{gathered}
\end{equation}	

We end the section with a Bogomolov-Gieseker type inequality for
{\it unstable} torsion free sheaves. This goes back to Miyaoka, see \cite{Mi2}.
\begin{lem}\label{lem:BGunst} Let $E$ be a torsion free sheaf on $X$ which is unstable
	with respect to some polarization $h$. We assume that its HN filtration
	$$
	E=E_m \supset E_{m-1} \supset \cdots \supset E_1 \supset E_0 =0
	$$
	with respect to $h$ has the property of positivity for its  quotients, that is,
	$
	Q_i=E_i /E_{i-1}
	$
	for every $i \in [1,m]$ is subject to
	$$
	c_1 (Q_i)\cdot h \geq 0.
	$$	Then the following inequality holds
	$$
	(ch_1(E)\cdot h)^2 -2h^2 ch_0(E_1)ch_2 (E)> \Delta^h (E_1).
	$$
	\end{lem}
\begin{pf}
	From the additivity of the Chern character it follows
	$$
	2ch_2 (E)=2ch_2 \Big(\sum^m_{i=1}Q_i \Big)=\sum^m_{i=1}2ch_2(Q_i) \leq 2ch_2(E_1)+\sum^m_{i=2}\frac{ch^2_1(Q_i)}{ch_0(Q_i)},
	$$
	where the inequality comes from the Bogomolov-Gieseker inequality
	\eqref{BG-ch} for quotient sheaves $Q_i$. The above and the Hodge Index inequality
	$$
	ch^2_1(Q_i)h^2 \leq (ch_1 (Q_i) \cdot h)^2
$$
give
$$
	2ch_2 (E)\leq 2ch_2 (E_1)+\sum^m_{i=2}\frac{(ch_1(Q_i)\cdot h)^2}{h^2 ch_0(Q_i)}.
$$
The assumption of positivity of the quotients tells us
$$
 	 \frac{(ch_1(Q_i)\cdot h)^2}{h^2ch_0(Q_i)}<\frac{ch_1(Q_1)\cdot h}{h^2ch_0(Q_1)} ch_1 (Q_i)\cdot h, 
 $$
 for all $i\geq 2$, and where we used the properties
 $$
 \frac{ch_1(Q_i)\cdot h}{ch_0(Q_i)} <\frac{ch_1(Q_1)\cdot h}{ch_0(Q_1)} =\frac{ch_1(E_1)\cdot h}{ch_0(E_1)} =,\,\forall i\geq 2,
 $$
 and $Q_1=E_1$ of the HN filtration. This gives
 $$
 \begin{gathered}
 	2ch_2 (E)\leq 2ch_2 (E_1)+\sum^m_{i=2}\frac{(ch_1(Q_i)\cdot h)^2}{h^2 ch_0(Q_i)} <2ch_2 (E_1) +\frac{ch_1(E_1)\cdot h}{h^2ch_0(E_1)}\sum^m_{i=2}{ch_1(Q_i)\cdot h}
 	\\
 	 =2ch_2 (E_1) +\frac{ch_1(E_1)\cdot h}{h^2ch_0(E_1)} (ch_1(E)\cdot h-ch_1 (E_1)\cdot h)=-\frac{1}{h^2ch_0(E_1)} \Delta^h(E_1) +\frac{(ch_1(E_1)\cdot h) (ch_1 (E)\cdot h)}{h^2ch_0 (E_1)}.
 	\end{gathered}
 $$ 
 This and the inequality $ch_1 (E_1)\cdot h \leq ch_1 (E)\cdot h$ (this again uses the positivity assumption) imply
 $$
 	2ch_2 (E) < -\frac{1}{h^2ch_0(E_1)} \Delta^h(E_1) +\frac{(ch_1 (E)\cdot h)^2}{h^2ch_0 (E_1)}
 	$$
 	or, equivalently, the inequality
 	$$
 	(ch_1 (E)\cdot h)^2 -2h^2ch_0 (E_1)ch_2 (E) >\Delta^h (E_1).
 	$$ 
	\end{pf}
 
\section{Preliminaries on Bridgeland stability}\label{S-BrS}

{\bf 3.1. Central charge.} As before $NS(X)$ is the N\'eron-Severi group of $X$ and
$NS(X)_{\RR} =NS(X)\otimes_{\ZZ} \RR $. We also remind the reader that $Amp(X)$ is the convex cone in $NS(X)_{\RR}$ generated by ample divisor classes in $NS(X)$.

For every $\omega \in Amp(X)$ and $B\in NS(X)_{\RR}$ one defines a linear complex valued function
$$
Z_{\omega , B} : K_0 (\DD) \longrightarrow \CC
$$
by the formula
\begin{equation}\label{cc}
\begin{gathered}
Z_{\omega , B} (E) =- \int_X ch(-B-i\omega) ch(E)=-\int_X\Big(1-(B+i\omega) +\frac{1}{2}(B+i\omega)^2 \Big)\cdot ch(E) \\
=-\int_X\Big(1-(B+i\omega) +\frac{1}{2}(B+i\omega)^2 \Big)\cdot (ch_0 (E) +ch_1 (E) +ch_2 (E)) \\
=-\Big(ch_2 (E) -ch_1 (E) \cdot B +\frac{1}{2} ch_0 (E) (B^2 -\omega^2) \Big) +
( ch_1 (E) \cdot \omega -ch_0 (E) \omega \cdot B) i,
\end{gathered} 
\end{equation}
where 
$$
ch(E)=ch_0(E) +ch_1(E) + ch_2(E)
$$
 is the Chern character of $E\in \DD$ and $K_0 (\DD)$ denotes the Grothendieck group of $\DD$.
This function is often called a {\it central charge function}; this is rooted in string theory literature, see \cite{Do1}.

The assignment $(\omega, B) \mapsto Z_{\omega , B}$ determines an inclusion
\begin{equation}\label{Z-incl}
Amp(X)\times NS(X)_{\RR} \hookrightarrow Hom_{\ZZ} ( K_0 (\DD),\CC)
\end{equation}
\vspace{0.2cm}
\noindent
{\bf 3.2. The $(\omega, B)$-slope function on $\AC$.} For a pair $(\omega, B) \in Amp(X)\times NS(X)_{\RR}$, one defines $(\omega, B)$-slope for every nonzero sheaf $E$:
\begin{equation}
\mu_{\omega , B} (E) =
\begin{cases}
\frac{c_1 (E) \cdot \omega}{rk (E)} -\omega \cdot B,& \text{if $E$ is not a torsion sheaf}, \\
+\infty,&\text{if $E$ is a torsion sheaf.}
\end{cases}
\end{equation}

A sheaf $E$ in $\AC$ is called $(\omega, B)$-(semi)stable if it is $\omega$-(semi)stable in the sense of Mumford-Takemoto, see Definition \ref{MT-stability}.

\vspace{0.2cm}
\noindent
{\bf 3.3. The tilting of $\AC$.} The slope function $\mu_{\omega , B} $ is used to define a pair of two full subcategories of $\AC$:
\begin{equation*}
\begin{matrix}
\bullet \,\,\TT_{\omega , B} :=&\!\!\!\!\! \text{the full extension closed subcategory of $\AC$ generated by $\mu_{\omega , B} $-stable sheaves $E$}\\
& \text{ of slope $\mu_{\omega , B} (E)>0 $},\\
\bullet \,\,\FF_{\omega , B} :=&\!\!\!\!\! \text{the full extension closed subcategory of $\AC$ generated by $\mu_{\omega , B} $-stable sheaves $E$}\\
&\text{ of slope $\mu_{\omega , B} (E)\leq 0 $.}
\end{matrix}
\end{equation*}

The following is well known, see \cite{HRS} or the overview \cite{Ma-S}.
\begin{pro}\label{pro:torsionpair}
The pair of full subcategories $(\TT_{\omega , B}, \FF_{\omega , B} )$ is a torsion pair of $\AC$, that is,
\begin{enumerate}
\item[(i)]
for every $T$ in $\TT_{\omega , B}$ and every $F$ in $\FF_{\omega , B}$, the group
of morphisms $Hom_{\AC} (T,F) =0$,
\item[(ii)]
every nonzero sheaf $E$ in $\AC$ admits an exact sequence
$$
\xymatrix@R=12pt@C=12pt{
0\ar[r]& T \ar[r]& E  \ar[r]& F \ar[r]&0,
}
$$
where $T$ and $F$ are respectively in $\TT_{\omega , B}$ and $\FF_{\omega , B}$. Such a sequence is unique up to a unique isomorphism. In the sequel we refer to the above exact sequence as the {\rm decomposition of $E$ with respect to a torsion pair $(\TT_{\omega , B}, \FF_{\omega , B})$}.
\end{enumerate}
\end{pro}

A torsion pair  $(\TT_{\omega , B}, \FF_{\omega , B} )$ gives rise to a new $t$-structure on $\DD$. Namely, one defines
\begin{equation}\label{t-struct}
\begin{matrix}
\DD^{\leq 0}_{\omega , B} =&\left\{ A\in \DD | \text{$\HH^0 (A) \in \TT_{\omega , B}$ and   $\HH^k (A) =0, \forall k\geq 1$} \right \}, \\
& \\
\DD^{\geq 0}_{\omega , B} =&\left\{ A\in \DD | \text{$\HH^{-1} (A) \in \FF_{\omega , B}$ and   $\HH^k (A) =0, \forall k\leq -2$} \right \}.
\end{matrix}
\end{equation}
\begin{pro}\label{pro:newt-struct}
The pair of subcategories $(\DD^{\leq 0}_{\omega , B}, \DD^{\geq 0}_{\omega , B} )$ is a $t$-structure of $\DD$ with the heart
$$
\AC_{\omega , B}=\DD^{\leq 0}_{\omega , B} \bigcap \DD^{\geq 0}_{\omega , B} = \left\{ E\in \DD | \text{$\HH^0 (E) \in \TT_{\omega , B}, \, \HH^{-1} (E) \in \FF_{\omega , B}$ and   $\HH^k (E) =0, \forall k\neq -1,0$} \right \}.
$$
 Every object $E\in \AC_{\omega , B}$ admits a unique distinguished triangle in $\DD$
\begin{equation}\label{dist-triang}
 \xymatrix@R=12pt@C=12pt{
E_{-1} \ar[r]& E \ar[r]& E_0  \ar[r]& E_{-1}[1]
}
\end{equation}
where $E_0 \cong \HH^0 (E) \in \TT_{\omega , B}$ and $E_{-1} \cong \HH^{-1} (E)[1] \in \FF_{\omega , B} [1]$. 

In the abelian category $\AC_{\omega , B}$ the above triangle gives rise to the exact sequence
$$
 \xymatrix@R=12pt@C=12pt{
0\ar[r]&\HH^{-1} (E)[1] \ar[r]& E \ar[r]& \HH^0 (E)  \ar[r]&0;
}
$$
the pair $( \FF_{\omega , B} [1], \TT_{\omega , B})$ is a torsion pair of the abelian category $\AC_{\omega , B}$, the above exact sequence is the decomposition of an object $E \in \AC_{\omega , B}$ with respect to the torsion pair $( \FF_{\omega , B} [1], \TT_{\omega , B})$.
\end{pro}

The main point of the above is that the central charge $Z_{\omega , B}$ in \eqref{cc} defines a
Bridgeland stability function on $\AC_{\omega , B}$. This means that we can define a notion of (semi)stable objects in  $\AC_{\omega , B}$ and every nonzero object in $\AC_{\omega , B}$ admits a Harder-Narasimhan filtration (HN for short), see \cite{Ma-S}.
Concretely, this works as follows.

\vspace{0.2cm}
1) The linear function
\begin{equation}\label{stab-func}
Z_{\omega , B} : K_0 (\AC_{\omega , B}) \longrightarrow \CC
\end{equation}
takes the isomorphism class of every
 nonzero object
$E$ in $\AC_{\omega , B}$,  into the completed half-plane
$$
{\bf H}=\{ z=r e^{i\pi\phi} | r>0, \, 0<\phi \leq 1 \}.
$$
 This property allows to define the argument $Arg(Z_{\omega , B} (E))$, for every
nonzero object $E$ in $\AC_{\omega , B}$, as the angle between $Z_{\omega , B} (E)$ and
the positive ray of the real axis in ${\bf H}$.

\vspace{0.2cm}
2) For every nonzero $E \in \AC_{\omega , B}$ one defines the {\it phase} of $E$ with respect to $Z_{\omega , B}$
\begin{equation}\label{phase}
\phi_{\omega , B} (E):= \frac{1}{\pi} Arg(Z_{\omega , B} (E)).
\end{equation}
From 1) it follows that the values of the phase of nonzero objects in $\AC_{\omega , B}$ lie in the interval $(0,1]$.

\vspace{0.2cm}
3) A nonzero object $E$ in $\AC_{\omega , B}$ is called Bridgeland semistable if any nonzero subobject $E'$ of $E$  has the phase smaller or equal to the phase of $E$:
$$
\phi_{\omega , B} (E') \leq \phi_{\omega , B} (E).
$$
 
\vspace{0.2cm}
4) Every nonzero object $E$ in $\AC_{\omega , B}$ admits a Harder-Narasimhan filtration in $\AC_{\omega , B}$, that is, there is a finite step filtration
$$
E=E_n \supset E_{n-1} \supset \cdots \supset E_1 \supset E_0=0
$$
of $E$ such that the quotients $A_i:= E_i /E_{i-1}$ are Bridgeland semistable, for every $i=1,\ldots, n$, and whose phases $\phi_{\omega , B} (A_i)$ form a strictly decreasing function of $i$:
$$
\phi_{\omega , B} (A_1)> \cdots > \phi_{\omega , B} (A_n).
$$

For every $\phi \in (0,1]$ one defines the full subcategory $\AC_{\omega , B}(\phi)$ of $\AC_{\omega , B}$ consisting of zero objects and nonzero objects of phase $\phi$. The following follows easily from the properties above.

\begin{pro}\label{pro:phi-slice}
1) $\AC_{\omega , B}(\phi)$ is the full abelian subcategory of $\AC_{\omega , B}$ whose nonzero objects are precisely the Bridgeland semistable objects of $\AC_{\omega , B}$ of phase $\phi$.

\vspace{0.2cm}
2) A nonzero object $E \in \AC_{\omega , B}(1)$ of phase one has the decomposition
$$
 \xymatrix@R=12pt@C=12pt{
0\ar[r]&\HH^{-1} (E)[1] \ar[r]& E \ar[r]& \HH^0 (E)  \ar[r]& 0
}
$$
with respect to the torsion pair $(\FF_{\omega , B}[1],\TT_{\omega , B})$, where $\HH^0 (E)$ is a torsion sheaf supported on at most $0$-dimensional subscheme of $X$ and $\HH^{-1} (E) \in \FF_{\omega , B}$ is $\mu_{\omega , B}$-semistable
with the slope $\mu_{\omega , B} (\HH^{-1} (E))=0$.
\end{pro}

As a result, the abelian category $\AC_{\omega , B}$ is sliced into the abelian subcategories with slices being $\AC_{\omega , B}(\phi)$, for $\phi \in (0,1]$. This can be extended to the slicing ${\cal P}_{\omega , B}:=\{{\cal P}_{\omega , B} (\phi)\}_{\phi \in \RR}$ of the whole derived category $\DD$ by setting
$$
{\cal P}_{\omega , B} (\phi) =\AC_{\omega , B}(\phi), \,\forall \phi \in (0,1]
$$
and postulating
$$
{\cal P}_{\omega , B} (\phi+1)={\cal P}_{\omega , B} (\phi)[1], \,\forall \phi \in \RR.
$$
 Observe that the central charge $Z_{\omega , B} :K_0 (\DD) \longrightarrow \CC$ defined in \eqref{cc} maps all nonzero objects of the slice ${\cal P}_{\omega , B} (\phi)$ into the ray $\RR_{>0}e^{i\pi\phi}$ of $\CC$. 

The pair $\sigma_{\omega , B}= ({\cal P}_{\omega,B},Z_{\omega , B} )$ is an example of a Bridgeland stability condition on the triangulated category $\DD$. We will call it {\it $(\omega , B)$-stability condition}. 

More generally, a Bridgeland stability condition on a small triangulated category ${\cal C}$ is a pair $\sigma=({\cal P}, Z)$, where 

1) ${\cal P}=\{{\cal P}(\phi)\}_{\phi\in \RR}$ is a slicing of ${\cal C}$, that is,
a family of full subcategories ${\cal P}(\phi)$, $\phi\in \RR$, of ${\cal C}$ subject to the following conditions:

1a) for all $\phi > \phi'$ one has
$$
Hom_{{\cal C}} (A,A')=0, \,\forall A\in {\cal P}(\phi), A'\in {\cal P}(\phi'),
$$

1b) for every nonzero object $E\in  {\cal C}$ there is a finite chain of morphisms
$$
0=E_0 \longrightarrow E_1 \longrightarrow \ldots \longrightarrow E_{n-1} \longrightarrow E_n =E
$$
 in 
 ${\cal C}$ such that for every $i=1,\ldots,n,$ the morphism $E_{i-1} \longrightarrow E_i$ can be completed to a distinguished triangle
$$
 \xymatrix@R=12pt@C=12pt{
E_{i-1} \ar[r]& E_{i} \ar[r]& A_i \ar[r]&E_{i-1} [1]
}
$$
in ${\cal C}$ with the object $A_i \in  {\cal P}(\phi_i)$; in addition the values $\{\phi_i\}$ form a strictly decreasing sequence with respect to $i$:
$$
\phi_1 > \cdots > \phi_{n-1} >\phi_{n};
$$

2) $Z: K_0({\cal C}) \longrightarrow \CC$ is a homomorphism of additive groups, the Grothendieck group $K_0({\cal C})$ of ${\cal C}$ and $\CC$, which  maps all nonzero objects of each slice ${\cal P}(\phi)$ into the ray $\RR_{>0}e^{i\pi\phi}$ of $\CC$. 

\vspace{0.2cm}
The set of Bridgeland stability conditions on the derived category $\DD$ is denoted $Stab_X$. By definition we have the forgetful map
$$
Stab_X \longrightarrow Hom_{\ZZ} (K_0 (\DD),\CC)
$$
sending  a stability condition $\sigma=({\cal P}, Z)$ to $Z \in Hom_{\ZZ} (K_0 (\DD),\CC)$. The Bridgeland's deformation theorem, \cite{Br}, says that  there is a natural topology on $Stab_X$ for which the forgetful map is continuous and locally it is a homeomorphism.

The construction assigning the  Bridgeland stability condition $\sigma_{\omega , B}= ({\cal P}_{\omega,B},Z_{\omega , B} )$ to $(\omega , B)$ in $Amp(X) \times NS_{\RR}$ gives an inclusion 
$$
  Amp(X) \times NS_{\RR} \subset Stab_X  
$$
for which the forgetful map $Stab_X \longrightarrow Hom_{\ZZ} (K_0 (\DD),\CC)$
induces the map
$$
Amp(X) \times NS_{\RR} \longrightarrow Hom_{\ZZ} (K_0 (\DD),\CC)
$$
which sends $(\omega, B) \in Amp(X) \times NS_{\RR}$ to the central charge
$Z_{\omega , B}$ in \eqref{cc}. 

By definition the central charge $Z_{\omega , B}$
factors through the Chern character map
$$
ch:K_0 (\DD) \longrightarrow H^0(X,\ZZ) \oplus NS(X) \oplus \HA H^4(X,\ZZ) \cong \ZZ\oplus NS(X) \oplus \HA \ZZ.
$$
The lattice on the right will be denoted $\Lambda_X$ and it is called Mukai lattice. Thus the forgetful map becomes the map
\begin{equation}\label{cc-fmap}
Amp(X) \times NS_{\RR} \longrightarrow Hom_{\ZZ} (\Lambda_X,\CC)
\end{equation}
sending $(\omega, B) \in Amp(X) \times NS_{\RR}$ to the central charge 
 $Z_{\omega , B}$ and that map is an embedding. 

\section{$(\omega , B)$-stability conditions and $\OO_X$ and $\TET[1]$} \label{s:omB}

{\bf \ref{s:omB}.1.} We fix $(\omega , B)$ in $Amp(X)\times \NSX$ and remind that the Bridgeland stability condition $\sigma_{\omega , B}= ({\cal P}_{\omega,B},Z_{\omega , B} )$ discussed in the previous section is called $(\omega , B)$-stability condition. Among $(\omega , B)$ we will distinguish those for which $\OO_X$ and $\TET[1]$ are objects of $\AC_{\omega , B}$. More precisely, since  $\AC_{\omega , B}$ is equipped with the torsion pair $(\FF_{\omega , B}[1], \TT_{\omega , B})$ we seek $(\omega , B)\in \AMX \times \NSX$ for which $\OO_X \in \TT_{\omega , B}$ and $\TET[1] \in \FF_{\omega , B}[1]$ or, equivalently,  $\TET \in \FF_{\omega , B}$.

\vspace{0.5cm}
\noindent
{\bf \ref{s:omB}.1.1.} As a line bundle, $\OO_X$ is $\omega $-stable, for any $\omega \in Amp(X) $. For $\OO_X $ to be in $\TT_{\omega , B}$ one needs, see {\bf 3.3}, 
$$
\mu_{\omega , B} (\OO_X )=(c_1(\OO_X)-B)\cdot \omega =-B\cdot \omega > 0.
$$
Thus we obtain
\begin{equation}\label{L-cond}
\OO_X \in \TT_{\omega , B} \subset \AC_{\omega , B}, \,\, \forall (\omega , B) \in Amp(X) \times NS_{\RR}, \,\,\rm{where}\,\, B\cdot \omega < 0.
\end{equation}

\noindent
{\bf \ref{s:omB}.1.2.} For $\TET$ we have
$$
\mu_{\omega , B} (\TET )=-\HA K_X \cdot \omega -B\cdot \omega.
$$
This implies that for $\TET$ to be in $\FF_{\omega , B}$ it is sufficient to have:

$\bullet$ $\TET$ is $\omega$-semistable,

$\bullet$ $B\cdot \omega \geq -\HA K_X \cdot \omega$.

Recall: we assume, unless explicitly stated otherwise, that $X$ is a a compact complex surface with $K_X$ ample and having the positive index, that is,
$$
\tau_X=K^2_X -2c_2 (X) >0.
$$

According to Proposition \ref{pro:stcone} (resp. \eqref{st-cone}), the cone $Amp^{\geq 0}(X)$ (resp. $Amp^{s}(X)$) parametrizes the polarizations for which $\TET$ is semistable (resp. stable). Thus we obtain
\begin{equation}\label{OMX-cond}
\TET[1]  \in \FF_{\omega , B}[1] \subset \AC_{\omega , B}, \,\, \forall (\omega , B) \in  Amp^{\geq 0}(X) \times NS_{\RR} \,\,\rm{where}\,\, B\cdot \omega \geq -\HA  K_X\cdot \omega .
\end{equation}
We summarize the above discussion in the following.
\begin{lem}\label{lem:OX-ThXin}
Let $X$ be a compact complex surface with $K_X$ ample and of positive index. Set
$$
Stab(\OO_X,\TET):= \{ (\omega , B)\in Amp^{\geq 0}(X) \times NS_{\RR} |\,\,  -\HA  K_X\cdot \omega  \leq  B\cdot \omega < 0 \}.
$$
This is a set of $(\omega , B)$-stability conditions for which $\OO_X$ and $\TET[1]$ belong to $\AC_{\omega , B}$. More precisely, for every $(\omega , B) \in Stab(\OO_X,\TET) $ one has:

$(i)$ $\OO_X  \in \TT_{\omega , B} \subset \AC_{\omega , B}$,

$(ii)$ $\TET[1] \in \FF_{\omega , B}[1] \subset \AC_{\omega , B}$.
\end{lem}

\vspace{0.5cm}
\noindent
{\bf \ref{s:omB}.2.} As we explained in \S1, we want to study $H^1(\TET)$ using the identification
$$
H^1(\TET)\cong Hom_{\DD} (\OO_X, \TET[1]).
$$
This enables to view cohomology classes in $H^1(\TET)$ as morphisms from $\OO_X$ to  $\TET[1]$ in the derived category $\DD$. We generalize this to linear subspaces. Namely, let $V$ be a nonzero linear subspace of $H^1(\TET)$. Then we have the identifications
\begin{equation}\label{V-ident}
Hom_{\CC} (V,H^1(\TET))\cong H^1(\TET)\otimes V^{\ast} \cong Hom_{\DD} (\OO_X, V^{\ast} \otimes \TET[1])\cong Hom_{\DD} (V\otimes \OO_X, \TET[1]).
\end{equation}
Thus the natural inclusion $i_V: V\hookrightarrow H^1(\TET)$ is identified with a unique nonzero morphism
\begin{equation}\label{F_V}
\Phi_V: V\otimes \OO_X  \longrightarrow  \TET[1]
\end{equation}
in $\DD$.

Choosing our stability condition $(\omega , B)\in Stab(\OO_X,\TET)$, see Lemma \ref{lem:OX-ThXin}, we obtain further identification
$$
Hom_{\CC} (V,H^1(\TET))\cong Hom_{\DD} (V\otimes \OO_X, \TET[1]) \cong Hom_{\AC_{\omega , B}} (V\otimes \OO_X, \TET[1]),
$$
since both objects, $V\otimes \OO_X$ and $\TET[1]$, are now in the full abelian subcategory 
$\AC_{\omega , B}$. To stress the fact that we are working in $\AC_{\omega , B}$, the morphism in \eqref{F_V} will be denoted
\begin{equation}\label{F_VomB}
\Phi^{\omega,B}_V: V\otimes \OO_X  \longrightarrow  \TET[1].
\end{equation}
Since we are in an abelian category, the morphism $\Phi^{\omega,B}_V$ has the kernel $ker (\Phi^{\omega,B}_V)$ (resp., the image, $im (\Phi^{\omega,B}_V)$, and the cokernel,  $coker (\Phi^{\omega,B}_V)$) in $\AC_{\omega , B}$. This gives us two exact sequences in $\AC_{\omega , B}$:
\begin{equation}\label{exseqst}
\xymatrix@R=12pt@C=12pt{
0\ar[r]&ker (\Phi^{\omega,B}_V) \ar[r]&V\otimes \OO_X \ar[r]&im (\Phi^{\omega,B}_V)\ar[r]&0, 
\\
0 \ar[r]&im (\Phi^{\omega,B}_V)\ar[r]& \TET[1] \ar[r]&coker (\Phi^{\omega,B}_V)\ar[r]&0.
}
\end{equation}

\begin{lem}\label{lem:ker-coker}
The kernel $ker (\Phi^{\omega,B}_V)$ of $\Phi^{\omega,B}_V$ is a sheaf in $\TT_{\omega , B}$ and the cokernel $coker (\Phi^{\omega,B}_V)$, if nonzero, is of the form $F_V[1]$, for some torsion free sheaf in $\FF_{\omega , B}$. The two exact sequences in \eqref{exseqst} give rise to two exact complexes in $\AC$
$$
\xymatrix@R=12pt@C=12pt{
0\ar[r]&\HH^{-1} (im(\Phi^{\omega,B}_V)) \ar[r]&ker (\Phi^{\omega,B}_V)\ar[r]&V\otimes \OO_X \ar[r]&\HH^0 (im (\Phi^{\omega,B}_V))\ar[r]&0 
\\
0 \ar[r]&\HH^{-1} (im (\Phi^{\omega,B}_V))\ar[r]& \TET \ar[r]&F_V  \ar[r]&\HH^0 (im (\Phi^{\omega,B}_V)) \ar[r]&0,
}
$$
where $\HH^0 (\bullet)$ is the homological functor associated with the standard $t$-structure on $\DD$, that is, the cohomology of a complex of sheaves in $\AC$. The sheaf $\HH^0 (im (\Phi^{\omega,B}_V)) \in \TT_{\omega , B}$ and $\HH^{-1} (im(\Phi^{\omega,B}_V))\in \FF_{\omega , B}$.
\end{lem}
\begin{pf}
The abelian category $\AC_{\omega , B}$ comes with the torsion pair 
$(\FF_{\omega , B}[1],  \TT_{\omega , B})$. So any nonzero object $E$ in $\AC_{\omega , B}$ admits the exact sequence
$$
\xymatrix@R=12pt@C=12pt{
0 \ar[r]&E^{-1} \ar[r]&E\ar[r]&E^0 \ar[r]&0
}
$$
in $\AC_{\omega , B}$ , where $E^{-1} =\HH^{-1} (E)[1] \in  \FF_{\omega , B}[1]$ and 
$E^0 =\HH^{0} (E)\in  \TT_{\omega , B}$. We may assume that $ker (\Phi^{\omega,B}_V)$ is nonzero, since all zero objects are in $\TT_{\omega , B}$, and apply the above decomposition to it:
\begin{equation}\label{kerPhi-seq}
\xymatrix@R=12pt@C=12pt{
0 \ar[r]&E^{-1} \ar[r]&ker (\Phi^{\omega,B}_V) \ar[r]&E^0 \ar[r]&0.
}
\end{equation}
Combining this with the monomorphism
$ker (\Phi^{\omega,B}_V) \longrightarrow V\otimes \OO_X $, we obtain from
$$
\xymatrix@R=12pt@C=12pt{
  E^{-1} \ar[r]&ker (\Phi^{\omega,B}_V) \ar[r]&V\otimes \OO_X
}
$$
the composite monomorphism
$$
\xymatrix@R=12pt@C=12pt{
  E^{-1}  \ar[r]&V\otimes \OO_X.
}
$$
But $E^{-1} \in \FF_{\omega , B}[1]$ and $V\otimes \OO_X \in \TT_{\omega , B}$ and by definition of the torsion pair $(\FF_{\omega , B}[1],  \TT_{\omega , B})$ no nonzero morphism can exist from $E^{-1}$ to $V\otimes \OO_X$. Hence $E^{-1}  \longrightarrow V\otimes \OO_X $ is zero and since it is a monomorphism, it follows $E^{-1}=0$.
Returning to \eqref{kerPhi-seq}, we deduce $ker (\Phi^{\omega,B}_V) \cong E^0$ and this means that $ker (\Phi^{\omega,B}_V)$ is a sheaf in  $\TT_{\omega , B}$.

Similarly, for $coker (\Phi^{\omega,B}_V)\neq 0$ we have the decomposition exact sequence
$$
\xymatrix@R=12pt@C=12pt{
0 \ar[r]&E^{-1} \ar[r]&coker (\Phi^{\omega,B}_V) \ar[r]&E^0 \ar[r]&0
}
$$
with respect to the torsion pair $(\FF_{\omega , B}[1],  \TT_{\omega , B})$. This gives the 
epimorphisms
$$
\xymatrix@R=12pt@C=12pt{
  \TET[1] \ar@{-{>>}}[r]&coker (\Phi^{\omega,B}_V)  \ar@{-{>>}}[r]&E^0.
  }
  $$
  Hence the composite epimorphism
$$
\xymatrix@R=12pt@C=12pt{
  \TET[1] \ar@{-{>>}}[r]  \ar@{-{>>}}[r]&E^0.
  }
  $$
  But this must be zero, since $\TET[1] \in \FF_{\omega , B}[1]$ and $E^0\in \TT_{\omega , B}$. Since in the composition above the epimorphism
  $$
  \xymatrix@R=12pt@C=12pt{
    \TET[1] \ar@{-{>>}}[r]&coker (\Phi^{\omega,B}_V)
  }
  $$
  is nonzero, it is the
  epimorphism
  $$
\xymatrix@R=12pt@C=12pt{
  coker (\Phi^{\omega,B}_V)  \ar@{-{>>}}[r]&E^0
  }
  $$
  that is zero and therefore $E^0=0$. This proves that  $coker (\Phi^{\omega,B}_V)$ is an object in $\FF_{\omega , B}[1]$. Hence it has the form
  $$
  coker (\Phi^{\omega,B}_V)=F_V[1]
  $$
  for a unique sheaf $F_V \in \FF_{\omega , B}$. By definition, see {\bf 3.3}, all nonzero sheaves in $\FF_{\omega , B}$ are torsion free.

  The two exact complexes of the lemma are obtain by applying the homological functor $\HH^0$ to the exact sequences \eqref{exseqst} or, more precisely, to the distinguished triangles
  $$
 \xymatrix@R=12pt@C=12pt{
ker (\Phi^{\omega,B}_V) \ar[r]&V\otimes \OO_X \ar[r]&im (\Phi^{\omega,B}_V)\ar[r]&ker (\Phi^{\omega,B}_V)[1] 
\\
im (\Phi^{\omega,B}_V)\ar[r]& \TET[1] \ar[r]&coker (\Phi^{\omega,B}_V)\ar[r]&im (\Phi^{\omega,B}_V)[1]
}
$$
in $\DD$.
\end{pf}

\vspace{0.5cm}
\noindent
{\bf \ref{s:omB}.3. (Co)kernel of $\Phi^{\omega,B}_V$.} We want to understand
the kernel $ker(\Phi^{\omega,B}_V)$ (resp. the cokernel $coker(\Phi^{\omega,B}_V)$) of $\Phi^{\omega,B}_V$. For this we specialize the choice of $(\omega,B)$ further. Namely, fix an ample divisor class $H$ in $NS(X)$ for which $\TET$ is {\it stable}; in the notation \eqref{st-cone} we take
$$
H\in Amp^{s}(X) \cap NS(X).
$$

Consider the half-plane
\begin{equation}\label{PiHplane}
  \Pi_H :=\{(cH,bH) \in Amp^{s}(X)\times \NSX |\, c\in \RR_{>0},\, b\in \RR\}
  \end{equation}
  in the space $Stab_X$ of Bridgeland stability conditions. The part of $\Pi_H$ contained in $Stab(\OO_X,\TET)$, see Lemma \ref{lem:OX-ThXin}, is the strip
  \begin{equation}\label{PiHstrip}
    \Pi'_H =\{(cH,bH) \in Amp^{s}(X)\times \NSX |\, c\in \RR_{>0},\, -\frac{K_X \cdot H}{2H^2} \leq b<0 \}.
 \end{equation}    
 The main result of this section is the following.
 \begin{pro}\label{pro:Phi_V-inj}
  The morphism
     $$
     \Phi^{cH,bH}_V= \Phi^{H,bH}_V: V\otimes \OO_X \longrightarrow \TET[1]
     $$
     is independent of $c$. Furthermore, it is injective in $\AC_{cH,bH}=\AC_{H,bH}$ for $(cH,bH)$ in $\Pi'_H$ with $b$ subject to
   $$
   -\frac{1}{H^2 dim(V)} \leq b <0.
     $$
    
     \end{pro}
\begin{pf}
The statement that $\Phi^{cH,bH}_V$ is independent of $c$ follows from the fact that the torsion pair $(\TT_{cH,bH}, \FF_{cH,bH})$ on $\AC$ is independent of $c$, see {\bf 3.3}. Hence neither $\AC_{cH,bH}$, nor the torsion pair $(\FF_{cH,bH}[1],\TT_{cH,bH})$ depend on $c$. Hence $\AC_{cH,bH}=\AC_{H,bH}$ for all $c>0$. In particular,
$ \Phi^{cH,bH}_V= \Phi^{H,bH}_V$.

Next we turn to the statement about the injectivity of $\Phi^{H,bH}_V$. Assume
$ker(\Phi^{H,bH}_V)\neq 0$. We will study the complexes in Lemma \ref{lem:ker-coker}.

If $\HH^0(im( \Phi^{H,bH}_V))=0$, then the second complex in Lemma \ref{lem:ker-coker} reads
\begin{equation}\label{sc}
\xymatrix@R=12pt@C=12pt{
0\ar[r]&\HH^{-1}(im( \Phi^{H,bH}_V))\ar[r]&\Theta_X \ar[r]&F_V \ar[r]&0.
}
\end{equation}
We know $\HH^{-1}(im( \Phi^{H,bH}_V))\neq 0$, since otherwise $im( \Phi^{H,bH}_V)=0$ and hence $\Phi^{H,bH}_V=0$ which is contrary to the definition of $\Phi^{H,bH}_V$. Furthermore $\HH^{-1}(im( \Phi^{H,bH}_V))$ is in $\FF_{H,bH}$. Hence it is torsion free and we proceed our considerations according to the rank of $\HH^{-1}(im( \Phi^{H,bH}_V))$. Since $\HH^{-1}(im( \Phi^{H,bH}_V))$ is a subsheaf of $\Theta_X$ that rank is $1$ or $2$. 

If the rank of $\HH^{-1}(im( \Phi^{H,bH}_V))$ is $2$, then $\HH^{-1}(im( \Phi^{H,bH}_V))=\Theta_X$, since the cokernel $F_V \in \FF_{H,bH}$ and hence must be either zero or a torsion free sheaf. Thus $F_V=0 $ and the first complex in Lemma \ref{lem:ker-coker} becomes
$$
\xymatrix@R=12pt@C=12pt{
0\ar[r]&\Theta_X\ar[r]&ker(\Phi^{H,bH}_V)\ar[r]&V\otimes \OO_X \ar[r]&0.
}
$$
From Lemma \ref{lem:ker-coker} we know that $ker(\Phi^{H,bH}_V)$ is a sheaf in 
$\TT_{H,bH}$. This and the above exact sequence imply
$$
0<\mu_{H,bH} (ker(\Phi^{H,bH}_V))=-\frac{K_X \cdot H}{(dim(V)+2)H^2} -b.
$$
Hence the inequality
$$
b<-\frac{K_X \cdot H}{(dim(V)+2)H^2}.
$$
But then the assumption $b\geq -\frac{1}{H^2dim(V)}$ implies
$$
-\frac{1}{H^2dim(V)} < -\frac{K_X \cdot H}{(dim(V)+2)H^2}
$$
and we obtain
$$
K_X \cdot H <1+\frac{2}{dim(V)}\leq 3.
$$
Hence $K_X \cdot H \leq 2$ and the Hodge Index inequality gives
$$
K^2_X H^2 \leq (K_X \cdot H)^2 \leq 4.
$$
From this it follows $K^2_X \leq 4$. But $X$ is a surface general type of positive index and hence $K^2_X >8\chi(\OO_X)\geq 8.$ Thus we conclude that the rank of $\HH^{-1}(im( \Phi^{H,bH}_V))$ is $1$. Then the exact sequence \eqref{sc}
has the form
$$
\xymatrix@R=12pt@C=12pt{
0\ar[r]&\OO_X(-M)\ar[r]&\Theta_X \ar[r]&{\cal J}_Z (-L) \ar[r]&0.
}
$$
with $\OO_X(-M)$ and ${\cal J}_Z (-L)$ in $\FF_{H,bH}$. 
This also tells us that the first complex in Lemma \ref{lem:ker-coker} is as follows
$$
\xymatrix@R=12pt@C=12pt{
0\ar[r]&\OO_X (-M)\ar[r]&ker(\Phi^{H,bH}_V)\ar[r]&V\otimes \OO_X \ar[r]&0.
}
$$
From this $c_1(ker(\Phi^{H,bH}_V))=-M$. Using again that $ker(\Phi^{H,bH}_V)$ is
a sheaf in $\TT_{H,bH}$, we obtain
$$
0<\mu_{H,bH} (ker(\Phi^{H,bH}_V))=-\frac{M\cdot H}{H^2(dim(V)+1)} -b.
$$
From this point on we continue as in the rank $2$ case. Namely, we have the double inequality
$$
-\frac{1}{H^2dim(V)} \leq b< -\frac{M\cdot H}{H^2(dim(V)+1)}
$$
and this gives
$$
{M\cdot H} <1+\frac{1}{dim(V)}\leq 2.
$$
The fact that $\Theta_X$ is $H$-stable also tells us that
$$
-M\cdot H< -\HA K_X\cdot H.
$$
This and the previous inequality give
$$
\HA K_X\cdot H < 2.
$$
Hence $K_X\cdot H \leq 3$. From the Hodge Index inequality we deduce $K^2_X \leq 9$. Since we have already seen that $K^2_X>8\chi(\OO_X)$, we deduce that we must have 
$K^2_X= 9$ and $\chi(\OO_X)=1$. But then $c_2(X)=12\chi(\OO_X)-K^2_X= 3$ and our surface is subject to $K^2_X=3c_2 (X)$. By a result of Yau, \cite{Y}, the surface $X$ is a compact quotient of the unit ball in $\CC^2$ and for such surfaces we have the vanishing
of $H^1(\TET)$, see \cite{CV}, and this is contrary to what we are assuming in the paper, see \eqref{H1non0}.

At this point we know that the sheaf $\HH^0(im( \Phi^{H,bH}_V))$ is nonzero. The first complex in Lemma \ref{lem:ker-coker} 
\begin{equation}\label{Vcpx}
\xymatrix@R=12pt@C=12pt{
0\ar[r]&\HH^{-1} (im(\Phi^{H,bH}_V)) \ar[r]&ker (\Phi^{H,bH}_V)\ar[r]&V\otimes \OO_X \ar[r]&\HH^0 (im (\Phi^{H,bH}_V))\ar[r]&0
} 
\end{equation}
also tells us that $\HH^0 (im (\Phi^{H,bH}_V))$ is generated by its global sections. Let $T$ (resp. ${\cal E}$) be the torsion (resp. torsion free) part of 
$\HH^0 (im (\Phi^{H,bH}_V))$. Then $\HH^0 (im (\Phi^{H,bH}_V))$ fits into the following exact sequence
$$
\xymatrix@R=12pt@C=12pt{
0\ar[r]&T\ar[r]&\HH^0 (im (\Phi^{H,bH}_V))\ar[r]&{\cal E}\ar[r]&0.
}
$$
From this we have
\begin{equation}\label{D}
D=c_1(\HH^0 (im (\Phi^{H,bH}_V)) =c_1(T) +c_1({\cal E})
\end{equation}
and from the global generation of $\HH^0 (im (\Phi^{H,bH}_V)$ this divisor is effective. Returning with this to the complex \eqref{Vcpx} we have
\begin{equation}\label{c1-1}
c_1 (ker (\Phi^{H,bH}_V))=-D +c_1(\HH^{-1} (im(\Phi^{H,bH}_V)) ).
\end{equation}
Since $ ker (\Phi^{H,bH}_V)\in \TT_{H,bH}$ we obtain
$$
c_1 (ker (\Phi^{H,bH}_V))\cdot H > rk( ker (\Phi^{H,bH}_V))bH^2=
\Big(rk(\HH^{-1} (im(\Phi^{H,bH}_V)))+dim(V)-rk({\cal E})\Big)bH^2.
$$
This and \eqref{c1-1} imply
$$
\begin{gathered}
(rk(\HH^{-1} (im(\Phi^{H,bH}_V)))+dim(V)-rk({\cal E}))bH^2 <-D\cdot H +c_1(\HH^{-1} (im(\Phi^{H,bH}_V)) )\cdot H \\
\leq
-D\cdot H +rk(\HH^{-1} (im(\Phi^{H,bH}_V)))bH^2,
\end{gathered}
$$
where the last inequality comes from $\HH^{-1} (im(\Phi^{H,bH}_V))\in \FF_{H,bH}$. The last inequality simplifies to give
$$
D\cdot H <(dim(V)-rk({\cal E}))H^2 (-b).
$$
This and the assumption $b\geq -\frac{1}{H^2dim(V)}$ imply
$$
D\cdot H <\frac{dim(V)-rk({\cal E})}{dim(V)}\leq 1.
$$
But $D$ is effective, so the above tells us that it must be zero. From the formula \eqref{D} it follows:

1) $c_1(T)=0$, meaning that the torsion part $T$ of $\HH^0 (im (\Phi^{H,bH}_V))$ is supported on at most $0$-dimensional subscheme of $X$,

\vspace{0.2cm}
2) $c_1({\cal E})=0$.

\vspace{0.2cm}
\noindent
We claim that ${\cal E}$ must be nonzero. Indeed, assume ${\cal E}=0$. Then
$\HH^0 (im (\Phi^{H,bH}_V))=T$ is torsion supported on a zero dimensional subscheme and the complex \eqref{Vcpx} has the form
\begin{equation}\label{Vcpx1}
\xymatrix@R=12pt@C=12pt{
0\ar[r]&\HH^{-1} (im(\Phi^{H,bH}_V)) \ar[r]&ker (\Phi^{H,bH}_V)\ar[r]&V\otimes \OO_X \ar[r]&T\ar[r]&0.
}
\end{equation}
Observe: $\HH^{-1} (im(\Phi^{H,bH}_V))$ must be nonzero, since otherwise the second complex in Lemma \ref{lem:ker-coker} becomes
$$
\xymatrix@R=12pt@C=12pt{
  0\ar[r]&\TET\ar[r]&F_V \ar[r]&T\ar[r]&0;
}
$$
since $F_V $ is torsion free and $T$ is supported on a $0$-dimensional subscheme it follows that $\TET \cong F_V$ and thus $T=0$ which is contrary to
$T=\HH^0 (im (\Phi^{H,bH}_V))$ being nonzero.

Once $\HH^{-1} (im(\Phi^{H,bH}_V))$ is nonzero, the second  complex in Lemma \ref{lem:ker-coker} tells us that it is a second syzygy sheaf and those are locally free on a surface, see \cite{OSS}. So $\HH^{-1} (im(\Phi^{H,bH}_V))$ is a line bundle and we write it as $\OO_X (-M)$. Substituting into \eqref{Vcpx1} we obtain
$$
\xymatrix@R=12pt@C=12pt{
0\ar[r]&\OO_X (-M) \ar[r]&ker (\Phi^{H,bH}_V)\ar[r]&{\cal E'}\ar[r]&0,
}
$$
where ${\cal E'}:=ker(V\otimes \OO_X \longrightarrow T)$ is a torsion free sheaf of rank $dim(V)$ and $c_1({\cal E'})=0$. From this it follows
$$
-M\cdot H=c_1(ker (\Phi^{H,bH}_V))\cdot H > (dim(V)+1)H^2 b,
$$
where the last inequality uses $ker (\Phi^{H,bH}_V)\in \TT_{H,bH}$. This together with the assumption $b\geq -\frac{1}{H^2 dim(V)}$ implies
$$
M\cdot H < 1 + \frac{1}{ dim(V)}\leq 2.
$$
Hence $M\cdot H \leq 1$. But $\OO_X (-M)=\HH^{-1} (im(\Phi^{H,bH}_V))$ is also a subsheaf of $\TET$. From $H$-stability of $\TET$ we obtain
$$
-M\cdot H < -\HA K_X \cdot H.
$$
Hence $K_X \cdot H\leq 1$. The Hodge Index inequality implies $K^2_X\leq 1$ which as we have seen before is impossible due to the positivity of the index of $X$.

Thus we have learned that the torsion free part ${\cal E}$ of
$\HH^0 (im (\Phi^{H,bH}_V))$ is nonzero. From \eqref{Vcpx} we obtain a surjective
morphism
$$
V\otimes \OO_X \longrightarrow {\cal E}.
$$
Hence ${\cal E}$ is also globally generated. Since $c_1 ({\cal E})=0$, we deduce that  ${\cal E}=V_0 \otimes \OO_X$ is trivial of rank $rk({\cal E})=dim(V_0)$.
Combining this with  \eqref{Vcpx} we obtain the following diagram:
\begin{equation}\label{Vcpx-diag}
\xymatrix@R=12pt@C=12pt{
  &&&&0\ar[d]&\\
  &&&&T\ar[d]&\\
  0\ar[r]&\HH^{-1} (im(\Phi^{H,bH}_V)) \ar[r]&ker (\Phi^{H,bH}_V)\ar[r]&V\otimes \OO_X \ar[r]\ar[dr]&\HH^0 (im (\Phi^{H,bH}_V))\ar[r]\ar[d]&0\\
  &&&&V_0 \otimes \OO_X\ar[d]&\\
  &&&&0&
}
\end{equation}
If $dim(V_0)=dim(V)$, then the slanted arrow is an isomorphism and we deduce
$$
\HH^{-1} (im(\Phi^{H,bH}_V)) \cong ker (\Phi^{H,bH}_V).
$$
Since $\HH^{-1} (im(\Phi^{H,bH}_V))\in \FF_{H,bH}$ and  $ker (\Phi^{H,bH}_V)\in \TT_{H,bH}$, both sheaves must be zero. Hence $\Phi^{H,bH}_V$ is injective as asserted.

Assume $dim(V_0)<dim(V)$. Then the kernel of the slanted arrow is the trivial
bundle $V_1 \otimes \OO_X$, where $V_1$ is the kernel of the homomorphism induced by the slanted arrow on the spaces of the global sections. So we obtain the following complex
$$
\xymatrix@R=12pt@C=12pt{
  0\ar[r]&\HH^{-1} (im(\Phi^{H,bH}_V)) \ar[r]&ker (\Phi^{H,bH}_V)\ar[r]&V_1\otimes \OO_X \ar[r]&T\ar[r]&0  
}
$$
But this is similar to the one we have encountered in \eqref{Vcpx1}. If $\HH^{-1} (im(\Phi^{H,bH}_V)) \neq 0$, we follow the same argument to obtain that
$\HH^{-1} (im(\Phi^{H,bH}_V))=\OO_X(-M)$ is a line bundle with
$$
M\cdot H < \frac{dim(V_1) +1}{dim(V)} \leq 1.
$$
But remember that $\HH^{-1} (im(\Phi^{H,bH}_V))=\OO_X(-M)$ is also a subsheaf of
$\TET$, see the second complex in Lemma \ref{lem:ker-coker}. The $H$-stability of $\TET$ imposes
$$
-M\cdot H < -\HA K_X \cdot H.
$$
Combining this with the previous inequality gives
$ K_X \cdot H\leq 0$ which contradicts that $X$ is of general type. Thus
$\HH^{-1} (im(\Phi^{H,bH}_V))=0$. With this $im(\Phi^{H,bH}_V )=\HH^{0} (im(\Phi^{H,bH}_V))$ and the second complex in Lemma \ref{lem:ker-coker} becomes
$$
\xymatrix@R=12pt@C=12pt{
  0\ar[r]&\TET\ar[r]&F_V\ar[r]&im(\Phi^{H,bH}_V )\ar[r]&0.  
}
$$
Furthermore, the vertical sequence in \eqref{Vcpx-diag} gives the diagram
$$
\xymatrix@R=12pt@C=12pt{
  &&&0\ar[d]&\\
&&&T\ar[d]&\\
0\ar[r]&\TET\ar[r]&F_V\ar[r]&im(\Phi^{H,bH}_V )\ar[r]\ar[d]&0\\
&&&V_0 \otimes \OO_X\ar[d] &\\
&&&0&
}
$$
This can be completed to the diagram
$$
\xymatrix@R=12pt@C=12pt{
  &&0\ar[d]&0\ar[d]&\\
0\ar[r]&\TET\ar[r]\ar@{=}[d]&F'\ar[d]\ar[r]&T\ar[d]\ar[r]&0\\
0\ar[r]&\TET\ar[r]&F_V\ar[r]\ar[d]&im(\Phi^{H,bH}_V )\ar[r]\ar[d]&0\\
&&V_0 \otimes \OO_X\ar@{=}[r]\ar[d]&V_0 \otimes \OO_X\ar[d] &\\
&&0&0&
}
$$
From the top row of the diagram we deduce $\TET\cong F'$. Hence
$T=0$ and $im(\Phi^{H,bH}_V )=V_0 \otimes \OO_X$, and $ker(\Phi^{H,bH}_V )=V_1 \otimes \OO_X$ . We now have a (split) exact sequence
$$
\xymatrix@R=12pt@C=12pt{
  0\ar[r]&V_1 \otimes \OO_X  \ar[r]&V \otimes \OO_X  \ar[r]&V_0 \otimes \OO_X  \ar[r]&0.
  }
  $$
 
It is a sequence in $\TT_{H,bH}$ and hence in $\AC_{H,bH}$. Applying to it the functor
  $Hom_{\AC_{H,bH}}(\bullet, \TET[1])$ we deduce the exact sequence
  {\small
  $$
  \xymatrix@R=12pt@C=12pt{
    0\ar[r]&Hom_{\AC_{H,bH}}(V_0 \otimes \OO_X,\TET[1]))  \ar[r]&Hom_{\AC_{H,bH}}(V \otimes \OO_X,\TET[1])\ar[r]&Hom_{\AC_{H,bH}}(V_1 \otimes \OO_X,\TET[1])\ar[r]&0.
    }
    $$
    }
    From what has been shown it follows that the morphism $\Phi^{H,bH}_V \in Hom_{\AC_{H,bH}}(V \otimes \OO_X,\TET[1])$ lies in the subspace $Hom_{\AC_{H,bH}}(V_0 \otimes \OO_X,\TET[1])$. But using our identification
    \begin{equation}\label{ident-pf}
    Hom_{\AC_{H,bH}}(V \otimes \OO_X,\TET[1])\cong Hom_{\CC}(V, H^1(\TET))
    \end{equation}
    and similarly, for the spaces
    $Hom_{\AC_{H,bH}}(V_i \otimes \OO_X,\TET[1]))\cong Hom_{\CC}(V_i, H^1(\TET))$, for $i=0,1$, we have the exact sequence
   {\small
  $$
  \xymatrix@R=12pt@C=12pt{
    0\ar[r]&Hom_{\CC}(V_0,H^1(\TET))  \ar[r]&Hom_{\CC}(V,H^1(\TET))\ar[r]&Hom_{\CC}(V_1, H^1(\TET))\ar[r]&0.
    }
    $$
  }
  With these identifications the morphism $\Phi^{H,bH}_V$ must correspond to a linear map 
$$
f:V \longrightarrow H^1(\TET)
$$
 which comes from $Hom_{\CC}(V_0,H^1(\TET))$ or, equivalently, which vanishes on the subspace $V_1 \subset V$. But under the identification \eqref{ident-pf} the morphism $\Phi^{H,bH}_V$ in $Hom_{\AC_{H,bH}}(V \otimes \OO_X,\TET[1])$ corresponds to the inclusion $i_V: V\hookrightarrow H^1(\TET)$. Hence  $V_1 =0$ or, equivalently, $V=V_0$ and this is contrary to the assumption $dim(V_0) < dim(V)$. 
\end{pf}

\vspace{0.5cm}
\noindent
{\bf \ref{s:omB}.4. The stability function $Z_{cH,bH}$ on $\AC_{H,bH}$.} So far we used only the tilted abelian category $\AC_{H,bH}$. We know that for an ample divisor class $H\in Amp^{s}(X)\cap NS(X)$, in the strip
\begin{equation}\label{strip-H}
 \Pi'_H=\{(cH,bH)| \,-\frac{K_X\cdot H}{2H^2}< b<0, c\in \RR_{>0}\}
\end{equation}
of the half-plane $\Pi_H$, see \eqref{PiHplane}, of the Bridgeland stability conditions, the objects $ \OO_X$ and $\TET[1]$ lie in $\TT_{H,bH}$ and $\FF_{H,bH}[1]$ respectively. Furthermore, the morphism
$$
\Phi_V: V\otimes \OO_X \longrightarrow \TET[1]
$$
in the derived category $\DD$, canonically associated to a nonzero linear subspace $V\subset H^1 (\TET)$, in the strip $ \Pi'_H $ becomes a morphism in the abelian category $\AC_{cH,bH}=\AC_{H,bH}$ and to stress this fact it is denoted
\begin{equation}\label{PhiVb}
\Phi^{H,bH}_V: V\otimes \OO_X \longrightarrow \TET[1].
\end{equation}
By Proposition \ref{pro:Phi_V-inj}, in a smaller strip
$$
\Pi^{V}_H=\{(cH,bH)| \,-\frac{1}{H^2dim(V)}\leq b<0, c\in \RR_{>0}\},
$$
that morphism is {\it injective} as a morphism in the abelian category $\AC_{H,bH}$.

We will now bring in the central charge
$$
Z_{cH,bH}: K_0 (\AC_{H,bH}) \longrightarrow \CC,
$$
where $K_0 (\AC_{H,bH})$ is the Grothendieck group of $\AC_{H,bH}$.
 As discussed in {\bf 3.3}, this gives a stability function on $\AC_{H,bH}$. So the domain and codomain objects of the morphism $\Phi^{H,bH}_V$ have their respective phases $\phi_{cH,bH} ( V\otimes \OO_X )$ and $\phi_{cH,bH} ( \TET[1] )$ with respect to $Z_{cH,bH}$. We seek a region in the strip  $\Pi'_H$ subject to the inequality
\begin{equation}\label{phase-ineq}
\phi_{cH,bH} ( V\otimes \OO_X )>\phi_{cH,bH} ( \TET[1] ).
\end{equation}
The motivation is obvious: in such a region either $V\otimes \OO_X$ or $\TET[1]$ (or both) is Bridgeland unstable. One expects then, that the HN filtration of unstable object(s) gives some interesting geometry. Before proceeding with calculations, we make an obvious remark:

\vspace{0.2cm}
\noindent
to determine the phase $\phi_{\omega,B} (E)$ of an object $E\in \AC_{\omega,B}$ it is enough to determine the value of a suitable trigonometric function of the argument $Arg(Z_{\omega,B}(E))=\pi \phi_{\omega,B} (E)$; suitable here means strictly monotone in the interval $(0,\pi]$.

\vspace{0.2cm}
\noindent
We will be using the cotangent function $cot(x)$. Since it is strictly decreasing on the interval $(0,\pi)$ the inequality \eqref{phase-ineq} is recast as follows
\begin{equation}\label{cot-ineq}
\frac{\Re(Z_{\omega,B}( \OO_X ))}{\Im(Z_{\omega,B}( \OO_X ))} =cot(\pi \phi_{\omega,B} ( V\otimes \OO_X ) ) <cot(\pi \phi_{\omega,B} ( \TET[1]) ) =\frac{\Re(Z_{\omega,B}( \TET[1] ))}{\Im(Z_{\omega,B}( \TET[1] ))},
\end{equation}
where $\Re(z)$ (resp. $\Im(z)$) denotes the real (resp. imaginary) part of a complex number $z$.

From the formula \eqref{cc} we obtain
\begin{equation}\label{ccOX}
  \begin{gathered}
    Z_{\omega,B}( \OO_X )=\HA(\omega^2-B^2) -(\omega \cdot B )i, \\
    Z_{\omega,B}( \TET[1]) = -Z_{\omega,B}( \TET)=\Big(\HA\tau_X +B\cdot K_X +B^2-\omega^2 \Big) +(K_X \cdot \omega +2B\cdot \omega)i.
  \end{gathered}
\end{equation}
The inequality \eqref{cot-ineq} now reads
$$
-\frac{\omega^2-B^2}{2\omega \cdot B}<\frac{\HA \tau_X +B\cdot K_X +B^2-\omega^2}{K_X \cdot \omega +2B\cdot \omega}.
$$
We now take $\omega=cH$ and $B=bH$ in the half-plane $\Pi_H$ to obtain
\begin{equation}\label{cb-ineq}
  -\frac{c^2-b^2}{2b}<\frac{\HA\tau_X +bH\cdot K_X +(b^2-c^2)H^2}{ 2b H^2+K_X \cdot H }.
\end{equation}
Specifying further to the values of $(cH,bH)$ in the strip $\Pi'_H$, that is, taking $b$ subject to
$$
-\frac{K_X \cdot H }{2H^2}<b<0,
$$
we obtain the inequality
$$
(b^2-c^2)(2bH^2 +K_X \cdot H )> 2b\Big(\HA\tau_X +bH\cdot K_X +(b^2-c^2)H^2 \Big).
$$
Simplifying gives
$$
b^2+c^2 +\frac{\tau_X}{H\cdot K_X} b <0,
$$
or, equivalently,
\begin{equation}\label{cb-ineq1}
  \Big(b+\frac{\tau_X}{2H\cdot K_X}\Big)^2 +c^2 <\Big(\frac{\tau_X}{2H\cdot K_X}\Big)^2.
\end{equation}
Thus, in the plane with coordinates $(b,c)$, we obtain the following picture:
$$
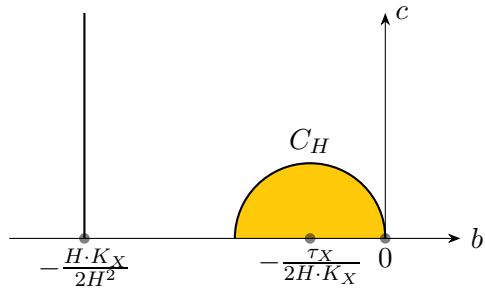
\begin{tikzpicture}[>=Stealth]
	\colorlet{semicircle}{yellow!80!red}
	\draw [->] (0,0)--(1,0)--(4,0)--(5,0)--(6,0);
	\coordinate[label=below:{$-\frac{H\cdot K_X}{2H^2}$}] (A) at (1,0);
	\coordinate [label=below:{$-\frac{\tau_X}{2H\cdot K_X}$}] (B) at (4,0);
	\coordinate [label=below:{$0$}] (C) at (5,0);
	\draw[->](5,0)--(5,3);
	\coordinate[label=right:{$c$}] (D) at (5,3);
	\coordinate[label=right:{$b$}] (E) at (6,0);
	\coordinate[label=above:{$C_H$}] (G) at (4,1);
	\draw[thick] (1,0)--(1,3);
	\draw[thick] (5,0) arc[start angle=0, end angle=180, radius=1];
	\begin{pgfonlayer}{background}
	\fill[semicircle!] (3,0)--(5,0)  arc[start angle=0, end angle=180, radius=1];
	\end{pgfonlayer}
		\foreach \point in {A,B,C}
		\fill[black,opacity=.5](\point)	circle(2pt);
\end{tikzpicture}
$$


\vspace{0.9cm}
\noindent
where the upper half $(b,c)$-plane is the half-plane $\Pi_H$, the vertical strip between the $c$-axis and the line $b=-\frac{H\cdot K_X}{2H^2}$ is the strip $\Pi'_H$ and the interior of the semicircle $C_H$ centered at $b_H=-\frac{\tau_X}{2H\cdot K_X}$ on the $b$-axis and of the radius $r_H=\frac{\tau_X}{2H\cdot K_X}$ is the region
subject to the inequality \eqref{cb-ineq1}. We summarize all of the above in the following.
\begin{pro}\label{pro:semicircH}
  Let $H$ be a divisor class in $Amp^s(X)\cap NS(X)$. Then the strip
  $$
  \Pi'_H=\Big\{(cH,bH)\Big| \,-\frac{K_X\cdot H}{2H^2}< b<0, c\in \RR_{>0} \Big\}
  $$
  is a set of Bridgeland stability conditions for which
  $\OO_X$ and $\TET[1]$ are objects of the abelian category $\AC_{cH,bH}=\AC_{H,bH}$ lying in $\TT_{H,bH}$ and $\FF_{H,bH}[1]$ respectively. Furthermore, for every
  $(b,c)$ inside the semicircle
  $$
  C_H:=\Big\{\Big(b+\frac{\tau_X}{2H\cdot K_X}\Big)^2 +c^2 =\Big(\frac{\tau_X}{2H\cdot K_X}\Big)^2 \,\Bigg| \, c> 0 \Big\} \subset  \Pi'_H,
  $$
  the phases $\phi_{cH,bH} ( \OO_X )$ and $\phi_{cH,bH} ( \TET[1] )$ of $\OO_X$ and $\TET[1] $ with respect to the central charge $Z_{cH,bH}$ are subject to the inequality
  $$
  \phi_{cH,bH} ( \OO_X )> \phi_{cH,bH} ( \TET[1] ).
  $$
  In particular, either $\OO_X$ or $\TET[1]$ is Bridgeland unstable, for every
  $(b,c)$ inside the semicircle $C_H$.
  \end{pro}
\begin{pf}
The statement about  $\OO_X$ and $\TET[1]$ lying in the subcategories $\TT_{H,bH}$ and $\FF_{H,bH}[1]$ respectively, for all $b$ as in the strip  $\Pi'_H$, is a consequence of Lemma \ref{lem:OX-ThXin}. The inequality between the phases of $\OO_X$ and $\TET[1]$ holding inside of the semicircle is the calculation leading to the inequality \eqref{cb-ineq1}. The fact that this inequality implies Bridgeland instability of either $\OO_X$ or $\TET[1]$ follows from the property of a slicing of a triangulated category: for two slices ${\cal P}(\phi)$, ${\cal P}(\phi')$ with $\phi>\phi'$, there are no nonzero morphisms from objects in ${\cal P}(\phi)$ to objects in ${\cal P}(\phi')$. But we know that 
$$
\Phi^{H,bH}_V: V\otimes\OO_X \longrightarrow \TET[1]
$$
is a nonzero morphism in $\AC_{H,bH}$ for all $b$ in the strip  $\Pi'_H$.

The only statement perhaps which one needs to check is that the semicircle 
$C_H$ is contained in the strip $\Pi'_H$. This comes down to verifying the inequality
$$
\frac{K_X\cdot H}{2H^2} > \frac{\tau_X}{H\cdot K_X}
$$
or, equivalently,
\begin{equation}\label{ineq-tau}
\frac{(K_X\cdot H)^2}{2H^2} > \tau_X.
\end{equation}
For this use the Hodge Index inequality $(K_X\cdot H)^2 \geq K^2_X H^2 $ to obtain
$$
\frac{(K_X\cdot H)^2}{2H^2} \geq \HA K^2_X.
$$
On the other hand $\tau_X =K^2_X -2c_2(X)=(1-2\alpha_X)K^2_X$, where
$\alpha_X=\frac{c_2(X)}{K^2_X}$ is the ratio of the Chern numbers, see \eqref{ratio-Chernnum}. Comparing with the right hand side of the inequality \eqref{ineq-tau} gives
$$
\HA K^2_X -\tau_X=(\HA-(1-2\alpha_X))K^2_X=(2\alpha_X-\HA)\alpha_X=2(\alpha_X-\frac{1}{4})K^2_X.
$$
The Bogomolov-Miyaoka-Yau inequality $3c_2(X)\geq K^2_X$ is equivalent to
$$
\alpha_X \geq \frac{1}{3}.
$$
Hence $2(\alpha_X-\frac{1}{4})K^2_X >0$. This proves the inequality \eqref{ineq-tau} and the assertion that the semicircle $C_H$ lies in the strip $\Pi'_H$.
\end{pf}
\begin{rem}\label{rem:ineq-tau=Hdisc}
	The inequality \eqref{ineq-tau} is equivalent to
	$$
	(K_X \cdot H)^2 -2H^2\tau_X >0.
	$$
	The expression on the left can be rewritten in the form
	$$
	(ch_1 (\TET[1])\cdot H )^2 -2H^2ch_0(\TET[1])ch_2 (\TET[1]);
	$$
	we remind the reader that this is $\Delta^H (\TET[1])$, the $H$-discriminant of $\TET[1]$, see \eqref{Deltah}. So the inequality
	\eqref{ineq-tau} can be restated as the positivity of the $H$-discriminant of $\TET[1]$:
	\begin{equation}\label{Hdisc-positive}
		\Delta^H (\TET[1]) > 0.
	\end{equation} 
\end{rem}

\begin{cor}\label{cor:TET-Brunst}
  Let $H$ and $C_H$ be as in Proposition \ref{pro:semicircH}. Then the object
  $\TET[1]$ is Bridgeland unstable in $\AC_{cH,bH} =\AC_{H,bH}$, for all $(b,c)$
  inside the semicircle $C_H$ where the morphism
  $$
  \Phi^{H,bH}_V : V\otimes \OO_X \longrightarrow \TET[1]
  $$
  in $\AC_{H,bH}$ is injective.
\end{cor}
\begin{pf}
  From Proposition \ref{pro:semicircH} we know that the phase $\phi_{cH,bH} (\OO_X)$ is bigger than the phase of $\phi_{cH,bH} (\TET[1])$ for all $(b,c)$ inside the semicircle $C_H$. But $V\otimes \OO_X$ has the same phase as $\OO_X$. So we have
  $$
  \phi_{cH,bH} (V\otimes\OO_X) > \phi_{cH,bH} (\TET[1]),
  $$
  for all $(b,c)$ inside the semicircle $C_H$. Assume $\TET[1]$ is Bridgeland semistable for some $(b,c)$ inside $C_H$. Then by Proposition \ref{pro:semicircH} the object $V\otimes\OO_X$ is unstable in $\AC_{cH,bH}$. Then the HN filtration of  $V\otimes\OO_X$ in $\AC_{cH,bH}$ gives a semistable subobject $E$ of $V\otimes\OO_X$ whose phase $\phi_{cH,bH} (E)$ is bigger than $\phi_{cH,bH} (V\otimes\OO_X)$. Composing the inclusion $E \longrightarrow V\otimes\OO_X$ with $\Phi^{H,bH}_V$
  gives a morphism
  $$
  E \longrightarrow \TET[1]
  $$
  between semistable objects in $\AC_{cH,bH}$ subject to $\phi_{cH,bH} (E) > \phi_{cH,bH} (\TET[1]) $. Hence that morphism must be zero. But this means that $\Phi^{H,bH}_V$ is not injective. Thus $\Phi^{H,bH}_V =\Phi^{cH,bH}_V$ fails to be injective for any $(b,c)$ inside $C_H$, where $\TET[1]$ is semistable. 
\end{pf}

From the above it follows that the region of Bridgeland instability of $\TET[1]$
is contained in the region of the injectivity of $\Phi^{H,bH}_V $. The latter was
partly explored in Proposition \ref{pro:Phi_V-inj}. We make it somewhat more precise here.

Set
\begin{equation}\label{inf-inj}
  i_H:=inf \Big\{ b\in \Big(-\frac{K_X \cdot H}{2H^2}, 0 \Big) |\, \Phi^{H,bH}_V \,\text{is injective} \Big\},
  \end{equation}
  the lower bound of values of $b$ in the strip $\Pi'_H$ for which the morphism $\Phi^{H,bH}_V$ is injective.
  \begin{lem}\label{lem:iH}
    1) $\displaystyle{-\frac{K_X \cdot H}{(dim(V)+2)H^2} \leq i_H \leq - \frac{1}{H^2 dim(V)}.}$

    2)  $\Phi^{H,bH}_V $ injective for all $b\in (-i_H,0)$.
  \end{lem}
  \begin{pf}
    The upper bound for $i_H$ in first assertion follows from Proposition \ref{pro:Phi_V-inj}. To see the lower bound, observe that the injectivity of $\Phi^{H,bH}_V$ gives the exact sequence
    $$
    \xymatrix{
    	0\ar[r]& V\otimes \OO_X \ar[r]^(.59){\Phi^{H,bH}_V}& \TET[1]\ar[r]& coker(\Phi^{H,bH}_V)\ar[r]&0    
}
$$ 
in $\AC_{H,b H}$. The cokernel must be nonzero and hence, according to Lemma \ref{lem:ker-coker}, it has the form $F_V[1]$ for a torsion free sheaf
$F_V \in \FF_{H,bH}$. So the above exact sequence takes the form
 $$
\xymatrix{
	0\ar[r]& V\otimes \OO_X \ar[r]^(.59){\Phi^{H,bH}_V}& \TET[1]\ar[r]& F_V[1]\ar[r]&0    
}
$$
Applying the homological functor $\HH^0$ gives
\begin{equation}\label{extVseq}
\xymatrix{
	0\ar[r]&\TET\ar[r]&F_V \ar[r]& V\otimes \OO_X \ar[r]&0.    
}	
\end{equation}
The condition for the sheaf $F_V$ to be in $\FF_{H,bH}$ reads
$$
-K_X \cdot H=c_1 (F_V)\cdot H \leq rk(F_V)H^2 b =(dim(V)+2)H^2 b.
$$
Hence we obtain
$$
b\geq -\frac{K_X \cdot H}{(dim(V)+2)H^2} 
$$
for all $b$ for which $\Phi^{H,bH}_V$ is injective in $\AC_{H,b H}$.
From this it follows
$$
i_H =inf \Big\{ b\in \Big(-\frac{K_X \cdot H}{2H^2}, 0 \Big) |\, \Phi^{H,bH}_V \,\text{is injective} \Big\} \geq -\frac{K_X \cdot H}{(dim(V)+2)H^2}.
$$

For the second assertion it is enough to show
    {\small
    $$
    \text{\it if $\Phi^{H,bH}_V $ is injective for $b\in \Big(-\frac{K_X \cdot H}{2H^2}, 0 \Big)$, then $\Phi^{H,b'H}_V $ is injective for all
      $b'\geq b$ in $\Big(-\frac{K_X \cdot H}{2H^2}, 0 \Big)$.}
    $$
    }
    But this follows from

    $(a)$ $ker(\Phi^{H,b'H}_V )$ is a sheaf in $\TT_{H,b'H}$, see Lemma \ref{lem:ker-coker},

    \vspace{0.2cm}

    $(b)$ $\TT_{H,b'H} \subset \TT_{H,bH}$, for all $b'>b$; this is immediate from the the definition of those subcategories, see {\bf 3.3}. 
  \end{pf}
  \begin{cor}\label{cor:ThXBrunst}
    Let $H$ and $C_H$ be as in Proposition \ref{pro:semicircH} an set
    $$
    \beta_H:=max \Big\{-\frac{\tau_X}{K_X \cdot H}, i_H\Big\},
    $$
    where $i_H$ is as in \eqref{inf-inj}. Then $\TET[1]$ is Bridgeland unstable
    in $\AC_{H,bH}=\AC_{cH,bH}$, for all $(b,c)$ inside the semicircle $C_H$ and
    $b \in (\beta_H,0)$.
  \end{cor} 
  \begin{pf}
    For all $(b,c)$ inside the semicircle $C_H$ and $b \in (\beta_H,0)$
    the morphism
  $$
  \Phi^{H,bH}_V : V\otimes \OO_X \longrightarrow \TET[1]
  $$
  in $\AC_{H,bH}$ is injective. So Corollary \ref{cor:TET-Brunst} applies and we deduce that the object
  $\TET[1]$ is Bridgeland unstable
    in $\AC_{H,bH}=\AC_{cH,bH}$ for all $(b,c)$ inside the semicircle $C_H$ and $b \in (\beta_H,0)$.
  \end{pf}

Let us summarize: 

\vspace{0.2cm}
1) for every ample divisor class $H\in Amp^{s} (X) \cap NS(X)$ we have the region
\begin{equation}\label{BH-unst}
B^{un}_H :=\{(cH,bH) | \text{$(b,c)$ is inside the semicircle $C_H$ and $b\in (\beta_H,0)$}\},
\end{equation}
where $C_H$ is as in Proposition \ref{pro:semicircH} and $\beta_H$ is defined in Corollary \ref{cor:ThXBrunst},

\vspace{0.2cm}
2) for every $(cH,bH)$-stability condition in $B^{un}_H$ the morphism
$$
\Phi^{cH,bH}_V=\Phi^{H,bH}_V : V\otimes \OO_X \longrightarrow \TET[1]
$$
intrinsically associated with a nonzero linear subspace $V\subset H^1(\TET)$ lies in the abelian subcategory $\AC_{cH,bH}=\AC_{H,bH}$ and that morphism is a monomorphism in $\AC_{H,bH}$,

\vspace{0.2cm}
3) the object $\TET[1]$ is Bridgeland unstable with respect to all 
$(cH,bH)$ in $B^{un}_H$. 
 
\vspace{0.2cm}
The usefulness of the above depends on the understanding of the Harder-Narasimhan
(HN) filtration of $\TET[1]$. This will be our next task.

\section{The HN filtration of $\TET[1]$}\label{s:HNfilt}

We fix a stability condition $(cH,bH)$ in $B^{un}_H$, the set defined in \eqref{BH-unst}. From the previous section we know that  the morphism
$$
\Phi^{cH,bH}_V=\Phi^{H,bH}_V : V\otimes \OO_X \longrightarrow \TET[1]
$$
is a monomorphism in the abelian category $\AC_{cH,bH}=\AC_{H,bH}$ and that the object $\TET[1]$ is Bridgeland unstable in $\AC_{H,bH}$. We consider
\begin{equation}\label{HN-Th}
 \TET[1]=E^{c,b}_n \hookleftarrow E^{c,b}_{n-1} \hookleftarrow \cdots \hookleftarrow  E^{c,b}_1 \hookleftarrow E^{c,b}_0 =0
\end{equation}
the HN filtration of $ \TET[1] $ in $\AC_{H,bH}$ with respect to the central charge $Z_{cH,bH}$, see \eqref{cc} for the definition. The inclusion
$E^{c,b}_{i-1} \hookrightarrow E^{c,b}_i$, for every $i=1,\ldots,n,$ defines the quotient object 
\begin{equation}\label{quot-ob}
A^{c,b}_i:=E^{c,b}_i /E^{c,b}_{i-1}.
\end{equation} 
Its phase with respect to $Z_{cH,bH}$ will be denoted $\phi^{c,b}_i$. These values
form a decreasing sequence
\begin{equation}\label{phases}
  \phi^{+}_{c,b} (\TET[1]):=\phi^{c,b}_1 > \cdots >\phi^{c,b}_n:= \phi^{-}_{c,b} (\TET[1]).
\end{equation}
We recall that by definition of $(\AC_{H,bH},Z_{cH,bH})$ those values lie in the interval $(0,1]$.

\vspace{0.5cm}
\noindent
{\bf \ref{s:HNfilt}.1.}
 We now begin investigating the properties of the filtration \eqref{HN-Th}.
\begin{lem}\label{lem:less1}
  $\phi^{+}_{c,b} (\TET[1])<1.$
\end{lem}
\begin{pf}
  Assume $\phi^{+}_{c,b} (\TET[1])=1$. Consider the monomorphism
  $$
  \xymatrix@R=12pt@C=12pt{
    0\ar[r]&E^{c,b}_1 \ar[r]&\TET[1]
  }
  $$
  in $\AC_{H,bH}$. Its cokernel is of the form $F_1[1]$ for a torsion free sheaf
  $F_1 \in \FF_{H,bH}$. So we have an exact sequence
   $$
  \xymatrix@R=12pt@C=12pt{
    0\ar[r]&E^{c,b}_1 \ar[r]&\TET[1]\ar[r]&F_1[1]\ar[r]&0
  }
  $$
  in $\AC_{H,bH}$, where $E^{c,b}_1 =A^{c,b}_1$ is a subobject of phase $1$.
  Applying the homological functor $\HH^0$ gives the exact complex of sheaves
  $$
  \xymatrix@R=12pt@C=12pt{
    0\ar[r]&\HH^{-1}(E^{c,b}_1 ) \ar[r]&\TET \ar[r]&F_1\ar[r]&\HH^{0}(E^{c,b}_1 \ar[r])&0.
  }
  $$
  Since $E^{c,b}_1 $ has phase $1$ we know from Proposition \ref{pro:phi-slice} that  $\HH^{0}(E^{c,b}_1)$, if nonzero, is a torsion sheaf supported on a $0$-dimensional subscheme of $X$, and $\HH^{-1}(E^{c,b}_1 )$, if nonzero, is in $\FF_{H,bH}$ and $c_1 (\HH^{-1}(E^{c,b}_1 )\cdot H=b \,rk(\HH^{-1}(E^{c,b}_1 ))H^2$.
  
  Assume $\HH^{-1}(E^{c,b}_1 )\neq 0$. Then it is a rank $1$ subsheaf of $\TET$ and, from the complex above, it is a second syzygy sheaf on $X$. Hence $\HH^{-1}(E^{c,b}_1 )$ is a line bundle and we denote it by $\OO_X (-L_1)$ with $L_1$ subject to
  $$
  -L_1 \cdot H=bH^2.
  $$
  But $\TET$ is $H$-stable, so we obtain
  $$
  -\HA K_X \cdot H >-L_1 \cdot H=bH^2
  $$
  and conclude that $b<-\frac{ K_X \cdot H}{2H^2}$. This is outside of the interval $(-\frac{ K_X \cdot H}{2H^2},0)$, where the values of $b$ belong. So $\HH^{-1}(E^{c,b}_1 )$ must be zero and the complex above becomes
   $$
  \xymatrix@R=12pt@C=12pt{
    0\ar[r]&\TET \ar[r]&F_1\ar[r]&E^{c,b}_1 \ar[r]&0,
  }
  $$
  where $E^{c,b}_1 $ is a sheaf supported on a $0$-dimensional subscheme of $X$.
   Since $F_1$ is torsion free, the exact sequence above implies that $\TET\cong F_1$. Hence $E^{c,b}_1=0$, a contradiction.
 \end{pf}

 Let $E=E^{c,b}_i$, for some $i=1,\ldots,n-1$, be one of the proper destabilizing subobjects of of $\TET[1]$ in the HN filtration \eqref{HN-Th}. The monomorphism
 $$
 \xymatrix@R=12pt@C=12pt{
   0\ar[r]&E \ar[r]&\TET[1]
 }
 $$
 gives rise to an exact sequence
 \begin{equation}\label{seqE}
 \xymatrix@R=12pt@C=12pt{
   0\ar[r]&E \ar[r]&\TET[1]\ar[r]&Q_E \ar[r]&0
 }
\end{equation}
 in $\AC_{H,bH}$, where $Q_E$, being a quotient of an object in $\FF_{H,bH}[1]$, must be also in $\FF_{H,bH}[1]$, so $Q_E=F_E[1]$ for some torsion free sheaf $F_E \in \FF_{H,bH}$.

 \begin{lem}\label{lem:FE-locfree}
   $F_E$ is a locally free sheaf.
 \end{lem}
 \begin{pf}
   We know that $F_E$ is torsion free. Let $F^{\ast \ast}_E$ be the double dual of $F_E$. This is a locally free sheaf on $X$ and we have an exact sequence
   $$
  \xymatrix@R=12pt@C=12pt{
   0\ar[r]&F_E \ar[r]&F^{\ast \ast}_E\ar[r]&{\cal S} \ar[r]&0,
 }
 $$
 in $\AC$, where ${\cal S} $ is a torsion sheaf supported on at most $0$-dimensional subscheme of $X$. In the derived category $\DD$ that sequence gives a distinguished triangle
 $$
 \xymatrix@R=12pt@C=12pt{
   F_E \ar[r]&F^{\ast \ast}_E\ar[r]&{\cal S}\ar[r]&F_E[1].
   }
 $$
 The axiom of turning triangles tells us that the triangle
 $$
 \xymatrix@R=12pt@C=12pt{
   {\cal S}\ar[r]&F_E[1]\ar[r]&F^{\ast \ast}_E[1] \ar[r]&{\cal S}[1]
   }
 $$
 is also distinguished. But the three terms starting from the left are in the abelian subcategory $\AC_{H,bH}$, therefore they form an exact sequence
 $$
 \xymatrix@R=12pt@C=12pt{
   0\ar[r]&{\cal S}\ar[r]&F_E[1]\ar[r]&F^{\ast \ast}_E[1] \ar[r]&0
   }
 $$
 in $\AC_{H,bH}$. If ${\cal S}$ is nonzero, it is a semistable object in $\AC_{H,bH}$ of phase $1$ admitting a monomorphism ${\cal S}\longrightarrow F_E[1]$. But $F_E[1] =E^{c,b}_n /E=E^{c,b}_n /E^{c,b}_i$ has semistable factors $A^{c,b}_l$, where $l=i+1, \ldots,n$. Hence the phases $\phi^{c,b}_l < \phi^{+}_{c,b} (\TET[1])$, for all $l\geq i+1$. From Lemma \ref{lem:less1} we know that $\phi^{+}_{c,b} (\TET[1])<1$. So all semistable factors of $F_E[1]$ have phase less than $1$. As a consequence we can not have nonzero morphisms ${\cal S}\longrightarrow F_E[1]$.
 This proves that ${\cal S}=0$ and hence  $F_E =F^{\ast \ast}_E$ is locally free.
\end{pf}

We go back to the exact sequence \eqref{seqE} and write it in the form
$$
\xymatrix@R=12pt@C=12pt{
  0\ar[r]&E \ar[r]&\TET[1]\ar[r]&F_E[1] \ar[r]&0,
  }
  $$
  where $F_E$ is in $\FF_{H,bH}$ and, according to Lemma \ref{lem:FE-locfree}, locally free. Applying the homological functor $\HH^0$ gives the exact complex
  \begin{equation}\label{E-FEcpx}
  \xymatrix@R=12pt@C=12pt{
  0\ar[r]&\HH^{-1}(E) \ar[r]&\TET\ar[r]&F_E \ar[r]&\HH^{0}(E)\ar[r]&0,
  }
\end{equation}
where $\HH^{0}(E)\in \TT_{H,bH}$ and $\HH^{-1}(E) \in \FF_{H,bH}$.
\begin{lem}\label{lem:H-1-E}
1) If $\HH^{-1}(E^{c,b}_i) =0$ for some $i\leq n-1$, then $\HH^{-1}(E^{c,b}_j)=0$, for all $0<j\leq i$.

2) Assume $i\leq n-2$ and $\HH^{-1}(E^{c,b}_i) \neq 0$, then $\HH^{-1}(E^{c,b}_i)$ is a line bundle in $\FF_{H,bH}$ and 
$$
\HH^{-1}(E^{c,b}_k) \cong \HH^{-1}(E^{c,b}_i),
$$
 for all $k$ subject to $i\leq k <n$.
\end{lem}
\begin{pf}
For each nonzero proper object $E^{c,b}_l$ of the HN filtration \eqref{HN-Th} we have the exact sequence
$$
\xymatrix@R=12pt@C=12pt{
  0\ar[r]&\HH^{-1}(E^{c,b}_{l})[1]\ar[r]& E^{c,b}_{l} \ar[r] &\HH^{0}(E^{c,b}_{l}) \ar[r]&0,
}
$$
the decomposition of $E^{c,b}_l$ with respect to the torsion pair
$(\FF_{H,bH}[1],\TT_{H,bH})$. This is functorial. In particular,
 the inclusion $E^{c,b}_{j} \hookrightarrow E^{c,b}_i$, for  $0<j\leq i$, gives rise to the morphism
of the exact sequences
$$
\xymatrix@R=12pt@C=12pt{
  0\ar[r]&\HH^{-1}(E^{c,b}_{j})[1] \ar[d] \ar[r]& E^{c,b}_{j} \ar[r] \ar@{^{(}->}[d]&\HH^{0}(E^{c,b}_{j}) \ar[r]\ar[d]&0\\
0\ar[r]&\HH^{-1}(E^{c,b}_{i})[1] \ar[r]& E^{c,b}_{i} \ar[r]&\HH^{0}(E^{c,b}_{j}) \ar[r]&0.
}
$$
From this it follows: if $\HH^{-1}(E^{c,b}_{i})=0$, then so is $\HH^{-1}(E^{c,b}_{j})$, for all $0<j\leq i$. This proves part 1) of the lemma.

Turning to the part 2), we deduce from the complex \eqref{E-FEcpx} that $\HH^{-1}(E^{c,b}_{i})$, whenever nonzero, is a subsheaf of rank one of $\TET$. Furthermore, from that complex it follows that it is a second syzygy sheaf and those are locally free on a surface. So if $\HH^{-1}(E^{c,b}_{i})$ is nonzero, it is a line bundle. From part 1) of the lemma it follows that $\HH^{-1}(E^{c,b}_k)$ is nonzero for all $k$ subject to $i\leq k<n$. Hence $\HH^{-1}(E^{c,b}_k)$ is a line bundle as well, for all $k$ subject to $i\leq k<n$. To compare $\HH^{-1}(E^{c,b}_k)$ with $\HH^{-1}(E^{c,b}_{i})$ we consider two exact sequences
$$
\xymatrix@R=12pt@C=12pt{
  0\ar[r]&E^{c,b}_{i}\ar[r]&\TET[1] \ar[r] &F_{E^{c,b}_{i}}[1]\ar[r]&0,\\
  0\ar[r]&E^{c,b}_{k}\ar[r]&\TET[1] \ar[r] &F_{E^{c,b}_{k}}[1]\ar[r]&0. 
}
$$
For all $k\geq i$ we have the monomorphism $E^{c,b}_{i} \hookrightarrow E^{c,b}_{k}$ which gives rise to a morphism of the sequences
$$
\xymatrix@R=12pt@C=12pt{
  0\ar[r]&E^{c,b}_{i}\ar[r]\ar[d]&\TET[1] \ar[r] \ar@{=}[d]&F_{E^{c,b}_{i}}[1]\ar[r]\ar[d]&0\\
  0\ar[r]&E^{c,b}_{k}\ar[r]&\TET[1] \ar[r] &F_{E^{c,b}_{k}}[1]\ar[r]&0 
}
$$
Applying the homological functor $\HH^0$ gives the morphism of two complexes in $\AC$:
$$
\xymatrix@R=12pt@C=12pt{
  0\ar[r]&\HH^{-1}(E^{c,b}_{i})\ar[r]\ar[d]&\TET \ar[r] \ar@{=}[d]&F_{E^{c,b}_{i}}\ar[r]\ar[d]&\HH^{0}(E^{c,b}_{i})\ar[r]\ar[d]&0\\
  0\ar[r]&\HH^{-1}(E^{c,b}_{k})\ar[r]&\TET \ar[r] &F_{E^{c,b}_{k}}\ar[r]&\HH^{0}(E^{c,b}_{k})\ar[r]&0
}
$$
The resulting morphism
\begin{equation}\label{H-1ik}
  \HH^{-1}(E^{c,b}_{i})\longrightarrow \HH^{-1}(E^{c,b}_{k})
  \end{equation}
is injective. We have just proved that both sheaves are line bundles, so the cokernel of that morphism must be a torsion sheaf. The above diagram tells us that the cokernel injects into the kernel of the arrow
$$
F_{E^{c,b}_{i}} \longrightarrow F_{E^{c,b}_{k}}.
$$
But the cokernel of the morphism in \eqref{H-1ik} is a torsion sheaf
while $F_{E^{c,b}_{i}}$ is torsion free. Hence
the cokernel of \eqref{H-1ik} is zero and we have an isomorphism
$ \HH^{-1}(E^{c,b}_{i}) \cong \HH^{-1}(E^{c,b}_{k})$.
\end{pf}

\begin{pro}\label{pro:i0}
 In the HN filtration \eqref{HN-Th} either all objects 
$E^{c,b}_i$, for $i\in [1,n-1]$, are sheaves in $\TT_{H,bH}$ or the set 
$$
\{ i\in [1,n-1] | \HH^{-1}(E^{c,b}_{i})\neq 0\}
$$
is nonempty and we set
$$
i_0:=min\{ i\in [1,n-1] | \HH^{-1}(E^{c,b}_{i})\neq 0\}
$$
and have  that $\HH^{-1}(E^{c,b}_{i_0})$ is a line bundle in $\FF_{H,bH}$ and 
$\HH^{-1}(E^{c,b}_{i}) \cong \HH^{-1}(E^{c,b}_{i_0})$, for all $i\in [i_0,n-1]$.

Assume the latter situation holds and denote the line bundle $\HH^{-1}(E^{c,b}_{i_0})$ by $\OO_X (-M_{cH,bH})$. Then it gives a (singular) foliation of $X$, that is, we have the following exact sequence
\begin{equation}\label{foli-cb}
\xymatrix@R=12pt@C=12pt{
  0\ar[r]&\OO_X (-M_{cH,bH})\ar[r]&\TET\ar[r]& {\cal J}_{W_{cH,bH}} (-L_{cH,bH})\ar[r]&0,
}
\end{equation}
where $ {\cal J}_{W_{cH,bH}} $ is the ideal sheaf of at most $0$-dimensional subscheme $W_{cH,bH}$ and $\OO_X (-L_{cH,bH}) \in \FF_{H,bH}$. One has
$$
\text{$L^2_{cH,bH} \leq 0$ and $M^2_{cH,bH} = \tau_X +2deg (W_{cH,bH}) -L^2_{cH,bH}$.}
$$
In particular, the divisor class $M_{cH,bH}$ lies in the positive cone $C^{+}(X)$.
\end{pro}

\begin{pf}
All assertions with the exception of the ones about $\OO_X (-M_{cH,bH})$ are contained in Lemma \ref{lem:H-1-E}. With the notation $\HH^{-1}(E^{c,b}_{i_0})=\OO_X (-M_{cH,bH})$ the complex in \eqref{E-FEcpx} reads as follows:
$$
\xymatrix@R=12pt@C=12pt{
  0\ar[r]&\OO_X (-M_{cH,bH})\ar[r]&\TET \ar[r] &F_{E^{c,b}_{i}}\ar[r]&\HH^{0}(E^{c,b}_{i})\ar[r]&0,
}
$$
for every $i_0 \leq i<n$. The cokernel of the monomorphism
$$
\OO_X (-M_{cH,bH}) \longrightarrow \TET
$$
is a rank $1$ subsheaf of a torsion free sheaf $F_{E^{c,b}_{i}}$. So it is torsion free and has the form ${\cal J}_{W_{cH,bH}} (-L_{cH,bH})$, where $ {\cal J}_{W_{cH,bH}} $ is the ideal sheaf of at most $0$-dimensional subscheme $W_{cH,bH}$. Thus we obtain an exact sequence
\begin{equation}\label{foli-cb1}
\xymatrix@R=12pt@C=12pt{
  0\ar[r]&\OO_X (-M_{cH,bH})\ar[r]&\TET \ar[r] &{\cal J}_{W_{cH,bH}} (-L_{cH,bH}).
}
\end{equation}
By definition $\OO_X (-M_{cH,bH})=\HH^{-1}(E^{c,b}_{i_0})$ is in $\FF_{H,bH}$. Since
${\cal J}_{W_{cH,bH}} (-L_{cH,bH})$ is defined as a subsheaf of $F_{E^{c,b}_{i}} \in \FF_{H,bH}$, it must be in $\FF_{H,bH}$ as well. Hence $\OO_X(-L_{cH,bH})\in \FF_{H,bH} $ and we have the inequality
$$
-L_{cH,bH} \cdot H \leq bH^2,
$$
or, equivalently, $L_{cH,bH} \cdot H \geq (-b)H^2$. Since $b$ is negative, we deduce $L_{cH,bH} \cdot H >0$. Dualizing the exact sequence \eqref{foli-cb1}, we see that $\OO_X(L_{cH,bH})$ is a subsheaf of $\Omega_X$. Hence $L_{cH,bH} $ can not be in the positive cone $C^{+}(X)$ of $X$. Thus $L^2_{cH,bH} \leq 0.$ 

Going back to the exact sequence \eqref{foli-cb1} we obtain
$$
\begin{gathered}
\tau_X=2ch_2(\TET) =2ch_2(\OO_X (-M_{cH,bH}))+2ch_2({\cal J}_{W_{cH,bH}} (-L_{cH,bH}))\\
=M^2_{cH,bH}+L^2_{cH,bH} -2deg(W_{cH,bH}).
\end{gathered}
$$
Expressing $M^2_{cH,bH}$ gives the equality
$$
M^2_{cH,bH}=\tau_X +2deg(W_{cH,bH}) -L^2_{cH,bH}.
$$
The positivity of $\tau_X$ and nonnegativity of two other terms imply
\begin{equation}\label{M2pos}
M^2_{cH,bH}\geq \tau_X >0.
\end{equation}
In addition, $\OO_X (-M_{cH,bH})\in \FF_{H,bH}$ and hence 
$$
-M_{cH,bH} \cdot H \leq b H^2.
$$
This inequality is equivalent to
$$
M_{cH,bH} \cdot H \geq (-b) H^2 >0
$$
and, together with the inequality \eqref{M2pos}, it tells us that $M_{cH,bH} \in C^{+}(X)$. 
\end{pf}

Each nonzero proper subobject $E^{c,b}_{i}$ of the HN filtration \eqref{HN-Th} is
Bridgeland destabilizing for $\TET[1]$ and the exact sequences
$$
\xymatrix@R=12pt@C=12pt{
  0\ar[r]&E^{c,b}_{i} \ar[r]&\TET[1]\ar[r]&F_{E^{c,b}_{i}}[1]\ar[r]&0,
  }
  $$
  for every $i\in[1,n-1]$, give us destabilizing sequences for $\TET[1]$ in the abelian category $\AC_{H,bH}$. This follows immediately from
  the formula
  $$
  Z_{cH,bH}(\TET[1])=Z_{cH,bH}(E^{c,b}_{i}) + Z_{cH,bH}(F_{E^{c,b}_{i}}[1])
  $$
  and the fact that the phase $\phi^{c,b} (E^{c,b}_{i})$ of $E^{c,b}_{i}$ is bigger than the phase $\phi^{c,b} (F_{E^{c,b}_{i}}[1])$ of $F_{E^{c,b}_{i}}[1]$; this in turn follows from the condition \eqref{phases} on phases of the semistable factors of the HN filtration \eqref{HN-Th}. Thus we have
  $$
  \phi^{c,b} (E^{c,b}_{i}) > \phi^{c,b} ( \TET[1]),
  $$
or, equivalently,
\begin{equation}\label{phase-ineq}
\phi^{c,b} ( \TET[1])>\phi^{c,b} (F_{E^{c,b}_{i}}[1]).
\end{equation} 
Similarly to the phase inequality in Proposition \ref{pro:semicircH}, these inequalities will produce walls in the half-plane 
$$
\Pi_H=\{(cH,bH) | c>0, b\in \RR \}
$$
 of stability conditions. In the next subsection we collect the facts about the structure of walls for $\TET[1]$ in $\Pi_H $.

\vspace{0.5cm}
\noindent
{\bf \ref{s:HNfilt}.2. The walls  for $\TET[1]$ in $\Pi_H $.} Let us recall the general set-up:

for an element $v$ in the Grothendieck group $K_0 (\DD)$ of the derived category
$D$ we denote by $ch(v)$ its Chern character lying in the Mukai lattice
$$
\Lambda_X =H^0(X,\ZZ)\oplus NS(X)\oplus H^4(X,\ZZ)[\HA]\cong \ZZ \oplus NS(X)\oplus \HA \ZZ.
$$
We write the vector 
$$
ch(v)=\begin{pmatrix} ch_0(v)\\ch_1(v)\\ch_2(v) \end{pmatrix} \in \Lambda_X.
$$
\begin{example}
$$
ch(\TET[1])=ch(-\TET)=-ch(\TET)=-\begin{pmatrix} 2\\-K_X\\\HA \tau_X \end{pmatrix},
$$
where the first equality comes from the identity
 $$
E[1]=-E
$$
in $K_0 (\DD)$, for any object $E \in \DD$.
\end{example}

Given two elements $v,w$ in $K_0 (\DD)$ whose Mukai vectors $ch(v)$ and $ch(w)$ in $\Lambda_X$ are $\QQ$-linearly independent, one defines the {\it numerical wall} 
$W_v(w)$ of $v$ with respect to $w$ as the subset of stability conditions $({\cal P}, Z) \in Stab_X$ for which the values of the central charge $Z(v)$ and $Z(w)$ are on the same ray in $\CC$. From Bridgeland deformation theorem it follows that $W_v(w)$ is a real codimension $1$ subvariety of $Stab_X$. The restriction of the wall $W_v(w)$ to a half-plane $\Pi_H$
$$
W^H_v (w):= W_v(w) \cap \Pi_H
$$
can be described quite explicitly.
\begin{lem}\label{lem:wallWvw}
 In the half-plane
$$
\Pi_H=\{(cH,bH)|c>0, b\in \RR\}
$$
the wall
$W^H_v (w)$ is described by the determinental equation
\begin{equation}\label{Wvw-eq}
\begin{vmatrix}
b^2 +c^2& ch_2(w)&ch_2(v)\\
2b& ch_1(w) \cdot H& ch_1(v) \cdot H \\
\frac{2}{H^2}& ch_0(w)&ch_0(v)
\end{vmatrix} 
=0,
\end{equation}
where $ch(v)=\begin{pmatrix} ch_0(v)\\ch_1(v)\\ch_2(v) \end{pmatrix} $ and 
 $ch(w)=\begin{pmatrix} ch_0(w)\\ch_1(w)\\ch_2(w) \end{pmatrix} $ are the Mukai vectors of $v$ and $w$.
\end{lem}
\begin{pf}
  Setting $\omega=cH$ and $B=bH$ in the formula \eqref{cc} we obtain
  \begin{equation}\label{Zcb(a)}
  Z_{cH,bH}(a)=-ch_2(a)+ch_1(a)\cdot H b -\frac{ch_0 (a)H^2}{2} (b^2-c^2) +(ch_1(a)\cdot H -ch_0 (a)H^2 b)c\, i,
 \end{equation} 
 for a Mukai vector $ch(a)=\begin{pmatrix} ch_0(a)\\ch_1(a)\\ch_2(a) \end{pmatrix} $. The condition that the values $Z_{cH,bH}(v)$ and $Z_{cH,bH}(w)$ of the central charge  $Z_{cH,bH}$ are aligned in the complex plane $\CC$ is given by the equation
 $$
 \begin{vmatrix}
   \Re (Z_{cH,bH}(w))& \Re (Z_{cH,bH}(v))\\
     \Im (Z_{cH,bH}(w))& \Im (Z_{cH,bH}(v))
   \end{vmatrix}=0,
   $$
   where $\Re (z)$ (resp. $\Im (z)$) denotes the real (resp. imaginary) part of a complex number $z$. From the formula in \eqref{Zcb(a)} we obtain
   $$
   \begin{vmatrix}
     -ch_2(w)+ch_1(w)\cdot H b -\frac{ch_0 (w)H^2}{2} (b^2-c^2) & -ch_2(v)+ch_1(v)\cdot H b -\frac{ch_0 (v)H^2}{2} (b^2-c^2)\\
     ch_1(w)\cdot H -ch_0 (w)H^2 b& ch_1(v)\cdot H -ch_0 (v)H^2 b
   \end{vmatrix}=0.
   $$
   From this it follows
   $$
   \begin{gathered}
 \begin{vmatrix}
     -ch_2(w)+ch_1(w)\cdot H b  & -ch_2(v)+ch_1(v)\cdot H b \\
     ch_1(w)\cdot H& ch_1(v)\cdot H 
   \end{vmatrix} +
   \begin{vmatrix}
     -ch_2(w)+ch_1(w)\cdot H b  &-\frac{ch_0 (v)H^2}{2} (b^2-c^2)\\
     ch_1(w)\cdot H&-ch_0 (v)H^2 b
   \end{vmatrix}\\
   +
   \begin{vmatrix}
     -\frac{ch_0 (w)H^2}{2} (b^2-c^2)&-ch_2(v)+ch_1(v)\cdot H b\\
     -ch_0 (w)H^2 b&ch_1(w)\cdot H
   \end{vmatrix}=0.
   \end{gathered}
   $$
   Simplifying we obtain
   $$
   \begin{gathered}
     -\begin{vmatrix}
    ch_2(w)  & ch_2(v) \\
     ch_1(w)\cdot H& ch_1(v)\cdot H 
   \end{vmatrix} +
   (ch_0 (v)ch_2(w) -ch_0 (w)ch_2(v))H^2b\\
   - \frac{H^2}{2}(ch_0(v)ch_1 (w)\cdot H -ch_0 (w)ch_1 (v)\cdot H)(b^2 +c^2) =0.
 \end{gathered}
 $$
 Multiplying by $(-\frac{2}{H^2})$ and writing all the coefficients in the determinental
 form gives
 $$
 \begin{vmatrix}
    ch_2(w)  & ch_2(v) \\
    ch_1(w)\cdot H& ch_1(v)\cdot H
  \end{vmatrix}\frac{2}{H^2}
  -
 2 \begin{vmatrix}
    ch_2(w)  & ch_2(v) \\
   ch_0 (w)& ch_0 (v)
 \end{vmatrix}b
 +
 \begin{vmatrix}
   ch_1(w)\cdot H& ch_1(v)\cdot H\\
   ch_0 (w)& ch_0 (v)
 \end{vmatrix}
 (b^2+c^2)
 =0.
 $$
 The expression on the left in the above equation is recognized as the expansion with respect to the first column of the $3\times 3$ determinant in the equation
 \eqref{Wvw-eq}.
\end{pf}

We now take $w=\TET[1]$ to obtain the equation
\begin{equation}\label{WvTh-eq}
\begin{vmatrix}
b^2 +c^2& \HA \tau_X&ch_2(v)\\
2b& -K_X \cdot H& ch_1(v) \cdot H \\
\frac{2}{H^2}& 2&ch_0(v)
\end{vmatrix} 
=0
\end{equation}
of the wall $W^H_v (\TET[1])$ of $v\in K_0 (\DD)$ with respect to $\TET[1]$ in the half-plane $\Pi_H$. This is nontrivial, unless the vector
  $\begin{pmatrix}ch_2(v)\\ch_1(v) \cdot H \\ch_0(v)\end{pmatrix}$ is a scalar multiple of $\begin{pmatrix}\HA \tau_X\\-K_X \cdot H \\2\end{pmatrix}$. In this case $\Pi_H$ is contained in $W_v (\TET[1])$.
\begin{lem}\label{lem:Th-walls}
  Assume the half-plane $\Pi_H$ is not contained in $W_v (\TET[1])$, that is, the vector
  $v_H=\begin{pmatrix}ch_2(v)\\ch_1(v) \cdot H \\ch_0(v)\end{pmatrix}$ is not a scalar multiple of $a_H=\begin{pmatrix}\HA \tau_X\\-K_X \cdot H \\2\end{pmatrix}$. Then
  the wall $W^H_v (\TET[1])$ is either the line with the equation
  $$
  b=-\frac{K_X \cdot H}{2H^2},
  $$
  and this happens if and only if in the $3\times 3$ determinant in \eqref{WvTh-eq} the minor 
  $$
M_{11}(v)= \begin{vmatrix}
 -K_X \cdot H& ch_1(v) \cdot H \\
 2&ch_0(v)
\end{vmatrix}=0,
$$
or $W^H_v (\TET[1])$ is the semicircle
\begin{equation}\label{semicirc-Cv}
W^H_v (\TET[1])=\{(b -b_{v,H})^2 + c^2 =r^2_{v,H}\},
\end{equation}
where 
$$
\begin{gathered}
b_{v,H} =\frac{\begin{vmatrix}\HA\tau_X&ch_2(v)\\
2&ch_0(v)
\end{vmatrix}}{M_{11}(v)}=-\frac{\HA ch_0(v)\tau_X -2ch_2(v)}{ch_0(v)K_X \cdot H+2ch_1(v)\cdot H},
\\
\\
r^2_{v,H} =b^2_{v,H} -\frac{2}{H^2}\frac{\begin{vmatrix}\HA\tau_X&ch_2(v)\\
 -K_X \cdot H&ch_1(v)\cdot H
\end{vmatrix}}{M_{11}(v)}=b^2_{v,H} +\frac{ch_1(v)\cdot H \tau_X +2ch_2(v)K_X \cdot H}{H^2(ch_0(v)K_X \cdot H+2ch_1(v)\cdot H)}.
\end{gathered}
$$
\end{lem}
\begin{pf} This follows directly from the equation \eqref{WvTh-eq}.
\end{pf}
Recall Proposition \ref{pro:semicircH}. The results there can be now understood from the perspective of numerical walls with respect to $\TET[1]$:

\vspace{0.2cm}
1) the vertical line
$ b=-\frac{K_X \cdot H}{2H^2}$ appears there as the left boundary of the strip $\Pi'_H$; now that line is understood as the only numerical wall with respect to $\TET[1]$ which is a line; it is usually referred to as the (unique){\it  vertical} wall for $\TET[1]$, since we choose to draw the walls in the $(b,c)$-plane;

\vspace{0.2cm}
2) the semicircle $C_H$ in Proposition \ref{pro:semicircH} is the wall $W^H_{\OO_X} (\TET[1])$ of $\OO_X$ with respect to $\TET[1]$.

\vspace{0.5cm}
The structure of the set of numerical walls in a half-plane $\Pi_H$ is well understood, see \cite{Ma-S} and references therein. For convenience of the reader in the next statement we give a summary of the properties of the numerical walls $W^H_{v} (\TET[1])$ for $\TET[1]$ .
\begin{pro}\label{pro:Th-walls}
  The set of numerical walls $\{W^H_{v} (\TET[1])\}$ for $\TET[1]$ has the following properties:

  \vspace{0.2cm}
  1) two distinct numerical walls do not intersect,
  
\vspace{0.2cm}
2) in $(b,c)$-plane, the tops of the semicircular walls $W^H_{v} (\TET[1])$ in \eqref{semicirc-Cv} are on the hyperbola
$$
\Gamma_H:=\{\Re(Z_{cH,bH}(\TET[1])=0\},
$$
where $\Re(Z_{cH,bH}(\TET[1])$ is the real part of $Z_{cH,bH}(\TET[1])$ and the explicit equation of $\Gamma_H$ is
$$
b^2-c^2 + \frac{K_X \cdot H}{H^2}b +\frac{\tau_X}{2H^2}=0.
$$
\end{pro}

Since it is easy in this case, we supply the proof (the general case can be done along the similar lines).

\vspace{0.2cm}
\noindent
{\it Proof.}
  Let $W^H_{v} (\TET[1])$ and $W^H_{v'} (\TET[1])$ be two distinct walls for
  $\TET[1]$ in the half-plane $\Pi_H$. Assume they intersect at a point $\begin{pmatrix}b_0\\c_0\end{pmatrix}$ of $\Pi_H$; that half-plane is identified with $(b,c)$-upper half-plane. This means that $\begin{pmatrix}b_0\\c_0\end{pmatrix}$ is a solution of the two determinental equations, see \eqref{WvTh-eq},
  $$
 \begin{vmatrix}
b^2_0 +c^2_0& \HA \tau_X&ch_2(v)\\
2b_0& -K_X \cdot H& ch_1(v) \cdot H \\
\frac{2}{H^2}& 2&ch_0(v)
\end{vmatrix} 
=
\begin{vmatrix}
  b^2_0 +c^2_0& \HA \tau_X&ch_2(v')\\
2b_0& -K_X \cdot H& ch_1(v') \cdot H \\
\frac{2}{H^2}& 2&ch_0(v')
\end{vmatrix}
=0.
$$
This in turn implies that the vector  $\begin{pmatrix}b^2_0 +c^2_0\\2b_0\\\frac{2}{H^2}\end{pmatrix}$ in the first column is a linear combination of the remaining two column-vectors. So from the first equation we have
\begin{equation}\label{lincmb}
\begin{pmatrix}b^2_0 +c^2_0\\2b_0\\\frac{2}{H^2}\end{pmatrix} =p a_H + q v_H,
\end{equation}
for some $p,q\in \RR$ and where $a_H$ (resp. $v_H$) is the second (resp. third)
column-vector in the first determinant.

If the scalar $q\neq 0$, substituting
the above linear combination in the second equation we deduce that the vectors
$v_H$, $a_H$ and $v'_H$ (the third column vector of the second determinant) are linearly dependent:
\begin{equation}\label{lindep}
sv_H + ta_H +s'v'_H=0,
\end{equation}
for some real constants $s,s',t$. Observe that both $s$ and $s'$ must be nonzero, since otherwise either $v_H$ or $v'_H$ is a multiple of $a_H$, but we know from Proposition \ref{lem:Th-walls} that
in this case either $W^H_{v} (\TET[1])$ or $W^H_{v'} (\TET[1])$ is not a wall in $\Pi_H$. This tells us that $v'_H$ is a linear combination of $v_H$ and $a_H$ and hence the walls $W^H_{v} (\TET[1])$ and $W^H_{v'} (\TET[1])$ coincide. This contradicts the assumption that the walls are distinct.

Thus we conclude that $q=0$ in \eqref{lincmb}. Then
$\begin{pmatrix}b^2_0 +c^2_0\\2b_0\\\frac{2}{H^2}\end{pmatrix} =p a_H$ and we have
$$
p=\frac{1}{H^2},\,\, b_0=-\frac{K_X \cdot H}{2H^2},\,\, b^2_0 +c^2_0=\frac{\tau_X}{2H^2}.
$$
Substituting for $b_0$ in the third equation gives
$$
c^2_0=\frac{\tau_X}{2H^2} -\Big(\frac{K_X \cdot H}{2H^2}\Big)^2 =\frac{K_X \cdot H}{2H^2} \Big(\frac{\tau_X}{K_X \cdot H} -\frac{K_X \cdot H}{2H^2}\Big) <0,
  $$
  where the last inequality comes from the inequality
  $$
  \frac{\tau_X}{K_X \cdot H}< \frac{K_X \cdot H}{2H^2}
  $$
  shown in the proof of Proposition \ref{pro:semicircH}. This completes the proof of the first assertion of the proposition.

  We now turn to the part 2) of the proposition. From the formula \eqref{Zcb(a)} it follows that the real part of $Z_{cH,bH}(\TET[1])$ has the form
  $$
  \Re(Z_{cH,bH}(\TET[1])=\HA \tau_X + K_X \cdot H b + H^2 (b^2 -c^2).
  $$
  Hence the set
  $$
  \Gamma_H=\{ \Re(Z_{cH,bH}(\TET[1])=0\}
  $$
  is determined by the equation
  \begin{equation}\label{Hyperb-eqn}
  b^2 -c^2 +\frac{K_X \cdot H}{H^2} b+\frac{ \tau_X}{2H^2}=0
  \end{equation}
  and this is the hyperbola appearing in 2).

  The `top' of the semicircle
  $C_{v,H}$ in $(b,c)$-plane is the point
  $$
  t_{v,H}=\begin{pmatrix}b_{v,H}\\r_{v,H} \end{pmatrix}.
  $$
  We need to check that this is a solution of the equation \eqref{Hyperb-eqn}.
  From the formula for $r^2_{v,H}$ in Lemma \ref{lem:Th-walls} it follows
  $$
 b^2_{v,H}-r^2_{v,H} =\frac{2}{H^2}\frac{\begin{vmatrix}\HA\tau_X&ch_2(v)\\
 -K_X \cdot H&ch_1(v)\cdot H
\end{vmatrix}}
{M_{11}(v)
},
$$
where $M_{11}(v)=\begin{vmatrix} -K_X \cdot H&ch_1(v)\cdot H\\
2&ch_0(v)
\end{vmatrix}$
is $(1,1)$-minor of the determinant in \eqref{WvTh-eq}.

From this and the formula for $b_{v,H}$ in Lemma \ref{lem:Th-walls}  the value of the left hand side in \eqref{Hyperb-eqn} becomes
$$
\begin{gathered}
\frac{2}{H^2}\frac{\begin{vmatrix}\HA\tau_X&ch_2(v)\\
 -K_X \cdot H&ch_1(v)\cdot H
\end{vmatrix}}
{M_{11}(v)}
+
\frac{K_X \cdot H}{H^2}\frac{\begin{vmatrix}\HA\tau_X&ch_2(v)\\
 2&ch_0(v)
\end{vmatrix}}
{M_{11} (v)}
+
\frac{ \tau_X}{2H^2}=\\
\\
\frac{1}{H^2M_{11} (v)}
\Big(
2\begin{vmatrix}\HA\tau_X&ch_2(v)\\
  -K_X \cdot H&ch_1(v)\cdot H
\end{vmatrix}
+K_X \cdot H \begin{vmatrix}\HA\tau_X&ch_2(v)\\
 2&ch_0(v)
\end{vmatrix}
+ \HA\tau_X M_{11} (v)
\Big)
\\
=\frac{1}{H^2M_{11} (v)}
\begin{vmatrix}\HA \tau_X&\HA\tau_X& ch_2(v)\\
  -K_X \cdot H&-K_X\cdot H& ch_1(v)\cdot H\\
  2&2&ch_0 (v)
 \end{vmatrix}=0. 
\end{gathered}
$$
This completes the proof of the proposition. $\Box$

\begin{rem}\label{rem:hyperbola}
	The equation of the hyperbola $\Gamma_H$, Proposition \ref{pro:Th-walls}, 2), which parametrizes the tops of semicircular numerical walls for $\TET[1]$ can be rewritten in the following form:
	$$
	\Big(b+\frac{K_X \cdot H}{2H^2}\Big)^2 -c^2 =\Big(\frac{K_X \cdot H}{2H^2}\Big)^2 -\frac{\tau_X}{2H^2}=\frac{\Delta^H (\TET[1])}{(2H^2)^2},
	$$
	where in the last expression $\Delta^H (\TET[1])$ is the $H$-discriminant of $\TET[1]$, see \eqref{Deltah}; we recall that
	$\Delta^H (\TET[1]) >0$, see Remark \ref{rem:ineq-tau=Hdisc}, the inequality \eqref{Hdisc-positive}. The part of the hyperbola $\Gamma_H$ in the upper half  $(b,c)$-plane is depicted below.
	$$
	
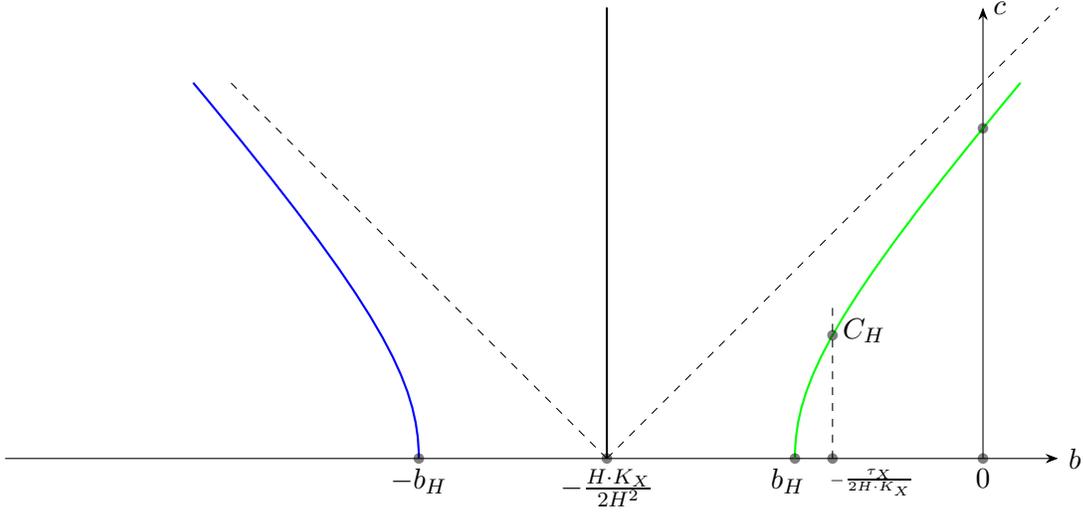
\begin{tikzpicture}[>=Stealth]
		\colorlet{semicircle}{yellow!80!red}
		\draw [->] (-7,0)--(1,0)--(4,0)--(6,0)--(7,0);
		\coordinate[label=below:{$-\frac{H\cdot K_X}{2H^2}$}] (A) at (1,0);
		\coordinate [label=below right:{$\text{\tiny $\!\!\!\!-\frac{\tau_X}{2H\cdot K_X}$}$}] (B) at (4,0);
		\coordinate [label=below:{$0$}] (C) at (6,0);
		\coordinate [label=below:{$\!\!\!b_H$}] (H) at (3.5,0);
		\coordinate [label=below:{$-b_H$}] (F) at (-1.5,0);
		\draw[->,name path=vertical axis](6,0)--(6,6);
		\coordinate[label=right:{$c$}] (D) at (6,6);
		\coordinate[label=right:{$b$}] (E) at (7,0);
		\draw[dashed] (1,0)--(5,4)--(7,6);
		\draw[dashed] (1,0)--(-4,5);
		\draw[name path=hyperbolar, thick,color=green] (3.5,0).. controls (3.5,1) and (4,2).. (6.5,5);
\draw[name path=hyperbolal, thick,color=blue] (-1.5,0).. controls (-1.5,1) and (-2,2).. (-4.5,5);
\draw[name intersections={of= vertical axis and hyperbolar, by=x}];					
\draw[name path=vertical t,dashed](4,0)--(4,2.1);
\draw[name intersections={of= vertical t and hyperbolar, by=y}];
		\coordinate[label=below right:{$C_H$}] (G) at (4,2);
		\draw[thick] (1,0)--(1,6);
		\foreach \point in {A,B,C,H,F,x,y}
		\fill[black,opacity=.5](\point)	circle(2pt);
	\end{tikzpicture}
	$$
	The numerical semicircular walls for $\TET[1]$ have their tops on the two upper half branches of the hyperbola which are symmetric with respect to the vertical wall $b=-\frac{H\cdot K_X}{2H^2}$. The points of intersection $\pm b_H$ of the hyperbola with the $b$-axis are
	$$
	\pm b_H= -\frac{H\cdot K_X}{2H^2} \pm \frac{\sqrt{\Delta^H (\TET[1])}}{2H^2} .
	$$
	The semicircular walls with the tops on each branch of the hyperbola, as we move from the point $\pm b_H$ to infinity, form an increasing family of semicircles: those semicircles are nested.
	
	The branch which we are concerned with is on the right side of the vertical wall. The point $C_H$ on the picture marks the top of the wall $W^H_{\OO_X} (\TET[1])$, the semicircle $C_H$ found in Proposition \ref{pro:semicircH}.  
\end{rem}

All of the above is summarized in the following statement.
\begin{cor}\label{cor:semicirwalls}
In the upper half $(b,c)$-plane the numerical walls for $\TET[1]$ are:

\vspace{0.2cm}
1) the vertical wall $b=-\frac{K\cdot H}{2H^2}$,

\vspace{0.2cm}
2) all other walls are semicircles with centers on the $b$-axis; they are placed either on the left or on the right of the vertical wall; they are disjoint from each other as well as from the vertical wall,

\vspace{0.2cm}
3) the walls on either side of the vertical wall are semicircles with their tops on the upper half of the hyperbola
$$
\Gamma_H =\{\Big(b+\frac{K_X \cdot H}{2H^2}\Big)^2 -c^2 =\frac{\Delta^H (\TET[1])}{(2H^2)^2}
\},
$$
where $\Delta^H (\TET[1])=\Big(K_X \cdot H \Big)^2 -2H^2{\tau_X}$ is the $H$-discriminant of $\TET[1]$. The hyperbola has two disjoint branches symmetrically placed with respect to the vertical wall. The semicircular numerical walls with the tops on the same branch of the hyperbola form a 
  sequence of nested semicircles.
\end{cor}

Among the numerical walls for $\TET[1]$ one distinguishes {\it actual} walls. Those 
are walls $W^H_v (\TET[1])$ in the half-plane $\Pi_H$ for which there exists
a stability condition $(c_0H,b_0H) \in W^H_v (\TET[1])$ such that the Mukai vector $v$ is realized by a subobject of $\TET[1]$ in $A_{H,b_0H}$. That is, the vector
$v=ch(E)$, for an object $E\in A_{H,b_0H}$  and we have a monomorphism
$$
\xymatrix@R=12pt@C=12pt{
  0\ar[r]&E\ar[r]&\TET[1] 
}
$$
in $A_{H,b_0H}$. This can be completed to the exact sequence
$$
\xymatrix@R=12pt@C=12pt{
  0\ar[r]&E\ar[r]&\TET[1]\ar[r]&F_E[1]\ar[r]&0  
}
$$
with $F_E \in \FF_{H,b_0H}$. The sequence is {\it strictly semistable for $\TET[1]$} with respect to the central charge $Z_{c_0H,b_0H}$. This is to say that all objects of the sequence have the same phase with respect to $Z_{c_0H,b_0H}$. To stress that a wall $W^H_v (\TET[1])$ is actual, we often write $W^H_E (\TET[1])$, where
$E$ is as just discussed and we call it a {\it subobject realizing $v$}. Observe that we have
\begin{equation}\label{wall-sub-quot}
W^H_E (\TET[1])= W^H_{F_E[1]} (\TET[1]),
\end{equation}  
so we can define an actual wall via the {\it quotient object ${F_E[1]}$} realizing $v$ as well. 
\begin{rem}\label{rem:actualwall}
 Let $W^H_v (\TET[1])$ be an actual wall for $\TET[1]$. It must be a semicircular wall on the right of the vertical wall $b=-\frac{K_X\cdot H}{2H^2}$, since $\TET[1]$ has no subobjects in $\AC_{H,b H}$ of phase $1$, see Lemma \ref{lem:less1}, and we are in $\AC_{H,bH}$ for $b\geq -\frac{K_X\cdot H}{2H^2}$. By definition there is some point $(c_0H,b_0H)$ on the semicircle $W^H_v (\TET[1])$, where $v$ is realized by a subobject $E$ of $\TET[1]$in $\AC_{H,b_0H}$. But in fact, the same holds for all points of $W^H_v (\TET[1])$. Indeed, we know that the quotient object
$F_E[1]$ realizes $v$ as a quotient object. So  $F_E\in \FF_{H,b_0H}$. From the inclusion
$\FF_{H,b_0H} \subset \FF_{H,bH} $ for all $b\geq b_0$, it follows that  $F_E\in \FF_{H,bH}$, for all $b\geq b_0$. Hence $F_E[1]$ realizes $v$ at all points of the 
semicircle $W^H_v (\TET[1])$ lying to the right of the point $(c_0H,b_0H)$.
A similar argument shows that $E$ lies in $\AC_{H,bH}$, for all $b\in (-\frac{K_X\cdot H}{2H^2},b_0]$. Hence,  $v$ is realized at all points of the semicircle
$W^H_v (\TET[1])$.
\end{rem}

The actual walls are subject to the following finiteness property.

\begin{lem}\label{lem:awall-finite}
Let $\beta_0$ be a rational number in the interval $(-\frac{K_X\cdot H}{2H^2},+\infty)$. Then there is at most finite number of actual walls for $\TET[1]$ intersecting the vertical line $b=\beta_0$.
\end{lem}
This is known to experts on Bridgeland stability, see \cite{Ma-S}, but we supply a proof here for
non-experts.

\vspace{0.2cm}
\noindent
{\it Proof.}
Let $W^H_v (\TET[1])$ an actual wall for $\TET[1]$. Then we know that it is a
semicircle
in the upper half of $(b,c)$-plane given
by the determinental equation
$$
\begin{vmatrix}
b^2 +c^2 &\HA \tau_X & ch_2 (v)\\
2b&-K\cdot H& ch_1 (v) \cdot H\\
\frac{2}{H^2}& 2&ch_0 (v)
\end{vmatrix}
=0.
$$
The condition of intersecting the vertical line $b=\beta_0$ says that
\begin{equation}\label{v-beta0-eq}
\begin{vmatrix}
\beta^2_0 +c^2_v &\HA \tau_X & ch_2 (v)\\
2\beta_0&-K\cdot H& ch_1 (v) \cdot H\\
\frac{2}{H^2}& 2&ch_0 (v)
\end{vmatrix}
=0
\end{equation}
where $(\beta_0, c_v)$ is the point of intersection of the semicircle  $W^H_v (\TET[1])$ with the line $b=\beta_0$ in the upper half of $(b,c)$-plane. We need to check that the set of column vectors
\begin{equation}\label{vvectors}
v_H=\begin{pmatrix}ch_2 (v)\\ch_1 (v) \cdot H\\ch_0 (v)\end{pmatrix}
\end{equation}
satisfying the above equation is at most finite for $v$ realized by objects
in $\AC_{H,\beta_0 H}$. To do this we use an exact sequence
\begin{equation}\label{E-F-seq}
\xymatrix@R=12pt@C=12pt{
  0\ar[r]&E\ar[r]&\TET[1] \ar[r]& F[1]\ar[r]&0
}
\end{equation}
in $\AC_{H,\beta_0 H}$ with $F[1]$ (resp. $E$) realizing $v$ as a quotient (resp. sub-) object of $\TET[1]$. From this exact sequence we obtain
$$
K_X \cdot H+2\beta_0 H^2=\Im(Z_{H,\beta_0 H}(\TET[1] ))=\Im(Z_{H,\beta_0 H}( F[1])+\Im(Z_{H,\beta_0 H}( E)).
$$
Both terms on the right side are positive since $E$ and $F[1]$ are in $\AC_{H,\beta_0 H}$ and the phase of $\TET[1]$ is less then $1$. The latter because
$K_X \cdot H+2\beta_0 H^2 >0$ according to the assumption on $\beta_0$. Thus we deduce
$$
0< ch_1(v)\cdot H -ch_0(v)\beta_0 H^2 <K_X \cdot H+2\beta_0 H^2,
$$
for all $v$ solving \eqref{v-beta0-eq} and realized in $\AC_{H,\beta_0 H}$.
Since the set of values $(ch_1(v)\cdot H -ch_0(v)\beta_0 H^2)$ is discrete in
$\RR$ because of the assumption that $\beta_0$ is rational, we deduce that the set 
\begin{equation}\label{Vbeta0}
	\begin{gathered}
V_{\beta_0}:=\{(ch_1 (v)\cdot H -ch_0(v)\beta_0 H^2)|\,\, \text{\it $v$ is a solution of \eqref{v-beta0-eq} and realized in $\AC_{H,\beta_0 H}$} \}
\\
\text{\it is finite. In addition, the same holds for $V_{\beta}$, for any rational $\beta \geq \beta_0$.}
\end{gathered}
\end{equation}

Take $v=F[1]$, where $F[1]$ is as in the exact sequence in \eqref{E-F-seq}
and consider the walls in $\Pi_H$ for that object. We know that the equation
$$
\Re(Z_{cH,bH} (F[1]))=0
$$
parametrizes the tops of the semicircular walls for $v=F[1]$. Writing out the left hand side explicitly we obtain 
$$
\HA ch_0(v)H^2 (b^2-c^2) -ch_1 (v)\cdot H b +ch_2(v) =0.
$$
Since $ch_0(v)=ch_0 (F[1])=-ch_0 (F)=-rk(F)$ is nonzero, the above is a hyperbola which we rewrite as follows
\begin{equation}\label{F-hyperbola}
\Big(b-\frac{ch_1(v)\cdot H}{ch_0(v)H^2}\Big)^2 -c^2 =\Big(\frac{ch_1(v)\cdot H}{ch_0(v)H^2}\Big)^2 -2\frac{ch_2(v)}{ch_0(v)H^2}.
\end{equation}
If the hyperbola intersects the $b$-axis, then the right hand side is nonnegative and we have
\begin{equation}\label{disc-ineq}
  \Delta^H (v)=(ch_1(v)\cdot H)^2 -2ch_0(v)H^2 ch_2(v)\geq 0.
\end{equation}
Since we are free to replace $v$ by $v(-\beta H)=v\otimes \OO_X (-\beta H)$, we deduce
$$
2ch_0(v)H^2ch_2 (v (-\beta H)) \leq (ch_1 (v (-\beta H))\cdot H)^2=(ch_1(v)\cdot H-ch_0 (v)\beta H^2)^2.
$$
From \eqref{Vbeta0} we know that for any given rational number $\beta\geq \beta_0$, the expression on the right takes at most finite number of values as $v$ runs trough all objects defining actual walls
intersecting the line $b=\beta_0$. This implies
\begin{equation}\label{ch0ch2ub}
	\begin{gathered}
	\text{\it for objects $v=F[1]$ subject to \eqref{disc-ineq} and defining actual walls for $\TET[1]$}
		\\
	\text{\it $2ch_0(v)H^2ch_2 (v (-\beta H)) $ is bounded from above, for every rational $\beta \geq \beta_0$.}
	\end{gathered}
\end{equation}
We now seek the lower bound for this expression. For this we relate it to the central charge $Z_{cH,\beta H}$.

From the formula of central charge $Z_{cH,bH}$ we have
\begin{equation}\label{ch2vb}
ch_2 (v (-bH)) =-\Re(Z_{cH,bH})+\HA ch_0 (v)H^2 c^2
\end{equation}
for all $c>0$ and all $b\in \RR$. We now take $b\geq\beta_0$ and $c_b$ the coordinates of the points lying on the wall $W^H_v (\TET[1])$. Remembering that for such pairs the central
charges $Z_{c_b H,\beta H} (v)$ and $Z_{c_b H,b H} (\TET[1])$ are aligned on the same ray in the upper half $(b,c)$-plane, see Remark \ref{rem:actualwall}, we obtain
\begin{equation}\label{aligned}
	Z_{c_b H,b H} (v)=\mu_v (b) Z_{c_b H,b H} (\TET[1]),
\end{equation}
for some $\mu_v(b) >0$. Furthermore, taking $v=F[1]$ and $v'=E$ in the exact sequence \eqref{E-F-seq} we deduce
$$
\begin{gathered}
	Z_{c_bH,b H} (\TET[1])=Z_{c_b H,b H} (E)+Z_{c_b H,b H} (F[1])=\mu_{v'}(b) Z_{c_b H,b H} (\TET[1])+\mu_{v} (b) Z_{c_b H,b H} (\TET[1])\\
	=(\mu_v (b) +\mu_{v'} (b))Z_{c_bH,b H} (\TET[1]),
\end{gathered}
$$
from which we deduce
$$
\mu_v (b) +\mu_{v'} (b)=1.
$$
Thus the scalars $\mu_v (b)$ (resp. $\mu_{v'} (b)$) in \eqref{aligned} are bounded
$$
0<\mu_v (b)<1\,\,\text{(resp. $0<\mu_{v'}(b) <1$)}.
$$
Putting together \eqref{aligned} and \eqref{ch2vb} we obtain
\begin{equation}\label{deform-formulas}
	\begin{gathered}
		ch_2(v (-bH)) =-\Re(Z_{c_b H,bH} (v)) + \HA  ch_0(v)H^2 c^2_b \\
		=-\mu_v (b)\Re(Z_{c_b H,bH} (\TET[1]))+ \HA  ch_0(v)H^2 c^2_b.
	\end{gathered}
\end{equation}
Substituting 
$$
-\Re(Z_{c_b H,bH} (\TET[1]))=ch_2(\TET[1]\otimes \OO_X (-b H))+H^2 c^2_b
$$
gives
$$
		ch_2(v (-bH)) =\mu_v (b) ch_2(\TET[1]\otimes \OO_X (-b H)) + \HA ( ch_0(v)+2\mu_v (b)) H^2 c^2_b.   
$$
From this it follows
$$
	ch_0 (v)ch_2(v(-bH)) =\mu_v (b)ch_0 (v) ch_2(\TET[1]\otimes \OO_X (-b H)) + \HA ch_0 (v)( ch_0(v)+2\mu_v (b)) H^2 c^2_b.
$$
Remember $v=F[1]$, so the above reads
\begin{equation}\label{ch0ch2-eq}
	ch_0 (v)ch_2(v (-bH)) =\mu_v (b) rk(F) ch_2(\TET\otimes \OO_X (-b H)) + \HA rk (F)( rk(F)-2\mu_v) H^2 c^2_b,
\end{equation}
for all values $b$ parametrizing the the semicircle $W^H_v (\TET[1])$.
There is a part of these values on the right of $\beta_0$ which lie under the right branch of the hyperbola
$$
\Big(b+\frac{K_X \cdot H}{2H^2}\Big)^2-c^2 -\Big(\Big(\frac{K_X \cdot H}{2H^2}\Big)^2 -\frac{\tau_X}{2H^2}\Big)=\Re(Z_{cH,bH}(\TET[1]))=0,
$$
see the drawing in Remark \ref{rem:hyperbola}.
For those values of $b$ we have
$$
 ch_2(\TET\otimes \OO_X (-b H)) \geq 0.
 $$ 
The expression on the right in \eqref{ch0ch2-eq} is positive for all such values of $b$, provided the rank $rk(F) \geq 2$, and we obtain
 $$
 	ch_0 (v)ch_2(v (-bH)) >0.
 	$$
 	Together with the upper bound \eqref{ch0ch2ub} we obtain
 	\begin{equation}\label{awallfinite1}
 		\begin{gathered}
 		\text{\it
 		the set of objects $v=F[1]$ with $rk(F) \geq 2$ and subject to
 		the inequality \eqref{disc-ineq}}
 	\\
 	\text{\it defining actual walls for $\TET[1]$ and intersecting the line $b=\beta_0$ is finite.}
 \end{gathered}
 	\end{equation}
 
 \vspace{0.2cm}
 We now turn to the case when the inequality \eqref{disc-ineq} fails.
  This means that the object $v=F[1]$ is Bogomolov unstable and 
 this in turn is equivalent to $F$ being Bogomolov unstable. In particular, it is
 $H$-unstable. Taking its HN filtration with respect to $H$
 $$
 F=F_m \supset F_{m-1}\supset \cdots \supset F_1 \supset F_0=0,
 $$
 we have its semistable factors $Q_i =F_i/F_{i-1}$ lying in the subcategory $\FF_{H,\beta_0 H}$. Hence all quotients $Q_i$ are subject to
 $$
 ch_1 (Q_i (-\beta H))\cdot H =c_1(Q_i)\cdot H - rk(Q_i)\beta H^2 \leq 0,\, \forall \beta \geq \beta_0.
 $$
  Thus the dual $F^{\ast} $ of $F$ is subject to the `twisted' positivity condition of Lemma \ref{lem:BGunst}. Namely, for every rational $\beta\geq \beta_0$, the $\QQ$-sheaf $F^{\ast} (\beta H)$ has the HN filtration with respect to $H$ whose $H$-semistable factors $Q^{\ast}_i (\beta H)$ are subject to
  $$
  \frac{ch_1 (Q^{\ast}_m (\beta H))\cdot H}{ch_0 (Q_m)} > \frac{ch_1 (Q^{\ast}_{m-1} (\beta H))\cdot H}{ch_0 (Q_{m-1})}>\cdots > \frac{ ch_1 (Q^{\ast}_1 (\beta H))\cdot H}{ch_0 (Q_1)} \geq 0.
  $$
  The proof of Lemma \ref{lem:BGunst} gives the following estimate
 \begin{equation}\label{Funst}
 	\begin{gathered}
 	2ch_2(F(-\beta H))=2ch_2(F^{\ast}(\beta H)) <
 	\\
 	 2ch_2 (Q^{\ast}_m (\beta H))+ \frac{ch_1 (Q^{\ast}_m (\beta H))\cdot H}{ch_0 (Q_m)H^2}
 	(ch_1 (F^{\ast}(\beta H))\cdot H)-ch_1 (Q^{\ast}_m (\beta H))\cdot H).
 	\end{gathered}
 \end{equation}
This implies 
$$
\begin{gathered}
2H^2ch_0(Q_m) ch_2(F(-\beta H))<  -\Delta^H (Q^{\ast}_m) +\Big( ch_1  (Q^{\ast}_m (\beta H))\cdot H \Big) \Big(ch_1 (F^{\ast}(\beta H))\cdot H \Big)
\\
\leq
\Big(ch_1(Q^{\ast}_m (\beta H))\cdot H \Big) \Big( ch_1 (F^{\ast}(\beta H))\cdot H \Big)=\Big(ch_1(Q_m (-\beta H))\cdot H \Big)\Big( ch_1 (F(-\beta H))\cdot H \Big),
\end{gathered}
$$  
 where the second inequality uses the $H$-Bogomolov-Gieseker inequality
 $$
 \Delta^H (Q^{\ast}_m)=\Big( ch_1  (Q^{\ast}_m)\cdot H\Big)^2  -2H^2ch_0 (Q_m)ch_2(Q^{\ast}_m)\geq 0
 $$
for $H$-semistable sheaf $Q^{\ast}_m$. But the factors 
$\Big( ch_1 (F(-\beta H))\cdot H \Big)$ and $\Big(ch_1(Q_m (-\beta H))\cdot H \Big)$ on the right hand side of the above inequality take finite number of values for every rational value $\beta \geq \beta_0$, see \eqref{Vbeta0}. Hence 
$$
\text{\it $ch_0(Q_m) ch_2(F(-\beta H))$ is bounded above for every rational $\beta\geq \beta_0$.}
$$ 
In addition, the failure of  \eqref{disc-ineq} means that the discriminant of the quadratic polynomial
$$
ch_2(F(-t H))=ch_2(F)-ch_1(F)\cdot H t +\HA H^2 ch_0 (F)t^2
$$
is negative. Hence it is always positive and we deduce the lower bound
$$
 ch_0(Q_m) ch_2(F(-\beta H))>0.
 $$
 Thus we obtain
 $$
 \text{\it $ch_0(Q_m) ch_2(F(-\beta H))$ takes finite number of values for every rational $\beta\geq \beta_0$.}
 $$
This implies 
$$
ch_2(F(-\beta H))=ch_2 (F) -ch_1 (F)\cdot H \beta +\HA H^2 ch_0 (F) \beta^2
$$
 takes finite number of values for every rational $\beta \geq \beta_0$. This is only possible if the set of coefficients 
 $$
 v_H=\begin{pmatrix}
 	ch_2(F)\\ch_1 (F)\cdot H \\ch_0 (F)
 \end{pmatrix}
$$
is finite. This proves the finiteness of the number of actual walls
$W^H_v (\TET[1])$ intersecting $b=\beta_0$ and realized by $v=F[1]$ Bogomolov unstable.

It remains to treat the case when $v=F[1]$ and $F$ has rank $1$, that is, $F$ is a rank $1$ torsion free sheaf on $X$ and we denote it by ${\cal J}_A (-L)$, where $L$ is a divisor and ${\cal J}_A$ is the ideal sheaf of at most $0$-dimensional subscheme $A$ of $X$.

The condition that ${\cal J}_A (-L)$ lies in $\FF_{H,\beta_0 H}$ means
\begin{equation}\label{LHlb}
L\cdot H \geq -\beta_0 H^2.
\end{equation}
On the other hand applying the homological functor $\HH^0$ to the exact sequence \eqref{E-F-seq} with $F={\cal J}_A (-L)$ we deduce
$$
\xymatrix@R=12pt@C=12pt{
	0\ar[r]&\HH^{-1}(E)\ar[r]&\TET \ar[r]& {\cal J}_A (-L) \ar[r]&\HH^{0}(E)\ar[r]&0.
}
$$
From this it follows that $\HH^{-1}(E)$ is a line bundle and we denote it $\OO_X (-M)$. Furthermore, the inclusion $\OO_X (-M)\hookrightarrow \TET$ can be completed to the exact sequence
$$
\xymatrix@R=12pt@C=12pt{
	0\ar[r]&\OO_X (-M)\ar[r]&\TET \ar[r]& {\cal J}_{A'} (-L-D) \ar[r]&0,
}
$$ 
where ${\cal J}_{A'}$ is the sheaf of ideals of at most $0$-dimensional subscheme $A'$ of $X$ and $D$ is an effective divisor. The $H$-stability of $\TET$ implies
$$
-\HA K_X \cdot H < -(L+D)\cdot H.
$$
This implies an upper bound
\begin{equation}\label{LHineq}
	L\cdot H < \HA K_X \cdot H.
\end{equation}
 Together with the lower bound in \eqref{LHlb} we deduce 
$$
\text{\it the set of values $L\cdot H$ is finite.}
$$
It remains to see that 
$$
2ch_2 ({\cal J}_A  (-L))=L^2 -deg(A)
$$ 
takes finite number of values. For the upper bound we have the Hodge Index inequality
$$
H^2 L^2 \leq (L\cdot H)^2:
$$
since $L\cdot H$ is bounded the above inequality tells us that $L^2$ and
hence $2ch_2 ({\cal J}_A  (-L))=(L^2 -2deg(A)) $ is bounded from above.

 For the lower bound on $2ch_2(v)=2ch_2 ({\cal J}_A  (-L))$, we use the fact that the wall $W^H_v (\TET[1])$ must lie on the right side of the vertical wall
 $b=-\frac{K_X \cdot H}{2H^2}$ in order to intersect the vertical line
$b=\beta_0 > -\frac{K_X \cdot H}{2H^2}$. In other words, we have the inequality
\begin{equation}\label{centerontheright}
-\frac{K_X \cdot H}{2H^2} <b_{v,H}=-\frac{\HA ch_0(v) \tau_X -2ch_2(v)}{ch_0 (v)K_X \cdot H +2ch_1(v)\cdot H},
\end{equation}
where $b_{v,H}$ is the $b$-coordinate of the center of $W^H_v (\TET[1])$ and the formula above comes from Lemma \ref{lem:Th-walls}. For $v= ({\cal J}_A  (-L))[1]$ the formula reads
$$
b_{v,H}=-\frac{\HA  \tau_X -2ch_2(v)}{K_X \cdot H -2L\cdot H}.
$$
 From \eqref{LHineq} we know that the denominator is positive, so the inequality \eqref{centerontheright}
  becomes
$$
\Big(K_X \cdot H -2L\cdot H \Big)\frac{K_X \cdot H}{2H^2} >\HA \tau_X -2ch_2 (v)
$$
or, equivalently, we obtain a lower bound
$$
2ch_2(v)>\HA  \tau_X -(K_X \cdot H -2L\cdot H)\frac{K_X \cdot H}{2H^2}.
$$
Since the values of $L\cdot H$ are bounded, the above gives a lower bound for $2ch_2 (v)$. This and the upper bound obtained above tells us that $ch_2 (v)$ is bounded and hence takes finite number of values.
So we conclude that there is at most finite set of $v= {\cal J}_A  (-L)[1]$'s solving 
\eqref{v-beta0-eq}.
This completes the proof of the lemma. $\Box$

\vspace{0.5cm}
From the general perspective of structure of walls for $\TET[1]$ obtained above, we can explain the results of Proposition \ref{pro:semicircH} as follows:

1) we choose to be on the right side of the vertical wall $b=-\frac{K\cdot H}{2H^2}$,
because it is for those stability conditions $(cH,bH)$ the object $\TET[1]$ belongs to $\FF_{H,bH}[1] \subset \AC_{H,bH}$; and we are in the strip
$\Pi'_H$ because  $\OO_X \in \TT_{H,bH}\subset \AC_{H,bH}$ for $b<0$; thus in the strip $\Pi'_H$ both $\OO_X$ and $\TET[1]$ are in $\AC_{H,bH}$,

2) the semicircle $C_H$ found in Proposition \ref{pro:semicircH} is the wall
$W^H_{\OO_X}$ of $\OO_X$ with respect $\TET[1]$; inside this wall we have the stability conditions $(cH,bH)$ for which the phase $\phi_{cH,bH} (\OO_X)$ of $\OO_X$
is bigger than the phase $\phi_{cH,bH} (\TET[1])$,

3) we found the region $B^{un}_H$, see \eqref{BH-unst}, inside the wall
$W^H_{\OO_X}$, where the morphism
$$
\Phi^V_{H,bH}: V\otimes \OO_X \longrightarrow \TET[1]
$$
is a monomorphism in $\AC_{H,bH}$;
as a consequence $\TET[1]$ is Bridgeland unstable for all stability conditions $(cH,bH)$ in $B^{un}_H$ and we obtain the HN filtration \eqref{HN-Th}; each intermediate object $E^{c,b}_i$ of the filtration is a destabilizing object for $\TET[1]$ with respect to the central charge $Z_{cH,bH}$.

Before going further with the investigation of the destabilizing objects 
$E^{c,b}_i$ we return to the cohomology $H^1 (\TET)$ to deduce the filtration stated in Theorem \ref{th:H1-filt-intro}.

\begin{cor}\label{cor:=thH1Thfilt}
	The HN filtration \eqref{HN-Th} induces the filtration of $H^1 (\TET)$
	$$
	\begin{gathered}
H^1 (\TET)=Hom_{\AC_{H,bH}}	(\OO_X, \TET[1])=Hom_{\AC_{H,bH}}	(\OO_X, E^{c,b}_n) \supset
\\
 Hom_{\AC_{H,bH}}	(\OO_X, E^{c,b}_{n-1})\supset \cdots \supset Hom_{\AC_{H,bH}}	(\OO_X, E^{c,b}_1) \supset Hom_{\AC_{H,bH}}	(\OO_X, E^{c,b}_0)=0,
\end{gathered}
$$
for every $(c,b)\in B^{un}_H$, see \eqref{BH-unst} for the definition of $B^{un}_H$.
\end{cor}
\begin{pf}
	For every nonzero intermediate subobject $E^{c,b}_i$ in the HN filtration \eqref{HN-Th} we have an exact sequence 
	$$
	\xymatrix{
	0\ar[r]& E^{c,b}_i \ar[r]& \TET[1] \ar[r]&F_i [1] \ar[r]&0
},
$$
where $F_i$ is a locally free sheaf in $\FF_{H,bH}$. Applying the homological functor $Hom_{\DD} (\OO_X, \bullet)$ gives the exact sequence
of groups
$$
\xymatrix{
Hom_{\DD}(\OO_X,F_i)\ar[r]&Hom_{\DD}(\OO_X,E^{c,b}_i  )\ar[r]& Hom_{\DD}(\OO_X,\TET[1]  )\ar[r]&Hom_{\DD}(\OO_X,F_i[1]).
}
$$
The first group from the left, $Hom_{\DD}(\OO_X,F_i)$, is zero:
$$
Hom_{\DD}(\OO_X,F_i)=Hom_{\AC}(\OO_X,F_i)=0,
$$
where the last equality because $\OO_X \in \TT_{H,bH}$ and $F_i \in \FF_{H,bH}$ and no nonzero morphisms are allowed from objects in $\TT_{H,bH}$ to objects in $\FF_{H,bH}$. Hence the exact sequence
\begin{equation}\label{i-incl-seq}
\xymatrix{
	0\ar[r]&Hom_{\DD}(\OO_X,E^{c,b}_i  )\ar[r]& Hom_{\DD}(\OO_X,\TET[1]  )\ar[r]&Hom_{\DD}(\OO_X,F_i[1])
}
\end{equation}
This gives the inclusions
\begin{equation}\label{i-incl}
		Hom_{\DD}(\OO_X,E^{c,b}_i  ) \hookrightarrow Hom_{\DD}(\OO_X,\TET[1]  )=H^1 (\TET),
\end{equation}
for every $i\in [1,n-1]$.

Given two consecutive subobjects $E^{c,b}_i$ and $E^{c,b}_{i+1}$ in the HN filtration \eqref{HN-Th} we have the commutative diagram in $\AC_{H,bH}$
$$
\xymatrix@R=12pt@C=12pt{
	&&&0\ar[d]&\\
	&0\ar[d]&&A^{c,b}_{i+1}\ar[d]&\\
	0\ar[r]& E^{c,b}_i \ar[r]\ar[d]& \TET[1] \ar@{=}[d]\ar[r]&F_i [1] \ar[d] \ar[r]&0\\
		0\ar[r]& E^{c,b}_{i+1} \ar[r]\ar[d]& \TET[1] \ar[r]&F_{i+1} [1]\ar[d] \ar[r]&0\\
		&A^{c,b}_{i+1}\ar[d]&&0&\\
		&0&&&	
}
$$
Applying the homological functor $Hom_{\DD}(\OO_X, \bullet)$ gives
$$
\xymatrix@R=12pt@C=12pt{
	0\ar[r]& Hom_{\DD}(\OO_X,E^{c,b}_i) \ar[r]\ar[d]& Hom_{\DD}(\OO_X,\TET[1] )\ar@{=}[d]\ar[r]&Hom_{\DD}(\OO_X,F_i [1]) \ar[d] \\
	0\ar[r]& Hom_{\DD}(\OO_X,E^{c,b}_{i+1}) \ar[r]& Hom_{\DD}(\OO_X,\TET[1]) \ar[r]&Hom_{\DD}(\OO_X,F_{i+1} [1]). 	
}
$$
The two horizontal arrows on the left in each row are monomorphisms. Hence the vertical arrow on the left is a monomorphism as well. This proves that
the inclusions in \eqref{i-incl} are compatible as $i$ varies in $[1,n-1]$, that is, we have the filtration asserted in the corollary.
\end{pf}

We now turn to the study of the destabilizing objects $E^{c,b}_i$ of the filtration \eqref{HN-Th} with an objective to obtain geometric/topological constraints they might impose on $X$ or/and the moduli space of $X$.

\vspace{0.5cm}
\noindent
{\bf \ref{s:HNfilt}.3. The walls $W^H_{E^{c,b}_i} (\TET[1])$.}
We let $E=E^{c,b}_i$ to be one of the intermediate objects of the HN filtration
\eqref{HN-Th}. It gives rise to the destabilizing exact sequence
\begin{equation}\label{E-FE-dest}
\xymatrix@R=12pt@C=12pt{
  0\ar[r]&E\ar[r]&\TET[1] \ar[r]& F_E[1]\ar[r]&0
}
\end{equation}
in $\AC_{H,bH}$, where $F_E$ is a locally free sheaf in $\FF_{H,bH}$, see Lemma \ref{lem:FE-locfree}. The destabilizing condition is recorded by the inequality of phases
$$
\phi_{cH,bH} (\TET[1]) >\phi_{cH,bH} (F_E[1]).
$$
This is equivalent to the inequality between the values of the cotangent:
$$
\frac{\Re(Z_{cH,bH}(\TET[1]))}{\Im(Z_{cH,bH}(\TET[1]))}=cot(\pi\phi_{cH,bH} (\TET[1])) <cot(\pi\phi_{cH,bH} (F_E[1]))=\frac{\Re(Z_{cH,bH}(F_E[1]))}{\Im(Z_{cH,bH}(F_E[1]))}.
$$
Using the formula \eqref{Zcb(a)} we obtain
$$
\frac{\HA\tau_X + K_X \cdot H b +(b^2-c^2)H^2}{K_X \cdot H +2H^2 b} <
\frac{-ch_2 (F_E) +c_1(F_E)\cdot H b -\HA rk(F_E)H^2(b^2-c^2)}{c_1(F_E)\cdot H - rk(F_E)H^2b}.
$$
This is of course similar to the considerations in the proof of Lemma \ref{lem:wallWvw} with the only difference that instead of equality we have an inequality.
Thus we obtain a `determinental' inequality
$$
\frac{1}{(K_X \cdot H +2H^2 b)(c_1(F_E)\cdot H - rk(F_E)H^2b)}
\begin{vmatrix}b^2+c^2& \HA\tau_X &ch_2 (F_E)\\
  2b& -K_X \cdot H& c_1 (F_E)\cdot H\\
  \frac{2}{H^2}& 2&rk(F_E)
\end{vmatrix} <0.
$$
Remembering that the multiple in front of the determinant is negative:
$(K_X \cdot H +2H^2 b)>0$ and $(c_1(F_E)\cdot H - rk(F_E)H^2b) <0$ because $\TET$ and $F_E$ are in $\FF_{H,bH}$. The above gives 
\begin{equation}\label{detineq}
\begin{vmatrix}b^2+c^2& \HA\tau_X &ch_2 (F_E)\\
  2b& -K_X \cdot H& c_1 (F_E)\cdot H\\
  \frac{2}{H^2}& 2&rk(F_E)
\end{vmatrix} >0.
\end{equation}
The coefficient of $(b^2 +c^2)$ in the expansion of the determinant with respect to the first column is the $(1,1)$-minor $M_{11} (F_E)$. It has the form
$$
M_{11} (F_E)=-(rk(F_E)K_X\cdot H+2c_1(F_E) \cdot H)
$$
and we need to know the sign of this expression.
\begin{lem}\label{lem:M11negative}
 $rk(F_E)K_X\cdot H+2c_1(F_E) \cdot H>0.$
\end{lem}
\begin{pf}
Write
\begin{equation}\label{-minor}
\begin{gathered}
  rk(F_E)K_X\cdot H+2c_1(F_E) \cdot H=rk(F_E)(K_X\cdot H +2bH^2) +2(c_1(F_E) \cdot H -brk(F_E) H^2)
  \\
=rk(F_E)\Im(Z_{H,bH}(\TET[1])) -2\Im(Z_{H,bH}(F_E[1]))\\
=
(rk(F_E)-2)\Im(Z_{H,bH}(\TET[1])) +2\Im(Z_{H,bH}(E)).
\end{gathered}
\end{equation}
If the rank $rk(F_E) \geq 2$, the last expression is positive since $\TET[1]$ and $E$ are objects in $\AC_{H,bH}$ and $\TET[1]$ has no subobjects of phase $1$.

For $rk(F_E)=1$ the expression we started with in \eqref{-minor} has the form
$$
K_X\cdot H+2c_1(F_E)\cdot H 
$$
and its positivity follows from $H$-stability of $\TET$. To see this apply the homological functor $\HH^0$ to the destabilizing sequence \eqref{E-FE-dest} to obtain the exact complex
$$
\xymatrix@R=12pt@C=12pt{
  0\ar[r]&\HH^{-1}(E)\ar[r]&\TET \ar[r]& F_E \ar[r]&\HH^{0}(E)\ar[r]&0.
}
$$
Since the rank of $F_E$ is $1$ the sheaf $\HH^{-1}(E)$ is a line subsheaf
of $\TET$ and its quotient injects into $F_E$. This implies that the divisor class
$$
c_1(F_E)-c_1(\TET /\HH^{-1}(E))
$$
is effective. Hence
$$
c_1(F_E)\cdot H \geq c_1(\TET /\HH^{-1}(E))\cdot H >\HA c_1(\TET)\cdot H =-\HA K_X \cdot H,
$$
where the second inequality uses $H$-stability of $\TET$. Thus we obtain
$$
c_1(F_E)\cdot H >-\HA K_X \cdot H
$$
or, equivalently, $(K_X \cdot H+2c_1(F_E)\cdot H )>0$.
\end{pf}

With the knowledge of sign of the minor $M_{11}(F_E)$ the inequality \eqref{detineq} becomes
$$
(b^2+c^2) +2b\frac{\HA rk(F_E)\tau_X -2ch_2(F_E)}{rk(F_E)K_X \cdot H+ 2c_1(F_E)\cdot H)} - \frac{ \tau_X c_1(F_E)\cdot H+2K_X \cdot H ch_2(F_E)}{H^2(rk(F_E)K_X \cdot H+ 2c_1(F_E)\cdot H)} <0.
$$
This can be formulated as follows.
\begin{lem}\label{lem:destwalls}
  Fix a stability condition $(c_0H,b_0H)$ in $B^{un}_H$. Every intermediate subobject $E=E^{b_0,c_0}_i$ in the HN filtration of $\TET[1]$ in \eqref{HN-Th} with respect to the central charge $Z_{c_0H,b_0H}$ gives rise to an exact sequence
  \begin{equation}\label{E-FE-seq-lem}
  \xymatrix@R=12pt@C=12pt{
    0\ar[r]&E\ar[r]&\TET[1] \ar[r]& F_E[1]\ar[r]&0
  }
  \end{equation}
  in $\AC_{H,b_0H}$ and to the semicircular wall $W_E(\TET[1])=W_{F_E[1]}(\TET[1])$
  in the half plane $\Pi_H$ with the equation
  $$
  (b -b_{E,H})^2+c^2=r^2_{E,H}
  $$
  where $b_{E,H}$, the $b$-coordinate of the center, and $r_{E,H}$, the radius of the semicircle, are determined by the formulas
  $$
  \begin{gathered}
    b_{E,H}=-\frac{\HA rk(F_E)\tau_X -2ch_2(F_E)}{rk(F_E)K_X \cdot H+ 2c_1(F_E)\cdot H},
    \\
    \\
    r^2_{E,H}= b^2_{E,H}  + \frac{ \tau_X c_1(F_E)\cdot H+2K_X \cdot H ch_2(F_E)}{H^2(rk(F_E)K_X \cdot H+ 2c_1(F_E)\cdot H)}.
  \end{gathered}
  $$
  Furthermore, the exact sequence \eqref{E-FE-seq-lem} is destabilizing for the stability conditions $(cH,bH)$ with $(b,c)$ lying inside the semicircle $W_E(\TET[1])$, that is, for $(b,c)$ subject to the inequality
  $$
  (b -b_{E,H})^2+c^2<r^2_{E,H} .
  $$
  In particular, the semicircles $W_E(\TET[1])$ and $W_{\OO_X}(\TET[1])$ are nested.
\end{lem}
\begin{pf}
  Only the last assertion perhaps needs a proof. From the structure of semicircular walls, see Corollary \ref{cor:semicirwalls}, 3), it follows that it is enough to check that both walls, $W_E(\TET[1])$ and $W_{\OO_X}(\TET[1])$, lie on the same side of the vertical wall for $\TET[1]$. This is obvious, since the point
  $(b_0,c_0)$ lies inside both of them.
\end{pf}

We remind that the region $B^{un}_H$ is the part of the half plane $\Pi_H$
defined by
$$
B^{un}_H =\{ (cH,bH)\in \Pi_H | \text{$(b,c)$ is inside the semicircle $W_{\OO_X}(\TET[1])$ and $b\in (\beta_H,0)$} \},
  $$
  where $\beta_H$ is the negative constant defined in Corollary \ref{cor:ThXBrunst}. In particular, the stability condition $(c_0 H, b_0 H)$ in the above lemma can be chosen so that $b_0$ comes as close to zero as we wish. Such a choice will guarantee that the wall $W_{\OO_X}(\TET[1])$ is nested {\it inside} the walls $W_E(\TET[1])$ of Lemma \ref{lem:destwalls}. 

  \begin{lem}\label{lem:beta}
 There exists $\beta \in (\beta_H,0)$ such that for any stability condition
    $(c_0 H, b_0 H)$ in $B^{un}_H$ with $b_0 >\beta$, the wall
    $W_{\OO_X}(\TET[1])$ is inside or equals the wall $W_E(\TET[1])$, for every intermediate subobject $E=E^{b_0,c_0}_i$ of the HN filtration of $\TET[1]$ with respect to the central charge $Z_{c_0H,b_0H}$.
  \end{lem}
  \begin{pf}
    From Lemma \ref{lem:destwalls} we know that the walls  $W_{\OO_X}(\TET[1])$ and $W_E(\TET[1])$ are nested. So the assertion of the lemma comes down to checking that for  $b_0$ sufficiently close to $0$ (from the left), the radius of the wall $W_E(\TET[1])$ is bigger (or equal to) than the radius of $W_{\OO_X}(\TET[1])$. This in turn means that the `top' of the semicircle $W_E(\TET[1])$ is at least as high as the `top' of $W_{\OO_X}(\TET[1])$. From Proposition \ref{pro:Th-walls}, 2), this is determined by the positions of their centers on the $b$-axis: the center $W_E(\TET[1])$ is to the right of the center of $W_{\OO_X}(\TET[1])$. We know that the center of $W_{\OO_X}(\TET[1])$ is at $\beta_0=-\frac{\tau_X}{2K_X \cdot H}$, see Proposition \ref{pro:semicircH}. This a rational number, so by Lemma \ref{lem:awall-finite} there is at most finite number of actual walls for $\TET[1]$ intersecting the vertical line
    $b=-\frac{\tau_X}{2K_X \cdot H}$. In particular, there is at most finite number of actual walls intersecting that line and nested inside of the semicircle
    $W_{\OO_X}(\TET[1])$. Let $b_m$ and $r_m$ be respectively the $b$-coordinate of the center and the radius of the largest among these `inside' semicircles. Set $\beta'=b_m +r_m$. Thus we have
    \begin{equation}\label{inside-awalls}
      \begin{gathered}
      \text{all actual walls for $\TET[1]$ with centers to the left of the vertical line $b=-\frac{\tau_X}{2K_X \cdot H}$}\\
      \text{lie to the left of the vertical line
        $ b=\beta'$.}
      \end{gathered}
        \end{equation}
        The constant
        $$
        \beta=max \{\beta_H,\beta'\}
        $$
        has the properties asserted in the lemma: for any $(c_0H,b_0H)$ in $B^{un}_H$ with $b_0 >\beta$, every wall
$W_E(\TET[1])$ under consideration is actual for $\TET[1]$ and the point $(b_0, c_0)$ is inside this wall; thus this wall is not to the left of  the vertical line
$ b=\beta'$ and hence by \eqref{inside-awalls}  the center of $W_E(\TET[1])$ is located on or to the right of the vertical line $b=-\frac{\tau_X}{2K_X \cdot H}$.
\end{pf}
\begin{lem}\label{lem:filt}
Let $\beta^{in}_H$ be the infimum of all values $\beta$ in Lemma \ref{lem:beta}.
Then for all stability conditions $(c_0 H,b_0H)$ in $B^{un}_H$ with $b_0 >\beta^{in}_H$ the HN filtration of $\TET[1]$
$$
\TET[1]=E^{c_0,b_0}_n \supset E^{c_0,b_0}_{n-1}\supset \cdots \supset E^{c_0,b_0}_{1}\supset E^{c_0,b_0}_{0}=0
$$
with respect to the central charge $Z_{c_0H,b_0H}$ has the following properties.

1) Every intermediate object $E^{c_0,b_0}_i$ of the filtration gives rise to the semicircular wall $W_{E^{c_0,b_0}_i} (\TET[1])$ with the wall $W_{\OO_X} (\TET[1])$
equals or nested inside it.

2) The radii $r_{E^{c_0,b_0}_i} $ and $r_{E^{c_0,b_0}_{i-1}}$ of the successive walls $W_{E^{c_0,b_0}_i} (\TET[1])$ and $W_{E^{c_0,b_0}_{i-1}} (\TET[1])$ are subject to
$$
r_{E^{c_0,b_0}_i} < r_{E^{c_0,b_0}_{i-1}},
$$
for $i\in [2,n-1]$, if and only if the slope inequality
$$
\frac{c_1(F_i)\cdot H}{rk(F_i)}>\frac{c_1(F_{i-1})\cdot H}{rk(F_{i-1})}
$$
holds. Here $F_i[1]$ (resp. $F_{i-1}[1]$) denotes the quotient object
$\TET[1]/E^{c_0,b_0}_i$ (resp. $\TET[1]/E^{c_0,b_0}_{i-1}$).
\end{lem}
\begin{pf}
  The assertion 1) is already known from Lemma \ref{lem:beta}. So only part 2) needs to be proved. Assume that for some $i\in [2,n-1]$ we have the inequality
  \begin{equation}\label{radii}
  r_{E^{c_0,b_0}_i} < r_{E^{c_0,b_0}_{i-1}}.
  \end{equation}
  To simplify the notation we omit the superscript and simply write $E_i$ instead of $E^{c_0,b_0}_i $; the corresponding quotient object $\TET[1]/E_i$ will be denoted $F_i [1]$, for some locally free sheaf $F_i \in \FF_{H,b_0H}$.
  
  The inequality \eqref{radii} means that the wall $W_{E_i} (\TET[1])$ is nested inside the wall $W_{E_{i-1}} (\TET[1])$. This tells us that for any
  point $(b,c)$ between those two walls we have the following relation of phases:

  $$
  \phi_{cH,bH} (F_i [1]) >\phi_{cH,bH} (\TET [1]) >\phi_{cH,bH} (F_{i-1} [1])
  $$
  with respect to the central charge $Z_{cH,bH}$. Taking $b=b_0$ gives the inequality
  \begin{equation}\label{i-i-1}
  \phi_{cH,b_0H} (F_i [1]) > \phi_{cH,b_0H} (F_{i-1} [1]),
  \end{equation}
  for all $(c,b_0)$ in the segment on the vertical line $b=b_0$ between the semicircles $W_{E_i} (\TET[1])$ and $W_{E_{i-1}} (\TET[1])$. On the other hand from the HN filtration of $\TET[1]$ the successive quotient objects $F_{i-1} [1]$ and $F_i [1]$ are related by the exact sequence
  $$
 \xymatrix@R=12pt@C=12pt{
   0\ar[r]&A_i \ar[r]&F_{i-1} [1]\ar[r]&F_i [1]\ar[r]&0,
 }
 $$
 where $A_i=E_i/E_{i-1}$. By the property of phases in the HN filtration it follows
\begin{equation}\label{i-i-1-c0} 
 \phi_{c_0H,b_0H} (F_i [1])<\phi_{c_0H,b_0H} (F_{i-1} [1]).
 \end{equation}
 Observe: all values $c$ in \eqref{i-i-1} are bigger than $c_0$.

 We will now rewrite the inequalities of phases in terms of the values of cotangents:
 $$
 \frac{\Re(Z_{cH,b_0H} (F_i [1]))}{\Im(Z_{cH,b_0H} (F_i [1]))}= cot(\pi \phi_{cH,b_0H} (F_i [1])) < cot(\pi \phi_{cH,b_0H} (F_{i-1} [1]))=\frac{\Re(Z_{cH,b_0H} (F_{i-1} [1]))}{\Im(Z_{cH,b_0H} (F_{i-1} [1]))}.
 $$
 The formula of the central charge gives the inequality:
 $$
 \HA H^2 \Big(\frac{ch_0(F_{i} [1])}{\Im(Z_{H,b_0H} (F_{i} [1]))} -\frac{ch_0(F_{i-1} [1])}{\Im(Z_{H,b_0H} (F_{i-1} [1]))} \Big)c^2 <T,
 $$
 where $T$ is an expression depending only on $b_0$ and the Chern characters of
 $F_{i-1} [1] $ and $F_{i} [1] $. Doing the same for the phase inequality \eqref{i-i-1-c0} gives the inequality
 $$
 T<\HA H^2 \Big(\frac{ch_0(F_{i} [1])}{\Im(Z_{H,b_0H} (F_{i} [1]))} -\frac{ch_0(F_{i-1} [1])}{\Im(Z_{H,b_0H} (F_{i-1} [1]))} \Big)c^2_0.
 $$
 Since $c>c_0 >0$ we deduce
 $$
 \frac{ch_0(F_{i-1} [1])}{\Im(Z_{H,b_0H} (F_{i-1} [1]))} -\frac{ch_0(F_{i} [1])}{\Im(Z_{H,b_0H} (F_{i} [1])) }>0.
 $$
 The objects involved in the inequality are both in $\AC_{H,b_0H}$ of phase less than $1$. So the denominators are positive and the above inequality simplifies to
 $$
ch_0(F_{i-1} [1])ch_1(F_{i} [1])\cdot H >ch_0(F_{i} [1])ch_1(F_{i-1} [1])\cdot H.
 $$
 This is equivalent to
  $$
rk(F_{i-1} )c_1(F_{i} )\cdot H >rk(F_{i})c_1(F_{i-1} )\cdot H
$$
or else to the inequality of slopes
$$
\frac{c_1(F_{i} )\cdot H }{rk(F_i)} > \frac{c_1(F_{i-1} )\cdot H }{rk(F_{i-1})}.
$$
\end{pf}
\begin{cor}\label{cor:tau-bound}
  Let $(c_0H,b_0H)$ be a stability condition in $B^{un}_H$ with $b_0 > \beta^{in}_H$, where $\beta^{in}_H$ is defined in Lemma \ref{lem:filt}. For every intermediate subobject $E^{c_0,b_0}_i$ of the HN filtration of $\TET[1]$ with respect to the central charge $Z_{c_0H,b_0H}$, the quotient object
  $F_i[1]=\TET[1]/E^{c_0,b_0}_i$ is subject to the inequality
  $$
      \tau_X 
      \leq \frac{2ch_2 (F_i)K_X \cdot H }{-c_1 (F_i)\cdot H}.
   $$
   In addition, one has the inequality
   $$
    \tau_X 
    \leq \frac{(-c_1 (F_i)\cdot H) K_X \cdot H}{rk(F_i)H^2},  
  $$
  unless the sheaf $F_i$ is Bogomolov unstable. If this is the case, then
  the HN filtration of $F_i$ with respect to $H$
  $$
  F_i =F^{l_i}_i \supset F^{l_i -1}_i \supset \cdots \supset F^{1}_i \supset F^{0}_i =0
  $$
  has all the semistable factors $Q^j_i=F^{j}_i / F^{j -1}_i$, for $j\in[1,l_i]$, in the subcategory $\FF_{H,b_0H}$ and one has the inequality
  \begin{equation}\label{taubound-cor}
  \tau_X 
  < \frac{(-c_1 (Q^{l_i}_i)\cdot H)}{rk(Q^{l_i}_i)H^2} K_X \cdot H <  \frac{(-c_1 (F_i)\cdot H)}{rk(Q^{l_i}_i)H^2} K_X \cdot H.  
\end{equation}
\end{cor}
\begin{pf}
  According to Lemma \ref{lem:filt} the wall $W^H_{\OO_X} (\TET[1])$ is nested inside or equals the wall $W^H_{E^{c_0,b_0}_i}$. In particular, the $b$-coordinate $\Big(-\frac{\tau_X}{2K_X\cdot H}\Big)$ of the center of $W^H_{\OO_X} (\TET[1])$ is to the left of the one of the center of $W^H_{E^{c_0,b_0}_i}$. In other words we have the inequality
  $$
  -\frac{\tau_X}{2K_X\cdot H}\leq b_{E^{c_0,b_0}_i,H} =-\frac{\HA rk(F_i)\tau_X -2ch_2 (F_i)}{rk(F_i)K_X \cdot H +2c_1 (F_i)\cdot H},
  $$
  where the last equality is the formula in Lemma \ref{lem:destwalls}. The denominator on the right is positive, see Lemma \ref{lem:M11negative}, so the above inequality simplifies to
  $$
  2K_X \cdot H ch_2(F_i)\geq -c_1(F_i)\cdot H \tau_X
  $$
  or, equivalently, since $(-c_1(F_i)\cdot H)$ is positive, we obtain the first  inequality of the corollary
  $$
  \tau_X \leq \frac{2ch_2(F_i)K_X \cdot H }{-c_1(F_i)\cdot H}.
  $$
  
   If, in addition, the quotient object $F_i[1]$ satisfies the inequality
  \begin{equation}\label{BGFi}
  (c_1 (F_i)\cdot H)^2 \geq 2H^2 rk(F_i)ch_2 (F_i),
\end{equation}
 then combining it with the one just proved gives
  $$
  \tau_X \leq \frac{(-c_1(F_i)\cdot H)K_X \cdot H}{rk(F_i)H^2}.
  $$
  Recall: the failure of the inequality \eqref{BGFi} means that $F_i$ is Bogomolov unstable. In particular, it is $H$-unstable and we have
   $$
  F_i =F^{l_i}_i \supset F^{l_i -1}_i \supset \cdots \supset F^{1}_i \supset F^{0}_i =0
  $$
  its HN filtration with respect to $H$. Since $F_i$ is in $\FF_{H,b_0H}$
  all its semistable factors $Q^j_i=F^{j}_i / F^{j -1}_i$, for $j\in[1,l_i]$, are in the subcategory $\FF_{H,b_0H}$ as well. Hence the inequality
  $$
  c_1(Q^j_i)\cdot H \leq rk(Q^j_i)H^2 b_0,
  $$
  for all $j\in[1,l_i]$. This gives
  $$
  -c_1(Q^j_i)\cdot H \geq rk(Q^j_i)H^2(- b_0)>0, \forall j\in [1,l_i],
  $$
  where the last inequality uses the fact that $b_0$ is negative. The positivity above allows to apply Lemma \ref{lem:BGunst} to the dual $F^{\ast}_i$. From the proof of that lemma we obtain
  $$
  2ch_2 (F_i)=2ch_2 (F^{\ast}_i)< 2ch_2 \Big((Q^{l_i}_i)^{\ast}\Big)+\frac{c_1\Big((Q^{l_i}_i)^{\ast}\Big)\cdot H}{rk(Q^{l_i}_i)H^2}(c_1 (F^{\ast}_i)\cdot H - c_1((Q^{l_i}_i)^{\ast})\cdot H).
  $$
  This gives 
  $$
  2rk(Q^{l_i}_i)H^2 ch_2 (F_i) <-\Delta^H((Q^{l_i}_i)^{\ast}) + \Big(c_1\Big((Q^{l_i}_i)^{\ast}\Big)\cdot H\Big) \Big(c_1 (F^{\ast}_i)\cdot H \Big) \leq \Big(c_1\Big((Q^{l_i}_i)^{\ast}\Big)\cdot H\Big) \Big(c_1 (F^{\ast}_i)\cdot H \Big),
  $$
  where the last inequality uses the $H$-Bogomolov-Gieseker inequality for $(Q^{l_i}_i)^{\ast}$. Thus we obtain
  $$
   2rk(Q^{l_i}_i)H^2 ch_2 (F_i) < \Big(c_1\Big((Q^{l_i}_i)^{\ast}\Big)\cdot H\Big) \Big(c_1 (F^{\ast}_i)\cdot H \Big).
   $$
   This together with the first inequality of the lemma give
   \begin{equation}\label{ubtauunst}
   \tau_X < \frac{c_1\Big((Q^{l_i}_i)^{\ast}\Big)\cdot H}{rk(Q^{l_i}_i)H^2} K_X \cdot H.
   \end{equation}
   From the relation
   $$
   -c_1(F_i)\cdot H =c_1(F^{\ast}_i)\cdot H=\sum^{l_i}_{j=1} c_1 \Big((Q^{j}_i)^{\ast}\Big)\cdot H
   $$
   and the positivity of the summands we deduce
   $$
  -c_1(F_i)\cdot H > c_1\Big((Q^{l_i}_i)^{\ast}\Big)\cdot H=-c_1(Q^{l_i}_i)\cdot H.
  $$
  Substituting into \eqref{ubtauunst} gives
  $$
  \tau_X < \frac{(-c_1(F_i)\cdot H)}{rk(Q^{l_i}_i)H^2} K_X \cdot H.
  $$
  \end{pf}
We can now start drawing further conclusions about properties of objects involved in the
$HN$ filtration of $\TET[1]$.
\begin{pro}\label{pro:no1rankbound}
Let $(c_0H,b_0H)$ be a stability condition chosen as in Corollary \ref{cor:tau-bound}. Then the quotient objects $F_i [1]=\TET[1]/E^{c_0,b_0}_i$, for all $i\in [1,n-1]$, have the following properties.

\vspace{0.2cm}
1) $rk(F_i) \geq 2$.

\vspace{0.2cm}
2) $(-c_1(F_i)\cdot H )\leq \begin{cases} \HA K_X \cdot H ,&\text{if $\HH^{-1} (E^{c_0,b_0}_i) \neq 0$,}
  \\
  \\
 K_X \cdot H ,&\text{if $\HH^{-1} (E^{c_0,b_0}_i) = 0$.}
\end{cases}
$

\vspace{0.2cm}
3) $ch_1(E^{c_0,b_0}_i)\cdot H  \geq 0$, for all $i\in[1,n-1]$.
\end{pro}
\begin{pf}
  For the first assertion, assume $rk(F_i) =1$. Then $F_i$ is a line bundle and we write it $\OO_X (-L)$ for some divisor $L$. Since $\OO_X (-L)$ lies in $\FF_{H,b_0H}$, we have
  $$
  -L\cdot H\leq b_0 H^2
  $$
  Hence the inequality
  \begin{equation}\label{LHpos}
    L\cdot H \geq -b_0 H^2 >0.
  \end{equation}
  From the exact sequence
  $$
 \xymatrix@R=12pt@C=12pt{
   0\ar[r]&E^{c_0,b_0}_i \ar[r]&\TET[1]\ar[r]&\OO_X (-L)[1]\ar[r]&0
   }
   $$
   in $\AC_{H,b_0H}$ we obtain, by applying the homological functor
   $\HH^0$, the exact complex of sheaves
   $$
 \xymatrix@R=12pt@C=12pt{
   0\ar[r]&\HH^{-1}(E^{c_0,b_0}_i )\ar[r]&\TET\ar[r]& \OO_X (-L)\ar[r]&\HH^{0}(E^{c_0,b_0}_i )\ar[r]&0
 }
 $$
 From this it follows that the morphism $\TET \longrightarrow \OO_X (-L)$
 is generically surjective. Dualizing gives a monomorphism
 $$
 \xymatrix@R=12pt@C=12pt{
   0\ar[r]&\OO_X (L)\ar[r]&\Omega_X
   }
   $$
   But the degree $L\cdot H$ is positive, see \eqref{LHpos}, and since $L$ can not lie in the positive cone of $X$, we deduce
   $$
   L^2 \leq 0.
   $$
   This and the first inequality in Corollary \ref{cor:tau-bound} applied to
   $F_i=\OO_X (-L)$ give
   $$
   \tau_X \leq \frac{L^2 (K_X \cdot H)}{L\cdot H} \leq 0.
   $$
   This is contrary to the assumption that $\tau_X$ is positive.
   
   \vspace{0.2cm}

   We now turn to the second assertion. For this we consider the destabilizing sequence
    $$
 \xymatrix@R=12pt@C=12pt{
   0\ar[r]&E^{c_0,b_0}_i \ar[r]&\TET[1]\ar[r]&F_i[1]\ar[r]&0
   }
   $$
   in $\AC_{H,b_0H}$, apply to it the homological functor $\HH^0$ to obtain the complex
   \begin{equation}\label{cpxFi}
\xymatrix@R=12pt@C=12pt{
   0\ar[r]&\HH^{-1}(E^{c_0,b_0}_i )\ar[r]&\TET\ar[r]& F_i \ar[r]&\HH^{0}(E^{c_0,b_0}_i )\ar[r]&0.
 }
\end{equation}
in $\AC$.

If $\HH^{-1}(E^{c_0,b_0}_i )$ is nonzero, it is a line bundle which we denote $\OO_X (-M)$ and the above complex breaks into two short exact sequences
$$
\xymatrix@R=12pt@C=12pt{
  0\ar[r]&\OO_X (-M)\ar[r]&\TET\ar[r]&{\cal J}_A (-L)\ar[r]&0,
  \\
 0\ar[r]& {\cal J}_A (-L)\ar[r]&  F_i \ar[r]&\HH^{0}(E^{c_0,b_0}_i )\ar[r]&0.
  }
  $$
  The second exact sequence implies
  $$
  -c_1( F_i )\cdot H=L\cdot H -c_1(\HH^{0}(E^{c_0,b_0}_i ))\cdot H <L\cdot H -
  rk (\HH^{0}(E^{c_0,b_0}_i )b_0 H^2=L\cdot H -(rk(F_i)-1)b_0H^2,
  $$
  where the inequality comes from the fact that $\HH^{0}(E^{c_0,b_0}_i )$ is a sheaf in $\TT_{H,b_0H}$. From the first exact sequence and $H$-stability of $\TET$ we deduce
  $$
  L\cdot H< \HA K_X \cdot H.
  $$
  Substituting into the previous inequality we obtain
  $$
  -c_1( F_i )\cdot H <\HA K_X \cdot H -(rk(F_i)-1)b_0H^2.
  $$

  If $\HH^{-1}(E^{c_0,b_0}_i ) =0$, then $E^{c_0,b_0}_i$ is a sheaf in $\TT_{H,b_0H}$ and the  complex \eqref{cpxFi} becomes the short exact sequence
  $$
 \xymatrix@R=12pt@C=12pt{
   0\ar[r]&\TET\ar[r]& F_i \ar[r]&E^{c_0,b_0}_i \ar[r]&0.
 }
 $$
 and we obtain
 $$
 -c_1( F_i )\cdot H=K_X \cdot H -c_1(\HH^{0}(E^{c_0,b_0}_i ))\cdot H <K_X \cdot H -
 rk (\HH^{0}(E^{c_0,b_0}_i )b_0 H^2=K_X\cdot H -(rk(F_i)-2)b_0H^2,
 $$
 where again the inequality uses that $\HH^{0}(E^{c_0,b_0}_i)$ is a sheaf in $\TT_{H,b_0H}$. 
 
 Putting the two cases together gives the inequality
 $$
 -c_1( F_i )\cdot H <\begin{cases}
 	\HA K_X \cdot H -(rk(F_i)-1)b_0H^2,& \text{if $\HH^{-1}(E^{c_0,b_0}_i )\neq 0$,}\\
 K_X \cdot H -(rk(F_i)-2)b_0H^2,& \text{if $\HH^{-1}(E^{c_0,b_0}_i )=0$.}	
 	\end{cases} 
 $$
 We now need to eliminate the dependence on $b_0$ in the above inequality. For this we recall that the subcategory
 $\FF_{H,b_0 H}$ is contained in $\FF_{H,bH}$, for all $b\geq b_0$. This implies that we can replace $b_0$ by any $b\in [b_0,0)$. Furthermore, by Lemma \ref{lem:awall-finite} the number of possible values for the rank $rk(F_i)$ is finite. Hence letting $b_0$ to converge to $0$, we obtain the asserted inequality
  $$
 -c_1( F_i )\cdot H \leq \begin{cases}
 	\HA K_X \cdot H ,& \text{if $\HH^{-1}(E^{c_0,b_0}_i )\neq 0$,}\\
 	K_X \cdot H,& \text{if $\HH^{-1}(E^{c_0,b_0}_i )=0$.}	
 \end{cases} 
 $$
 
 The inequality above implies the third assertion of the proposition -
 the non-negativity of $ch_1 (E^{c_0,b_0}_i)\cdot H$. Indeed, from the destabilizing sequence
 $$
 \xymatrix@R=12pt@C=12pt{
 	0\ar[r]&E^{c_0,b_0}_i \ar[r]&\TET[1]\ar[r]& F_i[1] \ar[r]&0
 }
 $$
 it follows
 $$
 K_X=ch_1 (\TET[1])=ch_1 (E^{c_0,b_0}_i) +ch_1(F_i[1])=ch_1 (E^{c_0,b_0}_i)+(-c_1(F_i)).
 $$
 Intersecting with $H$ and the inequality for $(-c_1(F_i)\cdot H)$ just proved give
 $$
 K_X \cdot H =ch_1 (E^{c_0,b_0}_i)\cdot H +(-c_1(F_i))\cdot H \leq
 ch_1 (E^{c_0,b_0}_i)\cdot H +K_X \cdot H.
 $$
 Hence the inequality
 $$
 ch_1 (E^{c_0,b_0}_i)\cdot H \geq 0.
 $$
\end{pf}
Combining the inequality in Proposition \ref{pro:no1rankbound}, 2), with the second
inequality in Corollary \ref{cor:tau-bound} gives the following.
\begin{cor}\label{cor:tau-bound1}
  With the notation of Proposition \ref{pro:no1rankbound}, for every $i \in [1,n-1]$ for which $F_i$ is Bogomolov semistable, one has the inequality
  $$
  \tau_X \leq \begin{cases} \frac{(K_X \cdot H)^2}{2rk(F_i)H^2},&\text{ if $\HH^{-1}(E^{c_0,b_0}_i ) \neq 0$,}
    \\
    \\
 \frac{(K_X \cdot H)^2}{rk(F_i)H^2}, &\text{ if $\HH^{-1}(E^{c_0,b_0}_i )= 0$.}   
\end{cases}
$$

For every $F_i$ which is Bogomolov unstable, let $Q^{l_i}_i$ be the last
$H$-semistable quotient of the HN filtration of $F_i$ with respect $H$,
see Proposition \ref{pro:no1rankbound} for notation. Then the following inequality holds
$$
\tau_X \leq  \begin{cases} \frac{(K_X \cdot H)^2}{2rk(Q^{l_i}_i)H^2},&\text{ if $\HH^{-1}(E^{c_0,b_0}_i ) \neq 0$,}
	\\
	\\
	\frac{(K_X \cdot H)^2}{rk(Q^{l_i}_i)H^2}, &\text{ if $\HH^{-1}(E^{c_0,b_0}_i )= 0$.}   
\end{cases}
$$ 
\end{cor}

 {\it Proof.} For $F_i$ Bogomolov semistable, we substitute the inequality in Proposition \ref{pro:no1rankbound}, 2), into the second inequality of Corollary \ref{cor:tau-bound} to obtain
  {\small
  $$
  \tau_X \leq \begin{cases} \HA \frac{(K_X\cdot H)^2}{rk(F_i)H^2}, & \text{ if $\HH^{-1}(E^{c_0,b_0}_i ) \neq 0$,}
    \\
    \\
    \frac{(K_X\cdot H)^2}{rk(F_i)H^2}, & \text{ if $\HH^{-1}(E^{c_0,b_0}_i ) = 0$.}
    \end{cases}
    $$
    \small}
 
 If $F_i$ is Bogomolov unstable, then we use the last inequality
  in Corollary \ref{cor:tau-bound}
  $$
  \tau_X <\frac{(-c_1 (F_i)\cdot H)}{rk(Q^{l_i})H^2} K_X \cdot H.
  $$
  From here on we argue as in the first part of the proof to obtain
  $$
  \tau_X <  \begin{cases} \frac{(K_X \cdot H)^2}{2rk(Q^{l_i}_i)H^2},&\text{ if $\HH^{-1}(E^{c_0,b_0}_i ) \neq 0$,}
  	\\
  	\\
  	\frac{(K_X \cdot H)^2}{rk(Q^{l_i}_i)H^2}, &\text{ if $\HH^{-1}(E^{c_0,b_0}_i )= 0$.}   
  \end{cases} 
  $$ 
\begin{flushright}
	$\Box$
\end{flushright} 
 
\begin{cor}\label{cor:-c1Fi}
With the notation of Corollary \ref{cor:tau-bound} the following holds.

\vspace{0.2cm}
1) For every object $F_i [1]$ which is Bogomolov semistable one has
$$
c^2_1 (F_i) > 0.
$$
In particular, the divisor class $(-c_1 (F_i))$ lies in the positive cone $C^+(X)$.

\vspace{0.2cm}
2) For every object $F_i [1]$ which is Bogomolov unstable either
$c^2_1 (Q^{l_i}_i) > 0$ and the divisor class $(-c_1 (Q^{l_i}_i))$ lies in the positive cone of $X$ and $rk(Q^{l_i}_i)\geq 2$, or $c^2_1 (Q^{l_i}_i) \leq 0$ and one has
 inequalities
$$
\begin{gathered}
ch_2 (F_i) < \frac{(c_1 (F_i)\cdot H)^2}{8rk(Q^{l_i}_i)H^2},
\\
\tau_X < \frac{(-c_1(F_i)\cdot H)K_X \cdot H}{4rk(Q^{l_i}_i)H^2}.
\end{gathered}
$$	
	\end{cor}
\begin{pf}
	If $F_i[1]$ is Bogomolov semistable, then we have
	$$
c^2_1 (F_i) \geq 2rk(F_i)ch_2 (F_i)
$$
the Bogomolov-Gieseker inequality for $F_i[1]$. This together with the first inequality in Corollary \ref{cor:tau-bound} imply
$c^2_1 (F_i) >0$. 

Since $F_i$ is a sheaf in $\FF_{H,b_0H}$, we also have
$$
c_1 (F_i)\cdot H \leq rk(F_i)H^2 b_0.
$$
Hence
$$
(-c_1 (F_i)\cdot H) \geq rk(F_i)H^2 (- b_0) >0,
$$
where the last inequality uses the fact that $b_0$ is negative. Thus
$(-c_1 (F_i))$ lies in the positive cone of $X$.

We now turn to the case when $F_i[1]$ is Bogomolov unstable. With the notation of Corollary \ref{cor:tau-bound}, we have the $H$-semistable factors $Q^j_i$ of the HN filtration of $F_i$ with respect to $H$.
They are all sheaves in $\FF_{H,b_0 H}$ and hence
$$
-c_1 (Q^j_i)\cdot H \geq rk(Q^j_i)H^2 (-b_0)> 0, \forall j\in [1,l_i].
$$
If, in addition, $c^2_1 (Q^{l_i}_i)> 0$, we deduce that the divisor class $(-c_1 (Q^{l_i}_i))$ lies in the positive cone of $X$.

We now consider the case $c^2_1 (Q^{l_i}_i)\leq 0$. The Bogomolov-Gieseker inequality for $Q^{l_i}_i$ gives
$$
2rk(Q^{l_i}_i)ch_2 (Q^{l_i}_i) \leq c^2_1 (Q^{l_i}_i)\leq 0.
$$
This and the estimate
$$
2ch_2 (F_i) < ch_2 (Q^{l_i}_i)+ \frac{(-c_1 (Q^{l_i}_i))\cdot H)}{rk(Q^{l_i}_i)H^2} (-c_1 (F_i)\cdot H + c_1 (Q^{l_i}_i) \cdot H)
$$
in the proof of Corollary \ref{cor:tau-bound} give
$$
2ch_2 (F_i) <  \frac{(-c_1 (Q^{l_i}_i))\cdot H)}{rk(Q^{l_i}_i)H^2} (-c_1 (F_i)\cdot H + c_1 (Q^{l_i}_i) \cdot H).
$$ 
This can be rewritten as quadratic relation in 	$t:=c_1 (Q^{l_i}_i) \cdot H$:
$$
t^2 -(c_1 (F_i)\cdot H)t +2rk(Q^{l_i}_i)H^2 ch_2 (F_i) <0.
$$
From this it follows that the discriminant of the quadratic polynomial on the left is positive. This gives the inequality
$$
(c_1 (F_i)\cdot H)^2 -8rk(Q^{l_i}_i)H^2 ch_2 (F_i) >0
$$
which is equivalent to the first inequality in the part 2) of the corollary.
Substituting into the first inequality of Corollary \ref{cor:tau-bound}
we obtain
$$
\tau_X <\frac{(-c_1 (F_i)\cdot H) K_X \cdot H}{4rk(Q^{l_i}_i)H^2}.
$$
\end{pf} 
  
  The inequalities of Corollary \ref{cor:tau-bound1} become particularly meaningful if we take the polarization $H=K_X$. Namely, for $F_i$'s which are Bogomolov semistable, we obtain
  \begin{equation}\label{tau-boundKX}
  \tau_X \leq \begin{cases} \frac{K^2_X}{2rk(F_i)}, & \text{ if $\HH^{-1}(E^{c_0,b_0}_i ) \neq 0$,}
     \\
     \\
    \frac{K^2_X}{rk(F_i)}, & \text{ if $\HH^{-1}(E^{c_0,b_0}_i ) = 0$.}
    \end{cases}
    \end{equation}  
    In terms of the ratio $\displaystyle{\alpha_X=\frac{c_2(X)}{K^2_X}}$ of the Chern numbers of $X$, the above inequality reads
   \begin{equation}\label{alphaX-bound}
1- 2\alpha_X \leq \begin{cases} \frac{1}{2rk(F_i)}, & \text{ if $\HH^{-1}(E^{c_0,b_0}_i ) \neq 0$,}
     \\
     \\
    \frac{1}{rk(F_i)}, & \text{ if $\HH^{-1}(E^{c_0,b_0}_i ) = 0$.}
    \end{cases}
  \end{equation}
For $F_i$'s which are Bogomolov unstable we have
\begin{equation}\label{alphaX-bound_Bogunst}
	1- 2\alpha_X \leq \begin{cases} \frac{1}{2rk(Q^{l_i}_i)}, & \text{ if $\HH^{-1}(E^{c_0,b_0}_i ) \neq 0$,}
		\\
		\\
		\frac{1}{rk(Q^{l_i)}_i}, & \text{ if $\HH^{-1}(E^{c_0,b_0}_i ) = 0$,}
	\end{cases}
\end{equation}
where $Q^{l_i}_i$ is as in Corollary \ref{cor:tau-bound1}.

  The point of the above inequalities is that the properties of HN filtration for the canonical polarization can be controlled by imposing additional restrictions on the ratio $\alpha_X$ of the Chern numbers of $X$.
  \begin{pro}\label{pro:HNalpha}
    Assume the ratio $\alpha_X <\frac{3}{8}$ and let $(c_0K_X,b_0 K_X)$ be a stability condition in $B^{un}_{K_X}$ with $b_0 >\beta^{in}_{K_X}$, where $\beta^{in}_{K_X}$ is defined in Lemma \ref{lem:filt}. Then the HN filtration
    \begin{equation}\label{HNfilt-pro}
    \TET[1]=E^{c_0,b_0}_n \supset E^{c_0,b_0}_{n-1}\supset \cdots \supset E^{c_0,b_0}_{1}\supset E^{c_0,b_0}_{0} =0
    \end{equation}
    of $\TET[1]$ with respect to the central charge $Z_{c_0K_X,b_0 K_X}$ has the following properties.

    1) All objects $E^{c_0,b_0}_{i}$, for $i\in[1,n-1]$, are sheaves.
    In addition, those having the quotient
    $F_i [1]=\TET[1]/E^{c_0,b_0}_{i}$ Bogomolov semistable, are sheaves of rank
    $0$ or $1$ in $\TT_{K_X,b_0 K_X}$.

    2) The sheaves $E^{c_0,b_0}_{i}$'s of rank $0$, that is, the torsion sheaves
    of the filtration are supported on subschemes of pure dimension $1$.

    3) In the quotient objects $F_i[1]=\TET[1]/E^{c_0,b_0}_{i}$, for  $i\in[1,n-1]$, the sheaves $F_i$'s are locally free and either Bogomolov unstable or have rank $2$ or $3$.

    4) All sheaves $F_i$
    of rank $2$ are $K_X$-stable and the divisor class $(-c_1(F_i))$ lies in the positive cone $C^+ (X)$.

    5) If a sheaf $F_i$ has rank $3$ and is Bogomolov semistable but $K_X$-unstable, then its
    HN filtration with respect to $K_X$ has the form
    $$
    F_i=F^2_i \supset F^1_i  \supset F^0_i=0,
    $$
    where $F^1_i$ is a line bundle and the divisor class
    $(-c_1( F_i/F^1_i ))$ lies in the positive cone $C^+ (X)$.
    
    6) If $F_i$ is Bogomolov unstable, then its HN filtration
    $$
    F_i=F^{l_i}_i \supset F^{l_i -1}_i \supset \cdots \supset F^{1}_i \supset F^{0}_i =0
    $$
    with respect to $K_X$ has its semistable factors $Q^j_i =F^j_i /F^{j-1}_i $ all lying in $\FF_{H,b_0H}$. In addition, the last factor $Q^{l_i}_i$ has rank $2$ or $3$ and the divisor class
    $(-c_1(Q^{l_i}_i))$ lies in the positive cone $C^+ (X)$. 
    \end{pro}
\begin{pf}
We begin by showing that all intermediate objects $E^{c_0,b_0}_{i}$ of the HN filtration \eqref{HNfilt-pro} are sheaves, provided the quotient object $F_i [1] =\TET[1]/E^{c_0,b_0}_{i}$ is Bogomolov semistable. This is equivalent to proving that
$\HH^{-1} (E^{c_0,b_0}_{i})$ is zero, for such $i\in [1,n-1]$. But if not, then
the first inequality in \eqref{alphaX-bound} holds and we have
$$
1-2\alpha_X \leq  \frac{1}{2rk(F_i)} \leq \frac{1}{4},
$$
where the second inequality uses the fact that the rank $rk(F_i)$ is at least $2$, see Proposition \ref{pro:no1rankbound}, 1). From this it follows
$$
\alpha_X \geq \frac{3}{8}
$$
and this contrary to the assumption that $\alpha_X < \frac{3}{8}$.

The second inequality in \eqref{alphaX-bound}
$$ 
1-2\alpha_X \leq  \frac{1}{rk(F_i)}
$$
also tells us that $rk(F_i)\leq 3$. This and Proposition \ref{pro:no1rankbound}, 1), imply the assertion that $rk(F_i)=2$ or $3$, for all $F_i$ which are Bogomolov semistable.

If $F_i$ is Bogomolov unstable, then the inequality \eqref{alphaX-bound_Bogunst} tells us that the rank $rk(Q^{l_i}_i)$ is at most $3$. In addition, Corollary \ref{cor:-c1Fi}, 2), implies that
$c^2_1(Q^{l_i}_i) >0$ and hence the divisor class $(-c_1(Q^{l_i}_i))$ lies in the positive cone of $X$ and the rank $rk(Q^{l_i}_i) \geq 2$. The latter and the first line in \eqref{alphaX-bound_Bogunst} imply
that $\HH^{-1}(E^{c_0,b_0}_{i})=0$. The statement that all $E^{c_0,b_0}_{i}$'s are sheaves is now proved. In the course of the proof we also obtained part 6) of the proposition.

For every $i\in [1,n-1]$ we have an exact sequence
$$
\xymatrix@R=12pt@C=12pt{
0\ar[r]&E^{c_0,b_0}_{i}\ar[r]&\TET[1]\ar[r]&F_i[1]\ar[r]&0
}
$$
in $\AC_{K_X,b_0K_X}$. Applying the homological functor $\HH^0$ and remembering that $E^{c_0,b_0}_{i}$ is a sheaf gives the exact sequence
\begin{equation}\label{Fi-Th}
\xymatrix@R=12pt@C=12pt{
0\ar[r]&\TET\ar[r]&F_i\ar[r]&E^{c_0,b_0}_{i}\ar[r]&0
}
\end{equation}
in $\AC$. In proving the first assertion in 1), we have seen that the rank of $F_i$ is $2$ or $3$, whenever that sheaf is Bogomolov semistable. From this and the exact sequence above it follows that $E^{c_0,b_0}_{i}$ is either a torsion sheaf or
of rank $1$ occurring respectively when the $rk(F_i)=2$ or $rk(F_i)=3$.

From the equality $E^{c_0,b_0}_{i} =\HH^0(E^{c_0,b_0}_{i})$ follows the assertion that
$E^{c_0,b_0}_{i}$ is a sheaf in the subcategory $\TT_{K_X,b_0K_X}$. Since $E^{c_0,b_0}_{i}$ is a subobject of $\TET[1]$ and the latter does not admit subobjects of phase $1$, it follows that $E^{c_0,b_0}_{i}$ can not have subsheaves supported on a $0$-dimensional subscheme. Hence the torsion part of $E^{c_0,b_0}_{i}$, whenever nonzero,
is supported on a subscheme of pure dimension $1$.

We now turn to the assertion 4) of the proposition. We need to check that 
for all quotient objects $F_i[1]$, with $F_i$ of rank $2$, that sheaf
must be $K_X$-stable. Assume the contrary. Then $F_i$ admits a (semi)destabilizing exact
sequence with respect to $K_X$
\begin{equation}\label{dest-F_i}
\xymatrix@R=12pt@C=12pt{
0\ar[r]&\OO_X (-M)\ar[r]&F_i\ar[r]&{\cal J}_A (-L)\ar[r]&0;
}
\end{equation}
where a (semi)destabilizing means
\begin{equation}\label{Mhalf}
  -M\cdot K_X \geq \HA c_1(F_i)\cdot K_X.
\end{equation}
We know that $F_i\in \FF_{K_X,b_0K_X}$. Hence the subsheaf $\OO_X (-M)$ must be in $\FF_{K_X,b_0K_X}$ as well. This means
$$
-M\cdot K_X \leq b_0K^2_X.
$$
From the inequality
\begin{equation}\label{M-L}
 -M\cdot K_X \geq -L\cdot K_X
\end{equation}
it follows that $\OO_X (-L)$ is in $\FF_{K_X,b_0K_X}$ as well. So we have
the inequality 
$$
L\cdot K_X \geq -b_0K^2_X >0,
$$
where the second inequality because $b_0$ is negative.
From the dual of \eqref{dest-F_i} we have that $\OO_X (L) $ injects into
$F^{\ast}_i$, the dual of $F_i$. Furthermore, dualizing \eqref{Fi-Th} gives an inclusion
$$
F^{\ast}_i \hookrightarrow \Omega_X.
$$
Hence $\OO_X (L) $ is a subsheaf of $\Omega_X$. We have just seen that it has positive degree with respect to $K_X$. Since it can not lie in the positive cone of $X$ we deduce
$$
L^2 \leq 0.
$$
From the exact sequence \eqref{dest-F_i} follows 
$$
2ch_2(F_i)=M^2 +L^2 -2deg(A)
$$
This together with the previous inequality gives
$$
2ch_2(F_i) \leq M^2.
$$
From  \eqref{dest-F_i} we also have
 the formula $-c_1(F_i)=L+M$. This gives
$$
-c_1(F_i)\cdot K_X =M\cdot K_X + L \cdot K_X \geq 2M\cdot K_X.
$$
Now we go to the first inequality in Corollary \ref{cor:tau-bound}
\begin{equation}\label{tau-bound-pf}
\tau_X \leq \frac{2ch_2(F_i)K^2_X}{-c_1(F_i)\cdot K_X}.
  \end{equation}
  Using it together with two previous ones we obtain
  $$
  \tau_X \leq \frac{M^2K^2_X}{2M\cdot K_X}.
  $$
  This and the Hodge Index inequality $M^2K^2_X \leq (M\cdot K_X)^2$ give
  $$
  \tau_X \leq \HA M\cdot K_X \leq \frac{1}{4} (-c_1(F_i)) \cdot K_X,
  $$
  where the second inequality uses \eqref{Mhalf}. Remembering that $F^{\ast}_i$ is a proper rank two subsheaf of $\Omega_X$, we deduce
  $$
  (-c_1(F_i)) \cdot K_X =c_1 (F^{\ast}_i)\cdot K_X < K^2_X.
  $$
  Hence the inequality
  $$
   \tau_X \leq \frac{1}{4} (-c_1(F_i)) \cdot K_X \leq \frac{1}{4} K^2_X.
   $$
  Substituting the formula
  $
  \tau_X =(1-2\alpha_X)K^2_X
  $
  implies
  $$
  \alpha_X \geq \frac{3}{8}
  $$
  contrary to our assumption on $\alpha_X$.

  Once we know that $F_i$ is $K_X$-stable, it is subject to the Bogomolov-Gieseker inequality
  $$
  c^2_1(F_i) \geq 2rk(F_i)ch_2 (F_i)=4ch_2 (F_i).
  $$
  The inequality \eqref{tau-bound-pf} used above also tells us that $ch_2 (F_i)$
  is positive. Hence
  $$
  c^2_1(F_i) >0.
  $$
  This and $(-c_1(F_i)\cdot K_X)>0$ imply that $(-c_1(F_i))$ is in the positive cone $C^+(X)$ of $X$.
  
\vspace{0.2cm}
  For the part 5) we write
  \begin{equation}\label{Fi-Kfilt-pf}
  F_i=F^{l_i}_i \supset \cdots \supset F^1_i \supset F^0_i=0,
\end{equation}
   the HN filtration of $F_i$ with respect to the canonical polarization.
  Set
  $$
  Q^j_i:=F^j_i /F^{j-1}_i,\hspace{0.2cm} L_j:=-c_1(Q^j_i),\hspace{0.2cm} g_j:=rk(Q^j).
  $$
  
  Remember we already proved 6), so the last quotient $Q^{l_i}_i$ must have rank $g_{l_i}\geq 2$.
  If $F_i$ is of rank $3$, then $Q^{l_i}_i$ must be of rank $2$ and the filtration \eqref{Fi-Kfilt-pf} has the form
$$
F_i=F^{2}_i  \supset F^1_i \supset F^0_i=0,
$$
where $ F^1_i $ is a line bundle. Furthermore,  from 6) of the proposition we also know that the divisor class $(-c_1(Q^2_i)) =(-c_1 (F_i /F^1_i))$ is in the positive cone of $X$.
\end{pf}
\begin{cor}\label{cor:Fi3}
 With the assumptions and notation of Proposition \ref{pro:HNalpha} the following holds.

1) $\displaystyle{\tau_X \leq \frac{-c_1(F_i)\cdot K_X}{3}}$, for every $i\in [1,m-1]$.

2) $0\leq c_1 (E^{c_0,b_0}_i)\cdot K_X \leq 2(3c_2(X)-K^2_X)$, for every $i\in [1,m-1]$.
\end{cor}
\begin{pf}
For $F_i$ Bogomolov semistable of rank $3$ the asserted inequality in 1) is the same as the second inequality in Corollary \ref{cor:tau-bound}.
Similarly, the asserted inequality holds if $F_i$ Bogomolov unstable and
its last semistable factor $Q^{l_i}_i$ has rank $3$, see the last inequality in  Corollary \ref{cor:tau-bound}. So only $rk(F_i)=2$ (resp. $rk(Q^{l_i}_i)=2$) is new. 

We recall the sequence
\begin{equation}\label{Fi-Th-Ei}
\xymatrix@R=12pt@C=12pt{
0\ar[r]&\TET\ar[r]&F_i \ar[r]&E^{c_0,b_0}_i \ar[r]&0,
}
\end{equation}
for all $i\in [1,m-1]$. In the case $rk(F_i)=2$ the quotient $E^{c_0,b_0}_i$ is a torsion sheaf. Dualizing we obtain
$$
F^{\ast}_i \hookrightarrow \Omega_X.
$$
In other words the dual $F^{\ast}_i $ is a rank $2$ subsheaf of $\Omega_X$. Furthermore, from Proposition \ref{pro:HNalpha} we know that 
$
c_1(F^{\ast}_i)=-c_1(F_i)
$
is in the positive cone of $X$. A result of Miyaoka tells us that the Chern classes of $F^{\ast}_i$ are subject to
$$
c^2_1(F_i)=c^2_1(F^{\ast}_i)\leq 3c_2(F^{\ast}_i)=3c_2(F_i),
$$
see \cite{Mi1}, Remark 4.18. Equivalently, this can be rewritten as
$$
2ch_2(F_i) \leq \frac{1}{3}c^2_1(F_i).
$$
Substituting this into the first inequality in Corollary \ref{cor:tau-bound}
we obtain
$$
\tau_X \leq \frac{2ch_2(F_i)K^2_X}{-c_1(F_i)\cdot K_X}\leq \frac{c^2_1(F_i)K^2}{-3c_1(F_i)\cdot K_X}.
$$
This and the Hodge Index inequality $c^2_1(F_i)K^2 \leq \Big(-c_1(F_i)\cdot K_X\Big)^2$ give
$$
\tau_X \leq \frac{-c_1(F_i)\cdot K_X}{3}.
$$

If $F_i$ is Bogomolov unstable with the destabilizing sequence
$$
\xymatrix@R=12pt@C=12pt{
	0\ar[r]&F^{l_i-1}_i\ar[r]&F_i \ar[r]&Q^{l_i}_i \ar[r]&0,
}
$$
where $Q^{l_i}_i$ is of rank $2$, we combine it with the exact sequence in \eqref{Fi-Th-Ei} to obtain the diagram
\begin{equation}\label{Fi-Theta-diag}
\xymatrix@R=12pt@C=12pt{\
	&&0\ar[d]&&\\
	& &F^{l_i-1}_i\ar[d]&&\\
	0\ar[r]&\TET\ar[r]\ar[dr]&F_i \ar[r]\ar[d]&E^{c_0,b_0}_i \ar[r]&0\\
	&&Q^{l_i}_i \ar[d]&&\\
	&&0&&
} 
\end{equation}
We claim that the resulting composition arrow, the slanted arrow of the diagram, is injective. Let us assume this for a moment and complete the argument: 

\noindent
dualize the slanted arrow to obtain that the
dual $(Q^{l_i}_i)^{\ast}$ is a subsheaf of $\Omega_X$ of rank $2$ and having $c_1( (Q^{l_i}_i)^{\ast})=-c_1(Q^{l_i}_i)$ in the positive cone of $X$; by the result of Miyaoka we used above
$$
c^2_1 (Q^{l_i}_i) \leq 3c_2 ( Q^{l_i}_i); 
$$ 
 this gives
$$
2ch_2 (Q^{l_i}_i) \leq \frac{1}{3} c^2_1 (Q^{l_i}_i);
$$
now use this in the estimate for $2ch_2 (F_i)$ in the proof of Corollary \ref{cor:tau-bound} to obtain
$$
\begin{gathered}
2ch_2 (F_i) < 2ch_2 (Q^{l_i}_i) +\frac{(-c_1 (Q^{l_i}_i)\cdot K_X)}{2K^2_X} (-c_1(F_i)\cdot K_X +c_1 (Q^{l_i}_i)\cdot K_X) \leq\\
\frac{1}{3} c^2_1 (Q^{l_i}_i) +\frac{(-c_1 (Q^{l_i}_i)\cdot K_X)}{2K^2_X} (-c_1(F_i)\cdot K_X +c_1 (Q^{l_i}_i)\cdot K_X);
\end{gathered}
	$$
	this and the Hodge Index inequality give
$$
2ch_2 (F_i) < -\frac{(c_1 (Q^{l_i}_i)\cdot K_X)^2}{6K^2} +
\frac{(c_1 (Q^{l_i}_i)\cdot K_X)}{2K^2_X} (c_1(F_i)\cdot K_X);
$$
setting $t= -c_1 (Q^{l_i}_i)\cdot K_X$, we obtain the quadratic relation
$$
t^2 +3 (c_1(F_i)\cdot K_X)t +12K^2_X ch_2 (F_i)<0.
$$
The quadratic polynomial on the left must have positive discriminant
$$
d:=9(c_1(F_i)\cdot K_X)^2 -48K^2 ch_2(F_i)>0
$$
 and the polynomial is non positive precisely when  $t$ is subject to the inequalities
$$
\frac{3}{2}(-c_1(F_i)\cdot K_X) -\HA \sqrt{d} \leq t \leq \frac{3}{2}(-c_1(F_i)\cdot K_X) +\HA \sqrt{d}.
$$
The value $t_0 = (-c_1 (Q^{l_i}_i)\cdot K_X)$ must be in this range and since 
$$
 -c_1 (Q^{l_i}_i)\cdot K_X \leq -c_1(F_i)\cdot K_X,
 $$ 
 we deduce
that $(-c_1(F_i)\cdot K_X)$ is also in that range. Hence the quadratic polynomial above is non positive for $t=(-c_1(F_i)\cdot K_X)$:
$$
-2(c_1(F_i)\cdot K_X)^2 +12K^2_X ch_2 (F_i)\leq 0.
$$ 
This gives
$$
2K^2ch_2 (F_i) \leq \frac{1}{3}(c_1(F_i)\cdot K_X)^2 .
$$
Substituting into the first inequality in Corollary \ref{cor:tau-bound} we obtain
$$
\tau_X \leq \frac{1}{3}(-c_1(F_i)\cdot K_X).
$$

We now return to the claim:

$$
\text{\it the slanted arrow in the diagram \eqref{Fi-Theta-diag} is injective.}
$$
For this we observe that the diagram comes from a similar one but in the abelian category $\AC_{K_X,b_0K_X}$. Namely, we have the exact sequence
\begin{equation}\label{intilted}
	\xymatrix@R=12pt@C=12pt{
0\ar[r]&E^{c_0,b_0}_i \ar[r]&\TET[1]\ar[r]&F_i[1]\ar[r]&0
}
\end{equation}
in $\AC_{K_X,b_0K_X}$ coming from Bridgeland instability of $\TET[1]$ and the exact sequence
$$
\xymatrix@R=12pt@C=12pt{
	0\ar[r]&F^{l_i -1}_i \ar[r]&F_i\ar[r]&Q^{l_i}_i\ar[r]&0
}
$$
in $\FF_{K_X,b_0K_X}$ resulting from the fact that $F_i$ is $K_X$-unstable. The last sequence gives rise to an exact sequence
$$
\xymatrix@R=12pt@C=12pt{
	0\ar[r]&F^{l_i -1}_i [1]\ar[r]&F_i[1]\ar[r]&Q^{l_i}_i[1]\ar[r]&0
}
$$
in $\FF_{K_X,b_0K_X}[1]$ and hence in $\AC_{K_X,b_0K_X}$. Putting it together with the exact sequence in \eqref{intilted} gives
$$
\xymatrix@R=12pt@C=12pt{
	&&&0\ar[d]&\\
	&&&F^{l_i -1}_i [1]\ar[d]&\\
	0\ar[r]&E^{c_0,b_0}_i \ar[r]&\TET[1]\ar[r]&F_i[1]\ar[r]\ar[d]&0\\
	&&&Q^{l_i}_i[1]\ar[d]&\\
	&&&0&
}
$$
which can be completed to the commutative diagram  		
\begin{equation}\label{Fi-Theta-diagtilt}
	\xymatrix@R=12pt@C=12pt{
		&&0\ar[d]&0\ar[d]&\\
		0\ar[r]&E^{c_0,b_0}_i \ar@{=}[d] \ar[r]&\widetilde{E_i}\ar[r]\ar[d]&F^{l_i-1}_i[1]\ar[d]\ar[r]&0\\
		0\ar[r]&E^{c_0,b_0}_i \ar[r]&\TET[1]\ar[r]\ar[d]&F_i[1]\ar[r]\ar[d]& 0\\
		&&Q^{l_i}_i[1] \ar@{=}[r]\ar[d]&Q^{l_i}_i[1]\ar[d]&\\
		&&0&0&
	} 
\end{equation}
in $\AC_{K_X,b_0K_X}$. The diagram \eqref{Fi-Theta-diag} we are concerned with is obtained by applying the homological functor $\HH^0$ to the diagram above:
$$
 	\xymatrix@R=12pt@C=12pt{
 	&0\ar[d]&0\ar[d]&&&\\ 0\ar[r]&\HH^{-1}(\widetilde{E_i})\ar[r]\ar[d]&F^{l_i-1}_i\ar[d]\ar[r]&E^{c_0,b_0}_i\ar[r]\ar@{=}[d]&\HH^0(\widetilde{E_i})\ar[r]&0\\
 	0\ar[r]&\TET\ar[r]\ar[d]&F_i\ar[r]\ar[d]&E^{c_0,b_0}_i \ar[r]& 0\\
 	&Q^{l_i}_i \ar@{=}[r]\ar[d]&Q^{l_i}_i \ar[d]&&\\
 	&\HH^0(\widetilde{E_i})\ar[d]&0&&\\
 	&0&&&
 } 
$$
 The slanted arrow in \eqref{Fi-Theta-diag} is the arrow 
 $
 \Theta_X \longrightarrow Q^{l_i}_i
 $
 in the column on the left in the above diagram. The claim of its injectivity comes down to showing that $\HH^{-1}(\widetilde{E_i})$ is zero.
 
 Assume $\HH^{-1}(\widetilde{E_i})\neq 0$. Then it is a line bundle which we call $\OO_X (-M)$ and the left column can be broken into two exact sequences
 $$
 \xymatrix@R=12pt@C=12pt{
 	0\ar[r]&\OO_X (-M)\ar[r]&\TET\ar[r] &{\cal J}_A (-L)\ar[r]&0,\\
0\ar[r] &{\cal J}_A (-L)\ar[r]& Q^{l_i}_i \ar[r]&\HH^0(\widetilde{E_i})\ar[r]&0, 	
 }
$$
where ${\cal J}_A $ is the ideal sheaf of at most $0$-dimensional subscheme $A$ and $\OO_X (-L)$ is a line bundle. From the second exact sequence it follows
\begin{equation}\label{c1Qli}
c_1 ( Q^{l_i}_i)\cdot K_X= -L\cdot K_X + c_1 (\HH^0(\widetilde{E_i})) \cdot K_X.
\end{equation}
The stability of $\TET$ with respect to $K_X$ gives
$$
 -L\cdot K_X > -\HA K^2_X,
$$
while $\HH^0(\widetilde{E_i}) \in \TT_{K_X,b_0 K_X}$ implies
$$
c_1 (\HH^0(\widetilde{E_i})) \cdot K_X > b_0 K^2_X.
$$
Substituting two estimates above into \eqref{c1Qli} we obtain
$$
c_1 ( Q^{l_i}_i)\cdot K_X > -\HA K^2_X +  b_0 K^2_X.
$$
Recall that we can replace $b_0$ in the above inequality by any value in the interval $[b_0,0)$, see proof of Proposition \ref{pro:no1rankbound}, 2). This gives the inequality
$$
c_1 ( Q^{l_i}_i)\cdot K_X \geq -\HA K^2_X.
$$ 
 Substituting into the first inequality in \eqref{taubound-cor} we deduce
 $$
 \tau_X <\frac{1}{4}K^2_X.
 $$
This is equivalent to $\alpha_X > \frac{3}{8}$ which contradicts the assumption $\alpha_X < \frac{3}{8}$. The proof of the part 1) of the corollary is now completed.
  
  The assertion 2) of the corollary follows immediately from 1) and the relation
$$
 -c_1(F_i)=K_X-c_1(E^{c_0,b_0}_i);
$$
the latter comes from the exact sequence \eqref{Fi-Th-Ei}. 
Substituting this into the inequality in 1) gives
$$
c_1(E^{c_0,b_0}_i)\cdot K_X \leq  K^2-3\tau_X=6c_2(X)-2K^2_X =2(3c_2(X)-K^2_X).
$$
The non-negativity of $c_1(E^{c_0,b_0}_i)\cdot K_X$ is Proposition \ref{pro:no1rankbound}, 3).
\end{pf}

We have seen that the quotient objects $F_i [1] =\TET[1]/E^{c_0,b_0}_i$ arising from the HN filtration \eqref{HNfilt-pro} are well controlled provided they are Bogomolov semistable. Otherwise we have no longer
in control of their ranks. However, from Proposition \ref{pro:HNalpha}, 6), we know that the ranks of their last quotients 
$Q^{l_i}_i$ are still of rank $2$ or $3$. We have already seen that we can perform a sort of reduction by replacing $F_i [1]$ by 
$Q^{l_i}_i [1]$. Namely, we have the diagram
\begin{equation}\label{diag-red}
 \xymatrix@R=12pt@C=12pt{
 	&&0\ar[d]&0\ar[d]&\\
 	0\ar[r]&E^{c_0,b_0}_i \ar[r]\ar@{=}[d]& \widetilde{E_i}\ar[d] \ar[r]&F^{l_i-1}_i[1] \ar[r]\ar[d]&0\\
 	0\ar[r]&E^{c_0,b_0}_i \ar[r]& \TET[1]\ar[r]\ar[d]&F_i[1] \ar[d]\ar[r]&0\\
 	&&Q^{l_i}_i [1]\ar@{=}[r]\ar[d]&Q^{l_i}_i [1]\ar[d]&\\
 	&&0&0&
 } 
\end{equation}
where the column in the middle replaces the quotient object $F_i[1]$ by
the `semistable' quotient $Q^{l_i}_i [1]$. The price to pay is that the
 new sequence is not necessarily {\it Bridgeland} destabilizing for $\TET[1]$. In particular, it is not clear whether the object
$\widetilde{E_i}$ enjoys the properties of $E^{c_0,b_0}_i$. The next result partially answers the question.

\begin{lem}\label{lem:tildEi}
	With the assumptions and notation of Proposition \ref{pro:HNalpha},
	consider a quotient object $F_i [1]$ with $F_i$ unstable with respect to $K_X$. Then it admits the reduction diagram \eqref{diag-red}, where $Q^{l_i}_i$ is $K_X$-semistable sheaf of rank $2$ or $3$ with $(-c_1 (Q^{l_i}_i))$ in the positive cone of $X$. In addition, $\widetilde{E_i}$ is a sheaf of rank
	$rk(\widetilde{E_i})=(rk(Q^{l_i}_i)-2)$ in $\TT_{K_X,b_0K_X}$ and subject to
	$$
	c_1 (\widetilde{E_i})\cdot K_X < (4-2rk(\widetilde{E_i}))(3c_2 (X)-K^2_X).
	$$    
\end{lem}
\begin{pf}
	The statements about $Q^{l_i}_i$ were proved in Prposition \ref{pro:HNalpha}. Only the part concerning $\widetilde{E_i}$ needs to be proved.
	Applying the homological functor $\HH^0$ to the middle column in \eqref{diag-red} we obtain the exact complex in $\AC$:
	$$
\xymatrix@R=12pt@C=12pt{
	0\ar[r]& \HH^{-1}(\widetilde{E_i})\ar[r]& \TET \ar[r]& Q^{l_i}_i \ar[r]& \HH^0 (\widetilde{E_i})\ar[r]&0.
}
	$$
	We need to check that $\HH^{-1}(\widetilde{E_i})$ is zero. Assume it is not. Then it is a line subsheaf of $\TET$ which is denoted
	$\OO_X (-M)	$. The exact complex can be broken into two short exact sequences
	$$
\xymatrix@R=12pt@C=12pt{
	0\ar[r]& \OO_X (-M)\ar[r]& \TET \ar[r]&{\cal J}_A (-L)\ar[r]&0,\\
	0\ar[r]& {\cal J}_A (-L)\ar[r]&Q^{l_i}_i \ar[r]& \HH^0 (\widetilde{E_i})\ar[r]&0.
}
$$
From the second sequence we deduce
$$
c_1 (Q^{l_i}_i)	=-L +c_1 (\HH^0 (\widetilde{E_i})).
$$
Intersecting with $K_X$ gives
$$
c_1 (Q^{l_i}_i)\cdot K_X	=-L\cdot K_X +c_1 (\HH^0 (\widetilde{E_i}))\cdot K_X > -L\cdot K_X +rk(\HH^0 (\widetilde{E_i}))(-b_0)K^2_X.
$$
where the inequality uses the fact that $\HH^0 (\widetilde{E_i})$ lies in $\TT_{K_X,b_0 K_X}$. Since $b_0$ can be replaced by any $b\in [b_0,0)$, the above inequality implies
$$
c_1 (Q^{l_i}_i)\cdot K_X	\geq -L\cdot K_X.
$$
From the first exact sequence above and $K_X$-stability of $\TET$ 
$$
L\cdot K_X < \HA K^2_X. 
$$
Substituting this into the previous inequality gives
$$
-c_1 (Q^{l_i}_i)\cdot K_X <	\HA K^2_X .
$$
Combining this with the first inequality in \eqref{taubound-cor} we obtain
$$
\tau_X <\frac{(-c_1 (Q^{l_i}_i)\cdot K_X) }{rk(Q^{l_i}_i)} <\frac{K^2_X}{2rk(Q^{l_i}_i)} \leq 	\frac{K^2_X}{4}.
$$
The resulting inequality
$$
\tau_X <\frac{K^2_X}{4}
$$
is equivalent to 
$$
\alpha_X > \frac{3}{8}
$$
which contradicts our assumption that $\alpha_X$ is less than $\frac{3}{8}$. Thus $\widetilde{E_i} =\HH^0 (\widetilde{E_i})$ is a sheaf in $\TT_{K_X,b_0 K_X}$.

We now turn to the upper bound on the degree $c_1 (\widetilde{E_i})\cdot K_X$. For this we use the estimate
$$
2ch_2 (F_i)< 2ch_2 (Q^{l_i}_i) + \frac{c_1 (Q^{l_i}_i)\cdot K_X}{rk(Q^{l_i}_i)} c_1 (F^{l_i -1}_i)\cdot K_X \leq \frac{1}{3}c^2_1 (Q^{l_i}_i)+ \frac{c_1 (Q^{l_i}_i)\cdot K_X}{rk(Q^{l_i}_i)K^2_X} c_1 (F^{l_i -1}_i)\cdot K_X .
$$
For $rk(Q^{l_i}_i)=3$, the above inequality gives
$$
\begin{gathered}
2ch_2 (F_i)<\frac{1}{3}c^2_1 (Q^{l_i}_i)+ \frac{c_1 (Q^{l_i}_i)\cdot K_X}{3} c_1 (F^{l_i -1}_i)\cdot K_X \leq \frac{\Big(c_1 (Q^{l_i}_i)\cdot K_X \Big)^2}{3K^2_X}+ \frac{c_1 (Q^{l_i}_i)\cdot K_X}{3K^2_X} c_1 (F^{l_i -1}_i)\cdot K_X
\\
= \frac{c_1 (Q^{l_i}_i)\cdot K_X}{3K^2_X} \Big(c_1 (Q^{l_i}_i) + c_1 (F^{l_i -1}_i)\Big)\cdot K_X = \frac{c_1 (Q^{l_i}_i)\cdot K_X}{3K^2_X}  c_1 (F_i)\cdot K_X,
\end{gathered}
$$
where the last equality uses the relation
$c_1 (F_i)=c_1 (Q^{l_i}_i) +c_1 (F^{l_i -1}_i)$. The above inequality implies
$$
\frac{-c_1 (Q^{l_i}_i)\cdot K_X}{3} > \frac{2ch_2 (F_i)K^2_X}{-
c_1(F_i)\cdot K_X}.
$$
This and the first inequality in Corollary \ref{cor:tau-bound}
imply
$$
\frac{-c_1 (Q^{l_i}_i)\cdot K_X}{3} > \tau_X.
$$
From the middle column in \eqref{diag-red} we have
$$
-c_1 (Q^{l_i}_i) =K_X -c_1 (\widetilde{E_i}).
$$
Substituting into the above inequality gives the estimate
$$
c_1 (\widetilde{E_i}) \cdot K_X <2(3c_2 (X)-K^2_X).
$$
Since $rk(\widetilde{E_i})=rk(Q^{l_i}_i)-2=3-2=1$, the above inequality is the one stated in the lemma.

We now turn to the case $rk(Q^{l_i}_i)=2$. The estimate for $2ch_2 (F_i)$
now reads
\begin{equation}\label{ch2Fi-2}
2ch_2 (F_i)<\frac{1}{3}c^2_1 (Q^{l_i}_i)+ \frac{c_1 (Q^{l_i}_i)\cdot K_X}{2K^2_X} c_1 (F^{l_i -1}_i)\cdot K_X
\end{equation}
Using the relations
$$
-c_1 (Q^{l_i}_i) =K_X -c_1 (\widetilde{E_i}),\hspace{0.2cm} 
-c_1 (F^{l_i-1}_i) =c_1 (\widetilde{E_i}) -c_1 (E^{c_0,b_0}_i),
$$
coming respectively from the middle column and the top row of the diagram \eqref{diag-red}, the inequality above becomes
$$
\begin{gathered}
2ch_2 (F_i)<\frac{1}{3}(K_X -c_1 (\widetilde{E_i}))^2+ \frac{ K^2_X -c_1 (\widetilde{E_i})\cdot K_X}{2K^2_X}  (c_1 (\widetilde{E_i})\cdot K_X -c_1 (E^{c_0,b_0}_i) \cdot K_X)
\\
=\frac{1}{3}K^2_X -\frac{1}{6} c_1 (\widetilde{E_i})\cdot K_X +\frac{1}{3}c^2_1 (\widetilde{E_i}) - \frac{\Big(c_1(\widetilde{E_i})\cdot K_X \Big)^2}{2K^2_X} -\Big(c_1 (E^{c_0,b_0}_i) \cdot K_X \Big) \frac{ K^2_X -c_1 (\widetilde{E_i})\cdot K_X}{2K^2_X}.
\end{gathered}
$$
From the first inequality of Corollary \ref{cor:tau-bound} we have
$$
\frac{(-c_1 (F_i)\cdot K_X)}{K^2_X}\tau_X \leq 2ch_2 (F_i)
$$
Using the relation
$$
-c_1 (F_i)=K_X -c_1 (E^{c_0,b_0}_i)
$$
from the middle row of the diagram \eqref{diag-red}, we obtain
$$
\tau_X -\frac{c_1 (E^{c_0,b_0}_i)\cdot K_X}{K^2_X}\tau_X \leq 2ch_2 (F_i).
$$
Combining with the upper bound for $2ch_2 (F_i)$ gives the inequality
$$
\begin{gathered}
\tau_X -\frac{c_1 (E^{c_0,b_0}_i)\cdot K_X}{K^2_X}\tau_X  
\\
<
\frac{1}{3}K^2_X -\frac{1}{6} c_1 (\widetilde{E_i})\cdot K_X +\frac{1}{3}c^2_1 (\widetilde{E_i}) - \frac{\Big(c_1(\widetilde{E_i})\cdot K_X \Big)^2}{2K^2_X} -\Big(c_1 (E^{c_0,b_0}_i) \cdot K_X \Big) \frac{ K^2_X -c_1 (\widetilde{E_i})\cdot K_X}{2K^2_X}.
\end{gathered} 
$$
This  is rewritten as follows
\begin{equation}\label{qin-lem}
\begin{gathered}
\frac{\Big(c_1(\widetilde{E_i})\cdot K_X \Big)^2}{2K^2_X} +\frac{1}{6} c_1(\widetilde{E_i})\cdot K_X -\frac{1}{3}c^2_1 (\widetilde{E_i})<\frac{1}{3}K^2_X -\tau_X 
-\frac{c_1 (E^{c_0,b_0}_i) \cdot K_X}{K^2}(\HA K^2_X -\tau_X -\HA c_1 (\widetilde{E_i})\cdot K_X)
\\
=
\frac{2}{3}(3c_2(X)-K^2_X)-\frac{c_1 (E^{c_0,b_0}_i) \cdot K_X}{K^2}(\HA K^2_X -\tau_X -\HA c_1 (\widetilde{E_i})\cdot K_X).
\end{gathered}
\end{equation}
We claim: 
\begin{equation}\label{posit-term}
\HA K^2_X -\tau_X -\HA c_1 (\widetilde{E_i})\cdot K_X) >0.
\end{equation}
With this in mind, let us  derive the asserted upper bound for 	
$c_1(\widetilde{E_i})\cdot K_X$: our claim \eqref{posit-term} and the non-negativity of $c_1 (E^{c_0,b_0}_i) \cdot K_X$, see Proposition \ref{pro:no1rankbound}, 3), allow to drop the last term in the inequality
\eqref{qin-lem} and thus to obtain
\begin{equation}\label{ineq-tildEi}
\frac{\Big(c_1(\widetilde{E_i})\cdot K_X \Big)^2}{2K^2_X} +\frac{1}{6} c_1(\widetilde{E_i})\cdot K_X -\frac{1}{3}c^2_1 (\widetilde{E_i})<\frac{2}{3}(3c_2(X)-K^2_X).
\end{equation}
The last term on the left can be bounded from above by using the Hodge Index inequality
$$
c^2_1 (\widetilde{E_i}) K^2_X \leq \Big(c_1(\widetilde{E_i})\cdot K_X \Big)^2
$$
This gives
$$
\frac{\Big(c_1(\widetilde{E_i})\cdot K_X \Big)^2}{6K^2_X} +\frac{1}{6} c_1(\widetilde{E_i})\cdot K_X <\frac{2}{3}(3c_2(X)-K^2_X).
$$ 
Setting $t=\displaystyle{\frac{c_1(\widetilde{E_i})\cdot K_X }{K^2_X}}$,
the above inequality becomes
$$
t^2 +t -4(3\alpha_X -1)<0.
$$
Solving for $t$ we obtain
\begin{equation}\label{solving-ineq}
\frac{c_1(\widetilde{E_i})\cdot K_X }{K^2_X} =t <\frac{-1+\sqrt{1+16(3\alpha_X -1)}}{2} =\frac{8(3\alpha_X -1)}{1+\sqrt{1+16(3\alpha_X -1)}}<4(3\alpha_X -1),
\end{equation}
where the last inequality uses the fact that $\alpha_X > \frac{1}{3}$.
Thus we obtain
$$
c_1(\widetilde{E_i})\cdot K_X <4(3\alpha_X -1)K^2_X =4(3c_2 (X)-K^2_X).
$$
Since the rank $rk(\widetilde{E_i})=rk(Q^{l_i}_i)-2=2-2=0$, the above inequality is the one asserted in the lemma.

We now return to the inequality claimed in \eqref{posit-term}. To prove it we go back to the estimate on $2ch_2 (F_i)$ in \eqref{ch2Fi-2} and write it as follows
$$
\begin{gathered}
2ch_2 (F_i)<\frac{1}{2}c^2_1 (Q^{l_i}_i) + \frac{c_1 (Q^{l_i}_i)\cdot K_X}{2K^2_X} c_1 (F^{l_i -1}_i)\cdot K_X -\frac{1}{6}c^2_1 (Q^{l_i}_i)
\\
\leq 
\frac{(c_1 (Q^{l_i}_i) \cdot K_X )^2}{2K^2_X} + \frac{c_1 (Q^{l_i}_i)\cdot K_X}{2K^2_X} c_1 (F^{l_i -1}_i)\cdot K_X -\frac{1}{6}c^2_1 (Q^{l_i}_i)
\\
=\frac{(c_1 (Q^{l_i}_i) \cdot K_X )}{2K^2_X} \Big(c_1 (Q^{l_i}_i) +c_1 (F^{l_i -1}_i) \Big)\cdot K_X -\frac{1}{6}c^2_1 (Q^{l_i}_i)=
\frac{(c_1 (Q^{l_i}_i) \cdot K_X )}{2K^2_X} c_1 (F_i)\cdot K_X -\frac{1}{6}c^2_1 (Q^{l_i}_i),
\end{gathered} 
$$
where the last equality uses the relation
$$
c_1 (F_i)=c_1 (Q^{l_i}_i) +c_1 (F^{l_i -1}_i) 
$$
coming from the right column in the diagram \eqref{diag-red}.
The resulting inequality
$$
2ch_2 (F_i)<\frac{(c_1 (Q^{l_i}_i) \cdot K_X )}{2K^2_X} c_1 (F_i)\cdot K_X -\frac{1}{6}c^2_1 (Q^{l_i}_i)
$$
together with the first inequality in Corollary \ref{cor:tau-bound}
give
$$
\tau_X \leq \frac{2ch_2 (F_i) K^2_X}{(-c_1 (F_i)\cdot K_X)} <
\frac{(-c_1 (Q^{l_i}_i)) \cdot K_X )}{2} - \frac{c^2_1 (Q^{l_i}_i) K^2_X}{6(-c_1 (F_i)\cdot K_X)}.
$$
Using the relation
$$
(-c_1 (Q^{l_i}_i)=K_X -c_1 (\widetilde{E_i})
$$
from the middle row of the diagram \eqref{diag-red} we obtain
$$
\frac{c^2_1 (Q^{l_i}_i) K^2_X}{6(-c_1 (F_i)\cdot K_X)} <\HA K^2_X -\tau_X -\HA c_1 (\widetilde{E_i})\cdot K_X.
$$
We know that the term on the left is positive. Hence the inequality
$$
0<\HA K^2_X -\tau_X -\HA c_1 (\widetilde{E_i})\cdot K_X
$$
asserted in \eqref{posit-term}. The proof of the lemma is now complete.
\end{pf}

We put a certain effort in obtaining the estimate in Lemma \ref{lem:tildEi}, however it is not very useful unless we know that the divisor class $c_1 (\widetilde{E_i})$ is effective. This is certainly  the case if $rk(\widetilde{E_i}) =0$: $\widetilde{E_i}$ is a torsion sheaf and the inequality gives the upper bound on the degree of the support of the sheaf in terms of the Chern numbers of $X$ (even better - in terms of $(3c_2(X)-K^2_X)$). But for $\widetilde{E_i}$ of rank $1$, one would not expect the effectiveness of its first Chern class to come
from general considerations of Bridgeland instability of $\TET[1]$.
In fact, so far we studied the HN filtration of $\TET[1]$ without using the fact that the Bridgeland instability of that object comes from the monomorphism
$$
\Phi^V_{H,bH}: V\otimes \OO_X \longrightarrow \TET[1]
$$
in the abelian category $\AC_{H,b H}$ for appropriate values of $b$.
We now discuss this aspect.

The monomorphism $\Phi^V_{H,bH}$ allows us to think of $V\otimes \OO_X$ as a {\it subobject} of $\TET[1]$. Furthermore, for the stability
conditions $(cH,bH) \in B^{un}_H$ this an unstable subobject of $\TET[1]$ and the monomorphism can be completed to the exact sequence
$$
\xymatrix@R=12pt@C=22pt{
	0\ar[r]&V\otimes \OO_X \ar[r]^{\Phi^V_{H,bH}}&\TET[1]\ar[r]& coker(\Phi^V_{H,bH}) \ar[r]& 0
}
	$$
	in the abelian category $\AC_{H,bH}$. We also know that the cokernel has the form $F_V [1]$, where $F_V$ is a sheaf in the subcategory $\FF_{H,bH}$. So we rewrite the above exact sequence as follows
	\begin{equation}\label{PhiV-extseq}
\xymatrix@R=12pt@C=22pt{
	0\ar[r]&V\otimes \OO_X \ar[r]^{\Phi^V_{H,bH}}&\TET[1]\ar[r]& F_V[1] \ar[r]& 0.
}
\end{equation}	
This an instance of a Bridgeland destabilizing sequence for $\TET[1]$ and it has the same aspect as the ones coming from the HN filtration of $\TET[1]$.
In addition, we can {\it identify} $F_V$: applying the homological functor $\HH^0$ to this sequence gives
\begin{equation}\label{Vext}
\xymatrix@R=12pt@C=22pt{
	0\ar[r]&\TET \ar[r]&F_V \ar[r]&V\otimes\OO_X \ar[r]& 0
}
\end{equation}
and this the extension sequence corresponding to the inclusion
$V\subset H^1 (\TET)$ via the identification
$$
Hom(V,H^1(\TET))\cong Ext^1 (V\otimes \OO_X, \TET).
$$
In particular, taking $H=K_X$, we apply Proposition \ref{pro:HNalpha} to obtain the following.
\begin{pro}\label{pro:FVunst}
	Assume $\displaystyle{\alpha_X<\frac{3}{8}}$ and $dim (V) \geq 2$.
	Let $(cK_X, bK_X)$ be a stability condition in $B^{un}_{K_X}$, with
	$b>\beta^{in}_{K_X}$, where $\beta^{in}_{K_X}$ is defined in Lemma \ref{lem:filt}.
	Then $F_V$ is Bogomolov unstable. Its HN filtration
	$$
	F_V=F^m_V \supset F^{m-1}_V \supset \cdots \supset F^{1}_V \supset 
	F^0_V =0
	$$
	with respect to $K_X$ has its semistable factors $Q^i_V =F^i_V / F^{i-1}_V$ all lying in $\FF_{K_X,bK_X}$. In addition, the last factor
	$Q^m_V$ has rank $2$ or $3$ and the divisor class $(-c_1 (Q^m_V))$
	lies in the positive cone $C^+(X)$.
\end{pro}
\begin{pf}
	From the extension sequence \eqref{Vext} the rank of $F_V$ is
	$(dim(V)+2)$ and this is at least $4$: the assumption
	$dim(V)\geq 2$. From Proposition \ref{pro:HNalpha}, 3), it follows that $F_V$ is Bogomolov unstable. The assertions about the HN filtration of $F_V$ are now obtained from Proposition \ref{pro:HNalpha}, 6).
\end{pf}

We now apply the reduction diagram \eqref{diag-red} to $F_V[1]$ to obtain

\begin{equation}\label{diag-red-FV}
	\xymatrix@R=12pt@C=12pt{
		&&0\ar[d]&0\ar[d]&\\
		0\ar[r]&V\otimes \OO_X \ar[r]\ar@{=}[d]& \widetilde{E_V}\ar[d] \ar[r]&F^{m-1}_V[1] \ar[r]\ar[d]&0\\
		0\ar[r]&V\otimes \OO_X \ar[r]& \TET[1]\ar[r]\ar[d]&F_V[1] \ar[d]\ar[r]&0\\
		&&Q^{m}_V [1]\ar@{=}[r]\ar[d]&Q^{m}_V [1]\ar[d]&\\
		&&0&0&
	} 
\end{equation}
The proof of Lemma \ref{lem:tildEi} goes through to give the following.
\begin{pro}\label{pro:tildEV}
	With the assumptions and notation of Proposition \ref{pro:FVunst},
	the object $F_V[1]$ admits the reduction diagram \eqref{diag-red-FV},
	where $Q^m_V$ is $K_X$-semistable of rank $2$ or $3$ with
	$(-c_1 (Q^m_V))$ in the positive cone of $X$. In addition, $\widetilde{E_V}$ in the diagram \eqref{diag-red-FV} is a sheaf of rank $rk(\widetilde{E_V})=(rk(Q^m_V)-2)$ generated by its global sections. The divisor class  $c_1 (\widetilde{E_V})$ is effective, nonzero and subject to
	$$
	c_1 (\widetilde{E_V})\cdot K_X <(4-2rk(\widetilde{E_V}))(3c_2(X)-K^2_X).
	$$
\end{pro}
\begin{pf}
	Only the statement of global generation is new and this follows from
	the top row of the diagram \eqref{diag-red-FV} by applying the homological functor $\HH^0$ to it:
	\begin{equation}\label{tildEV-V}
\xymatrix@R=12pt@C=12pt{
	0\ar[r]&F^{m-1}\ar[r]&V\otimes \OO_X \ar[r]&\widetilde{E_V} \ar[r]&0.
}
\end{equation}
The epimorphism in the sequence means that $\widetilde{E_V}$ is globally generated.

Once we know that $\widetilde{E_V}$ is globally generated, the effectivness of 
$c_1 (\widetilde{E_V})$ follows. The fact that it is nonzero is obvious in the case $rk(Q^m_V)=2$: $\widetilde{E_V}$ is torsion and must be nonzero, since otherwise the extension sequence \eqref{Vext} splits.
In the  case $rk(Q^m_V)=3$, the sheaf $\widetilde{E_V}$ is of rank $1$ and generated by its global sections. In these circumstances, $c_1(\widetilde{E_V})=0$ means that $\widetilde{E_V} =\OO_X$. But the cohomology sequence of \eqref{tildEV-V} tells that
$$
V\hookrightarrow H^0 (\widetilde{E_V});
$$
 with $\widetilde{E_V} =\OO_X$ we obtain $dim(V)=1$ and this contradicts the assumption that $dim(V)$ is at least $2$.	
\end{pf}

From the above perspective $V\otimes \OO_X$, viewed as a Bridgeland destabilizing subobject of $\TET[1]$, is `distinguished' because it gives rise to a {\it new} subobject of $\TET[1]$ which has interesting geometric properties: that subobject is $\widetilde{E_V}$ in Proposition \ref{pro:tildEV}.

But as we have seen Bridgeland stability conditions, a priori, give  rise to other destabilizing subobjects of $\TET[1]$. As we already said, see discussion after the proof of Lemma \ref{lem:tildEi}, one of the shortcomings of other destabilizing subobjects is a lack of effectiveness
property. After Proposition \ref{pro:tildEV}, one way to remedy this is
to compare them with the `effective' destabilizing sequence 
\eqref{PhiV-extseq}. Namely, we go back to the exact sequence
$$
\xymatrix@R=12pt@C=12pt{
	0\ar[r]& \widetilde{E_i}\ar[r]& \TET[1]\ar[r]& Q^{l_i}_i [1] \ar[r]&0,
	}
$$
the middle column in the diagram \eqref{diag-red}, and put it together with the exact sequence \eqref{PhiV-extseq}
$$
\xymatrix@R=12pt@C=12pt{
	&&0\ar[d]&&\\
	&&V\otimes \OO_X \ar[d]^{\Phi^V_{H,bH}}&&\\
	0\ar[r]& \widetilde{E_i}\ar[r]& \TET[1]\ar[r]\ar[d]& Q^{l_i}_i [1] \ar[r]&0\\
	&&F_V[1]\ar[d]&&\\
	&&0&&.
}
$$
From this diagram the effectiveness of $\widetilde{E_i}$ would be a consequence of the vertical arrow $\Phi^V_{H,bH}$ above factoring through the monomorphism of the horizontal sequence. This is now in the scope of Bridgeland stability:  the property of factoring through the monomorphism is equivalent to showing that the composite arrow
$$
\xymatrix@R=12pt@C=12pt{
	&&V\otimes \OO_X \ar[d]_{\Phi^V_{H,bH}} \ar[dr]&&\\
	0\ar[r]& \widetilde{E_i}\ar[r]& \TET[1]\ar[r]& Q^{l_i}_i [1] \ar[r]&0,
}
$$
the slanted arrow of the diagram, is zero and this can be achieved, if we know that $\OO_X$ is Bridgeland semistable.
This is the subject of the next section.

\section{Bridgeland semistability of $\OO_X$ in $\AC_{H,bH}$}
We fix $H\in Amp^s(X) \cap NS(X)$, that is, $H$ is an ample, integral divisor class with respect to which $\TET$ is stable. We have located the region $B^{un}_H$ in the half plane $\Pi_H$, where the morphism
$$
\Phi^V_{H,bH}: V \otimes \OO_X \longrightarrow \TET[1]
$$
is a monomorpism in $\AC_{cH,b H}=\AC_{H,b H}$, for all $(cH,bH)\in B^{un}_H$. We now investigate semistability of $\OO_X$ for stability conditions $(cH,bH)$ in $B^{un}_H$. 

We begin with an obvious remark:
\begin{equation}\label{subobOX}
	\text{any subobject $E'$ of $\OO_X$ in $\AC_{H,b H}$ is a sheaf in $\TT_{H,bH}$, for any $b<0$.}
\end{equation}

Let $E'$ be a subobject of $\OO_X$ in $\AC_{H,b H}$, where $b<0$. The inclusion $E'\hookrightarrow \OO_X$ gives rise to an exact sequence
$$
\xymatrix@R=12pt@C=12pt{
	0\ar[r]&E' \ar[r]&\OO_X \ar[r]& Q_{E'} \ar[r]&0
}
$$
in $\AC_{H,b H}$. Applying the homological functor $\HH^0$ gives an exact complex in $\AC$:
\begin{equation}\label{cpxE'OX}
\xymatrix@R=12pt@C=12pt{
	0\ar[r]&\HH^{-1}(Q_{E'})\ar[r]&E' \ar[r]&\OO_X \ar[r]&\HH^{0}(Q_{E'})\ar[r]&0.
}
\end{equation}
\begin{lem}\label{lem:E'subobj}
	Let $E'$ be a nonzero proper subobject of $\OO_X$ in $\AC_{H,b H}$.
	Then it is a torsion free sheaf in $\TT_{H,bH}$ with $c_1(E')\cdot H \leq 0$. Furthermore, the equality holds if and only if 
	$E'={\cal J}_A$, the ideal sheaf of a $0$-dimensional subscheme $A$ of $X$.
	\end{lem}
\begin{pf}
The fact that $E'$ is a sheaf in $\TT_{H,bH}$ is \eqref{subobOX}. For the rest we argue according to $\HH^{0}(Q_{E'})$.

	If the sheaf $\HH^{0}(Q_{E'})$ in the complex \eqref{cpxE'OX} is zero, then that complex becomes the exact sequence
	$$
	\xymatrix@R=12pt@C=12pt{
		0\ar[r]&\HH^{-1}(Q_{E'})\ar[r]&E' \ar[r]&\OO_X \ar[r]&0.
	}
$$
As an extension of torsion free sheaves $E'$ is torsion free. From that extension sequence its first Chern class
$$
c_1(E')=c_1(\HH^{-1}(Q_{E'})).
$$
Intersecting with $H$ gives
$$
c_1(E')\cdot H=c_1(\HH^{-1}(Q_{E'}))\cdot H \leq b\, rk(\HH^{-1}(Q_{E'}))H^2,
$$
where the inequality comes from the fact that the sheaf $\HH^{-1}(Q_{E'})$ lies in $\FF_{H,bH}$. Since that sheaf is nonzero (otherwise $E'=\OO_X$, contradicting the assumption that $E'$ is a proper subobject of $\OO_X$), we deduce that $c_1(E')\cdot H $ is negative.

We now turn to the case  $\HH^{0}(Q_{E'})$ is nonzero. We claim that it must be a torsion sheaf. Indeed, consider the standard sequence
$$
\xymatrix@R=12pt@C=12pt{
	0\ar[r]&T\ar[r]&\HH^{0}(Q_{E'})\ar[r]&E'' \ar[r]&0.
}
$$
where $T$ (resp. $E''$) is the torsion (resp. torsion free) part of 
$\HH^{0}(Q_{E'})$. We need to check that $E''$ is zero. Assume it is not zero. Then the epimorphism in the above sequence together with the epimorphism $\OO_X \longrightarrow \HH^{0}(Q_{E'})$ in \eqref{cpxE'OX} give an epimorphism
$$
\OO_X \longrightarrow E''.
$$
Since $E''$ is torsion free, this is an isomorphism. But then $E'\cong \HH^{-1}(Q_{E'})$. Since $E'$ is in $\TT_{H,bH}$ and $\HH^{-1}(Q_{E'})$ is in $\FF_{H,bH}$, this can only occur if both sheaves are zeros and this is contrary to our assumption that $E'$ is nonzero.

Once we know that $\HH^{0}(Q_{E'})$ is a torsion sheaf, the kernel
of $\OO_X \longrightarrow \HH^{0}(Q_{E'})$ in \eqref{cpxE'OX} has the form ${\cal J}_A (-D)$, where $D$ is an effective divisor and ${\cal J}_A$ is the sheaf of ideals of at most $0$-dimensional subscheme $A$ of $X$.
From \eqref{cpxE'OX} we deduce that $E'$ fits into the following exact sequence
\begin{equation}\label{E'ext}
	\xymatrix@R=12pt@C=12pt{
	0\ar[r]&\HH^{-1}(Q_{E'})\ar[r]&E' \ar[r]&{\cal J}_A (-D) \ar[r]&0.
}
\end{equation}

Again $E'$ is an extension of torsion free sheaves and hence is torsion free. The above extension also tells us
$$
c_1(E')=-D+c_1(\HH^{-1}(Q_{E'})).
$$
Intersecting with $H$ gives
$$
c_1(E')\cdot H=-D\cdot H+c_1(\HH^{-1}(Q_{E'}))\cdot H
$$
with the conclusion $c_1(E')\cdot H \leq 0$. Furthermore, the equality holds if and only if $D=0$ and $\HH^{-1}(Q_{E'})=0$. From \eqref{E'ext} we conclude that $E'={\cal J}_A$.
\end{pf}
\begin{pro}\label{pro:OXBrsemist}
	There is a nonempty open subset ${}'\!B^{un}_H$ in $B^{un}_H$ (determined in the proof, see \eqref{B'uns}) such
	that for any stability condition $(cH,bH)\in {}'\!B^{un}_H$  
 the object $\OO_X$ in $\AC_{H,b H}$ is semistable with respect to the central charge $Z_{cH,bH}$.
\end{pro}
\begin{pf}
Let $(cH,bH)$ be a stability condition in $B^{un}_H$ and let $E'$ be a destabilizing subobject of $\OO_X$ with respect to the central charge
$Z_{cH,bH}$.
This means that the phase $\phi_{cH,bH} (E')$ of $E'$ is bigger than the phase $\phi_{cH,bH} (\OO_X)$ of $\OO_X$. Hence the inequality of  values of the cotangent:
$$
\frac{\Re(Z_{cH,bH}(E'))}{\Im(Z_{cH,bH}(E'))}=cot(\pi \phi_{cH,bH} (E')) < cot(\pi\phi_{cH,bH} (\OO_X) )=\frac{\Re(Z_{cH,bH}(\OO_X))}{\Im(Z_{cH,bH}(\OO_X))}
$$
 From the formula for $Z_{cH,bH}$ we deduce
$$
-\frac{ch_2(E')-c_1(E')\cdot H b +\HA rk(E')H^2(b^2-c^2) }{c_1(E')\cdot H-brk(E')H^2} < \frac{b^2-c^2}{2b}.
$$
Since $b<0$ and $(c_1(E')\cdot H-brk(E')H^2) >0$, we obtain
\begin{equation}\label{wallE'ineq}
(c_1(E')\cdot H) (b^2 +c^2) -2ch_2 (E')b >0.
\end{equation}
Observe:
\begin{equation}\label{c1E'negative}
	c_1(E')\cdot H <0.
\end{equation}
Indeed, in view of Lemma \ref{lem:E'subobj} we have to rule out the case
$c_1(E')\cdot H=0$. According to that lemma this occurs only for
$E'={\cal J}_A$, the sheaf of ideals of a $0$-dimensional subscheme $A$ of $X$. But then $ch_2(E')=-deg(A)$ and the inequality \eqref{wallE'ineq} becomes
$$
2deg(A)b>0
$$
and this is impossible, since $b$ is negative. 

With the inequality \eqref{c1E'negative} we rewrite \eqref{wallE'ineq} as follows:
$$
\Big(b-\frac{ch_2(E')}{c_1(E')\cdot H}\Big)^2 +c^2 <\Big(\frac{ch_2(E')}{c_1(E')\cdot H}\Big)^2 
$$
Thus $E'$ gives the semicircular wall 
\begin{equation}\label{WE'OXwall}
W^H_{E'} (\OO_X)=\Big\{ \Big(b-\frac{ch_2(E')}{c_1(E')\cdot H}\Big)^2 +c^2 =\Big(\frac{ch_2(E')}{c_1(E')\cdot H}\Big)^2 \Big\}
\end{equation}
for $\OO_X$ and $E'$ is a destabilizing subobject for $\OO_X$ for all stability conditions $(cH,bH)$ inside the semicircle 
$W^H_{E'} (\OO_X)$.  Observe that $E'$ is also a destabilizing subobject
of $\TET[1]$. This is because we have an inclusion
$$
\OO_X \hookrightarrow V\otimes \OO_X
$$
in $\TT_{H,bH}$ and hence in $\AC_{H,b H}$; composing with
the monomorphism
$$
\Phi^V_{H,bH}:V\otimes \OO_X \longrightarrow \TET[1] 
$$
in $\AC_{H,b H}$ gives the composite monomorphism
$$
\OO_X \longrightarrow \TET[1]
$$
and via this monomorphism we have:
{\small
\begin{equation}\label{unstinOX=unstTh}
	\begin{gathered}
\text{\it for any stability condition $(cH,bH)\in B^{un}_H$, a subobject $E'$ of $\OO_X$ which is destabilizing}
\\
\text{\it with respect to the central charge $Z_{cH,bH} $ is also a destabilizing subobject of $\TET[1]$.}
\end{gathered}
\end{equation}
}

From the structure of walls for $\OO_X$ we know that the walls $W^H_{E'} (\OO_X)$ and $W^H_{\TET[1]} (\OO_X)=W^H_{\OO_X} (\TET[1])$ are nested. We claim that there is a semicircle
\begin{equation}\label{E'wall}
	\begin{gathered}
	W_{max}=\{ (b-b_m)^2 +c^2=b^2_m\}
\\
\text{for some $b_m \in (-\frac{\tau_X}{2K_X \cdot H},0)$ such that the wall $W^H_{\OO_X} (\TET[1])$ is nested}
\\
\text{{\it inside or equal to} the wall $W^H_{E'} (\OO_X)$,}
\\
	\text{for all stability conditions $(cH,bH)\in B^{un}_H$ and lying outside the semicircle $W_{max}$.}
	\end{gathered} 
\end{equation}
To see this consider the actual walls for $\OO_X$ nested {\it inside} the wall
$W^H_{\OO_X} (\TET[1])$. Such walls 
 either 
 
 (i) intersect the vertical line
$b=-\frac{\tau_X}{2K_X \cdot H}$, the line passing through the center of the semicircle $W^H_{\OO_X} (\TET[1])$,

 or
 
 \vspace{0.2cm}
 (ii) lie on the right of that line.
 
 \vspace{0.2cm}
In the latter case the walls are inside or equal to the semicircle
\begin{equation}\label{halfradiuswall}
\Big(b+\frac{\tau_X}{4K_X \cdot H}\Big)^2 +c^2 =\Big(\frac{\tau_X}{4K_X \cdot H}\Big)^2 
\end{equation}
of half of the radius of $W^H_{\OO_X} (\TET[1])$ and centered on the $b$-axis at the point $b=-\frac{\tau_X}{4K_X \cdot H}$.

In the case of walls in (i), that is, the walls intersecting the vertical line
$b=-\frac{\tau_X}{2K_X \cdot H}$, there are only finite number of them in view of \eqref{unstinOX=unstTh} and Lemma \ref{lem:awall-finite}.  We take $W_0$ to be the largest of them. Since it is an actual wall for $\OO_X$, from \eqref{WE'OXwall} it has the form
$$
W_0=\{(b-b_0)^2+c^2=b^2_0\},
$$
for some $b_0 \in (-\frac{\tau_X}{2K_X \cdot H},-\frac{\tau_X}{4K_X \cdot H})$. Thus in either case
all actual walls for $\OO_X$ which are nested inside the wall
$W^H_{\OO_X} (\TET[1])$ are bounded by the maximal semicircle
\begin{equation}\label{Wmax}
W_{max}=\{(b-b_m)^2+c^2=b^2_m\}, \hspace{0.2cm} b_m=min \Big\{b_0, -\frac{\tau_X}{4K_X \cdot H} \Big\}.
\end{equation}
 By definition that maximal semicircle lies inside the wall 
$W^H_{\OO_X} (\TET[1])$. Hence for any stability condition $(cH,bH) \in B^{un}_H$ outside $W_{max}$, the subobjects of $\OO_X$ defining  walls inside $W_{max}$ will no longer be destabilizing for $\OO_X$.
Therefore, for all $(cH,bH) \in B^{un}_H $ lying outside $W_{max}$, the actual walls $W_{E'} (\OO_X)$ must be larger or equal to the wall
 $W_{\OO_X} (\TET[1])$. This means that the center of $W_{E'} (\OO_X)$
 is to the left or equal to the center of $W_{\OO_X} (\TET[1])$. In other words, we have the inequality
 \begin{equation}\label{center-ineq}
 -\frac{\tau_X}{2K_X \cdot H} \geq \frac{ch_2(E')}{c_1(E')\cdot H},
 \end{equation}
 where the right (resp. left) side of the inequality is the $b$-coordinate of the center of $W_{E'} (\OO_X)$ (resp. $W_{\OO_X} (\TET[1])$), see the equation \eqref{WE'OXwall} of the wall $W_{E'} (\OO_X)$ (resp. see Proposition \ref{pro:semicircH} for the equation of the wall $W_{\OO_X} (\TET[1])$). From now on we will be working with the stability conditions in
 \begin{equation}\label{B'uns}
 	{}'B^{un}_H:=\{(cH,bH)\in B^{un}_H | \text{$(cH,bH)$ is outside $W_{max}$} \},
 \end{equation}
where $W_{max}$ is the semicircle in \eqref{Wmax}.
 With this in mind we go back to the inequality \eqref{center-ineq}. Remembering that $(c_1(E')\cdot H)$ is negative, see \eqref{c1E'negative}, that inequality reads
 \begin{equation}\label{ch2E'}
 ch_2(E')\geq (-c_1(E')\cdot H)\frac{\tau_X}{2K_X \cdot H}.
\end{equation} 
 This inequality tells us:
 \begin{equation}\label{b-Bunst}
 	\begin{gathered}
 	\text{\it for stability conditions $(cH,bH)$ in ${}'B^{un}_H$, the region defined in \eqref{B'uns},}
 	\\
 	\text{\it the sheaf $E'$ is Bogomolov unstable.}
 \end{gathered}
 \end{equation}

 Indeed, the Bogomolov semistability of $E'$ says
 $$
 2ch_2(E')\leq \frac{c^2_1(E')}{rk(E')} .
 $$
 This and the Hodge Index inequality 
 $$
c^2_1(E') H^2 \leq {(c_1(E') \cdot H)^2}
$$
give
$$
2ch_2(E')\leq \frac{(c_1(E') \cdot H)^2}{rk(E')H^2}.
$$
Combining with the inequality \eqref{ch2E'} we obtain
$$
\frac{(-c_1(E')\cdot H)}{rk(E')H^2}\geq \frac{\tau_X}{K_X \cdot H}.
$$
But $E'$ is a sheaf in $\TT_{H,bH}$ and hence subject to
$$
c_1(E')\cdot H >rk(E')H^2b.
$$
This and the previous inequality give
$$
b <-\frac{\tau_X}{K_X \cdot H}.
$$
But this is outside $B^{un}_H$.
 This proves \eqref{b-Bunst}.

Once $E'$ is Bogomolov unstable, it is $H$-unstable and we take its
HN filtration
$$
E'=E'_l\supset E'_{l-1} \supset \cdots \supset E'_1 \supset E'_0 =0
$$
with respect to $H$. Since $E'\in \TT_{H,bH}$, the factors of the filtration
$$
Q'_i =E'_i/E'_{i-1}
$$
are also in $\TT_{H,bH}$ and we have
\begin{equation}\label{c1Q'ib}
c_1 (Q'_i)\cdot H >brk(Q'_i)H^2
\end{equation}
for all $i\in[1,l]$. On the other hand we have
\begin{equation}\label{c1Q'inegative}
c_1 (Q'_i)\cdot H \leq 0,\,\,\forall i\in [1,l].
\end{equation}
This is seen as follows:
from the HN filtration of $E'$ we have
\begin{equation}\label{E'-HN}
\frac{c_1 (Q'_1)\cdot H}{rk(Q'_1)}>\cdots >\frac{c_l (Q'_l)\cdot H}{rk(Q'_l)}.
\end{equation}
So it is enough to check that $c_1 (Q'_1)\cdot H \leq 0$. By definition
$Q'_1=E'_1$ is the maximal $H$-destabilizing subsheaf of $E'$ and it defines the destabilizing exact sequence 
$$
\xymatrix@R=12pt@C=12pt{
	0\ar[r]&E'_1 \ar[r]&E'\ar[r]&E'/E'_1 \ar[r]&0.
}
$$
On the other hand from the proof of Lemma \ref{lem:E'subobj} we know that
$E'$ fits into the following exact sequence
$$
\xymatrix@R=12pt@C=12pt{
	0\ar[r]&\HH^{-1}(Q_{E'}) \ar[r]&E'\ar[r]&{\cal J}_A (-D) \ar[r]&0.
}
$$
Putting the two exact sequences together gives the diagram
$$
\xymatrix@R=12pt@C=12pt{
	&&0\ar[d]&&\\
	&&E'_1\ar[d]&&\\
	0\ar[r]&\HH^{-1}(Q_{E'})\ar[r]&E'\ar[r]\ar[d]&{\cal J}_A (-D)\ar[r]&0\\
	&&E'/E'_1\ar[d]&&\\
	&&0&&
}
$$
 The resulting composite arrow
$$
s:E'_1 \longrightarrow {\cal J}_A (-D)
$$
is nonzero, since otherwise the vertical arrow factors through
$\HH^{-1}(Q_{E'})$ giving a nonzero morphism
$$
 E'_1 \longrightarrow \HH^{-1}(Q_{E'});
 $$
 but $E'_1$ is in $\TT_{H,bH}$ and $\HH^{-1}(Q_{E'}) \in \FF_{H,bH}$ and no nonzero morphism is allowed between such sheaves. 
 
 Once the composite arrow $s$ is nonzero, its image, being a subsheaf of ${\cal J}_A (-D)$ , must be of the form ${\cal J}_{A'} (-D')$, for some effective divisor $D'$ and the ideal sheaf ${\cal J}_{A'} $ of a subscheme $A'$ of dimension at most $0$.
 Thus $E'_1$ fits into the following exact sequence
 $$
 \xymatrix@R=12pt@C=12pt{
 	0\ar[r]&ker(s) \ar[r]&E'_1 \ar[r]&{\cal J}_{A'} (-D') \ar[r]&0.
 }
 	$$
 	In addition, $ker(s)$ from the above diagram is a subsheaf of $\HH^{-1}(Q_{E'}) \in \FF_{H,bH}$, hence it is also in $\FF_{H,bH}$. From the exact sequence above we deduce
 	$$
 	c_1(E'_1)=c_1 (ker(s)) -D'.
 	$$
 	The intersection with $H$ gives
 	$$
 	c_1(E'_1)\cdot H=c_1 (ker(s))\cdot H -D'\cdot H \leq c_1 (ker(s))\cdot H \leq rk(ker(s))bH^2 \leq 0.
 	$$
 	
 	We now use the HN filtration of $E'$ to estimate $ch_2 (E')$. Namely, the additivity of the Chern character gives
 	$$
 	2ch_2(E')=\sum^{l}_{i=1}2ch_2 (Q'_i) \leq \sum^{l}_{i=1}\frac{c^2_1 (Q'_i)}{rk(Q'_i)},
 	$$
 	where the inequality comes from Bogomolov-Gieseker inequality for 
 	$Q'_i$'s. This and the Hodge Index inequality
 	$c^2_1(Q'_i)H^2 \leq (c_1(Q'_i)\cdot H)^2$ give
 	$$
 		2ch_2(E') \leq \sum^{l}_{i=1}\frac{c^2_1 (Q'_i)}{rk(Q'_i)}\leq 
 		\sum^{l}_{i=1}\frac{(c_1(Q'_i)\cdot H)^2}{rk(Q'_i)H^2},
 	$$
 Using \eqref{c1Q'inegative} and \eqref{E'-HN} we deduce
 $$
 \frac{(c_1(Q'_i)\cdot H)^2}{rk(Q'_i)H^2} \leq \frac{(-c_1(Q'_l)\cdot H)}{rk(Q'_l)H^2} (-c_1(Q'_i)\cdot H) ,
 $$
 for all $i$ and where the inequality is strict for all $i<l$. Substituting into the previous inequality we obtain
 $$
 	2ch_2(E') < \frac{(-c_1(Q'_l)\cdot H)}{rk(Q'_l)H^2} (-c_1(E')\cdot H).
 	$$
 	Combining with the inequality in \eqref{ch2E'} gives
 	$$
 	\frac{(-c_1(Q'_l)\cdot H)}{rk(Q'_l)H^2} >\frac{\tau_X}{K_X \cdot H}.
 	$$	

 Recalling the inequality in \eqref{c1Q'ib} expressing that $Q'_l$ is in $\TT_{H,bH}$, we again deduce
$$
b <-\frac{\tau_X}{K_X \cdot H}.
$$
But this takes us outside of the wall
$$
W^H_{\OO_X}  (\TET[1])=\Big \{ \Big(b+\frac{\tau_X}{2K_X \cdot H}\Big)^2 +c^2 =\Big(\frac{\tau_X}{2K_X \cdot H}\Big)^2 \Big\}
$$
and hence outside of $B^{un}_H$. 
Thus we proved that for all stability conditions in ${}'\!B^{un}_H$, see \eqref{B'uns}, the object $\OO_X$ is Bridgeland semistable.
\end{pf}

We now have the region $ {}'\!B^{un}_H$ of the stability conditions $(cH,bH)$ for which

\vspace{0.2cm}
$\bullet$ the phase $\phi_{cH,bH} (\OO_X)$ of $\OO_X$ is bigger than the phase $\phi_{cH,bH} (\TET[1])$,

\vspace{0.2cm}
$\bullet$ le morphism $\Phi^V_{H,bH}: V\otimes \OO_X \longrightarrow \TET[1]	$ is injective in $\AC_{H,b H}$ and hence $\TET[1]$ is unstable with respect to $(cH,bH) \in  {}'\!B^{un}_H$ and we have the HN filtration
\begin{equation}\label{Brfilt6}
\TET[1]=E^{c,b}_n \supset E^{c,b}_{n-1} \supset \cdots \supset E^{c,b}_1 \supset E^{c,b}_0 =0   
\end{equation}
of $\TET[1]$ in $\AC_{H,bH}$,

\vspace{0.2cm}
$\bullet$ $V\otimes \OO_X$ is a semistable object	 in $\AC_{H,bH}$ with respect to $(cH,bH)$ in $ {}'\!B^{un}_H$.

\vspace{0.2cm}
These properties imply that the morphism $\Phi^V_{H,bH}$ must factor through the intermediate objects of the HN filtration in the second item. More precisely, let 
$$
A^{c,b}_i:=E^{c,b}_{i}/E^{c,b}_{i-1} 
$$
be the $i$-th semistable factor, for $i\in [1,n]$, of the HN filtration in the second item and denote its phase with respect to the central charge $Z_{cH,bH}$  by $\phi^{c,b}_i$. These values form a strictly decreasing function of $i$, that is, we have
$$
\phi^+_{c,b}(\TET[1]):=\phi^{c,b}_1 >\cdots> \phi^{c,b}_n:=\phi^-_{c,b}(\TET[1]).
$$
From the first item  the phase $\phi^{c,b}(\TET[1])$ is less than the phase $\phi^{c,b}(\OO_X)=\phi^{c,b} (V\otimes \OO_X)$. Hence the inequality
$$
\phi^{c,b} (V\otimes \OO_X) >\phi^{c,b} (\TET[1]) >\phi^-_{c,b}(\TET[1]) =\phi^{c,b}_n.
$$
From this it follows that the composition morphism, the slanted arrow in the following diagram,
$$
\xymatrix@R=18pt@C=12pt{
	&&V\otimes \OO_X \ar[d]_{\Phi^V_{H,bH}}\ar[rd]&&\\
0\ar[r]	&E^{c,b}_{n-1}\ar[r]&\TET[1]\ar[r]&A^{c,b}_n \ar[r]&0
}
$$
is zero. Thus the vertical arrow, the morphism  $\Phi^V_{H,bH}$, factors through $E^{c,b}_{n-1}$. This motivates the following.

Set
\begin{equation}\label{iV}
	\begin{gathered}
E^{c,b}_{i_{{}_V}}:=\text{the lowest subobject in the filtration \eqref{Brfilt6}}
	\\
	\text{through which factors the morphism $\Phi^V_{H,bH}$,}
\end{gathered}
\end{equation}
that is, $\Phi^V_{H,bH}$ factors through $E^{c,b}_{i_{{}_V}}$ but does not factor through $E^{c,b}_{i_{{}_V}-1}$. In particular, $\Phi^V_{H,bH}$ induces a nonzero morphism
$$
\Big(\Phi^V_{H,bH}\Big)_{i_{{}_V}} :V\otimes \OO_X  \longrightarrow A^{c,b}_{i_{{}_V}}
$$
and this implies
\begin{equation}\label{phaseOXbound}
	\phi^{c,b} (\OO_X)\leq  \phi^{c,b}_{i_{{}_V}}.
\end{equation}

The monomorphism  
$$
\Phi^V_{H,bH} :V\otimes \OO_X \longrightarrow \TET[1]
$$
is completed to the exact sequence
$$
\xymatrix@R=18pt@C=25pt{
0\ar[r]&V\otimes \OO_X \ar[r]^{\text{\tiny$\Phi^V_{H,bH}$}}&\TET[1]\ar[r]&F_V [1]\ar[r]&0
}
$$
in $\AC_{H,bH}$ which we already encountered in \eqref{PhiV-extseq}. We also have the destabilizing sequence
$$
\xymatrix@R=12pt@C=12pt{
0\ar[r]&E^{c,b}_{i_{{}_V}}\ar[r]&\TET[1]\ar[r]& F_{i_{{}_V}}[1]\ar[r]&0.
}
$$
Putting the above two sequences together gives the diagram
$$
\xymatrix@R=12pt@C=12pt{
	&&0\ar[d]&&\\
	&&V\otimes \OO_X \ar[d]&&\\
	0\ar[r]&E^{c,b}_{i_{{}_V}}\ar[r]&\TET[1]\ar[r] \ar[d]& F_{i_{{}_V}}[1]\ar[r]&0\\
	&&F_V[1]\ar[d]&&\\
	&&0&&
}
$$
In addition, the monomorphism in the column is our $\Phi^V_{H,bH}$ and we know that it factors through $E^{c,b}_{i_{{}_V}}$, so the above diagram can be completed as follows.
\begin{equation}\label{iV-diag}
\xymatrix@R=12pt@C=12pt{
	&0\ar[d]&0\ar[d]&&\\
	&V\otimes \OO_X \ar[d]\ar@{=}[r]&V\otimes \OO_X \ar[d]&&\\
	0\ar[r]&E^{c,b}_{i_{{}_V}}\ar[r]\ar[d]&\TET[1]\ar[r] \ar[d]& F_{i_{{}_V}}[1]\ar[r]\ar@{=}[d]&0\\
	0\ar[r]&E'_{i_{{}_V}}\ar[r]\ar[d]&F_V[1]\ar[d]\ar[r]&F_{i_{{}_V}}[1]\ar[r]&0\\
	&0&0&&
}
\end{equation}	
where $E'_{i_{{}_V}}$ is the cokernel of the monomorphism
$$
V\otimes \OO_X \longrightarrow E^{c,b}_{i_{{}_V}}
$$
induced by $\Phi^V_{H,bH}$. In other words $V\otimes \OO_X$ is a subobject of $E^{c,b}_{i_{{}_V}}$.

 The above diagram establishes a relation of the `effective' destabilizing subobject $V\otimes \OO_X$ of $\TET[1]$ and  destabilizing subobjects of $\TET[1]$ coming from the HN filtration of $\TET[1]$.
We will use that diagram to investigate further properties of $E^{c,b}_{i_{{}_V}}$. 

\section{Surfaces with $\displaystyle{\alpha_X < \frac{3}{8} }$}
 We recall that the assumption $\displaystyle{\alpha_X < \frac{3}{8} }$
 on the ratio of the Chern numbers allows to control the HN filtration of $\TET[1]$ for the canonical polarization. From now on we assume that $H=K_X$. Then Proposition \ref{pro:HNalpha} tells us:
 
 a) all subobjects $E^{c,b}_i$ of the HN filtration in \eqref{Brfilt6} are sheaves in $\TT_{K_X,bK_X}$,
 
 b) $E^{c,b}_i$ has rank $0$ or $1$, unless the quotient object
 $F_i[1]=\TET[1]/E^{c,b}_i$ is Bogomolov unstable,
 
 c) for $F_i[1]=\TET[1]/E^{c,b}_i$ Bogomolov unstable, the HN filtration
 of $F_i$
 $$
 F_i=F^{l_i}_i \supset \cdots \supset  F^{1}_i \supset F^0_i =0
 $$
 with respect to $K_X$ has the last quotient $Q^{l_i}_i :=F^{l_i}_i /F^{l_i-1}_i $ of rank $2$ or $3$.
 
 \vspace{0.2cm}
 We go back to $E^{c,b}_{i_{{}_V}}$, the lowest degree subobject through which
 the monomorphism $\Phi^V_{K_X,bK_X} : V\otimes \OO_X \longrightarrow \TET[1]$ factors. According to item b) above, the sheaf $E^{c,b}_{i_{{}_V}}$ has rank $0$ or $1$, unless the quotient object $F_{i_{{}_V}}[1]$ in the diagram \eqref{iV-diag} is Bogomolov unstable. In the latter case we can replace $F_{i_{{}_V}}$ by its last quotient
 $Q^{l_{i_{{}_V}}}_{i_{{}_V}}$, see the item c) above, and have the diagram similar to the one in \eqref{iV-diag}. More precisely, we combine the
 destabilizing sequence
 $$
 \xymatrix@R=12pt@C=12pt{
 	0\ar[r]&E^{c,b}_{i_{{}_V}}\ar[r]&\TET[1]\ar[r]& F_{i_{{}_V}}[1]\ar[r]&0.
 }
 $$
 with the exact sequence
 $$
  \xymatrix@R=12pt@C=12pt{
 	0\ar[r]&F^{l_{i_{{}_V}} -1}_{i_{{}_V}}[1]\ar[r]& F_{i_{{}_V}}[1]\ar[r]& Q^{l_{i_{{}_V}}}_{i_{{}_V}}[1] \ar[r]&0
 }
$$
to obtain the diagram which we encountered in \eqref{diag-red} and which we reproduce below
\begin{equation}\label{EiVtild-diag}
\xymatrix@R=12pt@C=12pt{
	&&0\ar[d]&0\ar[d]&\\
0\ar[r]	&E_{i_{{}_V}}\ar[r]\ar@{=}[d]&\widetilde{E_{i_{{}_V}}}\ar[r]\ar[d]&F^{l_{i_{{}_V}} -1}_{i_{{}_V}}[1]\ar[d]\ar[r]&0\\
	0\ar[r]&E_{i_{{}_V}}\ar[r]&\TET[1]\ar[r]\ar[d]& F_{i_{{}_V}}[1]\ar[r]\ar[d]&0\\
	&&Q^{l_{i_{{}_V}}}_{i_{{}_V}}[1] \ar@{=}[r]\ar[d]&Q^{l_{i_{{}_V}}}_{i_{{}_V}}[1] \ar[d]&\\
	&&0&0&	
}
\end{equation}
We also recall that $\widetilde{E_{i_{{}_V}}}$ is a sheaf in $\TT_{K_X,bK_X}$, see Lemma \ref{lem:tildEi}.

Applying the functor $Hom_{\AC_{K_X,bK_X}}(V\otimes \OO_X,\bullet)$ we
obtain the commutative square
$$
\xymatrix@R=12pt@C=12pt{
Hom_{\AC_{K_X,bK_X}}(V\otimes \OO_X,\TET[1])\ar[r]\ar[d]& Hom_{\AC_{K_X,bK_X}}(V\otimes \OO_X,F_{i_{{}_V}}[1])\ar[d]\\
Hom_{\AC_{K_X,bK_X}}(V\otimes \OO_X,Q^{l_{i_{{}_V}}}_{i_{{}_V}}[1])\ar@{=}[r]& Hom_{\AC_{K_X,bK_X}}(V\otimes \OO_X,Q^{l_{i_{{}_V}}}_{i_{{}_V}}[1])
}
$$
from which it follows that the monomorphism
$\Phi^V_{K_X,bK_X} \in Hom_{\AC_{K_X,bK_X}}(V\otimes \OO_X,\TET[1])$ 
maps to zero in $Hom_{\AC_{K_X,bK_X}}(V\otimes \OO_X,Q^{l_{i_{{}_V}}}_{i_{{}_V}}[1])$. This means that $\Phi^V_{K_X,bK_X}$ factors through $\widetilde{E_{i_{{}_V}}}$, see the middle column in the diagram, and gives the diagram similar to the one
in \eqref{iV-diag}. Namely, we have

\begin{equation}\label{iV-tilde-diag}
	\xymatrix@R=12pt@C=12pt{
		&0\ar[d]&0\ar[d]&&\\
		&V\otimes \OO_X \ar[d]\ar@{=}[r]&V\otimes \OO_X \ar[d]&&\\
		0\ar[r]&\widetilde{E_{i_{{}_V}}}\ar[r]\ar[d]&\TET[1]\ar[r] \ar[d]& Q^{l_{i_{{}_V}}}_{i_{{}_V}}[1]\ar[r]\ar@{=}[d]&0\\
		0\ar[r]&\widetilde{E_{i_{{}_V}}}'\ar[r]\ar[d]&F_V[1]\ar[d]\ar[r]&Q^{l_{i_{{}_V}}}_{i_{{}_V}}[1]\ar[r]&0\\
		&0&0&&
	}
\end{equation}

We summarize the above discussion in the following statement.
\begin{pro}\label{pro:QV-23}
	Assume $\displaystyle{\alpha_X < \frac{3}{8} }$. Then 
	for every stability condition $(cK_X,bK_X)\in {}'\!B^{un}_{K_X}$
	there is a distinguished quotient object $Q^{l_{i_{{}_V}}}_{i_{{}_V}}[1]$ of $\TET[1]$ with
$Q^{l_{i_{{}_V}}}_{i_{{}_V}}$ a $K_X$-semistable sheaf of rank $2$ or $3$ lying in $\FF_{K_X,bK_X}$.
It fits into the commutative diagram \eqref{iV-tilde-diag}, where
$\widetilde{E_{i_{{}_V}}}$ is a sheaf of rank $0$ or $1$ in $\TT_{K_X,bK_X}$.
	\end{pro}

The left column of the diagram \eqref{iV-tilde-diag} will give us the
effectiveness of $c_1 (\widetilde{E_{i_{{}_V}}})$ missing in Lemma \ref{lem:tildEi}.
In addition, the diagram will give a better geometric and cohomological meaning of those divisor classes.

 \vspace{1.5cm}
 \noindent
 {\bf 7.1. From the diagram \eqref{iV-tilde-diag} to the category $\AC$.}
We begin by translating the diagram \eqref{iV-tilde-diag} in  $\AC_{K_X,bK_X}$ into the one in $\AC$ by applying the homological functor $\HH^0$. This gives the diagram of coherent sheaves

\begin{equation}\label{QV-diag-sh}
	\xymatrix@R=12pt@C=12pt{
		&&0\ar[d]&0\ar[d]&\\
		&&\HH^{-1}(\widetilde{E_{i_{{}_V}}}')\ar@{=}[r]\ar[d]&\HH^{-1}(\widetilde{E_{i_{{}_V}}}')\ar[d] &\\
		0\ar[r]&\TET \ar[r]\ar@{=}[d]&F_V \ar[r]\ar[d]& V\otimes \OO_X \ar[r]\ar[d]&0\\
		0\ar[r]&\TET \ar[r]&Q^{l_{i_{{}_V}}}_{i_{{}_V}} \ar[r]\ar[d]& \widetilde{E_{i_{{}_V}}} \ar[r]\ar[d]&0\\
			&&\HH^0(\widetilde{E_{i_{{}_V}}}')\ar@{=}[r]\ar[d]&\HH^0(\widetilde{E_{i_{{}_V}}}')\ar[d]&&\\		
		&&0&0&
	}
\end{equation}
In this diagram 

-the row in the middle is the functor $\HH^0$ applied to 
the column in the middle of \eqref{iV-tilde-diag}; 

-the bottom row is the functor $\HH^0$ applied to 
the row in the middle of \eqref{iV-tilde-diag};

- the right and the middle columns are obtained by applying  the functor $\HH^0$ respectively to the left column and bottom row of \eqref{iV-tilde-diag}.

Observe that the middle row is the extension sequence
corresponding to the inclusion $V\subset H^1 (\TET)$ under the identifications
$$
Hom_{\DD}(V\otimes \OO_X, \TET[1])\cong Ext^1 (V\otimes \OO_X, \TET)\cong Hom(V,H^1(\TET)).
$$
So the above diagram is a translation of Bridgeland instability of $\TET[1]$ into instability properties of the vector bundle $F_V$, the middle term of the extension sequence in the middle row of  \eqref{QV-diag-sh}.

We proceed to investigate the geometric and cohomological consequences of this diagram according to the rank of $Q^{l_{i_{{}_V}}}_{i_{{}_V}}$.

\vspace{0.3cm}
\noindent
{\bf 7.2. Case: $rk(Q^{l_{i_{{}_V}}}_{i_{{}_V}})=2$.} To simplify the notation we set
$$
F'_{i_{{}_V}}:=\HH^{-1}(\widetilde{E_{i_{{}_V}}}').
$$
From the middle column of the diagram 	\eqref{QV-diag-sh} it follows that  it is a subsheaf of $F_V$ of rank $rk(F'_{i_{{}_V}})=dim(V)$. Hence it is nonzero. Furthermore, it is a second syzygy sheaf and therefore locally free. We simplify the diagram to
\begin{equation}\label{rkQV=2diag}
\xymatrix@R=12pt@C=12pt{
	&&0\ar[d]&0\ar[d]&\\
&&F'_{i_{{}_V}}\ar[d]\ar@{=}[r]&F'_{i_{{}_V}} \ar[d]&\\	
0\ar[r]&\TET\ar[r]\ar@{=}[d]&F_V \ar[r]\ar[d]&V\otimes \OO_X \ar[r]\ar[d]&0\\
0\ar[r]&\TET \ar[r]&Q'_{i_{{}_V}}\ar[r]\ar[d]&{S_{i_{{}_V}}}\ar[r]\ar[d]&0\\
&&0&0&
}
\end{equation}
where $Q'_{i_{{}_V}}$ is the quotient sheaf $F_V /F'_{i_{{}_V}}$. 

\begin{lem}\label{lem:DV}
	1) ${S_{i_{{}_V}}}$ is a torsion sheaf whose support $D_{i_{{}_V}}=supp(S_{i_{{}_V}})$ is a subscheme of pure dimension $1$. It is generated by its global sections: the vector space $V$ is identified with a subspace of global sections of $S_{i_{{}_V}}$ and that subspace globally generates $S_{i_{{}_V}}$; in addition, this identification fits into the commutative diagram
	$$
\xymatrix@R=12pt@C=12pt{
	V \ar@{^{(}->}[r]  \ar@{^{(}->}[d]&H^1(\TET)\ar@{=}[d]	\\
	H^0 (S_{i_{{}_V}})\ar@{^{(}->}[r]&H^1(\TET)
}
$$
	
	2) The quotient sheaf $Q'_{i_{{}_V}}$ in \eqref{rkQV=2diag} is torsion free.
	
	3) Let $\delta_{i_{{}_V}}$ be a global section of $\OO_X (D_{i_{{}_V}})$ defining the divisor $D_{i_{{}_V}}$, that is, $D_{i_{{}_V}}=(\delta_{i_{{}_V}}=0)$ is the zero locus of $\delta_{i_{{}_V}}$. Then $\delta_{i_{{}_V}}$ annihilates the subspace $V\subset H^1 (\TET)$, this is to say that $V$ is contained in the kernel of the
	homomorphism
	$$
	H^1 (\TET)\stackrel{\delta_{i_{{}_V}}}{\longrightarrow} H^1 (\TET(D_{i_{{}_V}}))
	$$
	defined by the multiplication morphism
	$\xymatrix@R=12pt@C=12pt{
		\TET \ar^(.39){\delta_{i_{{}_V}}}[r]&\TET(D_{i_{{}_V}}).}
		$
		 Furthermore, it is a minimal divisor having this property, that is, no proper subdivisor $D'$ of $D_{i_{{}_V}}$ annihilates $V$.
		 
		4) The degree of $D_{i_{{}_V}}$ with respect to $K_X$ is subject to
		the inequality
		$$
		D_{i_{{}_V}} \cdot K_X <\frac{\sqrt{1+16(3\alpha_X-1)} -1}{2}K^2_X=
		\frac{8(3\alpha_X-1)K^2_X}{\sqrt{1+16(3\alpha_X-1)} +1} <4(3c_2(X)-K^2_X).
		$$
	\end{lem}
\begin{pf}
	By definition $S_{i_{{}_V}}$ is a subsheaf of $\widetilde{E_{i_{{}_V}}}$ in the diagram \eqref{QV-diag-sh} and this one has no torsion subsheaves supported on a $0$-dimensional subscheme. So $S_{i_{{}_V}}$ has no subsheaves of supported on a $0$-dimensional subscheme as well. Hence the support of $S_{i_{{}_V}}$ has no $0$-dimensional components and this means that $D_{i_{{}_V}}$, the support of $S_{i_{{}_V}}$, is a subscheme of pure dimension $1$.
	
	The epimorphism $V\otimes \OO_X \longrightarrow S_{i_{{}_V}}$ in the right 
	column of \eqref{rkQV=2diag} tells us that $V$ is identified with a subspace of $H^0 (S_{i_{{}_V}})$ and $S_{i_{{}_V}}$ is globally generated by that subspace. 
	
	The commutative diagram in 1) of the lemma comes from cohomology sequences of the two rows of the diagram \eqref{rkQV=2diag}:
	$$
\xymatrix@R=12pt@C=12pt{
	&0\ar[d]&\\
	0\ar[r]&V\ar[r]\ar[d]& H^1 (\TET)\ar@{=}[d]\\
		0\ar[r]&H^0 (S_{i_{{}_V}})\ar[r]& H^1 (\TET)
	}
$$	
The injectivity of arrows at $V$ and $H^0(S_{i_{{}_V}})$ in the above diagram is the vanishing of $H^0 (F_V)$, $H^0(F'_{i_{{}_V}})$ and $H^0(Q'_{i_{{}_V}})$ and this is because the sheaves involved are in $\FF_{K_X,bK_X}$.
	
	The part 2), saying that $Q'_{i_{{}_V}}$ is torsion free, follows from the fact that $Q'_{i_{{}_V}}$ is a subsheaf of the torsion free sheaf $Q^{l_{i_{{}_V}}}_{i_{{}_V}}$ in \eqref{QV-diag-sh}.

	We now turn to the part 3) of the lemma. For this we observe that
	the diagram \eqref{rkQV=2diag} gives
	$$
	c_1 (F'_{i_{{}_V}})=-K_X -c_1 (Q'_{i_{{}_V}})=-c_1 (S_{i_{{}_V}})=-D_{i_{{}_V}}.
	$$
	Taking the exterior power of degree $v:=rk(F'_{i_{{}_V}})=dim(V)$ of the diagram \eqref{rkQV=2diag} we obtain
	\begin{equation}\label{deltaV-diag}
	\xymatrix@R=12pt@C=12pt{
		&&0\ar[d]&&\\
		&&\OO_X (-D_{i_{{}_V}})\ar[d] \ar[dr]^(.65){\delta_{i_{{}_V}}}&&\\ 
	0\ar[r]& P_V	\ar[r]& \bigwedge^{v} F_V \ar[r]& \OO_X \ar[r]& 0
	}
\end{equation}
where the slanted arrow is the multiplication by a global section $\delta_{i_{{}_V}}$ defining the divisor $D_{i_{{}_V}}$. The sheaf $P_V$
in the above diagram fits into the following exact sequence
$$
\xymatrix@R=12pt@C=12pt{
0\ar[r]&\OO_X (-K_X)\ar[r]& P_V	\ar[r]&\bigwedge^{v-1}V \otimes \TET \ar[r]&0.
}
$$
Observe: $1\in H^0 (\OO_X)$ in the diagram \eqref{deltaV-diag} goes under the coboundary map
$$
\mbox{$
H^0 (\OO_X)\longrightarrow H^1 (P_V) \hookrightarrow \bigwedge^{v-1}V \otimes H^1(\TET) \cong V^{\ast}\otimes H^1 (\TET)$
}
$$
to the cohomology class
$\xi_V \in Ext^1(V\otimes\OO_X,\TET) \cong V^{\ast}\otimes H^1 (\TET)$ corresponding to the extension sequence
$$
\xymatrix@R=12pt@C=12pt{
0\ar[r]&\TET\ar[r]&F_V \ar[r]&V\otimes \OO_X\ar[r]&0,
}
$$
the middle row of the diagram \eqref{rkQV=2diag}.

Tensoring the diagram \eqref{deltaV-diag} with $\OO_X (D_{i_{{}_V}})$ and using an identification $\bigwedge^{v-1}V \cong V^{\ast}$, we obtain
that the extension class $\xi_V \in Ext^1 (V\otimes\OO_X, \TET)\cong V^{\ast}\otimes H^1 (\TET)$ goes to zero in $V^{\ast}\otimes H^1 (\TET(D_{i_{{}_V}}))$ under the homomorphism
$$
\xymatrix@R=12pt@C=38pt{
V^{\ast}\otimes H^1 (\TET)\ar^(.45){id_{V^{\ast}} \otimes \delta_{i_{{}_V}}}[r]&  V^{\ast}\otimes H^1 (\TET(D_{i_{{}_V}})).
}
$$
The above can be restated as
$$
V \subset ker\Big(H^1(\TET)\stackrel{\delta_{i_{{}_V}}}{\longrightarrow}H^1 (\TET(D_{i_{{}_V}}))\Big).
$$
This proves the statement that $\delta_{i_{{}_V}}$ annihilates the subspace $V\subset H^1 (\TET)$.

Next we need to check that no proper subdivisor $D'\subset D_{i_{{}_V}}$ annihilates $V$: assume such a subdivisor $D'$ exists; then the
vertical arrow in the diagram \eqref{deltaV-diag} factors through
$\OO_X (-D')$ and this means that the quotient sheaf
$coker(\OO_X (-D_{i_{{}_V}})\longrightarrow \bigwedge^v F_V)$ has a torsion part 
supported on $(D_{i_{{}_V}}-D')$; but that quotient is locally free outside of the singularity locus of $Q'_{i_{{}_V}}$ which is, according to 2) of the lemma, at most a $0$-dimensional subscheme. So the support of $(D_{i_{{}_V}} -D')$ is at most $0$-dimensional and this contradicts the assumption that $D'$ is a proper subdivisor of $D_{i_{{}_V}}$.

The part 4) of the lemma asserting an upper bound for $D_{i_{{}_V}}\cdot K_X$ comes from	the relation of $S_{i_{{}_V}}$ to the sheaf $\widetilde{E_{i_{{}_V}}} $ in the diagram \eqref{QV-diag-sh}. Namely, we have the exact sequence
$$
\xymatrix{
	0\ar[r]&S_{i_{{}_V}} \ar[r]& \widetilde{E_{i_{{}_V}}} \ar[r]&\HH^0(\widetilde{E_{i_{{}_V}}}')\ar[r]&0
}
$$
implying 
$$
c_1 (\widetilde{E_{i_{{}_V}}}) -c_1 (S_{i_{{}_V}})=c_1 (\widetilde{E_{i_{{}_V}}}) -D_{i_{{}_V}}=\text{an effective divisor.}
$$
Hence the inequality
$$
D_{i_{{}_V}}\cdot K_X \leq c_1 (\widetilde{E_{i_{{}_V}}}) \cdot K_X.
$$
The asserted upper bound now follows from the inequality
$$
\frac{c_1 (\widetilde{E_{i_{{}_V}}}) \cdot K_X}{K^2_X}< \frac{\sqrt{16(3\alpha_X -1)+1}-1}{2}
$$
in the proof of Lemma \ref{lem:tildEi}, see \eqref{solving-ineq}.
\end{pf}

We continue to be in the setting of Lemma \ref{lem:DV}. The lemma tells us that
all cohomology classes in the subspace $V$ are annihilated by $D_{i_{{}_V}}$.
In particular, taking a nonzero class $\xi$ in $V$, we have the extension
sequence 
$$
\xymatrix@R=12pt@C=12pt{
	0\ar[r]& \TET \ar[r]& F_{[\xi]} \ar[r]& \OO_X \ar[r]& 0
}
$$
corresponding to $[\xi]$, the one dimensional subspace of $H^1(\TET)$ spanned by $\xi$. Tensoring with $\OO_X (D_{i_{{}_V}})$ gives us the cohomology sequence
$$
\xymatrix@R=12pt@C=12pt{
 H^0 (F_{[\xi]} (D_{i_{{}_V}}))\ar[r]& H^0 (\OO_X (D_{i_{{}_V}})) \ar[r]^{\xi}& H^1 (\TET(D_{i_{{}_V}})),
}
$$
where the last map on the right is the cup-product with the cohomology class $\xi$. The fact that $\xi$ is annihilated by $D_{i_{{}_V}}$ means that
the global section $\delta_{i_{{}_V}}$ of $\OO_X (D_{i_{{}_V}})$ defining $D_{i_{{}_V}}$, see Lemma \ref{lem:DV}, 3), lies in the kernel of that cup-product. Hence it comes from a global section $s_{\xi}$ of $F_{[\xi]} (D_{i_{{}_V}})$. Equivalently, we think of this section as a morphism
$$
s_{\xi}:\OO_X (-D_{i_{{}_V}})\longrightarrow  F_{[\xi]} 
$$
which fits into the commutative diagram
\begin{equation}\label{diag-sxi}
\xymatrix@
R=12pt@C=12pt{
	&&0\ar[d]&&\\
	&&\OO_X (-D_{i_{{}_V}})\ar[d]_{s_{\xi}} \ar[dr]^(.55){\delta_{i_{{}_V}}}&&\\
	0\ar[r]& \TET \ar[r]& F_{[\xi]} \ar[r]& \OO_X \ar[r]& 0.
}
\end{equation}
The minimality property in Lemma \ref{lem:DV}, 3), can be rephrased as
follows.
\begin{lem}\label{lem:cokersxi}
	There is a dense Zariski open subset $U_{D_{i_{{}_V}}}$ of $\PP(V)$ parametrizing those $[\xi]$ in $\PP(V)$ for which the divisor $D_{i_{{}_V}}$ is precisely the divisor annihilating $\xi$, that is, no proper subdivisor $D$ of $D_{i_{{}_V}}$ annihilates $\xi$. Equivalently, the cokernel 
	$$
	Q_{[\xi]} :=coker(s_{\xi})
	$$
	of the morphism $s_{\xi}$ in the diagram \eqref{diag-sxi} is torsion free and fits into the exact sequence
	$$
	\xymatrix@
	R=12pt@C=12pt{
		0\ar[r]&\TET \ar[r]&Q_{[\xi]} \ar[r]&\OO_{D_{i_{{}_V}}} \ar[r]&0.
	}
		$$
 \end{lem}
\begin{pf}
	Consider the set
	$$
	\text{${\mathfrak{C}}_{D_{i_{{}_V}}}:=$ the set of proper nonzero components of $D_{i_{{}_V}}$.}
	$$
	This is a finite set. For every $D\in  {\mathfrak{C}}_{D_{i_{{}_V}}}$ we let
	$$
	V_D :=ker \Big(V\stackrel{\delta_D}{\longrightarrow} H^1(\TET(D))\Big),
	$$
	where $\delta_D$ is a defining section of $D$. This is a linear subspace of $V$ and Lemma \ref{lem:DV}, 3), tells us that it is a proper subspace of $V$. This and the finiteness of the set ${\mathfrak{C}}_{D_{i_{{}_V}}}$ imply
	$$
V':=	\bigcup_{D\in {\mathfrak{C}}_{D_{i_{{}_V}}}} V_D
	$$
	is a proper closed subvariety of $V$. Hence the complement
	$$
	\widetilde{U}_{D_{i_{{}_V}}}=V\setminus V'
	$$
	is the Zariski dense open subset of $V$ parametrizing nonzero cohomology classes in $V$ which are annihilated by $D_{i_{{}_V}}$ and by no  proper nonzero component of $D_{i_{{}_V}}$. 
	
	The set of lines
	$$
	 U_{D_{i_{{}_V}}}=\{[\xi] \in \PP(V)|\,\, \xi \in \widetilde{U}_{D_{i_{{}_V}}}\}
	 $$
	 gives the Zariski dense open subset of the projective space
	 $\PP(V)$ parametrizing the diagrams in \eqref{diag-sxi} for which
	 the cokernel of $s_{\xi}$ is torsion free. Indeed, the diagram \eqref{diag-sxi} can be expanded as follows
	 \begin{equation}\label{diag-sxi-exp}
	 	\xymatrix@
	 	R=12pt@C=12pt{
	 		&&0\ar[d]&0\ar[d]&\\
	 		&&\OO_X (-D_{i_{{}_V}})\ar[d]_{s_{\xi}} \ar@{=}[r]&\OO_X (-D_{i_{{}_V}})\ar[d]^{\delta_{{i_{{}_V}}}}&\\
	 		0\ar[r]& \TET \ar[r]\ar@{=}[d]& F_{[\xi]} \ar[r]\ar[d]& \OO_X \ar[r]\ar[d]& 0\\
	 		0\ar[r]&\TET\ar[r]&Q_{[\xi]}\ar[r]\ar[d]&\OO_{D_{i_{{}_V}}} \ar[r]\ar[d]&0\\
	 		&&0&0&
	 	}
	 \end{equation}
where $Q_{[\xi]}$ denotes the cokernel of $s_{\xi}$. Assume the torsion part $Tor(Q_{[\xi]})$ of $Q_{[\xi]}$ is nonzero. It must inject into $\OO_{D_{i_{{}_V}}}$, see the bottom exact sequence of \eqref{diag-sxi-exp}. This means that the support $T_{\xi}$ of $Tor(Q_{[\xi]})$ is a component of $D_{i_{{}_V}}$. Factoring out the torsion part gives the diagram

\begin{equation}\label{diag-sxi-exp1}
	\xymatrix@
	R=12pt@C=12pt{
		&&0\ar[d]&0\ar[d]&\\
		&&\OO_X (-D_{\xi})\ar[d]_{s'_{\xi}} \ar@{=}[r]&\OO_X (-D_{\xi})\ar[d]^{\delta_{D_{\xi}}}&\\
		0\ar[r]& \TET \ar[r]\ar@{=}[d]& F_{[\xi]} \ar[r]\ar[d]& \OO_X \ar[r]\ar[d]& 0\\
		0\ar[r]&\TET\ar[r]&Q'_{[\xi]}\ar[r]\ar[d]&\OO_{D_{\xi}} \ar[r]\ar[d]&0\\
		&&0&0&
	}
\end{equation}
where $D_{\xi}=D_{i_{{}_V}} -T_{\xi}	$ and $Q'_{[\xi]}=Q_{[\xi]} /Tor(Q_{[\xi]})$.
This implies that $\xi$ is annihilated by the component $D_{\xi}$ of $D_{i_{{}_V}}$. From the first part of the proof
$D_{\xi}=D_{i_{{}_V}}$, for all $[\xi] \in U_{D_{i_{{}_V}}}$. In other words, $T_{\xi}=0$ or, equivalently, $Q_{[\xi]}$ is torsion free, for all $[\xi] \in U_{D_{i_{{}_V}}}$. 

The exact sequence in the last assertion of the lemma is the one on the bottom of the diagram \eqref{diag-sxi-exp}.
\end{pf}

We will now spell out the relation of the one dimensional extension corresponding to a line $[\xi]$, the middle row in the diagram \eqref{diag-sxi-exp}, and the extension corresponding to the inclusion $V\subset H^1 (\TET)$, the middle row of
\eqref{rkQV=2diag}.

The inclusion of the line $[\xi]$ in $V$ gives rise to the inclusion of the extensions

\begin{equation}\label{Fxi-FV}
	\xymatrix@
R=12pt@C=12pt{
	&&0\ar[d]&0\ar[d]&\\
	0\ar[r]& \TET \ar[r]\ar@{=}[d]& F_{[\xi]} \ar[r]\ar[d]& \OO_X \ar[r]\ar[d]& 0\\
		0\ar[r]& \TET \ar[r]& F_{V} \ar[r]\ar[d]& V\otimes \OO_X \ar[r]\ar[d]& 0\\	
		&&V/[\xi] \otimes \OO_X \ar@{=}[r]\ar[d]&V/[\xi] \otimes \OO_X \ar[d]&\\
		&&0&0
	}
\end{equation} 
The exact sequence in the middle column combined with the vertical column in \eqref{rkQV=2diag} give the diagram
$$
\xymatrix@R=12pt@C=12pt{
	&&0\ar[d]&&\\
	&&F'_{i_{{}_V}} \ar[d]\ar@{=}[r]&F'_{i_{{}_V}} \ar[d]&\\
	 0\ar[r] &F_{[\xi]} \ar[r]\ar@{=}[d]&F_{V} \ar[r]\ar[d]& V/[\xi] \otimes \OO_X \ar[r] \ar[d]&0\\
		&F_{[\xi]} \ar[r]&Q'_{i_{{}_V}} \ar[d] \ar[r]&S_{[\xi]}\ar[r]\ar[d]&0\\
	&&0&0&
}
$$

The resulting morphism 
\begin{equation}\label{FxiQ'V}
\xymatrix@R=12pt@C=12pt{
F_{[\xi]} \ar[r]&Q'_{i_{{}_V}}
}
\end{equation}
fits into the diagram

$$
\xymatrix@R=12pt@C=12pt{
	0\ar[d]&0\ar[d]\\
	\TET\ar[d]\ar@{=}[r]&\TET \ar[d]&\\	
F_{[\xi]} \ar[r]\ar[d]&Q'_{i_{{}_V}} \ar[d]\ar[r]&S_{[\xi]}\ar[r]&0\\
\OO_X \ar[d]	&S_{i_{{}_V}} \ar[d]&\\
0&0&
}
$$

\vspace{1.2cm}
\noindent
From this it follows that the kernel of the morphism in \eqref{FxiQ'V} is a line bundle denoted $\OO_X (-D)$ and the above diagram can be completed as follows

\begin{equation}\label{FxiQ'V-diag}
\xymatrix@R=12pt@C=12pt{
	&&0\ar[d]&0\ar[d]\\
	&&\TET\ar[d]\ar@{=}[r]&\TET \ar[d]&\\	
	0\ar[r]&\OO(-D)\ar[r]\ar@{=}[d]&F_{[\xi]} \ar[r]\ar[d]&Q'_{i_{{}_V}} \ar[d]\ar[r]&S_{[\xi]}\ar[r]\ar@{=}[d]&0\\
	0\ar[r]&\OO_X (-D)\ar[r]&\OO_X \ar[d]\ar[r]	&S_{i_{{}_V}} \ar[d]\ar[r]&S_{[\xi]}\ar[r]&0\\
	&&0&0&
}
\end{equation}
From this we deduce that $D$ is a nonzero component of $D_{i_{{}_V}}$ and $D$  annihilates $\xi$. We also can interpret geometrically the sheaf $S_{[\xi]}$: 
from Lemma \ref{lem:DV}, 1), the vector space $V$ is identified with the subspace of global sections of $S_{i_{{}_V}}$; from this perspective, the morphism $\OO_X \longrightarrow S_{i_{{}_V}}$ in the bottom row of the above diagram is defined by $\xi$, viewed as a corresponding global section of $S_{i_{{}_V}}$, and
the sheaf $S_{[\xi]}$ is supported on the zero locus of that global section. With those remarks in mind and Lemma \ref{lem:cokersxi} we can characterize the Zariski open subset $U_{D_{i_{{}_V}}}$ in that lemma as follows.

\begin{lem}\label{lem:V-secSV}
	With the identification of $V$ as a subspace of global sections of
	$S_{i_{{}_V}}$ in Lemma \ref{lem:DV}, 1), the Zariski open subset $U_{D_{i_{{}_V}}}$ in Lemma \ref{lem:cokersxi} parametrizes global sections of
	$S_{i_{{}_V}}$ in $V$ whose zero locus is at most $0$-dimensional. 
	
	For every $[\xi]\in U_{D_{i_{{}_V}}}$, the sheaf $S_{i_{{}_V}}$ fits into the exact sequence
	$$
	\xymatrix@R=12pt@C=12pt{
			0\ar[r]&\OO_{D_{i_{{}_V}}}\ar[r]^{\xi}&S_{i_{{}_V}} \ar[r]&S_{[\xi]}\ar[r]&0,
	}
	$$
	where $\xi$ is viewed as a global section of $S_{i_{{}_V}}$ and 	the quotient sheaf $S_{[\xi]}$ is supported on the zero locus of that global section and that zero locus is at most $0$-dimensional subscheme of $X$.
\end{lem}

\begin{pf}
	From Lemma \ref{lem:cokersxi} it follows that $[\xi]\in U_{D_{i_{{}_V}}}$ if and only if $D$ in the diagram \eqref{FxiQ'V-diag} equals $D_{i_{{}_V}}$:
	$$
\xymatrix@R=12pt@C=12pt{
	&&0\ar[d]&0\ar[d]\\
	&&\TET\ar[d]\ar@{=}[r]&\TET \ar[d]&\\	
	0\ar[r]&\OO(-D_{i_{{}_V}})\ar[r]\ar@{=}[d]&F_{[\xi]} \ar[r]\ar[d]&Q'_{i_{{}_V}} \ar[d]\ar[r]&S_{[\xi]}\ar[r]\ar@{=}[d]&0\\
	0\ar[r]&\OO_X (-D_{i_{{}_V}})\ar[r]&\OO_X \ar[d]\ar[r]	&S_{i_{{}_V}} \ar[d]\ar[r]&S_{[\xi]}\ar[r]&0\\
	&&0&0&
}
$$
This implies that the bottom row can be broken into the following two exact sequences
\begin{equation}\label{SV-DV-seq-pf}
\xymatrix@R=12pt@C=12pt{
	0\ar[r]&\OO_X (-D_{i_{{}_V}})\ar[r]&\OO_X \ar[r]&\OO_{D_{i_{{}_V}}} \ar[r]&0,\\
	0\ar[r]&\OO_{D_{i_{{}_V}}}\ar[r]&S_{i_{{}_V}} \ar[r]&S_{[\xi]}\ar[r]&0.
}
\end{equation}
From $c_1 (S_{i_{{}_V}})=c_1 (\OO_{D_{i_{{}_V}}})=D_{i_{{}_V}}$ and the second exact sequence above we deduce that $S_{[\xi]}$ is supported on at most $0$-dimensional subscheme of $D_{i_{{}_V}}$. That exact sequence is the one in the last assertion of the lemma. 		
\end{pf}

\begin{lem}\label{lem:Q'Vstable}
	1) The sheaf $Q'_{i_{{}_V}}$ in \eqref{rkQV=2diag} is $K_X$-stable and
	$(-c_1 (Q'_{i_{{}_V}}))$ lies in the positive cone of $X$.
	
	2) The Chern numbers of $Q'_{i_{{}_V}}$ are subject to the inequality of Miyaoka
	$$
	3c_2 (Q'_{i_{{}_V}}) -c^2_1 (Q'_{i_{{}_V}})\geq -\frac{1}{4}(c^{-}_1 (Q'_{i_{{}_V}}) )^2 +3deg(Sing(Q'_{i_{{}_V}})),
	$$
	where $(-c_1 (Q'_{i_{{}_V}}))=c^{+}_1 (Q'_{i_{{}_V}}) + c^{-}_1 (Q'_{i_{{}_V}})$ is the Zariski decomposition of the divisor class $(-c_1 (Q'_{i_{{}_V}}))$ with
	$c^{\pm}_1 (Q'_{i_{{}_V}})$ denoting respectively the positive and negative parts of the decomposition and $Sing(Q'_{i_{{}_V}})$ is the singularity sheaf of $Q'_{i_{{}_V}}$.
	
	3) The Chern numbers of $Q'_{i_{{}_V}}$ and $X$ are related as follows
	$$
	3c_2 (X)-K^2_X = 	(3c_2 (Q'_{i_{{}_V}}) -c^2_1 (Q'_{i_{{}_V}})) +K_X \cdot D_{i_{{}_V}} -2D^2_{i_{{}_V}} +3deg(S_{[\xi]}),
	$$
	where $S_{[\xi]} $ as described in Lemma \ref{lem:V-secSV}.
\end{lem}
\begin{pf} We begin with the first assertion of the lemma.
	Assume $Q'_{i_{{}_V}}$ is not $K_X$-stable. This means that it admits a (semi)destabilizing with respect to $K_X$ exact sequence
	\begin{equation}\label{Q'V-destseq}
\xymatrix@R=12pt@C=12pt{
	0\ar[r]&Q'_1 \ar[r]& Q'_{i_{{}_V}} \ar[r]& Q'_2 \ar[r]&0,
}
\end{equation}
where $Q'_i$ ($i=1,2$) are sheaves of rank $1$ in $\FF_{K_X,bK_X}$ and
the (semi)destabilizing means
$$
c_1 (Q'_1)\cdot K_X \geq \HA c_1 (Q'_{i_{{}_V}})\cdot K_X.
$$
This and the Hodge Index inequality imply
\begin{equation}\label{Q'1-bound}
	c^2_1 (Q'_1)\leq \frac{(-c_1 (Q'_1)\cdot K_X)^2}{K^2}\leq \frac{(-c_1 (Q'_{i_{{}_V}})\cdot K_X)^2}{4K^2}.
\end{equation}
For $ Q'_2 $ we have
\begin{equation}\label{Q'2-bound}
	c^2_1 (Q'_2)\leq 0.
\end{equation}
This is because the dual of the exact sequence \eqref{Q'V-destseq} gives the inclusion
$(Q'_2)^{\ast} \hookrightarrow  (Q'_{i_{{}_V}})^{\ast}$ and the dual of the bottom row in the diagram \eqref{rkQV=2diag} tell us that
$(Q'_{i_{{}_V}})^{\ast}$ is a subsheaf of $\Omega_X$; hence
$(Q'_2)^{\ast}$ is a line subsheaf of $\Omega_X$; furthermore
$$
c_1 ((Q'_2)^{\ast}) \cdot K_X =-c_1(Q'_2)\cdot K_X \geq (-b)K^2_X >0,
$$
where the first inequality comes from the fact that $Q'_2$ is a sheaf in $\FF_{K_X,bK_X}$. Since $\Omega_X$ has no line subsheaves lying in the positive cone, the inequality \eqref{Q'2-bound} follows.
	
From the exact sequence \eqref{Q'V-destseq} we deduce
$$
2ch_2 (Q'_{i_{{}_V}})=2ch_2 (Q'_1)+2ch_2 (Q'_2)\leq c^2_1 (Q'_1) + c^2_1 (Q'_2).
$$
This together with \eqref{Q'1-bound} and \eqref{Q'2-bound} give the estimate
\begin{equation}\label{ch2Q'V-ub-pf}
2ch_2 (Q'_{i_{{}_V}}) \leq \frac{(-c_1 (Q'_{i_{{}_V}})\cdot K_X)^2}{4K^2}.
\end{equation}
From the bottom sequence in \eqref{rkQV=2diag} we have
$$
2ch_2 (Q'_{i_{{}_V}})=\tau_X +2ch_2 (S_{i_{{}_V}}).
$$
From the exact sequences in \eqref{SV-DV-seq-pf} 
$$
2ch_2 (S_{i_{{}_V}})=-D^2_{i_{{}_V}} +2deg(S_{[\xi]})\geq -D^2_{i_{{}_V}}
$$
and we obtain the lower bound
$$
2ch_2 (Q'_{i_{{}_V}})=\tau_X +2ch_2 (S_{i_{{}_V}}) \geq \tau_X -D^2_{i_{{}_V}}.
$$
Combining this with the upper bound on $2ch_2 (Q'_{i_{{}_V}}) $ in \eqref{ch2Q'V-ub-pf} gives
$$
\tau_X \leq \frac{(-c_1 (Q'_{i_{{}_V}})\cdot K_X)^2}{4K^2} +D^2_{i_{{}_V}}.
$$
From the bottom sequence in \eqref{rkQV=2diag} we also have
$$
c_1(Q'_{i_{{}_V}})=-K_X +D_{i_{{}_V}}.
$$
Substituting for $c_1 (Q'_{i_{{}_V}})$ into the previous inequality gives
$$
\tau_X \leq \frac{(K^2_X -D_{i_{{}_V}}\cdot K_X)^2}{4K^2} +D^2_{i_{{}_V}} =\frac{1}{4}K^2_X -\HA D_{i_{{}_V}} \cdot K_X +\frac{(D_{i_{{}_V}} \cdot K_X)^2}{4K^2_X} +D^2_{i_{{}_V}}.
$$
This and the Hodge Index inequality
$$
D^2_{i_{{}_V}} \leq \frac{(D_{i_{{}_V}} \cdot K_X)^2}{K^2_X}
$$
imply
$$
\frac{5}{4} \frac{(D_{i_{{}_V}} \cdot K_X)^2}{K^2_X} -\HA D_{i_{{}_V}} \cdot K_X +(\frac{1}{4}K^2_X -\tau_X) \geq 0.
$$

Setting $t=\frac{D_{i_{{}_V}} \cdot K_X}{K^2_X}$, the above inequality becomes
$$
5t^2 -2t -(3-8\alpha_X)\geq 0.
$$
Solving for $t$ we obtain
$$
\frac{D_{i_{{}_V}} \cdot K_X}{K^2_X}=t\geq \frac{1+\sqrt{1+5(3-8\alpha_X)}}{5}.
$$
On the other hand we have the upper bound
$$
\frac{D_{i_{{}_V}} \cdot K_X}{K^2_X} < \frac{\sqrt{1+16(3\alpha_X -1)}-1}{2},
$$
see Lemma \ref{lem:DV}, 4).
Combining the two inequalities gives
$$
 \frac{\sqrt{1+16(3\alpha_X -1)}-1}{2} > \frac{1+\sqrt{1+5(3-8\alpha_X)}}{5}
 $$
 or, equivalently,
 \begin{equation}\label{alpha-function}
 \frac{\sqrt{1+16(3\alpha_X -1)}-1}{2} - \frac{1+\sqrt{1+5(3-8\alpha_X)}}{5} > 0.
\end{equation}
 The expression on the left is an increasing function of $\alpha_X$, for   all positive  $\alpha_X$, where the function is defined, and we are concerned with the values 
 $$
 \frac{1}{3}<\alpha_X <\frac{3}{8}.
 $$
  Evaluating that function at $\alpha_X=\frac{3}{8}$ gives
 $$
 \frac{\sqrt{3} -1}{2} -\frac{2}{5}=
 \frac{5\sqrt{3} -9}{10} <0.
 	$$
 Thus the inequality \eqref{alpha-function} is invalid in the interval	
 $(\frac{1}{3},\frac{3}{8})$. Hence $Q'_V$ must be $K_X$-stable.
 
 \vspace{0.2cm} 
 We turn to the part 3) of the lemma.
 From the bottom sequence of the diagram \eqref{rkQV=2diag} it follows
 $$
 2ch_2(Q'_{i_{{}_V}})=\tau_X +2ch_2 (S_{i_{{}_V}}),
 $$
 while the exact sequence 
 $$
 \xymatrix@R=12pt@C=12pt{
 	0\ar[r]& \OO_{D_{i_{{}_V}}} \ar[r]& S_{i_{{}_V}} \ar[r]& S_{[\xi]}\ar[r]&0
 }
$$
in Lemma \ref{lem:V-secSV} gives
$$
2ch_2(S_{i_{{}_V}})=-D^2_{i_{{}_V}} -2c_2(S_{[\xi]})=-D^2+2deg(S_{[\xi]}).
$$
Substituting into the previous equation gives
\begin{equation}\label{ch2Q'iV-tauX}
2ch_2(Q'_{i_{{}_V}})=\tau_X +2ch_2 (S_{i_{{}_V}}) = \tau_X-D^2_{i_{{}_V}} +2deg(S_{[\xi]}).
\end{equation}
This is rewritten as follows
$$
\begin{gathered}
2c_2(X)-K^2_X =2c_2(Q'_{i_{{}_V}})-c^2_1 (Q'_{i_{{}_V}})-D^2_{i_{{}_V}} +2deg(S_{[\xi]})
\\
=
2c_2(Q'_{i_{{}_V}})-\frac{2}{3}c^2_1 (Q'_{i_{{}_V}}) -\frac{1}{3}c^2_1 (Q'_{i_{{}_V}})-D^2_{i_{{}_V}} +2deg(S_{[\xi]})
\\
=\frac{2}{3} (3c_2(Q'_{i_{{}_V}})-c^2_1 (Q'_{i_{{}_V}})) -\frac{1}{3}c^2_1 (Q'_{i_{{}_V}})-D^2_{i_{{}_V}} +2deg(S_{[\xi]}).
\end{gathered}
$$
Substituting $c_1 (Q'_{i_{{}_V}})=-K_X +D_{i_{{}_V}}$, see the bottom sequence \eqref{rkQV=2diag}, in the second term of last expression we obtain
$$
2c_2(X)-K^2_X =\frac{2}{3} (3c_2(Q'_{i_{{}_V}})-c^2_1 (Q'_{i_{{}_V}})) -\frac{1}{3}(K_X-D_{i_{{}_V}})^2-D^2_{i_{{}_V}} +2deg(S_{[\xi]}).
$$
This is rewritten as 
$$
\frac{2}{3} (3c_2 (X)-K^2_X)=\frac{2}{3} (3c_2(Q'_{i_{{}_V}})-c^2_1 (Q'_{i_{{}_V}})) +\frac{2}{3}K_X\cdot D_{i_{{}_V}}-\frac{4}{3} D^2_{i_{{}_V}} +2deg(S_{[\xi]}).
$$
or, equivalently,
$$
 (3c_2 (X)-K^2_X)= (3c_2(Q'_{i_{{}_V}})-c^2_1 (Q'_{i_{{}_V}})) +K_X\cdot D_{i_{{}_V}}-2 D^2_{i_{{}_V}} +3deg(S_{[\xi]})
 $$
 as asserted in 3) of the lemma.
 
 The above also implies the assertion in 1) that $(-c_1 (Q'_{i_{{}_V}}))$ in the positive cone of $X$: the $K_X$-stability implies the Bogomolov-Gieseker inequality
 $$
 c^2_1(Q'_{i_{{}_V}}) \geq 4ch_2 (Q'_{i_{{}_V}}),
 $$
 while the positivity of $ch_2 (Q'_{i_{{}_V}})$ follows from \eqref{ch2Q'iV-tauX} and the upper bound
 $$
 D_{i_{{}_V}} \cdot K_X <4(3c_2 (X)-K^2_X)
 $$
 in Lemma \ref{lem:DV}, 4); this upper bound also implies the inequality
 $$
 (-c_1 (Q'_{i_{{}_V}}))\cdot K_X =(K_X -D_{i_{{}_V}})\cdot K_X >\HA K^2_X.
 $$
 The proof of part 1) of the lemma is now complete.
 
 For part 2) of the lemma, we dualize the bottom sequence \eqref{rkQV=2diag} to deduce
 that the dual
 $Q^{'\ast}_{i_{{}_V}}$ of $Q'_{i_{{}_V}}$ is a subsheaf of rank $2$ of $\Omega_X$.
 Furthermore, from part 1) the divisor class $(-c_1 (Q'_{i_{{}_V}}))=c_1(Q^{'\ast}_{i_{{}_V}})$ lies  
 in the positive cone of $X$ and hence admits the Zariski decomposition
 $$
 (-c_1 (Q'_{i_{{}_V}}))=c_1(Q^{'\ast}_{i_{{}_V}})=c^{+}_1 (Q'_{i_{{}_V}}) +c^{-}_1 (Q'_{i_{{}_V}}),
 $$
 where $c^{\pm}_1 (Q'_{i_{{}_V}})$ denotes respectively the positive and negative parts of the Zariski decomposition of $(-c_1 (Q'_{i_{{}_V}}))$.
 A result of Miyaoka tells us that the Chern numbers of $Q^{'\ast}_{i_{{}_V}}$ are subject to the inequality
 \begin{equation}\label{chMiyaoka}
 3c_2(Q^{'\ast}_{i_{{}_V}})-c^2_1(Q^{'\ast}_{i_{{}_V}}))\geq -\frac{1}{4}(c^{-}_1 (Q'_{i_{{}_V}}))^2,
\end{equation}
see \cite{Mi1}, Remark 4.18. If $Q'_{i_{{}_V}}$ is locally free, then the above inequality is the one asserted in 2) of the lemma, since
 $$
 3c_2(Q^{'\ast}_{i_{{}_V}})-c^2_1(Q^{'\ast}_{i_{{}_V}}))=3c_2(Q'_{i_{{}_V}})-c^2_1(Q'_{i_{{}_V}}))
 $$
 when $Q'_{i_{{}_V}}$ is locally free. In general, $Q'_{i_{{}_V}}$ is torsion free
 and its double dual $(Q'_{i_{{}_V}})^{\ast \ast}$ is locally free and
 we have the exact sequence
 \begin{equation}\label{Q'V-dd}
 \xymatrix@R=12pt@C=12pt{
 0\ar[r]& Q'_{i_{{}_V}} \ar[r]& (Q'_{i_{{}_V}})^{\ast \ast} \ar[r]& Sing(Q'_{i_{{}_V}}) \ar[r]&0,
}
\end{equation}
where $Sing(Q'_{i_{{}_V}})$ is the sheaf supported on the $0$-dimensional subscheme where $Q'_{i_{{}_V}}$ fails to be locally free. Dualizing this sequence we obtain
$$
Q^{'\ast}_{i_{{}_V}} \cong {\cal H}om ((Q'_{i_{{}_V}})^{\ast \ast}, \OO_X),
$$
since ${\cal E}xt^i (Sing(Q'_{i_{{}_V}}), \OO_X)=0$, for $i=0,1$, because the support of $Sing(Q'_{i_{{}_V}})$ is $0$-dimensional.	Hence we have
$$
3c_2(Q^{'\ast}_{i_{{}_V}})-c^2_1(Q^{'\ast}_{i_{{}_V}}))=3c_2(Q^{'\ast\ast}_{i_{{}_V}})-c^2_1(Q^{'\ast \ast}_{i_{{}_V}}))
$$
and the inequality \eqref{chMiyaoka} can be rewritten as
$$
 3c_2(Q^{'\ast\ast}_{i_{{}_V}})-c^2_1(Q^{'\ast\ast}_{i_{{}_V}}))\geq -\frac{1}{4}(c^{-}_1 (Q'_{i_{{}_V}}))^2.
 $$
 From the exact sequence \eqref{Q'V-dd} it follows
 $$
 c_1(Q^{'\ast \ast}_{i_{{}_V}})=c_1(Q'_{i_{{}_V}}), \hspace{0.2cm}
 c_2(Q^{'\ast \ast}_{i_{{}_V}})=c_2(Q'_{i_{{}_V}}) +c_2(Sing(Q'_{i_{{}_V}}))=c_2(Q'_{i_{{}_V}}) -deg(Sing(Q'_{i_{{}_V}})).
 $$
 Substituting into the previous inequality gives
 $$
 3c_2(Q^{'}_{i_{{}_V}})-c^2_1(Q^{'}_{i_{{}_V}}))\geq -\frac{1}{4}(c^{-}_1 (Q'_{i_{{}_V}}))^2 +3deg(Sing(Q'_{i_{{}_V}})),
 $$
 as asserted in 2) of the lemma.
\end{pf}

The exact sequence
	\begin{equation}\label{FV-dest-seq-lem}
	\xymatrix@R=12pt@C=12pt{
		0\ar[r]&F'_{i_{{}_V}} \ar[r]& F_V \ar[r]&Q'_{i_{{}_V}} \ar[r]&0,
	}
\end{equation}
the middle column of the diagram \eqref{rkQV=2diag} behaves as a destabilizing sequence of the sheaf $F_V$ in the middle. The following statement confirms it.

\begin{lem}\label{lem:Q'Vdual-maxdist}
	Assume the dimension of the vector space $V$ is at least $2$. Then
	the exact sequence \eqref{FV-dest-seq-lem}
 is a $K_X$-destabilizing
	sequence for the vector bundle $F_V$. Furthermore, the quotient sheaf $Q'_{i_{{}_V}}$ is the last (semi)stable factor of the HN filtration of $F_V$ with respect to $K_X$, see Proposition \ref{pro:FVunst}, unless 
	$F'_{i_{{}_V}}$ is $K_X$-unstable and its last $K_X$-semistable factor has rank $1$. If this happens, then
	$$
	D_{i_{{}_V}} \cdot K_X > \frac{K^2_X}{3}
	$$
	and the ratio $\alpha_X$ of the Chern numbers of $X$ must be subject to
	$$
	\frac{1}{9} <3\alpha_X -1 <\frac{1}{8}.
	$$
	\end{lem}
\begin{pf}
To see that the sequence \eqref{FV-dest-seq-lem} is $K_X$-destabilizing for 
$F_V$ we need to check the inequality
$$
-D_{i_{{}_V}} \cdot K_X =c_1(F'_{i_{{}_V}}) \cdot K_X	> -\frac{dim(V)}{dim(V) +2}K^2_X.
$$
Assume the inequality fails. Then we have
$$
D_{i_{{}_V}} \cdot K_X =-c_1(F'_{i_{{}_V}}) \cdot K_X	\geq \frac{dim(V)}{dim(V) +2}K^2_X \geq \HA K^2,
$$
where the last inequality uses the facts that $\frac{dim(V)}{dim(V) +2}$ is an increasing function of $dim(V)$ and the values of $dim(V)\geq 2$.
On the other hand the upper bound on $D_{i_{{}_V}} \cdot K_X $ from Lemma \ref{lem:DV}, 4), tells us
$$
4(3c_2(X)-K^2_X) >D_{i_{{}_V}} \cdot K_X.
$$
Combining with the previous inequality gives
$$
4(3c_2(X)-K^2_X) >\HA K^2.
$$
Rewriting this in terms of the ratio $\alpha_X=\frac{c_2(X)}{K^2_X}$ of the Chern numbers of $X$ we obtain
$$
\alpha_X > \frac{3}{8}.
$$
This contradicts the assumption that $\alpha_X < \frac{3}{8}$. This proves that the sequence \eqref{FV-dest-seq-lem} is $K_X$-destabilizing for $F_V$.

The above argument also says that
 $$
F_V \supset F'_{i_{{}_V}} \supset 0
$$
is the HN filtration of $F_V$ with respect to $K_X$, unless $F'_{i_{{}_V}}$ is $K_X$-unstable. So
if $Q'_{i_{{}_V}}$ fails to be the last $K_X$-semistable factor of $F_V$, then $F'_{i_{{}_V}}$ must be $K_X$-unstable and we take its
 HN filtration with respect to $K_X$:
 $$
 F'_{i_{{}_V}}=F'_m \supset F'_{m-1}\supset \cdots \supset F'_1 \supset F'_0=0.
 $$
  The condition that $Q'_{i_{{}_V}}$ is not the last $K_X$-semistable factor of $F_V$ comes down to the following inequality of slopes 
 \begin{equation}\label{slopeQ'V-m}
 \frac{c_1 (Q'_{i_{{}_V}})\cdot K_X}{2} =\mu_{K_X} (F_V/F'_m) \geq \mu_{K_X} (F'_m /F'_{m-1})=\frac{c_1 (F'_m/F'_{m-1})\cdot K_X}{r_m}
\end{equation}
 where $r_i =rk(F'_i /F'_{i-1})$, for $i\in [1,m]$. 
 
 From the first part of the proof we know
 $$
 \frac{c_1 (Q'_{i_{{}_V}})\cdot K_X}{2} <-\frac{D_{i_{{}_V}} \cdot K_X}{2}<\frac{c_1 (F'_m /F'_{m-1})\cdot K_X}{2},
 $$
where the second inequality comes from the fact that
\begin{equation}\label{DiV-sum}
-D_{i_{{}_V}}=c_1(F'_{i_{{}_V}}) =\sum^m_{i=1} c_1 (F'_i/F'_{i-1})
\end{equation}
and all summands have negative degree with respect to $K_X$; this last property holds because $F'_V$ is a sheaf in $\FF_{K_X,bK_X}$. This and the inequality \eqref{slopeQ'V-m} imply $r_m =1$. The inequality
\eqref{slopeQ'V-m}  now reads
$$
\frac{c_1 (Q'_{i_{{}_V}})\cdot K_X}{2} \geq -L_m \cdot K_X,
$$
where $(-L_m):=c_1 (F'_m/F'_{m-1})$.
 Hence the inequality
\begin{equation}\label{c1Q'VLm}
-c_1 (Q'_{i_{{}_V}})\cdot K_X \leq 2L_m \cdot K_X.
\end{equation}
This and the formula \eqref{DiV-sum} imply
$$
-c_1 (Q'_{i_{{}_V}})\cdot K_X \leq 2L_m \cdot K_X<2D_{i_{{}_V}} \cdot K_X.
$$
Combining together with the formula
$$
K_X =-c_1 (Q'_{i_{{}_V}}) + D_{i_{{}_V}},
$$
we deduce
$$
K^2_X <3D_{i_{{}_V}} \cdot K_X.
$$
Hence the inequality for $D_{i_{{}_V}} \cdot K_X$ asserted in the lemma.
 
 The restrictions on $(3\alpha_X -1)$ are seen as follows. The upper bound is equivalent to the inequality
 $$
 \alpha_X <\frac{3}{8};
 $$
 the lower bound is obtained by combining the inequality
 $$
 D_{i_{{}_V}} \cdot K_X > \frac{1}{3}K^2_X
 $$
 obtained above with the upper bound 
 $$
  D_{i_{{}_V}} \cdot K_X < \frac{\sqrt{1+16(3\alpha_X -1)}-1}{2} K^2_X
  $$
  in Lemma \ref{lem:DV}, 4):
  $$
  \frac{1}{3} <\frac{\sqrt{1+16(3\alpha_X -1)}-1}{2} \Leftrightarrow
  3\alpha_X -1 >\frac{1}{9}.
  $$
\end{pf}

	We summarize all the properties of rank $2$ case in the following statement.
\begin{pro}\label{pro:QVrk2}
 Assume $\displaystyle{\alpha_X < \frac{3}{8} }$ and the rank of $Q^{l_{i_{{}_V}}}_{i_{{}_V}}$ in Proposition \ref{pro:QV-23} is $2$. Then the extension
 $$
 \xymatrix@R=12pt@C=12pt{
 	0\ar[r]& \TET \ar[r]& F_V \ar[r]& V\otimes \OO_X \ar[r]& 0
 }
$$
corresponding to the natural inclusion $V\subset H^1 (\TET)$ fits into the diagram \eqref{rkQV=2diag} as the middle row of that diagram. Furthermore, the following properties hold.

1) The sequence
$$
	\xymatrix@R=12pt@C=12pt{
	0\ar[r]&F'_{i_{{}_V}} \ar[r]& F_V \ar[r]&Q'_{i_{{}_V}} \ar[r]&0,
}
$$
the middle column of the diagram \eqref{rkQV=2diag}, gives rise to
 the decomposition of the canonical divisor $K_X$
$$
K_X =D_{i_{{}_V}} + (-c_1 (Q'_{i_{{}_V}})),
$$
where $(-c_1 (Q'_{i_{{}_V}}))$ lies in the positive cone of $X$ and $D_{i_{{}_V}}$ is an effective, nonzero divisor subject to the properties 3) - 4) of Lemma \ref{lem:DV}.

2) The sheaf $Q'_{i_{{}_V}}$ is torsion free and fits into the exact sequence
$$
\xymatrix@R=12pt@C=12pt{
	0\ar[r]&\TET \ar[r]&Q'_{i_{{}_V}}\ar[r]& S_{i_{{}_V}} \ar[r]&0,
}
$$  
the bottom row of the diagram \eqref{rkQV=2diag}; the sheaf $S_{i_{{}_V}}$ is a torsion sheaf with
$c_1 (S_{i_{{}_V}})=D_{i_{{}_V}}$ and the subspace $V$ injects into
$H^0 (S_{i_{{}_V}})$; in addition, the sheaf $S_{i_{{}_V}}$ is globally generated by that subspace.

3) For every $[\xi] \in \PP(V)$ the corresponding extension
$$
\xymatrix@R=12pt@C=12pt{
	0\ar[r]& \TET \ar[r]& F_{[\xi]} \ar[r]&  \OO_X \ar[r]& 0
}
$$
 fits into the commutative diagram
$$
\xymatrix@R=12pt@C=12pt{
	&&0\ar[d]&0\ar[d]&\\
	&&\OO_X (-D_{i_{{}_V}}) \ar[d]\ar@{=}[r]&\OO_X (-D_{i_{{}_V}}) \ar[d] &\\
		0\ar[r]& \TET \ar[r]\ar@{=}[d]& F_{[\xi]} \ar[r]\ar[d]&  \OO_X \ar[r]\ar[d]& 0\\
	0\ar[r]& \TET \ar[r]& Q_{[\xi]} \ar[r]\ar[d]&  \OO_{D_{i_{{}_V}}} \ar[r]\ar[d]& 0\\
	&&0&0		
}
$$
Furthermore, for general $[\xi]\in \PP(V)$, the quotient sheaf $Q_{[\xi]}$ 
is torsion free and the bottom row of the above diagram is related to  the exact sequence in 2) via the diagram
$$
\xymatrix@R=12pt@C=12pt{
	&&0\ar[d]&0\ar[d]&\\
0\ar[r]& \TET \ar[r]\ar@{=}[d]& Q_{[\xi]} \ar[r]\ar[d]&  \OO_{D_{i_{{}_V}}} \ar[r]\ar[d]& 0\\
	0\ar[r]&\TET \ar[r]&Q'_{i_{{}_V}}\ar[d]\ar[r]& S_{i_{{}_V}}\ar[d] \ar[r]&0\\
	&&S_{[\xi]}\ar@{=}[r]\ar[d]&S_{[\xi]}\ar[d]&\\
	&&0&0
}
$$
where $S_{[\xi]}$ is a sheaf supported on the zero locus of $\xi$, viewed as a global section of $S_{i_{{}_V}}$ via the inclusion
$V\subset H^0(S_{i_{{}_V}})$; the zero locus is at most $0$-dimensional.

4) If $dim(V)\geq 2$, then the sheaf $F_V$ is Bogomolov unstable and the exact sequence in 1) is a $K_X$-destabilizing sequence for $F_V$.
Furthermore, if the ratio $\alpha_X$ of the Chern numbers of $X$ is subject to the inequality
$$
3\alpha_X -1 \leq \frac{1}{9} \Leftrightarrow \alpha_X \leq \frac{10}{27},
$$
then the sheaf $Q'_{i_{{}_V}}$ is the last $K_X$-semistable factor in the HN filtration of $F_V$ with respect to $K_X$; in other words
$Q'_{i_{{}_V}}$ coincides with the sheaf $Q^m_V$ in Proposition \ref{pro:FVunst}. In addition, for every $[\xi]\in \PP(V)$, the vector bundle
$F_{[\xi]}$ defined by the extension corresponding to the one dimensional subspace $[\xi]\subset V$, see
the extension sequence in 3), is $K_X$-unstable; for a general 
$[\xi] \in \PP(V)$ the exact sequence
\begin{equation}\label{Fxiseq-pro}
\xymatrix@R=12pt@C=12pt{
	0\ar[r]&\OO_X (-D_{i_{{}_V}})\ar[r]&F_{[\xi]} \ar[r]&Q_{[\xi]}\ar[r]&0,
}
\end{equation}
the middle column of the diagram in 3), is the maximal $K_X$-destabilizing sequence of $F_{[\xi]}$, that is, the inclusion
$$
F_{[\xi]}\supset \OO_X (-D_{i_{{}_V}}) \supset 0
$$
is the HN filtration of $F_{[\xi]}$ with respect to $K_X$.
\end{pro} 
\begin{pf}
	Parts 1) and 2) of the proposition are proved in Lemma \ref{lem:DV};
	part 3) is proved in Lemma \ref{lem:V-secSV}.
	
	Part 4) is Lemma \ref{lem:Q'Vdual-maxdist}. Furthermore, imposing the condition
	$$
	3\alpha_X -1 \leq \frac{1}{9},
	$$
	means that the divisor $D_{i_{{}_V}}$ is subject to
	\begin{equation}\label{DiV-third}
	D_{i_{{}_V}} \cdot K_X < \frac{1}{3}K^2_X.
\end{equation}
	This tells us that the subsheaf $\OO_X (-D_{i_{{}_V}})$ in the exact sequence \eqref{Fxiseq-pro} is $K_X$-destabilizing for $F_{[\xi]}$.
	Thus $F_{[\xi]}$ is $K_X$-unstable and it remains to check that  $\OO_X (-D_{i_{{}_V}})$ is its maximal $K_X$-destabilizing subsheaf.
	
	We know that the HN filtration of $F_{[\xi]}$ has the form
	$$
	F_{[\xi]} \supset \OO_X (-D) \supset 0,
	$$
	where  $\OO_X (-D)$ is the maximal $K_X$-destabilizing subsheaf of  $F_{[\xi]}$. Combining with the exact sequence \eqref{Fxiseq-pro}
	gives the diagram
	$$
\xymatrix@R=12pt@C=12pt{
	&&0\ar[d]&&\\
	&&\OO_X (-D)\ar[d]\ar[dr]&&\\
	0\ar[r]&\OO_X (-D_{i_{{}_V}})\ar[r]&F_{[\xi]} \ar[r]&Q_{[\xi]}\ar[r]&0,
}
$$
where the slanted arrow is the composition of the inclusion with the epimorphism of the horizontal sequence. The proof amounts to showing that the slanted arrow is zero: once that arrow is zero, the vertical arrow in the above diagram factors through $\OO_X (-D_{i_{{}_V}})$ and this must be an isomorphism, because the quotient
$F_{[\xi]} /\OO_X (-D) $ is torsion free.	

Let us assume that the slanted arrow in the above diagram is not zero. Then we have a monomorphism
$$
\xymatrix@R=12pt@C=12pt{
	0\ar[r]& \OO_X (-D)\ar[r]& Q_{[\xi]}
}
	$$
	which can be completed to the exact sequence
	$$
\xymatrix@R=12pt@C=12pt{
	0\ar[r]& \OO_X (-D)\ar[r]& Q_{[\xi]}\ar[r]& Q'\ar[r]&0.
}
$$
Factoring out by the torsion part of $Q'$, if nonzero, we obtain an exact sequence
\begin{equation}\label{Qxiseq}
\xymatrix@R=12pt@C=12pt{
	0\ar[r]& {\cal J}_A (-D+T)\ar[r]& Q_{[\xi]}\ar[r]& {\cal J}_{A'} (-L) \ar[r]&0,
}
\end{equation}
where $L=K_X -D_{i_{{}_V}} -(D-T)$ and ${\cal J}_A$ (resp. ${\cal J}_{A'}$) is the ideal sheaf of at most $0$-dimensional subscheme $A$ (resp. $A'$). From this it follows
\begin{equation}\label{ch2Qxi}
2ch_2 (Q_{[\xi]})\leq (D-T)^2 + L^2
\end{equation}
We claim that $L^2 \leq 0$. Indeed, dualizing \eqref{Qxiseq} we obtain
the inclusion
$$
\OO_X(L)\hookrightarrow Q^{\ast}_{[\xi]}.
$$
On the other hand $ Q^{\ast}_{[\xi]}$ is a subsheaf of $\Omega_X$, see the dual of the bottom sequence of the diagram in 3) of the proposition. Thus an inclusion
$$
\OO_X(L)\hookrightarrow \Omega_X.
$$
By definition of $L$ we have
$$
L\cdot K_X =K^2_X -D_{i_{{}_V}}\cdot K_X -(D-T)\cdot K_X \geq K^2_X -D_{i_{{}_V}}\cdot K_X -D\cdot K_X \geq   K^2_X -2D_{i_{{}_V}}\cdot K_X,	
$$
where the last inequality follows from the fact that
$\OO_X (-D)$ is the maximal $K_X$-destabilizing subsheaf of $F_{[\xi]}$.
This and the upper bound for $D_{i_{{}_V}}\cdot K_X$ in \eqref{DiV-third}
give
$$
L\cdot K_X \geq   K^2_X -2D_{i_{{}_V}}\cdot >  K^2_X -\frac{2}{3}K^2_X =\frac{1}{3}K^2_X >0.
$$
Since $L$ can not lie in the positive cone of $X$, we deduce
$$
L^2\leq 0.
$$
This and the estimate \eqref{ch2Qxi} imply
$$
2ch_2(Q_{[\xi]}) \leq (D-T)^2.
$$
From the bottom sequence of the first diagram in 3) of the proposition it follows
$$
 2ch_2(Q_{[\xi]}) =\tau_X -D^2_{i_{{}_V}}.
 $$
 Combining with the previous inequality we obtain
 $$
 \tau_X \leq (D-T)^2 +D^2_{i_{{}_V}}.
 $$
 This and the Hodge Index inequality for $(D-T)^2$ and $D^2_{i_{{}_V}}$
 imply
 $$
 \tau_X \leq (D-T)^2 +D^2_{i_{{}_V}} \leq \frac{((D-T)\cdot K_X)^2}{K^2_X} +\frac{(D_{i_{{}_V}}\cdot K_X)^2}{K^2_X}\leq \frac{2(D_{i_{{}_V}}\cdot K_X)^2}{K^2_X},
 $$
 where the last inequality uses
 $$
 (D-T)\cdot K_X \leq D\cdot K_X \leq D_{i_{{}_V}}\cdot K_X.
 $$
Using the upper bound \eqref{DiV-third} once again we obtain
$$
(1-2\alpha_X)K^2_X <\frac{2}{9}K^2_X
$$
or, equivalently,
$$
\alpha_X >\frac{7}{18}.
$$
But this is in contradiction with the assumption
$\displaystyle{\alpha_X \leq \frac{10}{27}}$.
\end{pf}

\vspace{0.2cm}
\noindent
{\bf 7.2.1. Geometric properties of $D_{i_{{}_V}}$.}
The divisor $D_{i_{{}_V}}$ is a geometric output of our considerations.
We recall that it comes from the exact sequence
\begin{equation}\label{ThetaX-QiV-SiV}
\xymatrix@R=12pt@C=12pt{
	0\ar[r]& \TET\ar[r]&Q'_{i_{{}_V}} \ar[r]& S_{i_{{}_V}}\ar[r]&0
}
\end{equation}
see Proposition \ref{pro:QVrk2}, 2), where $S_{i_{{}_V}}$ is a torsion sheaf with $c_1(S_{i_{{}_V}})=D_{i_{{}_V}}$. The above exact sequence
could be viewed as a homological algebra analog of a ramified covering
map 
$$
f:X \longrightarrow X'
$$
from $X$ to some other surface $X'$. The differential of $f$ gives rise
to 
$$
\xymatrix@R=12pt@C=12pt{
	0\ar[r]& \TET\ar[r]&f^{\ast} (\Theta_{X'})\ar[r]& {\cal R}_f\ar[r]&0,
}
$$ 
the geometric counterpart of the exact sequence \eqref{ThetaX-QiV-SiV}. From this point of view
the divisor $D_{i_{{}_V}}$ becomes an abstract version of the ramification divisor
of a ramified covering $f:X \longrightarrow X'$.

In this subsection we give some geometric properties of $D_{i_{{}_V}}$.
In particular, we identify the curves on $X$ intersecting the divisor 
negatively. These curves turn out to be smooth rational curves and they can be used to produce particular filtrations of $V$.
Most of the results here have appeared in \cite{R}, so if no new assertions are made, we merely
state the relevant properties without proofs.

The technical tool to analyze the geometry of $D_{i_{{}_V}}$ is the diagram
\begin{equation}\label{DiV-Fxi-diag}
\xymatrix@R=12pt@C=12pt{
	&&0\ar[d]&0\ar[d]&\\
	&&\OO_X (-D_{i_{{}_V}}) \ar[d]\ar@{=}[r]&\OO_X (-D_{i_{{}_V}}) \ar[d] &\\
	0\ar[r]& \TET \ar[r]\ar@{=}[d]& F_{[\xi]} \ar[r]\ar[d]&  \OO_X \ar[r]\ar[d]& 0\\
	0\ar[r]& \TET \ar[r]& Q_{[\xi]} \ar[r]\ar[d]&  \OO_{D_{i_{{}_V}}} \ar[r]\ar[d]& 0\\
	&&0&0		
}
\end{equation}
from Proposition \ref{pro:QVrk2}, 3). We choose $[\xi]$ in the Zariski dense open subset where the quotient sheaf $Q_{[\xi]}$ is torsion free.
The main observation is that the monomorphism
in the middle column of the diagram restricted to the divisor $D_{i_{{}_V}}$ must factor through $\TET \otimes \OO_{D_{i_{{}_V}}}$. In other words, the above diagram gives rise to a morphism
\begin{equation}\label{sigmaxi}
\sigma_{[\xi]}:\OO_{D_{i_{{}_V}}} (-D_{i_{{}_V}})\longrightarrow \TET \otimes \OO_{D_{i_{{}_V}}}.
\end{equation}
In addition, the zero locus of this morphism is the zero scheme of the
global section $\widetilde{\sigma_{[\xi]}}$ of
$ F_{[\xi]}(D_{i_{{}_V}})$ corresponding to the monomorphism in the middle column of \eqref{DiV-Fxi-diag} under the identification
$$
Hom_{\AC} (\OO_X (-D_{i_{{}_V}}), F_{[\xi]})=H^0 (F_{[\xi]}(D_{i_{{}_V}})).
$$
The zero scheme of $\widetilde{\sigma_{[\xi]}}$ is precisely the singularity scheme of the quotient sheaf $Q_{[\xi]}$ and according to our choice of $[\xi]$ that scheme is at most $0$-dimensional.
So the restriction $\sigma^D_{[\xi]}$ of $\sigma_{[\xi]}$ to any component $D$ of $D_{i_{{}_V}}$ is also at most $0$-dimensional.
In particular, for any reduced irreducible component $C$ of $D_{i_{{}_V}}$ we obtain the monomorphism
$$
 \sigma^C_{[\xi]}:\OO_C (-D_{i_{{}_V}})\longrightarrow \TET \otimes \OO_{C}.
$$
To extract the geometry we combine this with the normal sequence of 
$C \subset X$:
\begin{equation}\label{sigmaC-normalC}
	\xymatrix@R=25pt@C=12pt{
		&&\OO_C (-D_{i_{{}_V}})\ar^{\sigma^C_{[\xi]}}[d]\ar^(.55){\phi_{{}_C}}[rd]&&\\
		0\ar[r]&\Theta_C \ar[r]& \TET \otimes \OO_C \ar[r]& \OO_C (C)
	}
\end{equation}
where $\Theta_C$ is the tangent sheaf of $C$ and the slanted arrow
$$
\phi_{{}_C}: \OO_C (-D_{i_{{}_V}}) \longrightarrow  \OO_C (C)
$$
is the composition morphism. We distinguish two types of irreducible components of $D_{i_{{}_V}}$ according to whether  the morphism $\phi_{{}_C}$ is zero or not.

\vspace{0.2cm}
{\bf Type I: $\phi_{{}_C} \neq 0$.} In this case the zero locus
$A^C_{[\xi]}$
of  $\sigma^C_{[\xi]}$ is contained in the zero locus of $\phi_{{}_C}$ and we have the estimate
\begin{equation}\label{type1-est}
	deg(A^C_{[\xi]}) \leq deg(\{\phi_{{}_C}=0\}) =C^2 +C\cdot D_{i_{{}_V}}.
\end{equation}

\vspace{0.2cm}
{\bf Type II: $\phi_{{}_C} = 0$.} 
In this case $\sigma^C_{[\xi]}$ factors through $\Theta_C$. To estimate
$deg(A^C_{[\xi]})$ we take
\begin{equation}\label{etaC}
\eta_C: \widetilde{C} \longrightarrow C\subset X
\end{equation}
the normalization of $C$ and obtain:
\begin{equation}\label{Cnorm}
	\begin{gathered}
	\text{the pull back $\eta^{\ast}\Big(\OO_C ( -D_{i_{{}_V}})\Big) \otimes \OO_{\widetilde{C}} \Big(\eta^{\ast}(A^C_{[\xi]})\Big)$ is isomorphic}\\ 
	\text{to the tangent bundle $\Theta_{\widetilde{C}}$ of $\widetilde{C}$:
$\eta^{\ast}\Big(\OO_C ( -D_{i_{{}_V}})\Big) \otimes \OO_{\widetilde{C}} \Big(\eta^{\ast}(A^C_{[\xi]})\Big) \cong \Theta_{\widetilde{C}}$.}	
\end{gathered}
\end{equation}
The details of the proof of this statement can be found in \cite{R}, p.433.
From the isomorphism in \eqref{Cnorm} it follows
\begin{equation}\label{type2-estim}
	\begin{aligned}
\bullet\,\, &	deg(A^C_{[\xi]}) \leq deg(\eta^{\ast}(A^C_{[\xi]}))=C\cdot D_{i_{{}_V}} +(2-2g_{\widetilde{C}}),\\
\bullet \,\,&\text{the morphism $\eta_C$ in \eqref{etaC} is an immersion,}	
	\end{aligned}
\end{equation}
where $g_{\widetilde{C}}$ in the first item is the genus of $\widetilde{C}$, for the second item see \cite{R}, Remark 4.3.

\vspace{0.2cm}
We now decompose $D_{i_{{}_V}}$ into the sum of two parts
\begin{equation}\label{type-decomp}
	D_{i_{{}_V}}=D'_{i_{{}_V}}+D''_{i_{{}_V}},
\end{equation}
where $D'_{i_{{}_V}}$ (resp. $D''_{i_{{}_V}}$) is the part of $D_{i_{{}_V}}$ composed of reduced irreducible components of type I (resp. type II). This decomposition fits with thinking of $D_{i_{{}_V}}$ as an abstract version of the ramification divisor suggested at the beginning:

 the part $D'_{i_{{}_V}}$ should be envisaged as the part of the ramification divisor on which the (hypothetical) ramified covering
is finite, while the part $D''_{i_{{}_V}}$ is composed of the irreducible
components of $D_{i_{{}_V}}$ which are contracted to points.
 
We will see in a moment that the latter part contains all curves 
 intersecting  the divisor $D_{i_{{}_V}}$ negatively.
\begin{lem}\label{lem:Cneg}
	Let $C$ be a reduced irreducible curve on $X$ intersecting
	$D_{i_{{}_V}}$ negatively. Then $C$ is a smooth rational curve, $C\cdot D_{i_{{}_V}} =-1$ or $-2$ and $C$ is a component of $D''_{i_{{}_V}}$.	
\end{lem}
\begin{pf}
	From the assumption $C\cdot D_{i_{{}_V}} <0$ it follows that $C$
	is an irreducible component of $D_{i_{{}_V}}$ and its self-intersection $C^2 <0$. From \eqref{type1-est} it follows that
	$C$ is an irreducible component of the second type subject to
	$$
	deg(\eta^{\ast}(A^C_{[\xi]}))=C\cdot D_{i_{{}_V}} +(2-2g_{\widetilde{C}}),
	$$
	see \eqref{type2-estim}. Since the first term is negative we deduce
	that $g_{\widetilde{C}} =0$ and $C\cdot D_{i_{{}_V}} =-1$ or $-2$.
	Furthermore, the above formula implies
	$$
	deg(\eta^{\ast}(A^C_{[\xi]})) \leq 1.
	$$
	Hence $A^C_{[\xi]}$ is contained in the smooth locus of $C$ and the diagram \eqref{sigmaC-normalC} gives an isomorphism
	$$
	\OO_C ( A^C_{[\xi]}-D_{i_{{}_V}})  \cong \Theta_{C}.
	$$
	From this it follows that $\Theta_{C}$ is locally free $\OO_C$-module. By Lipman's criterion of smoothness, the dual of the Jacobian criterion, see \cite{Li},
	we deduce that $C$ is smooth.
\end{pf}

The isomorphism in \eqref{Cnorm} also establishes the relation between
the singular locus of the reduced part $(D''_{i_{{}_V}})_{red}$ of $D''_{i_{{}_V}}$ and the zeros
of $\widetilde{\sigma_{[\xi]}}$ on $(D''_{i_{{}_V}})_{red}$.
\begin{lem}\label{lem:singD''}
	Let $(D''_{i_{{}_V}})_{red}$ be the reduced part of $(D''_{i_{{}_V}})$. Then the singular locus of $(D''_{i_{{}_V}})_{red}$ is contained in the scheme 
	$\{\widetilde{\sigma_{[\xi]}}=0\} \bigcap D''_{i_{{}_V}}$.
\end{lem}
\begin{pf}
	 The restriction of $\sigma_{[\xi]}$ in \eqref{sigmaxi} to
	$(D''_{i_{{}_V}})_{red}$ factors trough the tangent sheaf
	$\Theta_{(D''_{i_{{}_V}})_{red}}$ of $(D''_{i_{{}_V}})_{red}$:
	$$
	\sigma^{(D''_{i_{{}_V}})_{red}}_{[\xi]}:\OO_{(D''_{i_{{}_V}})_{red}} (-D_{i_{{}_V}}) \longrightarrow \Theta_{(D''_{i_{{}_V}})_{red}}
	$$
	and this is an isomorphism away from the zero locus of
	$\widetilde{\sigma_{[\xi]}}$. This means that
	$\Theta_{(D''_{i_{{}_V}})_{red},p}$ is a locally free $\OO_{(D''_{i_{{}_V}})_{red},p}$-module for any point
	$p \in (D''_{i_{{}_V}})_{red} \setminus \{\widetilde{\sigma_{[\xi]}}=0\}
	$. By the Lipman's criterion of smoothness, \cite{Li}, all points
	in $(D''_{i_{{}_V}})_{red} \setminus \{\widetilde{\sigma_{[\xi]}}=0\}$ are smooth points of
	$(D''_{i_{{}_V}})_{red}$.	
\end{pf}

From Lemma \ref{lem:Cneg} the obstruction for $D_{i_{{}_V}}$ to be nef
is the set of smooth rational curves
\begin{equation}\label{Rset}
\mathfrak{R}:=\{\text{the irreducible components $C$ of $D_{i_{{}_V}}$ subject to $C\cdot D_{i_{{}_V}}=-2$ or $-1$} \}.	
\end{equation}
We divide it into two parts according to the value of the intersection with $ D_{i_{{}_V}}$:
\begin{equation}\label{Rsetdeg}
\mathfrak{R}^d:=\{C\in \mathfrak{R} |   C\cdot D_{i_{{}_V}}=d \},	
\end{equation}
where $d=-2$ or $(-1)$. The properties of these sets are as  follows.

\begin{pro}\label{pro:CinR}
	1) Every curve $C\in \mathfrak{R}^{-2}$ is disjoint from any other irreducible component of $D''_{i_{{}_V}}$ and intersects some components in $D'_{i_{{}_V}}$; more precisely, the intersection
	$$
	C \cdot D'_{i_{{}_V}}=-m_C C^2 -2, 
	$$
	where $m_C$ is the multiplicity of $C$ in $D_{i_{{}_V}}$.
	
	2) Every curve $C\in \mathfrak{R}^{-1}$ passes exactly through one point, denoted $p^{[\xi]}_C$, in the scheme $\{\widetilde{\sigma_{[\xi]}}=0\}$. 
		Any other irreducible component $C'$ of $D''_{i_{{}_V}}$ intersecting $C$ has the set-theoretic intersection
	$$
	C'\cap C=\{p^{[\xi]}_C\}.
	$$
	In particular, if such components $C'$ exist, then, set-theoretically, we have  
	$$
p^{[\xi]}_C = C\bigcap  \big(D''_{i_{{}_V}}\big)^C,
$$
 where $\big(D''_{i_{{}_V}}\big)^C$ is the part of $D''_{i_{{}_V}}$ formed by the irreducible components different from $C$.
 In this situation the point $p^{[\xi]}_C$ is independent of $[\xi]$ and it will be denoted $p_C$.
\end{pro}
\begin{pf}
	For $C\in \mathfrak{R}^{-2}$ we have
	$$
	C\cdot D_{i_{{}_V}}=-2.
	$$
	This and the formula in \eqref{type2-estim} imply that
	$\widetilde{\sigma_{[\xi]}}$ has no zeros on $C$. From Lemma \ref{lem:singD''} it follows that $C$ is contained in the smooth part of $(D''_{i_{{}_V}})_{red}$. Hence no other irreducible component of $D''_{i_{{}_V}}$ intersects $C$.
	
	From this we have
	$$
		-2=C\cdot D_{i_{{}_V}}=C\cdot D'_{i_{{}_V}} +C\cdot D''_{i_{{}_V}}=C\cdot D'_{i_{{}_V}} +m_C C^2,
	$$
where $m_C$ is the multiplicity of $C$ in $D_{i_{{}_V}}$. Hence the formula
$$
C\cdot D'_{i_{{}_V}} =-m_C C^2 -2.
$$
Observe that the expression on the right is positive: we are assuming that $K_X$ is ample and this implies by adjunction formula that all smooth rational curves on $X$ have self-intersection at most $(-3)$.
	
	For $C\in \mathfrak{R}^{-1}$ we have
	$$
		C\cdot D_{i_{{}_V}}=-1.
	$$
	This and the formula in \eqref{type2-estim} imply that
	$\widetilde{\sigma_{[\xi]}}$ vanishes precisely at one point on $C$. Call this point $p^{[\xi]}_C$. From Lemma \ref{lem:singD''} it follows that
	any other irreducible component of $D''_{i_{{}_V}}$
	intersects $C$ along a subscheme supported on $p^{[\xi]}_C$. 
\end{pf}

\begin{rem}\label{rem:Rset-marked}
All curves in $\mathfrak{R}$ come with marked points on them:

1) if $C\in {\mathfrak{R}^{-2}}$, then the marked points on $C$ are the points of the intersection
$$
P_C:=C\cdot D'_{i_{{}_V}};
$$

2) if $C\in  {\mathfrak{R}^{-1}}$, then we have
$$
P_C:=C\cdot D'_{i_{{}_V}} + p^{[\xi]}_C,
$$
where $p^{[\xi]}_C$ is as in Proposition \ref{pro:CinR}.
\end{rem}

\vspace{0.2cm}
\noindent
{\bf 7.2.2. $D_{i_{{}_V}}$ and filtrations of $V$.}
The understanding of the subspace $V$ of $H^1(\TET)$ annihilated by the divisor $D_{i_{{}_V}}$ essentially reduces to the case when this divisor  is nef and/or enumeration of certain smooth rational curves
on $X$. Namely,
using the curves in $\mathfrak{R}$ introduced in \eqref{Rset}, we can define geometrically meaningful filtrations of $V$. This relies on the following observation.
\begin{lem}\label{lem:VC}
	Assume the set $\mathfrak{R}$ in \eqref{Rset} is nonempty and let $C$ be a curve in $\mathfrak{R}$. Then there is a linear map
	$$
\lambda_C:	V\longrightarrow H^0(\TET \otimes \OO_C (D_{i_{{}_V}}))
	$$
	whose kernel $V_C$ is a subspace of $V$ of codimension at most $2$. Furthermore, the codimension $2$ may only occur if
	$C\in \mathfrak{R}^{-1}$ and no other irreducible component
	of $D''_{i_{{}_V}}$ intersects $C$. 
\end{lem}
\begin{pf}
	We know that the subspace $V$ is annihilated by $D_{i_{{}_V}}$, see
	Lemma \ref{lem:DV}. This gives an injective linear map
	$$
	V \longrightarrow H^0(\TET \otimes \OO_{D_{i_{{}_V}}} (D_{i_{{}_V}})).
	$$
	Furthermore, $D_{i_{{}_V}}$ is minimal with respect to the property
	of annihilating $V$. This means that for any proper component
	$D$ of $D_{i_{{}_V}}$ the composition map,
	$$
	\xymatrix@R=12pt@C=12pt{
	V \ar[r] \ar[dr]& H^0(\TET \otimes \OO_{D_{i_{{}_V}}}(D_{i_{{}_V}}))\ar[d]\\
	& H^0(\TET \otimes \OO_D (D_{i_{{}_V}}))
}
$$
the slanted arrow of the diagram, is nonzero. Applying this to $D=C$, a curve in $\mathfrak{R}$, gives a nonzero homomorphism
$$
\lambda_C:V\longrightarrow  H^0(\TET \otimes \OO_C (D_{i_{{}_V}})).
$$
For $C$ in $\mathfrak{R}^{d}$, the set defined in \eqref{Rsetdeg}, the space on the right is
\begin{equation}\label{Theta-on-C}
 H^0(\TET \otimes \OO_C (D_{i_{{}_V}})) \cong H^0 (\OO_{\PP^1} (2+d))\oplus  H^0 (\OO_{\PP^1} (C^2+d)) =H^0 (\OO_{\PP^1} (2+d)).
\end{equation}
 
 For $d=-2$, the last space is $1$-dimensional, so the kernel
 $$
 V_C:=ker(\lambda_C)
 $$
 is a codimension one subspace of $V$.
 
 For $d=-1$, according to \eqref{Theta-on-C} the space $H^0(\TET \otimes \OO_C (D_{i_{{}_V}})) \cong H^0 (\OO_{\PP^1} (1))$ is two dimensional. In addition, from Proposition
 \ref{pro:CinR}, 2), for a general  $\xi \in V$, the global section
  $\lambda_C (\xi)$ of $\TET \otimes \OO_C (D_{i_{{}_V}})$ must vanish at a single point $p^{[\xi]}_C$ on $C$. Under the identification
  $$
  H^0(\TET \otimes \OO_C (D_{i_{{}_V}})) \cong H^0 (\OO_{\PP^1} (1))
  $$
 that point is the zero divisor of the global section of $\OO_{\PP^1} (1)$ under the above isomorphism.
 From Proposition
 \ref{pro:CinR}, 2), we also know that the point $p^{[\xi]}_C$ is 
 independent of $\xi$, when there are other irreducible components of $D''_{i_{{}_V}}$ intersecting $C$. 
  From this it follows that the image of $\lambda_C$ is one dimensional, with the possible exception, when
   no other irreducible component of $D''_{i_{{}_V}}$ intersects $C$.   
\end{pf}

In view of the previous results it will be useful to introduce the following.
\begin{defi}\label{def:V-D-pair}
	Let $W$ be a nonzero subspace of $H^1 (\TET)$ and
	let $D$ be a nonzero effective divisor annihilating $W$ and it is minimal with respect to this property, this is to say that no proper component of $D$ annihilates $W$.
	The pair $(W,D)$ is called {\rm numerically effective (nef pair)} if the divisor $D$ is nef. Otherwise we say that $(W,D)$ has negative part; this, according to Lemma \ref{lem:Cneg}, is the collection
	of smooth rational components $C$ of $D$ having the intersection
	$C\cdot D =-2$ or $-1$. This collection from now on will be denoted ${\mathfrak{R}}_{W,D}$ and decomposed into two disjoint parts
	$$
	{\mathfrak{R}}_{W,D}={\mathfrak{R}}^{-2}_{W,D} \coprod {\mathfrak{R}}^{-1}_{W,D},
	$$
	where ${\mathfrak{R}}^{d}_{W,D}$ consists of curves $C$ in ${\mathfrak{R}}_{W,D}$ with $C\cdot D=d$, for $d=-2, -1$. 	
\end{defi}

We can now describe a stratification of the pair 
$(V, D_{i_{{}_V}})$ which in essence `peels off' the negative part of $D_{i_{{}_V}}$:

\vspace{0.2cm}
1) if $D_{i_{{}_V}}$ is nef, then $(V,D_{i_{{}_V}})$ is already a nef pair, so the process terminates;

\vspace{0.2cm}
2) if the negative part $\mathfrak{R}_{V,D_{i_{{}_V}}}$ is nonempty, choose $C_0$ in $\mathfrak{R}_{V,D_{i_{{}_V}}}$ so that the codimension of
$V_{C_0}$ in Lemma \ref{lem:VC} is the smallest possible; if $V_{C_0} =0$, the process terminates, if not, we set $V^{(1)}:=V_{C_0}$ and continue: the divisor $(D_{i_{{}_V}} -C_0)$
annihilates the subspace $V^{(1)}$, we choose a component which is {\it minimal}
with respect to this property and denote it by $D^{(1)}$; thus we obtain 
a {\it subpair}
\begin{equation}\label{subpair}
	(V^{(1)},D^{(1)}) \subset (V,D_{i_{{}_V}})
\end{equation}
 meaning that $V^{(1)} \subset V$ and $D^{(1)}\subset D_{i_{{}_V}}$;
 furthermore, by construction we have
\begin{equation}\label{subpair1} 
	\begin{gathered}
 dim(V/V^{(1)})=dim(V/V_{C_0})=\min \{codim(V_C) | C\in \mathfrak{R}_{V,D_{i_{{}_V}}}\}=\text{$1$ or $2$},
 \\
 \text{the divisor $(D_{i_{{}_V}} -D^{(1)})$ contains $C_0$};
\end{gathered}
\end{equation}

\vspace{0.2cm}
3) with the pair $(V^{(1)}, D^{(1)})$, return to 1) and repeat the process.

\vspace{0.2cm}
 Using the above process we obtain the following.
 
 \begin{pro}\label{pro:VDfilt}
 	Let the pair $(V,D_{i_{{}_V}})$ be as in Proposition \ref{pro:QVrk2}.
 		Then either the pair is nef in the sense of Definition \ref{def:V-D-pair} or the set $\mathfrak{R}_{V,D_{i_{{}_V}}}$ is nonempty and there is a filtration
 		$$
 		(V,D_{i_{{}_V}})=(V^{(0)},D^{(0)}) \supset (V^{(1)},D^{(1)}) \supset \cdots \supset (V^{(l)},D^{(l)}),
 		$$ 
 		subject to the following properties:
 		
 		\vspace{0.2cm}
 		(i) for every $m\in [0,l-1]$, the quotient space $V^{(m)}/V^{(m+1)}$ has dimension
 		$$
 		dim(V^{(m)}/V^{(m+1)})=\min \{codim(V^{(m)}_C) | C \in {\mathfrak{R}_{V^{(m)},D^{(m)}}}\} =\text{$1$ or $2$},
 		$$
 		and the divisor $(D^{(m)} -D^{(m+1)})$ contains at least
 		one curve $C \in \mathfrak{R}_{V^{(m)},D^{(m)}}$ with
 		$V^{(m)}_C=V^{(m+1)}$,
 		
 			\vspace{0.2cm}
 		(ii) if  
 		$
 		dim(V^{(m)}/V^{(m+1)})=2 
 		$, 
 		for some $m \in [0,l-1]$, then
 		 
 		 $(1ii)_m$
 		$ \mathfrak{R}_{V^{(m)},D^{(m)}}=\mathfrak{R}^{-1}_{V^{(m)},D^{(m)}},
 		 $
 		   
 		$(2ii)_m$ every curve $C\in\mathfrak{R}_{V^{(m)},D^{(m)}}$ is disjoint from any other irreducible component of $\big(D^{(m)}\big)''$,
 		
 		$(3ii)_m$ the divisors
 		$$
 		\begin{aligned}
 		\big(D^{(m)}\big)^-:=&\sum_{C\in \mathfrak{R}_{V^{(m)},D^{(m)}}}
 		\frac{C}{(-C^2)},\\
 		\big(D^{(m)}\big)^+ :=&D^{(m)}-\big(D^{(m)}\big)^-=D^{(m)}+\sum_{C\in \mathfrak{R}_{V^{(m)},D^{(m)}}}
 		\frac{C}{C^2},
 	\end{aligned}
 			$$
 			are respectively the negative and positive parts of
 			 the Zariski decomposition of $D^{(m)}$,
 		
 		\vspace{0.2cm}	
 		(iii) the last pair 	$(V^{(l)},D^{(l)})$ of the filtration is either nef or for every $C\in \mathfrak{R}_{V^{(l)},D^{(l)}}$
 		$V^{(l)}_C=0$; in the latter case $dim(V^{(l)})\leq 2$ and if the equality holds, then the properties $(1ii)_m-(3ii)_m$ hold for $m=l$.
 \end{pro}
\begin{pf}
	Starting from $(V,D_{i_{{}_V}}):=(V^{(0)},D^{(0)})$ the filtration is constructed inductively using the procedure described in 1) - 3) prior to the statement of the proposition.
	
	The item (i) follows from \eqref{subpair1}. We turn now to the item
	(ii) of the proposition. It is assumed that
	$$
		dim(V^{(m)}/V^{(m+1)})=\min \{codim(V^{(m)}_C) | C \in {\mathfrak{R}_{V^{(m)},D^{(m)}}}\} =2.
		$$
		This and Lemma \ref{lem:VC} imply that $codim(V^{(m)}_C)=2$, for
		every $C \in {\mathfrak{R}_{V^{(m)},D^{(m)}}}$ and this in turn
		tell us
		
		$(1ii)_m$: $\mathfrak{R}_{V^{(m)},D^{(m)}}=\mathfrak{R}^{-1}_{V^{(m)},D^{(m)}}$,
		
		$(2ii)_m$: every curve $C$ in $\mathfrak{R}_{V^{(m)},D^{(m)}}$
		is disjoint from any other irreducible component of
		$\big(D^{(m)}\big)''$.
		
		Next we investigate the Zariski decomposition of $D^{(m)}$.
			For this we set
			$$
			\big(D^{(m)}\big)^-:=\sum_{C\in \mathfrak{R}_{V^{(m)},D^{(m)}}}
			\frac{C}{(-C^2)}.
			$$
			From $(2ii)_m$ it follows that the support of this divisor consists of disjoint curves and we have
			$$
			\big(D^{(m)}\big)^- \cdot C=C\cdot \Big(\sum_{C'\in \mathfrak{R}_{V^{(m)},D^{(m)}}}
			\frac{C'}{(-{C'}^2)}\Big) =	\frac{C\cdot C}{(-{C}^2)}=-1,
			$$
			for every $C\in \mathfrak{R}_{V^{(m)},D^{(m)}}$. This implies
			$$
			(D^{(m)}-\big(D^{(m)}\big)^-)\cdot C=0, \,\forall C\in \mathfrak{R}_{V^{(m)},D^{(m)}}.
			$$
			Thus setting						
		$$
			\big(D^{(m)}\big)^+ :=D^{(m)}-\big(D^{(m)}\big)^-=D^{(m)}+\sum_{C\in \mathfrak{R}_{V^{(m)},D^{(m)}}}
			\frac{C}{C^2}
$$
gives the orthogonal decomposition
$$
D^{(m)}=\big(D^{(m)}\big)^+ + \big(D^{(m)}\big)^-.
$$
To see that this is the Zariski decomposition of $D^{(m)}$ it remains to check that 	$\big(D^{(m)}\big)^+$ is nef.
For this we need to verify that
$$
\big(D^{(m)}\big)^+ \cdot \Gamma \geq 0,
$$
for every irreducible component $\Gamma$ of $D^{(m)}$. From $(2ii)_m$ we know that
this is the case for all irreducible components of $\big(D^{(m)}\big)''$. So we take $\Gamma \in \big(D^{(m)}\big)'$ and write
\begin{equation}\label{DmGamma-est}
D^{(m)} \cdot \Gamma \geq m_{\Gamma}\Gamma^2 + \sum_{C\in \mathfrak{R}_{V^{(m)},D^{(m)}}} m_C C\cdot \Gamma,
\end{equation}
where $m_{\Gamma}$ (resp. $m_C$) is the multiplicity of $\Gamma$ (resp. $C$) in $D^{(m)}$. Since $\Gamma$ is an irreducible component of type I
in $D^{(m)}$, we have the estimate
$$
\Gamma^2 + D^{(m)} \cdot \Gamma \geq 0,
$$
see \eqref{type1-est}. We rewrite it as $\Gamma^2 \geq -D^{(m)} \cdot \Gamma$ and substitute into \eqref{DmGamma-est} to deduce
$$
(m_{\Gamma} +1)D^{(m)} \cdot \Gamma \geq \sum_{C\in \mathfrak{R}_{V^{(m)},D^{(m)}}} m_C C\cdot \Gamma
$$
From this it follows
$$
\begin{gathered}
\big(D^{(m)}\big)^+ \cdot \Gamma=(D^{(m)} -\big(D^{(m)}\big)^-) \cdot \Gamma=
D^{(m)} \cdot \Gamma +\sum_{C\in \mathfrak{R}_{V^{(m)},D^{(m)}}} \frac{ C\cdot \Gamma}{C^2} \geq
\\ \frac{1}{m_{\Gamma}+1}\sum_{C\in \mathfrak{R}_{V^{(m)},D^{(m)}}} m_C C\cdot \Gamma
+\sum_{C\in \mathfrak{R}_{V^{(m)},D^{(m)}}} \frac{ C\cdot \Gamma}{C^2}
=\sum_{C\in \mathfrak{R}_{V^{(m)},D^{(m)}}} \Big(\frac{m_C}{m_{\Gamma}+1}+ \frac{1}{C^2} \Big) C\cdot \Gamma.
\end{gathered}
$$
Thus we obtain
\begin{equation}\label{Dm+Gamma-est}
	\big(D^{(m)}\big)^+ \cdot \Gamma \geq \sum_{C\in \mathfrak{R}_{V^{(m)},D^{(m)}}} \Big(\frac{m_C}{m_{\Gamma}+1}+ \frac{1}{C^2} \Big) C\cdot \Gamma. 
\end{equation}
 Recall the decomposition
$$
D^{(m)} =\big(D^{(m)}\big)' + \big(D^{(m)}\big)''
$$
according to the types of the irreducible components of $D^{(m)}$, see \eqref{type-decomp}. Intersecting with $C \in \mathfrak{R}_{V^{(m)},D^{(m)}}$ and remembering that $C\cdot D^{(m)} =-1$, we obtain
$$
-1=C\cdot D^{(m)}=C\cdot \Big( \big(D^{(m)}\big)' + \big(D^{(m)}\big)'' \Big)=C\cdot  \big(D^{(m)}\big)' + m_C C^2 \geq m_{\Gamma} C\cdot \Gamma + m_C C^2.
$$
The self-intersection $C^2$ is negative, so above inequality can be rewritten as follows
$$
m_C \geq  \frac{ m_{\Gamma} C\cdot \Gamma +1}{(-C^2)}.
$$
Substituting for $m_C$ in \eqref{Dm+Gamma-est} gives
$$
\begin{gathered}
 	\big(D^{(m)}\big)^+ \cdot \Gamma \geq \sum_{C\in \mathfrak{R}_{V^{(m)},D^{(m)}}} \Big(\frac{m_C}{m_{\Gamma}+1}+ \frac{1}{C^2} \Big) C\cdot \Gamma \geq 
 	 \sum_{C\in \mathfrak{R}_{V^{(m)},D^{(m)}}} \Big(\frac{m_{\Gamma}C\cdot \Gamma +1}{(m_{\Gamma}+1)(-C^2)}+ \frac{1}{C^2} \Big) C\cdot \Gamma
 	 \\
 	 = \sum_{C\in \mathfrak{R}_{V^{(m)},D^{(m)}}} \frac{m_{\Gamma}(C\cdot \Gamma -1)}{m_{\Gamma}+1}\frac{ C\cdot \Gamma}{(-C^2)} \geq 0,
 	\end{gathered}
 	$$
 	for every irreducible component $\Gamma$ in $\big(D^{(m)}\big)'$. This completes the proof of $(3ii)_m$.
 	
 	The part (iii) follows directly from the conditions in 1) and 2) for the process of constructing a smaller pair to terminate.	
\end{pf}

It should be clear that for $(V,D_{i_{{}_V}})$ to be a nef pair in the sense of Definition \ref{def:V-D-pair} is a `good' property to have.
Not only it tells us that the divisor $D_{i_{{}_V}}$ is nef, but it also
gives a good estimate of the dimension of $V$.

\begin{pro}\label{pro:nefpair}
	Assume $(V,D_{i_{{}_V}})$ is as in Proposition \ref{pro:QVrk2} and assume $(V,D_{i_{{}_V}})$ is a nef pair. Then either $D_{i_{{}_V}}$ is nef and big and the dimension $V$ is subject to
	$$
	dim(V)\leq \frac{1}{3}(3c_2(X)-K^2_X),
	$$
	or $D_{i_{{}_V}}$ is nef and $D^2_{i_{{}_V}}=0$. In this case the 
	following holds.
	
	1) the parts $D'_{i_{{}_V}}$ and $D''_{i_{{}_V}}$
	of the decomposition
	$$
	D_{i_{{}_V}}=D'_{i_{{}_V}}+D''_{i_{{}_V}}
	$$
	are subject to
	$$
	D'_{i_{{}_V}} \cdot D''_{i_{{}_V}} =\big(D'_{i_{{}_V}}\big)^2 =\big(D''_{i_{{}_V}}\big)^2 =0.
	$$
	In particular, the parts $D'_{i_{{}_V}} $ and $D''_{i_{{}_V}}$ are disjoint and each part is a nef divisor.
	
	2) If  $D'_{i_{{}_V}}$ is nonzero, then it consists of mutually disjoint irreducible components; each irreducible component $C'$  of $D'_{i_{{}_V}}$ is a smooth curve with the split normal sequence
	$$
		\TET \otimes \OO_{C'} \cong \OO_{C'} (C')\oplus \Theta_{C'},
		$$
where $ \Theta_{C'}$ is the tangent bundle of $C'$, and $\OO_{C'} ((m_{C'} +1) C') =\OO_{C'}$, where $m_{C'}$ is the multiplicity of $C'$ in  $D'_{i_{{}_V}}$. The latter property means that $\OO_{C'} (C')$ is a torsion point of the Jacobian variety $Pic^{\circ} (C')$ of $C'$; letting $t_{C'}$ to be the order of torsion of $\OO_{C'} (C')$ in $Pic^{\circ} (C')$, gives 
$$
(m_{C'}+1)\equiv 0 \text{(mod $t_{C'}$)}.
$$

3) If $D''_{i_{{}_V}}$ is nonzero, then each of its irreducible components $C''$ is a rational curve. It is either smooth or has
a single double point.

4) If $	D_{i_{{}_V}}$ has at least two connected components, then 
$X$ is a fibered surface with the divisor $	D_{i_{{}_V}}$ contained 
in the fibres, that is, there is a surjective morphism
$$
f:X\longrightarrow B
$$
onto a smooth projective curve $B$ with connected fibres and an effective nonzero divisor $A$ on $B$ such that
$$
T:=f^{\ast}(A)- 	D_{i_{{}_V}}
$$
is an effective divisor. Furthermore, denote by $F$ the class of a fibre of $f$

$\bullet$ then $F=t_{C'}C'$,
for any irreducible component $C'$ of $	D'_{i_{{}_V}}$ and where $t_{C'}$ is as in 2),

$\bullet$ let $\displaystyle{D''_{i_{{}_V}}=\sum^{N}_{k=1}D''_{i_{{}_V},k}}$ be the decomposition of $D''_{i_{{}_V}}$ into its distinct connected components, then
$$
F=r_k D''_{i_{{}_V},k},
$$
for some positive rational number $r_k$. 

5) $dim(V) \leq \frac{1}{6} K_X \cdot 	D_{i_{{}_V}} +\frac{1}{3}(3c_2 (X)-K^2_X) \leq \HA \big(3c_2 (X)-K^2_X -deg(S_{[\xi]})\big)$, where
$S_{[\xi]}$ as in Proposition \ref{pro:QVrk2}, 3). 
\end{pro}
\begin{pf}
	To estimate the dimension of $V$ we recall the inclusion
	$$
	V\subset H^0 (S_{i_{{}_V}})
		$$
in Proposition \ref{pro:QVrk2}, 2). Hence the inequality
\begin{equation}\label{dimVesth0SiV}
dim(V) \leq h^0 (S_{i_{{}_V}}).
\end{equation}
We estimate the right hand side by choosing a general $\xi \in V$ and using the exact sequence
$$
\xymatrix{
	0\ar[r]& \OO_{D_{i_{{}_V}}} \ar[r]& S_{i_{{}_V}} \ar[r]& S_{[\xi]}\ar[r]&0,
}
$$
see the column on the right in the second diagram of
Proposition \ref{pro:QVrk2}, 3):
$$
h^0 (S_{i_{{}_V}}) \leq h^0 (\OO_{D_{i_{{}_V}}}) +h^0(S_{[\xi]})=h^0 (\OO_{D_{i_{{}_V}}}) +deg(S_{[\xi]}).
$$	
From Lemma \ref{lem:Q'Vstable}, 3), we deduce
$$
deg(S_{[\xi]}) \leq \frac{1}{3}(3c_2 (X)-K^2_X) -\frac{1}{3}(K_X \cdot D_{i_{{}_V}}-2D^2_{i_{{}_V}}).
$$
Substituting into the previous inequality gives
\begin{equation}\label{h0SiV-estim}
h^0 (S_{i_{{}_V}}) \leq h^0 (\OO_{D_{i_{{}_V}}}) +\frac{1}{3}(3c_2 (X)-K^2_X) -\frac{1}{3}(K_X \cdot D_{i_{{}_V}}-2D^2_{i_{{}_V}}).
\end{equation}
If the divisor $D_{i_{{}_V}}$ is nef and big
$$
h^0 (\OO_{D_{i_{{}_V}}}) =1
$$
and we obtain
$$
h^0 (S_{i_{{}_V}}) \leq \frac{1}{3}(3c_2 (X)-K^2_X) +1-\frac{1}{3}(K_X \cdot D_{i_{{}_V}}-2D^2_{i_{{}_V}}).
$$ 
Combining with \eqref{dimVesth0SiV} gives the inequality
$$
dim(V) \leq \frac{1}{3}(3c_2 (X)-K^2_X)+1-\frac{1}{3}(K_X \cdot D_{i_{{}_V}}-2D^2_{i_{{}_V}}).
$$
To obtain the inequality asserted at the beginning of the proposition it is enough to show
\begin{equation}\label{bigger3}
	K_X \cdot D_{i_{{}_V}}-2D^2_{i_{{}_V}} \geq 3.
\end{equation}
Assume we have
$$
	K_X \cdot D_{i_{{}_V}}-2D^2_{i_{{}_V}} \leq 2.
	$$
	This can be rewritten as
	\begin{equation}\label{D2HA}
	D^2_{i_{{}_V}} \geq \HA 	K_X \cdot D_{i_{{}_V}} -1.
\end{equation}
	We claim that this implies
	\begin{equation}\label{KDiValmostHA}
		K_X \cdot D_{i_{{}_V}} \geq \HA K^2_X -2.
	\end{equation}
Let us assume this for a moment and complete the argument: compare the above lower bound with the upper bound
$$
	K_X \cdot D_{i_{{}_V}} < \frac{\sqrt{1+16(3\alpha_X -1)} -1}{2} K^2_X
	$$
	in Lemma \ref{lem:DV}, 4), to obtain the inequality
	$$
\frac{\sqrt{1+16(3\alpha_X -1)} -1}{2}  > \HA  -\frac{2}{K^2_X}.
$$
This is rewritten as follows
$$
\sqrt{1+16(3\alpha_X -1)} >2-	\frac{4}{K^2_X}
	$$
	Squaring and simplifying gives
	$$
	16(3\alpha_X -1)> 3 -\frac{16}{K^2_X} +\frac{16}{(K^2_X)^2}.
	$$
	Remember that we are working under the assumption
	$$
	\alpha_X < \frac{3}{8}.
	$$
	Substituting into the previous inequality gives
	$$
	0> 1 -\frac{16}{K^2_X} +\frac{16}{(K^2_X)^2}.
	$$ 
	Thus $K^2_X$ is subject to 
	$$
	(K^2_X)^2 -16 K^2_X +16 <0
	$$
	and this occurs only if $K^2_X \leq 14$. But the positivity of the index
	$\tau_X$ imposes $K^2_X > 8\chi(\OO_X)$. Hence $\chi(\OO_X)=1$ and
	$c_2 (X)=12-K^2_X$. The only admissible pair subject to the above constraints is
	$$
	(c_2 (X),K^2_X)=(3,9)
	$$
	lying on the BMY line which is excluded from our considerations because
	those surfaces are infinitesimally rigid, while we are assuming 
	$H^1 (\TET)\neq 0$.
	 
	 We now return to the proof of \eqref{KDiValmostHA}. Assume the contrary
	 $$
	 K_X \cdot D_{i_{{}_V}} < \HA K^2_X -2
	 $$
	or, equivalently,
	\begin{equation}\label{K2lb}
	K^2_X \geq 2 K_X \cdot D_{i_{{}_V}} +5.
\end{equation}
	Combining this with the inequality \eqref{D2HA} and the Hodge Index inequality we obtain
	$$
	 (K_X \cdot D_{i_{{}_V}})^2 \geq  K^2_X  D^2_{i_{{}_V}} \geq 
	 (2 K_X \cdot D_{i_{{}_V}} +5)\big(\HA  K_X \cdot D_{i_{{}_V}} -1\big).
	 $$
	 This simplifies to
	 $$
	 K_X \cdot D_{i_{{}_V}} \leq 10.
	 $$
	 The equality holds if we have equality everywhere. In particular,
	 $$
	 K^2_X = 2 K_X \cdot D_{i_{{}_V}} +5=25.
	 $$
	 Taken together with $K^2_X > 8\chi(\OO_X)$ gives
	 $$
	 \chi(\OO_X) \leq 3.
	 $$
	 From Noether's formula
	 $$
	 c_2 (X)=12\chi(\OO_X) -K^2_X =12\chi(\OO_X) -25
	 $$
	 and the positivity of $c_2 (X)$ we deduce
	 $$
	 (c_2 (X),K^2_X)=(11,25).
	 $$
	 However the ratio
	 $$
	 \alpha_X =\frac{c_2(X)}{K^2_X}=\frac{11}{25}>\frac{3}{8}.
	 $$
	 This implies that
	 $$
	  K_X \cdot D_{i_{{}_V}} \leq 9.
	  $$
	  On the other hand from the positivity of $\tau_X$ we know that $K^2_X \geq 17$. This and the Hodge Index inequality implies
	  $$
	  D^2_{i_{{}_V}} \leq \frac{(K_X \cdot D_{i_{{}_V}})^2}{K^2_X} \leq \frac{81}{17} <5.
	  $$
	  This and the lower bound \eqref{D2HA} imply that for the value
	  $K_X \cdot D_{i_{{}_V}} = 9$ we must have $D^2_{i_{{}_V}}=4$  and hence $K^2_X \leq 20$. But the lower bound in \eqref{K2lb} gives
	  $K^2_X \geq 23$. This rules out the value $K_X \cdot D_{i_{{}_V}} = 9$
	  and we are left with
	   $$
	  K_X \cdot D_{i_{{}_V}} \leq 8.
	  $$
	  Repeating the above argument we deduce that for $K_X \cdot D_{i_{{}_V}} = 8$ the value $D^2_{i_{{}_V}}$ is $3$. From the Hodge Index
	  inequality
	  $$
	  K^2_X \leq \frac{(K_X \cdot D_{i_{{}_V}})^2}{D^2_{i_{{}_V}}}=\frac{64}{3} <22
	  $$
	  This and the lower bound in \eqref{K2lb} imply
	  $K^2_X=21$ and $\chi(\OO_X) \leq 2$ contradicting the BMY inequality
	  $$
	  21=K^2_X \leq 9\chi(\OO_X) \leq 18.
	  $$
	  Hence we now have
	  $$
	   K_X \cdot D_{i_{{}_V}} \leq 7.
	   $$
	   Arguing as before we deduce $ D^2_{i_{{}_V}}\leq 2$ and this, in view of the lower bound \eqref{D2HA}, rules out $K_X \cdot D_{i_{{}_V}} = 7$.
	   
	   For $K_X \cdot D_{i_{{}_V}} =6$ the value $D^2_{i_{{}_V}} =2$ and the Hodge Index inequality gives
	   $$
	   K^2_X \leq \frac{(K_X \cdot D_{i_{{}_V}})^2}{D^2_{i_{{}_V}}}=\frac{36}{2} =18.
	   $$
	   and the only numerical possibility is the pair
	   $$
	   (c_2 (X), K^2_X)=(6,18)
	   $$
	   on the BMY line and those are not allowed in our considerations.
	   
The value  $K_X \cdot D_{i_{{}_V}} =5$ is the last possibility to examine, since the smaller values will contradict the assumption that the divisor $D_{i_{{}_V}}$ is big. The Hodge Index inequality implies
$$
D^2_{i_{{}_V}} \leq \frac{(K_X \cdot D_{i_{{}_V}})^2}{K^2_X}\leq \frac{25}{17}<2,
$$
while the lower bound \eqref{D2HA} tells us that $D^2_{i_{{}_V}} \geq 2$. This completes the proof of the upper bound for the dimension of $V$ when the divisor $D_{i_{{}_V}}$ is nef and big.

We now turn to the case: $D_{i_{{}_V}}$ is nef and $D^2_{i_{{}_V}}=0$.
The type decomposition
$$
D_{i_{{}_V}}=D'_{i_{{}_V}}+D''_{i_{{}_V}}
$$
defined in \eqref{type-decomp} is subject to
$$
\text{$D_{i_{{}_V}} \cdot D'_{i_{{}_V}} \geq 0$ and  $D_{i_{{}_V}} \cdot D''_{i_{{}_V}} \geq 0$.}
	$$
	This together with
	$$
0=D^2_{i_{{}_V}}=	D_{i_{{}_V}} \cdot D'_{i_{{}_V}} + D_{i_{{}_V}} \cdot D''_{i_{{}_V}} 
$$
imply the equality
$$
	D_{i_{{}_V}} \cdot D'_{i_{{}_V}} =D_{i_{{}_V}} \cdot D''_{i_{{}_V}}=0.
	$$ 
From this it follows
\begin{equation}\label{D'2=D''2}
	\big(D'_{i_{{}_V}}\big)^2=\big(D''_{i_{{}_V}}\big)^2.
\end{equation}
The argument above shows
\begin{equation}\label{CintDiV=0}
	D_{i_{{}_V}} \cdot C=0,
\end{equation}
for every irreducible component $C$ of $D_{i_{{}_V}}$.  Hence
$$
C^2 \leq 0,
$$
for every irreducible component $C$ of $D_{i_{{}_V}}$.

We now consider the component $D'_{i_{{}_V}}$. 
For an irreducible component $C'$ of $D'_{i_{{}_V}}$ we have a nonzero morphism
$$
\phi_{C'}:\OO_{C'} (-D_{i_{{}_V}}) \longrightarrow \OO_{C'} (C'), 
$$
see \eqref{sigmaC-normalC}. From the above it follows
\begin{equation}\label{C'comp}
(C')^2=C'\cdot D_{i_{{}_V}}=0,
\end{equation}
for every irreducible component $C'$ of $D'_{i_{{}_V}}$. This implies that

$\bullet$ $D'_{i_{{}_V}}$ is composed of mutually disjoint irreducible components of
self-intersection $0$,

$\bullet$  $(D'_{i_{{}_V}})^2=D'_{i_{{}_V}} \cdot D''_{i_{{}_V}} =0$. 

The second item and \eqref{D'2=D''2} give
$$
(D''_{i_{{}_V}})^2=(D'_{i_{{}_V}})^2=D'_{i_{{}_V}} \cdot D''_{i_{{}_V}} =0.
$$
This in turn implies that the divisor classes 
$D'_{i_{{}_V}}$, $D''_{i_{{}_V}}$ and $D_{i_{{}_V}}$ are colinear in the 
vector space 
$NS(X)_{\QQ}$. Since $D_{i_{{}_V}}$ is nef, so are  $D'_{i_{{}_V}}$ and $D''_{i_{{}_V}}$. In addition, $D'_{i_{{}_V}}$ and $D''_{i_{{}_V}}$ have no component in common, so  $D'_{i_{{}_V}} \cdot D''_{i_{{}_V}} =0$ implies that $D'_{i_{{}_V}}$ and $D''_{i_{{}_V}}$ are disjoint. This proves the part 1) of the proposition.

 Next we complete the part 2) of the proposition. For this we return to the
 morphism $\phi_{C'}$ above: the formulas \eqref{C'comp} imply that it is an isomorphism. More precisely, we have
$$
\OO_{C'} (-m_{C'}C')=\OO_{C'} (-D_{i_{{}_V}}) \cong \OO_{C'} (C'),
$$
where $m_{C'}$ is the multiplicity of $C'$ in $D_{i_{{}_V}}$. Hence the assertion
\begin{equation}\label{C'-tor}
	\OO_{C'} ((m_{C'}+1)C')=\OO_{C'}.
\end{equation}
From the fact that $\phi_{C'}$ is an isomorphism it also follows
 that the normal sequence of $C'$ splits, see the diagram \eqref{sigmaC-normalC}. In particular, $C'$ is smooth, the line bundle
 $\OO_{C'} (C')\in Pic^{\circ} (C')$ and the equation \eqref{C'-tor} tells us that it is a torsion point of $Pic^{\circ} (C')$. Setting
 $$
 t_{C'}:=\text{the order of torsion of $\OO_{C'} (C')$ in $Pic^{\circ} (C')$,}
 $$
 we deduce from \eqref{C'-tor} that $(m_C'+1)$ is a multiple of $t_{C'}$.
 This completes the proof of part 2) of the proposition.
 
 We now turn to part 3), the study of the component $D''_{i_{{}_V}}$. Let $C''$ be its irreducible component and let
 $$
 \eta: \widetilde{C''} \longrightarrow C''
 $$
 be the normalization of $C''$. Then for a general choice of $\xi \in V$
 we have the formula
 \begin{equation}\label{C''comp}
 deg(\eta^{\ast}A^{C''}_{[\xi]})=C''\cdot D_{i_{{}_V}} +2-2g_{\widetilde{C''}}=2-2g_{\widetilde{C''}},
\end{equation}
 see \eqref{type2-estim}; here $g_{\widetilde{C''}}$ is the genus of $\widetilde{C''}$ and the last equality uses \eqref{CintDiV=0}. The equation
 implies
 $$
 g_{\widetilde{C''}}=0:
 $$
 it is obvious that the equation implies $g_{\widetilde{C''}} \leq 1$; the equality means that $A^{C''}_{[\xi]}=0$ and this, according to Lemma \ref{lem:singD''}, tells us that $C''$ is smooth and disjoint from any other irreducible component of $D''_{i_{{}_V}}$; but then
 $$
 0=C''\cdot D_{i_{{}_V}}=C'' \cdot D''_{i_{{}_V}}=m_{C''}(C'')^{2},
 $$
 in other words $(C'')^{2}=0$ and by the adjunction formula
 $$
 0=K_X \cdot C'' +(C'')^2=K_X \cdot C'',
 $$
  contradicting that we are on a minimal surface of general type.
 
 Once we know that $g_{\widetilde{C''}}=0$, the formula in \eqref{C''comp}
 reads
 $$
 deg(\eta^{\ast}A^{C''}_{[\xi]})=2.
 $$
 This implies:
 
 - either $deg(A^{C''}_{[\xi]})=2$ and then $C''$ is smooth rational curve,
  
 -  or $deg(A^{C''}_{[\xi]})=1$ and then $A^{C''}_{[\xi]}$ is a
 (unique) double point of $C''$.
 
 \vspace{0.2cm}
 For the part 4) of the proposition, we refer the reader to \cite{R}, Proposition 5.5: the proof there goes through without changes.
 
 We now turn to the proof of 5): the upper bound for the dimension of $V$.
 Arguing as in the case of $D_{i_{{}_V}}$ big we obtain the estimate
 \begin{equation}\label{dimV-isot}
 dim(V) \leq h^0(S_{i_{{}_V}})\leq h^0 (\OO_{D_{i_{{}_V}}}) +\frac{1}{3}(3c_2(X)-K^2_X) -\frac{1}{3} K_X \cdot D_{i_{{}_V}},
\end{equation}
 see \eqref{h0SiV-estim}. So we need to estimate $h^0 (\OO_{D_{i_{{}_V}}})$.
 This is done by decomposing $D_{i_{{}_V}}$ into its distinct connected components
 $$
 D_{i_{{}_V}} =\sum^N_{k=1}D_k.
 $$
From $D^2_k=D_{i_{{}_V}} \cdot D_k=0$ it follows that each component is on the ray $\QQ^{+}D_{i_{{}_V}}$ in $NS(X)_{\QQ}$
 generated by $ D_{i_{{}_V}}$. Hence each $D_k$ is nef and $D^2_k=0$ and we have
 $$
 h^0 (\OO_{D_{i_{{}_V}}}) = \sum^N_{k=1}h^0 (\OO_{D_k}).
 $$ 
 For every $D_k$ with $h^0 (\OO_{D_k})>1$ we can continue breaking $D_k$. This is well-known, but we recall the procedure:
 there is a nonzero global section $s$ of $\OO_{D_k}$ which vanishes on
 a proper component of ${D_k}$ and furthermore, does not vanish on any large component of $D_k$; denote that component by $D_k(s)$; this gives the exact sequence
 $$
 \xymatrix{
 0\ar[r]& \OO_{D^c_k (s)} (-D_k (s)) \ar[r]& \OO_{D_k} \ar[r]& \OO_{D_k(s)} \ar[r]&0,
}
 $$
 where $D^c_k (s):=D_k -D_k(s)$. By definition the global section $s$ comes from
 the global section of $\OO_{D^c_k (s)} (-D_k (s))$ and that section vanishes on at most $0$-dimensional subscheme of $D^c_k (s)$ (otherwise $s$ would vanish on a strictly larger component of $D_k$); from this it follows
 \begin{equation}\label{Dks}
D^c_k (s) \cdot D_k (s) \leq 0.
\end{equation}
But $D_k$ is nef, so we have
$$
D_k \cdot D_k(s),\, D_k \cdot D^c_k(s) \geq 0.
$$
This together with $D^2_k =0$ give
$$
0=D^2_k=D_k \cdot D_k(s)+ D_k \cdot D^c_k(s)
$$
and hence
$$
0=D_k \cdot D_k(s)= D_k \cdot D^c_k(s).
$$
This implies
$$
 D^2_k(s)= ( D^c_k(s))^2.
 $$
 Both are at most zero. This, \eqref{Dks} and the equation
 $$
 0=D_k \cdot D_k(s)=  D^2_k(s) +D^c_k (s) \cdot D_k (s) 
 $$
 imply
 $$
 ( D^c_k(s))^2= D^2_k(s)= D^c_k (s) \cdot D_k (s) =0.
 $$
 This in turn implies that the components $D^c_k (s)$, $D_k (s)$ and $D_k$
 are colinear in $NS(X)_{\QQ}$. Hence $D^c_k (s)$ and  $D_k (s)$ are both nef and isotropic.
 
 The above procedure implies that $D_{i_{{}_V}}$ can be decomposed
 $$
 D_{i_{{}_V}}=\sum_{k,j} D_{k,j}
 $$
 so that all components $D_{k,j}$ are effective divisors lying on the ray
 $\QQ^+ D_{i_{{}_V}}$ and subject to
 $$
 h^0(\OO_{D_{k,j}}) =1,
 $$
 for all $(k,j)$. Thus we deduce
 $$
 h^0 (\OO_{D_{i_{{}_V}}}) \leq \sum_{k,j}h^0(\OO_{D_{k,j}})=\sum_{k,j} 1 \leq
 \sum_{k,j} \HA K_X \cdot D_{k,j},
 $$
 where the last inequality comes from $K_X \cdot D_{k,j} \geq 2$, for all
 $(k,j)$. Hence the inequality
 $$
 h^0 (\OO_{D_{i_{{}_V}}}) \leq \sum_{k,j} \HA K_X \cdot D_{k,j}=\HA K_X \cdot D_{i_{{}_V}}.
 $$
 Substituting into \eqref{dimV-isot} gives
 $$
 dim(V) \leq \HA K_X \cdot D_{i_{{}_V}} +\frac{1}{3}(3c_2(X)-K^2_X)-\frac{1}{3}K_X \cdot D_{i_{{}_V}} =
 \frac{1}{6}K_X \cdot D_{i_{{}_V}} +\frac{1}{3}(3c_2(X)-K^2_X).
 $$
 This is the first inequality asserted in 5). The second is deduced
 from the above by using the upper bound
 $$
 K_X \cdot D_{i_{{}_V}} \leq 3c_2(X)-K^2_X -3deg(S_{[\xi]})
 $$
 in Lemma \ref{lem:Q'Vstable}, 3).
\end{pf}
	
\vspace{0.2cm}
\noindent
{\bf 7.3. Case:} $rk(Q^{l_{i_{{}_V}}}_{i_{{}_V}})=3$. We return to the diagram \eqref{QV-diag-sh}. The morphism
\begin{equation}\label{morphFVtoQiV}
F_V \longrightarrow Q^{l_{i_{{}_V}}}_{i_{{}_V}}
\end{equation}
in that diagram is of rank $2$ or $3$. If it has rank $2$, then we can reduce the diagram to the rank $2$ case considered in the previous subsection. So we assume the morphism above is of rank $3$. 

As in the rank $2$ case  we set
$$
F'_{i_{{}_V}}:=\HH^{-1}(\widetilde{E_{i_{{}_V}}}').
$$
The assumption on the rank of the morphism in \eqref{morphFVtoQiV} tells us that it is a sheaf of rank 
$$
rk(F'_{i_{{}_V}})=dim(V)-1.
$$
\begin{lem}\label{lem:F'Vrk3}
	Either $F'_{i_{{}_V}}=0$ and then $dim(V)=1$ and the diagram \eqref{QV-diag-sh} takes the form
	$$
		\xymatrix@R=12pt@C=12pt{
			&&0\ar[d]&0\ar[d]&\\
	0\ar[r]&\TET\ar@{=}[d]\ar[r]&F_V\ar[r]\ar@{=}[d]&\OO_X\ar[r]\ar[d]&0
	\\		
	0\ar[r]&\TET\ar[r]&Q^{l_{i_{{}_V}}}_{i_{{}_V}}\ar[r]\ar[d]&\widetilde{E_{i_{{}_V}}}\ar[r]\ar[d]&0\\
			&&\widetilde{E_{i_{{}_V}}}'\ar@{=}[r]\ar[d]&
			\widetilde{E_{i_{{}_V}}}'\ar[d]&\\		
			&&0&0&\\
		}
	$$
	or $dim(V)\geq 2$ and $F'_{i_{{}_V}}$ is a locally free subsheaf of $F_V$ of rank  $rk(F'_{i_{{}_V}})=dim(V)-1$. In the latter case the diagram \eqref{QV-diag-sh}  simplifies to	

\begin{equation}\label{rkQV=3diag}
	\xymatrix@R=12pt@C=12pt{
		&&0\ar[d]&0\ar[d]&\\
		&&F'_{i_{{}_V}}\ar[d]\ar@{=}[r]&F'_{i_{{}_V}} \ar[d]&\\	
		0\ar[r]&\TET\ar[r]\ar@{=}[d]&F_V \ar[r]\ar[d]&V\otimes \OO_X \ar[r]\ar[d]&0\\
		0\ar[r]&\TET \ar[r]&Q'_{i_{{}_V}}\ar[r]\ar[d]&{S_{i_{{}_V}}}\ar[r]\ar[d]&0\\
		&&0&0&
	}
\end{equation}
where $Q'_{i_{{}_V}}$ is the quotient sheaf $F_V /F'_{i_{{}_V}}$ and it is a part of the following commutative diagram 
\begin{equation}\label{Q'iVSiVdiag3}
\xymatrix@R=12pt@C=12pt{
	&0\ar[d]&0\ar[d]&&\\
	&\TET\ar[d]\ar@{=}[r]&\TET\ar[d]&&\\
	0\ar[r]&Q'_{i_{{}_V}} \ar[d]\ar[r]&Q^{l_{i_{{}_V}}}_{i_{{}_V}}\ar[d] \ar[r]&\HH^0 (\widetilde{E_{i_{{}_V}}}')\ar[r]\ar@{=}[d]&0\\
0\ar[r]&S_{i_{{}_V}}\ar[d]\ar[r]&\widetilde{E_{i_{{}_V}}}\ar[d]\ar[r]
&\HH^0 (\widetilde{E_{i_{{}_V}}}')\ar[r]&0\\
	&0&0&&
}	 
\end{equation}
In particular, $Q'_{i_{{}_V}}$ is a torsion free sheaf of rank $3$ in $\FF_{K_X,bK_X}$, $S_{i_{{}_V}}$ is globally generated and $c_1 ( \widetilde{E_{i_{{}_V}}})$ is an effective nonzero divisor subject to
\begin{equation}\label{tildEiV-formula}
c_1 (\widetilde{E_{i_{{}_V}}}) =c_1 (S_{i_{{}_V}}) +c_1 (\HH^0 (\widetilde{E_{i_{{}_V}}}'))
\end{equation}
and the upper bound
\begin{equation}\label{tildEiV-degbound}
c_1 (S_{i_{{}_V}})\cdot K_X \leq 	c_1 (\widetilde{E_{i_{{}_V}}})\cdot K_X \leq 2(3c_2(X)-K_X).
\end{equation}
\end{lem}
\begin{pf}
	Only the statement that $c_1 (\widetilde{E_{i_{{}_V}}})$ is effective and nonzero perhaps needs a proof. The formula
	$$
c_1 (\widetilde{E_{i_{{}_V}}}) =c_1 (S_{i_{{}_V}}) +c_1 (\HH^0 (\widetilde{E_{i_{{}_V}}}'))
$$
certainly tells us that $c_1 (\widetilde{E_{i_{{}_V}}})$ is effective, since both summands are: the first because $S_{i_{{}_V}}$ is globally generated and the second because $\HH^0 (\widetilde{E_{i_{{}_V}}}')$ is a torsion sheaf. Furthermore, according to the above formula the divisor
$c_1 (\widetilde{E_{i_{{}_V}}})$ is zero if and only if
$\HH^0 (\widetilde{E_{i_{{}_V}}}')$ is supported on at most $0$-dimensional subscheme of $X$ and $S_{i_{{}_V}}=\OO_X$. But then the right column in the diagram \eqref{rkQV=3diag} implies that $V$ is one dimensional which contrary to the situation under consideration.

The first inequality in \eqref{tildEiV-degbound} follows from the formula	\eqref{tildEiV-formula} and the effectiveness of $c_1 (\HH^0 (\widetilde{E_{i_{{}_V}}}'))$ and the second comes from Lemma \ref{lem:tildEi}.
\end{pf}

From now on we assume that the vector space $V$ is at least of dimension $2$ and we wish to understand stability properties of $Q'_{i_{{}_V}}$ with respect to $K_X$.

\begin{lem}\label{lem:Q'iV-stab}
	Assume $dim(V) \geq 2$. Then the exact sequence
	\begin{equation}\label{FiVdest}
	\xymatrix@R=12pt@C=12pt{
		0\ar[r]& F'_{i_{{}_V}} \ar[r]& F_V \ar[r]& Q'_{i_{{}_V}}\ar[r]&0,
	}
\end{equation}
the middle column of the diagram \eqref{rkQV=3diag}, is a $K_X$-destabilizing sequence of $F_V$.

	If $Q'_{i_{{}_V}}$ is $K_X$-semistable, then it is the last
	$K_X$-semistable factor of the HN filtration of $F_V$ with respect to
	$K_X$, in other words $Q'_{i_{{}_V}}=Q^m_V$, the sheaf defined in Proposition \ref{pro:FVunst}; dually, the sheaf $(Q'_{i_{{}_V}})^{\ast}$, the dual of $Q'_{i_{{}_V}}$, is the maximal $K_X$-destabilizing subsheaf of
	$F^{\ast}_V$, the dual of $F_V$.

	If $Q'_{i_{{}_V}}$ is $K_X$-unstable, then it fits into the exact sequence
	$$
\xymatrix@R=12pt@C=12pt{
	0\ar[r]& {\cal J}_A (-M) \ar[r]& Q'_{i_{{}_V}}\ar[r]&Q^m_V \ar[r]&0,
}
$$
where ${\cal J}_A (-M) $ is the maximal $K_X$-destabilizing subsheaf of
$Q'_{i_{{}_V}}$ and $Q^m_V$ is the last semistable factor in the HN filtration of $F_V$ with respect to $K_X$, 
 the sheaf defined in Proposition \ref{pro:FVunst}. 
\end{lem}
\begin{pf}
To see that the sequence \eqref{FiVdest} is $K_X$-destabilizing for 
$F_V$ amounts to showing that $F'_{i_{{}_V}} $ is a $K_X$-destabilizing subsheaf of $F_V$ and this comes down to checking the inequality of slopes
$$
\frac{c_1(F'_{i_{{}_V}})\cdot K_X}{dim(V)-1} > \frac{c_1 (F_V)\cdot K_X}{dim(V)+2}=-\frac{ K^2_X}{dim(V)+2}
$$
or, equivalently,
$$
-c_1(F'_{i_{{}_V}})\cdot K_X <\frac{dim(V)-1}{dim(V)+2}K^2_X.
$$
This is seen as follows: from the right column in the diagram \eqref{rkQV=3diag} we obtain
$$
-c_1(F'_{i_{{}_V}})\cdot K_X=c_1 (S_{i_{{}_V}})	\cdot K_X;
$$
now use the formula \eqref{tildEiV-formula} to deduce
$$
-c_1(F'_{i_{{}_V}})\cdot K_X=c_1 (S_{i_{{}_V}})	\cdot K_X \leq c_1(\widetilde{E_{i_{{}_V}}})\cdot K_X;
$$
next use the the estimate in Lemma \ref{lem:tildEi} to obtain
$$
-c_1(F'_{i_{{}_V}})\cdot K_X \leq c_1(\widetilde{E_{i_{{}_V}}})\cdot K_X
<2(3c_2(X)-K^2_X)=2(3\alpha_X -1)K^2_X< \frac{1}{4}K^2_X,
$$
where the last inequality uses the assumption
$\alpha_X <\frac{3}{8}$. It remains to observe that
$$
\frac{1}{4} \leq \frac{dim(V)-1}{dim(V)+2},
$$
for any value $dim(V)\geq 2$. Thus the inequality of slopes
\begin{equation}\label{c1FiVquater}
-c_1(F'_{i_{{}_V}})\cdot K_X < \frac{1}{4}K^2_X \leq \frac{dim(V)-1}{dim(V)+2}K^2_X
\end{equation}
is proved.

If $Q'_{i_{{}_V}}$ is $K_X$-semistable, the above estimate also
shows that $Q'_{i_{{}_V}}$ is the last semistable factor of the HN filtration of $F_V$ with respect to $K_X$. Indeed, if $F'_{i_{{}_V}}$
is $K_X$-semistable, then we are done since
$$
F_V \supset F'_{i_{{}_V}} \supset 0
$$
is the HN filtration of $F_V$ with respect to $K_X$.
If $F'_{i_{{}_V}}$ is $K_X$-unstable we take its HN filtration with respect to $K_X$:
$$
F'_{i_{{}_V}}=F'_l \supset\cdots \supset F'_1 \supset F'_0=0.
$$
To show that $Q'_{i_{{}_V}}$ is the last $K_X$-semistable factor of the 
HN filtration of $F_V$ come down to checking the inequality of slopes
$$
\frac{c_1(Q'_{i_{{}_V}})\cdot K_X}{3} <\frac{c_1 (F'_l /F'_{l-1})\cdot K_X}{r_l},
$$
where $r_l$ is the rank of $F'_l /F'_{l-1}$.

Assume the opposite inequality holds
$$
\frac{c_1(Q'_{i_{{}_V}})\cdot K_X}{3} \geq \frac{c_1 (F'_l /F'_{l-1})\cdot K_X}{r_l}.
$$
 This gives
$$
-\frac{c_1(Q'_{i_{{}_V}})\cdot K_X}{3} \leq -\frac{c_1 (F'_l /F'_{l-1})\cdot K_X}{r_l}
$$
Now recall: $F'_{i_{{}_V}}$ is a sheaf in $\FF_{K_X,bK_X}$, so all its
semistable factors $F'_i/F'_{i-1}$ are in $\FF_{K_X,bK_X}$ and hence
$$
-c_1(F'_i/F'_{i-1})\cdot K_X >0
$$
for all $i\in[1,l]$. Thus previous inequality implies
$$
-\frac{c_1(Q'_{i_{{}_V}})\cdot K_X}{3} < -c_1(F'_{i_{{}_V}})\cdot K_X=c_1
(S_{i_{{}_V}})\cdot K_X.
$$
This and the relation
$$
K_X =-c_1(Q'_{i_{{}_V}}) +c_1
(S_{i_{{}_V}})
$$
coming from the bottom row of the diagram \eqref{rkQV=3diag} give the inequality
$$
K^2_X <4c_1
(S_{i_{{}_V}})\cdot K_X=-4c_1(F'_{i_{{}_V}})\cdot K_X
$$
which is in contradiction with the first inequality in \eqref{c1FiVquater}.

We now turn to the case of $Q'_{i_{{}_V}}$ being $K_X$-unstable. Denote by $Q'_1$ its maximal $K_X$-destabilizing subsheaf. Our first assertion says that $Q'_1$ has rank $1$. 

Assume it is not the case. Then the rank of $Q'_1$ is $2$ and it gives the destabilizing exact sequence
$$
\xymatrix@R=12pt@C=12pt{
	0\ar[r]&Q'_1 \ar[r]&Q'_{i_{{}_V}}\ar[r]&{\cal J}_A (-L) \ar[r]&0,
}
$$
where ${\cal J}_A $ is the ideal sheaf of at most $0$-dimensional subscheme $A$ of $X$ and $L$ is a divisor subject to
\begin{equation}\label{slopedestQiV}
	-L\cdot K_X <\frac{c_1 (Q'_1)\cdot K_X}{2}.
\end{equation}
To rule out this situation it is enough to show that ${\cal J}_A (-L)$ is the last semistable factor in the HN filtration of $F_V$ with respect to $K_X$: recall, this 
 this is impossible in view of Proposition \ref{pro:FVunst} which says that the last $K_X$-semistable factor of  the HN filtration of $F_V$ with respect to $K_X$ has rank $2$ or $3$.
 
 Consider the composite epimorphism
 $$
 \xymatrix@R=12pt@C=12pt{
 F_V \ar[r] \ar@/_1.5pc/[rr]_h &Q'_{i_{{}_V}} \ar[r]& {\cal J}_A (-L)
}
$$
and set
$$
F':=ker(h).
$$
This gives the filtration
$$
F_V \supset F' \supset F'_{i_{{}_V}}.
$$
We complete this by the HN filtration of $F'_{i_{{}_V}}$ with respect to is $K_X$:
$$
F'_{i_{{}_V}}=F'_l \supset \cdots \supset F'_1 \supset F'_0=0
$$
For the resulting filtration of $F_V$ to be its HN filtration we need to
check the inequality
of slopes
$$
\frac{c_1 (Q'_1)\cdot K_X}{2} < \frac{c_1(F'_l /F'_{l-1})\cdot K_X}{r_l}
$$
where $r_l$ is the rank of $F'_l/F'_{l-1}$.
Since all semistable factors of $F'_{i_{{}_V}}$ have negative degrees with respect to $K_X$ (recall: all those sheaves are in $\FF_{K_X,bK_X}$)
it will be enough to show
$$
\frac{c_1 (Q'_1)\cdot K_X}{2} < c_1(F'_{i_{{}_V}})\cdot K_X =-c_1 (S_{i_{{}_V}})\cdot K_X.
$$
Assume the opposite inequality holds
$$
\frac{c_1 (Q'_1)\cdot K_X}{2} \geq -c_1 (S_{i_{{}_V}})\cdot K_X.
$$
 Observe that we have inclusions
$$
Q'_1 \hookrightarrow Q'_{i_{{}_V}} \hookrightarrow  Q^{l_{i_{{}_V}} }_{i_{{}_V}}.
$$
Since the rightmost sheaf is $K_X$-semistable we deduce
$$
 \frac{c_1 ( Q^{l_{i_{{}_V}} }_{i_{{}_V}})\cdot K_X}{3} \geq  \frac{c_1 (Q'_1)\cdot K_X}{2} \geq -c_1 (S_{i_{{}_V}})\cdot K_X.
 $$
 This and the middle column of the diagram \eqref{Q'iVSiVdiag3} give
 $$
 c_1(\widetilde{E_{i_{{}_V}}})\cdot K_X -K^2_X \geq -3c_1 (S_{i_{{}_V}})\cdot K_X
 $$
 or, equivalently,
 $$
 c_1(\widetilde{E_{i_{{}_V}}})\cdot K_X +3c_1 (S_{i_{{}_V}})\cdot K_X
 \geq K^2_X.
 $$
 This and the formula for $ c_1(\widetilde{E_{i_{{}_V}}})$ in \eqref{tildEiV-formula} imply
 $$
 4c_1(\widetilde{E_{i_{{}_V}}})\cdot K_X \geq K^2_X.
 $$
 Combining this with the upper bound in Lemma \ref{lem:tildEi} we obtain
 $$
 8(3c_2(X)-K^2_X) >K^2_X.
 $$
 Rewriting in terms of the ratio $\displaystyle{\alpha_X =\frac{c_2(X)}{K^2_X}}$ of the Chern numbers of $X$ gives
 $$
 3\alpha_X -1 > \frac{1}{8}
 $$
 and this contradicts the assumption $\alpha_X <\frac{3}{8}$.
 
 We know now that the maximal $K_X$-destabilizing subsheaf $Q'_1$ of 
 $Q'_{i_{{}_V}}$ has rank $1$. Being torsion free it has the form
$$
Q'_1={\cal J}_A (-M)
$$
and gives rise to the exact sequence
\begin{equation}\label{Q'iVrk1dest}
	\xymatrix@R=12pt@C=12pt{
		0\ar[r]&{\cal J}_A (-M)\ar[r]&Q'_{i_{{}_V}} \ar[r]& Q''_{i_{{}_V}} \ar[r]&0
	}
		\end{equation}
		with the slope relation
		$$
		-M\cdot K_X > \frac{c_1 (Q''_{i_{{}_V}})\cdot K_X}{2}.
		$$
		Furthermore, the argument used to prove that the rank of $Q'_1$ is one also gives the inequality
		$$
	c_1(F'_{i_{{}_V}})\cdot K_X =	-c_1(S_{i_{{}_V}})\cdot K_X >-M\cdot K_X.
	$$
	Setting as before
	$$
	h:F_V \longrightarrow Q'_{i_{{}_V}} \longrightarrow Q''_{i_{{}_V}},
	$$
	the composition of two epimorphisms, gives the filtration
	$$
	F_V \supset F'' \supset F'_{i_{{}_V}},
	$$
	where $F''=ker(h)$. Completing this with the HN filtration of 
	$ F'_{i_{{}_V}}$ with respect to $K_X$, we will obtain the HN filtration of $F_V$ as soon as we know that $Q''_{i_{{}_V}}$ is
	$K_X$-semistable. So we need to prove the $K_X$-semistability of
	that sheaf.
	
	For this assume on the contrary that $Q''_{i_{{}_V}}$ is $K_X$-unstable and let
	$$
\xymatrix@R=12pt@C=12pt{
	0\ar[r]&{\cal J}_{A'} (-L_1) \ar[r]&Q''_{i_{{}_V}} \ar[r]&{\cal J}_{A''} (-L_2) \ar[r]&0
}
$$
to be its destabilizing sequence, that is,
$$
-L_1 \cdot K_X > -L_2 \cdot K_X.
$$
Furthermore, we must have
\begin{equation}\label{L1-M}
-L_1\cdot K_X \geq -M\cdot K_X,
\end{equation}
since otherwise ${\cal J}_{A''} (-L_2)$ is the last $K_X$-semistable factor of the HN filtration of $F_V$ and we know from Proposition \ref{pro:FVunst} that this is impossible.

We now combine the destabilizing sequence \eqref{Q'iVrk1dest} for $Q'_{i_{{}_V}}$ with the one for $Q''_{i_{{}_V}} $:
$$
\xymatrix@R=12pt@C=12pt{
	&&&0\ar[d]&\\
	&&&{\cal J}_{A'} (-L_1)\ar[d]&\\
	0\ar[r]&{\cal J}_A (-M)\ar[r]&Q'_{i_{{}_V}} \ar[r]& Q''_{i_{{}_V}} \ar[r]\ar[d]&0\\
	&&&{\cal J}_{A''} (-L_2)\ar[d]\\
	&&&0
}
$$
This can be completed to the following commutative diagram
$$
\xymatrix@R=12pt@C=12pt{
	&&0\ar[d]&0\ar[d]&\\
	0\ar[r]&{\cal J}_A (-M)\ar[r]\ar@{=}[d]&Q\ar[r]\ar[d]&{\cal J}_{A'} (-L_1)\ar[d]\ar[r]&0\\
	0\ar[r]&{\cal J}_A (-M)\ar[r]&Q'_{i_{{}_V}} \ar[r]\ar[d]& Q''_{i_{{}_V}} \ar[r]\ar[d]&0\\
	&&{\cal J}_{A''} (-L_2)\ar@{=}[r]\ar[d]&{\cal J}_{A''} (-L_2)\ar[d]\\
	&&0&0
}
$$
The subsheaf $Q$ of $Q'_{i_{{}_V}}$ resulting from this diagram has rank $2$ an degree 
$$
c_1 (Q)\cdot K_X =-M\cdot K_X -L_1 \cdot K_X \geq -2M\cdot K_X,
$$
where the last inequality comes from \eqref{L1-M}. This tells us that 
$$
\frac{c_1 (Q)\cdot K_X }{2} \geq -M\cdot K_X.
$$
From the maximality of ${\cal J}_A (-M)$ it follows that the above inequality must be equality and this means
$$
-L_1 \cdot K_X =-M\cdot K_X.
$$
 But then $Q$ is $K_X$-semistable and it is the maximal $K_X$-destabilizing subsheaf of $Q'_{i_{{}_V}}$ and this is impossible since we already proved that the maximal $K_X$-destabilizing subsheaf of $Q'_{i_{{}_V}}$ has rank $1$.	 
\end{pf}

The above result tells us that when $Q'_{i_{{}_V}}$ is $K_X$-unstable
we can perform a further reduction of the bundle $F_V$  to its last $K_X$-semistable factor $Q^m_V$ of rank $2$
$$
\xymatrix@R=12pt@C=12pt{
	0\ar[r]&F^{m-1}_V \ar[r]&F_V\ar[r]&Q^m_V \ar[r]&0
}
$$
and this reduction is related to the rank $3$ reduction, the middle column in \eqref{Q'iVSiVdiag3}, by the following commutative diagram
\begin{equation}\label{QmV-diag}
\xymatrix@R=12pt@C=12pt{
	&&0\ar[d]&0\ar[d]&&\\
&0\ar[r]&F'_{i_{{}_V}}\ar[d]\ar[r]&F^{m-1}_V \ar[r]\ar[d]&{\cal J}_A (-M)\ar[r]&0\\	
&&F_V \ar[d]\ar@{=}[r]&F_V\ar[d]&&\\
0\ar[r]&{\cal J}_A (-M) \ar[r]&Q'_{i_{{}_V}}\ar[r]\ar[d]&Q^m_V\ar[r]\ar[d]&0&\\
&&0&0&&
}		
\end{equation}
In addition, the two columns in the above diagram fit together
with the extension sequence
$$
\xymatrix@R=12pt@C=12pt{
	0\ar[r]&\TET \ar[r]&F_V\ar[r]&V\otimes \OO_X \ar[r]&0
}
$$
in the sense that the monomorphism $\TET \longrightarrow F_V$ composes with epimorphisms of the two columns to give the commutative diagram
\begin{equation}\label{Q'iVQmV-diag}
\xymatrix@R=12pt@C=12pt{
	&&0\ar[d]&0\ar[d]&\\	
	&& \TET \ar@{=}[r]\ar[d] & \TET \ar[d] & \\
0\ar[r]&{\cal J}_A (-M)\ar[r]\ar@{=}[d]&Q'_{i_{{}_V}}\ar[r]\ar[d]&
Q^m_V \ar[r]\ar[d]&0\\
0\ar[r]&{\cal J}_A (-M)\ar[r]&S_{i_{{}_V}}\ar[r]\ar[d]&\widetilde{E_V}\ar[r]\ar[d]&0\\
&&0&0&
}
\end{equation}
Similarly, the monomorphisms of two columns in \eqref{QmV-diag}	compose with the epimorphism 
$$
F_V \longrightarrow V\otimes \OO_X
$$ 
of the extension sequence to give the commutative diagram
\begin{equation}\label{F'iVFmV-diag}
	\xymatrix@R=12pt@C=12pt{
	&&0\ar[d]&0\ar[d]&&\\
	&0\ar[r]&F'_{i_{{}_V}}\ar[d]\ar[r]&F^{m-1}_V \ar[r]\ar[d]&{\cal J}_A (-M)\ar[r]&0\\			
	&	& V\otimes \OO_X\ar@{=}[r]\ar[d] & V\otimes \OO_X \ar[d] & &\\
		0\ar[r]&{\cal J}_A (-M)\ar[r]&S_{i_{{}_V}}\ar[r]\ar[d]&
		\widetilde{E_V} \ar[r]\ar[d]&0&\\
	&	&0&0&&
	}
\end{equation}
All of the above can be summarized as follows
\begin{cor}
	If we are in the situation of Lemma \ref{lem:Q'iV-stab} and
	the quotient sheaf $Q'_{i_{{}_V}}$ is $K_X$-unstable then the
	extension sequence
	$$
\xymatrix@R=12pt@C=12pt{
	0\ar[r]&\TET \ar[r]&F_V\ar[r]&V\otimes \OO_X \ar[r]&0
}
$$
reduces to the case of rank $2$ considered in the subsection {\bf 7.2}.
More precisely, the extension sequence
fits into the diagram
$$
\xymatrix@R=12pt@C=12pt{
	&&0\ar[d]&0\ar[d]&\\
	&&F^{m-1}_V \ar[d]\ar@{=}[r]&F^{m-1}_V \ar[d]&\\
	0\ar[r]&\TET \ar[r]\ar@{=}[d]&F_V\ar[r]\ar[d]&V\otimes \OO_X \ar[r]\ar[d]&0\\
		0\ar[r]&\TET \ar[r]&Q^m_V \ar[r]\ar[d]&\widetilde{E_V}\ar[r]\ar[d]&0\\
		&&0&0
				}
			$$
			where the sheaf $Q^m_V$ is of rank $2$ and it is the last 
			$K_X$-semistable factor in the HN filtration of $F_V$ with respect to $K_X$. All the properties of Proposition \ref{pro:QVrk2} are valid for this diagram with the sheaf
			$Q^m_V$ and the divisor class $D_V:=c_1 (\widetilde{E_V})$
			playing respectively the roles of $Q'_{i_{{}_V}}$ and $D_{i_{{}_V}}$ of {\rm that} proposition.
\end{cor}

The moral of the above discussion is that in treating the case of rank
$3$ we may as well assume that the rank $3$ sheaf $Q'_{i_{{}_V}}$ in the diagram \eqref{rkQV=3diag} is $K_X$-semistable (otherwise we are back in the rank $2$ case.)
	
\vspace{0.5cm}
\noindent
{\bf 7.3.1. The sheaf $S_{i_{{}_V}}$.} The geometric output of our considerations is partly contained in the sheaf $S_{i_{{}_V}}$. We examine it in more details.

From the diagram \eqref{rkQV=3diag} we deduce the following.
\begin{lem}\label{lem:SiV-Vgen}
	1) $S_{i_{{}_V}}$ is a sheaf of rank $1$.
	
	2) $V$ injects into $H^0 (S_{i_{{}_V}})$.
	
	3) $S_{i_{{}_V}}$ is generated by its global sections. More precisely, the sheaf $S_{i_{{}_V}}$ is globally generated by the space  $V$, viewed as a subspace of global sections of $S_{i_{{}_V}}$ via the inclusion in 2).
\end{lem}
\begin{pf}
	The rank of $Q'_{i_{{}_V}}$ is $3$, so the bottom row in \eqref{rkQV=3diag} implies that the rank of  $S_{i_{{}_V}}$ is $1$.
	
	The second assertion follows from the right column in \eqref{rkQV=3diag}: the cohomology sequence gives
	$$
	\xymatrix@R=12pt@C=12pt{
		0\ar[r]& H^0 (F'_{i_{{}_V}})  \ar[r]&V \ar[r]&H^0(S_{i_{{}_V}});
	}
$$ 
the asserted injectivity comes from the vanishing
$H^0 (F'_{i_{{}_V}})=Hom (\OO_X,F'_{i_{{}_V}})=0$ and this is because $F'_{i_{{}_V}}$ is a sheaf in
$\FF_{K_X,bK_X}$ and $\OO_X$ is in $\TT_{K_X,bK_X}$ and no nonzero morphisms
$$
\OO_X \longrightarrow F'_{i_{{}_V}}
$$
is allowed. 

The epimorphism $V\otimes \OO_X \longrightarrow S_{i_{{}_V}}$ in the right column of \eqref{rkQV=3diag} means that $S_{i_{{}_V}}$ is globally generated by $V$, viewed as a subspace of $H^0(S_{i_{{}_V}})$.
\end{pf}
	
The sheaf $S_{i_{{}_V}}$ may have torsion. This is recorded by the standard exact sequence
$$
	\xymatrix@R=12pt@C=12pt{
		0\ar[r]& Tor(S_{i_{{}_V}})\ar[r]&S_{i_{{}_V}}\ar[r]&S^{tf}_{i_{{}_V}} \ar[r]&0,
	} 
$$
where $Tor(S_{i_{{}_V}})$ and $S^{tf}_{i_{{}_V}}$ are torsion and torsion free parts of $S_{i_{{}_V}}$ respectively. We want to investigate how the global sections of $S_{i_{{}_V}}$ coming from $V$ are distributed between the parts of $S_{i_{{}_V}}$.

Being of rank $1$, the sheaf $S^{tf}_{i_{{}_V}}$ has the form
$$
S^{tf}_{i_{{}_V}}:={\cal J}_{A_{i_{{}_V}}} (L_{i_{{}_V}}),
$$
where $L_{i_{{}_V}}$ is a divisor and ${\cal J}_{A_{i_{{}_V}}}$ is the ideal sheaf of at most $0$-dimensional subscheme of $X$.

The surjective morphism in \eqref{rkQV=3diag}
$$
V\otimes \OO_X \longrightarrow S_{i_{{}_V}}
$$
composed with the epimorphism in the standard exact sequence gives
an epimorphism
\begin{equation}\label{epiLV}
V\otimes \OO_X \longrightarrow {\cal J}_{A_{i_{{}_V}}} (L_{i_{{}_V}}).
\end{equation}
Thus ${\cal J}_{A_{i_{{}_V}}} (L_{i_{{}_V}})$ is generated by its global sections and we deduce
\begin{lem}\label{lem:LV}
	The linear system $|L_{i_{{}_V}}|$ has at most $0$-dimensional base locus. In particular, the divisor $L_{i_{{}_V}}$ is nef and 
	$L^2_{i_{{}_V}} \geq deg(A_{i_{{}_V}})$.
	\end{lem}

Set $F_{L_{i_{{}_V}}}$ to be the kernel of the epimorphism in \eqref{epiLV}. This fits into the following exact sequences
\begin{equation}\label{FLV}
	\xymatrix@R=12pt@C=12pt{
		0\ar[r]&F_{L_{i_{{}_V}}}\ar[r]& V\otimes \OO_X \ar[r]& {\cal J}_{A_{i_{{}_V}}} (L_{i_{{}_V}})\ar[r]&0,\\
0\ar[r]&F'_{i_{{}_V}} \ar[r]&F_{L_{i_{{}_V}}}\ar[r]& Tor(S_{i_{{}_V}})	\ar[r]&0.
}
\end{equation}
In addition, the epimorphism \eqref{epiLV} fits into the commutative diagram
$$
	\xymatrix@R=12pt@C=12pt{
		&&0\ar[d]&0\ar[d]&&\\
	&0\ar[r]	&F'_{i_{{}_V}} \ar[r]\ar[d]&F_{L_{i_{{}_V}}} \ar[r]\ar[d]&Tor(S_{i_{{}_V}})\ar[r]&0\\
	0\ar[r]&\TET\ar[r]&	F_V \ar[r]\ar[d]&V\otimes \OO_X \ar[r]\ar[d]&0&\\
	&&	Q'_{i_{{}_V}} \ar[r]\ar[d]&{\cal J}_{A_{i_{{}_V}}} (L_{i_{{}_V}}) \ar[r]\ar[d]&0&\\
	&&0&0&
	}
$$
where the middle column is the one in the middle of the diagram \eqref{rkQV=3diag}. Setting
$$
Q''_{i_{{}_V}}:=ker(Q'_{i_{{}_V}} \longrightarrow {\cal J}_{A_{i_{{}_V}}} (L_{i_{{}_V}})),
$$
the kernel of the bottom row of the above diagram, we obtain the exact sequence
\begin{equation}\label{Q''V}
	\xymatrix@R=12pt@C=12pt{
		0\ar[r]&\TET \ar[r]&Q''_{i_{{}_V}} \ar[r]&Tor(S_{i_{{}_V}}) \ar[r]&0
	} 
\end{equation}
\begin{lem}\label{lem:H0FLV}
	Either $H^0 (F_{L_{i_{{}_V}}})=0$ and then there is an injective homomorphism
	$$
	V\hookrightarrow H^0 ({\cal J}_{A_{i_{{}_V}}} (L_{i_{{}_V}}));
	$$
	in particular, $dim(V)\leq h^0 ({\cal J}_{A_{i_{{}_V}}} (L_{i_{{}_V}}))$,
	
	 or 	$H^0 (F_{L_{i_{{}_V}}})\neq 0$ and the sheaf $F_{L_{i_{{}_V}}}$ splits as follows:
	$$
	F_{L_{i_{{}_V}}}=H^0 (F_{L_{i_{{}_V}}})\otimes \OO_X \oplus F'_{L_{i_{{}_V}}}.
	$$
	\end{lem}
\begin{pf}
	From the first exact sequence in \eqref{FLV} we have
	$$
	\xymatrix@R=12pt@C=12pt{
		0\ar[r]& H^0 (F_{L_{i_{{}_V}}})	\ar[r]& V \ar[r]&H^0 ({\cal J}_{A_{i_{{}_V}}} (L_{i_{{}_V}})). 
	}
$$ 
Hence the first assertion of the lemma. We now assume 
	$H^0 (F_{L_{i_{{}_V}}})\neq 0$ and consider a nonzero global section $s$ of
$F_{L_{i_{{}_V}}}$ as a monomorphism
$$
s: \OO_X \longrightarrow F_{L_{i_{{}_V}}}.
$$
On the other hand the dual of the first exact sequence in \eqref{FLV}
tells us that $ F^{\ast}_{L_{i_{{}_V}}}$, the dual of $ F_{L_{i_{{}_V}}}$, is generated by its global sections outside of the subscheme $A_{i_{{}_V}}$. So for a general
global section $t^{\ast}$ of $ F^{\ast}_{L_{i_{{}_V}}}$ the contraction morphism
$$
t^{\ast}: F_{L_{i_{{}_V}}} \longrightarrow \OO_X
$$
composed with $s$ must be nonzero. In other words in the diagram
$$
\xymatrix@R=12pt@C=12pt{
	\OO_X \ar[d]^s \ar[rd]&\\
	  F^{\ast}_{L_{i_{{}_V}}}\ar[r]^{t^{\ast}}&\OO_X
	}
$$
the composition, the slanted arrow of the diagram, is nonzero and hence an isomorphism giving a splitting
$$
 F^{\ast}_{L_{i_{{}_V}}} =\OO_X \oplus ker(t^{\ast}).
$$
If the summand $ker(t^{\ast})$ has a nonzero global section we repeat the
above argument until we arrive to the asserted splitting.
\end{pf}

Assume $H^0(F_{L_{i_{{}_V}}})$ is nonzero. From the first exact sequence in \eqref{FLV} it follows that
$H^0(F_{L_{i_{{}_V}}})$ injects into $V$. We denote the resulting subspace of $V$ by
$V^0$. In addition, the commutative diagram 
$$
\xymatrix@R=12pt@C=12pt{
0\ar[r]&V^0 \otimes \OO_X\ar[r]&V\otimes \OO_X \ar[r]\ar[d]&\Big(V/V_0 \Big) \otimes \OO_X \ar[r]\ar[d]&0\\
0\ar[r]& Tor(S_{i_{{}_V}})\ar[r]& S_{i_{{}_V}} \ar[r]\ar[d]& {\cal J}_{A_{i_{{}_V}}} (L_{i_{{}_V}})\ar[d]\ar[r]&0\\
	&&0&0&
}
$$
can be completed in a unique way to a morphism of exact sequences
\begin{equation}\label{V0V-diag}
\xymatrix@R=12pt@C=12pt{
	0\ar[r]&V^0 \otimes \OO_X\ar[r] \ar[d]&V\otimes \OO_X \ar[r]\ar[d]&\Big(V/V_0 \Big) \otimes \OO_X \ar[r]\ar[d]&0\\
	0\ar[r]& Tor(S_{i_{{}_V}})\ar[r]& S_{i_{{}_V}} \ar[d]\ar[r]& {\cal J}_{A_{i_{{}_V}}} (L_{i_{{}_V}})\ar[r]\ar[d]&0\\
	&&0&0&
}
\end{equation}
\begin{lem}\label{lem:V0}
	The morphism $V^0 \otimes \OO_X \longrightarrow Tor(S_{i_{{}_V}})$
in the above diagram is surjective and gives an inclusion
	$$
	V^0 \hookrightarrow H^0 (Tor(S_{i_{{}_V}})).
	$$
In particular, 	$Tor(S_{i_{{}_V}})$ is globally generated by the subspace $V^0$ of its global sections. 

The inclusion above fits into the following commutative diagram
	$$
\xymatrix@R=12pt@C=12pt{
	0\ar[r]&V^0 \ar[r]\ar@{^{(}->}[d]& V\ar@{^{(}->}[d]\\
	0 \ar[r]& H^0 (Tor(S_{i_{{}_V}})) \ar[r]&H^1 (\TET)
}
$$
where the monomorphism at the bottom comes from the cohomology sequence of the exact sequence \eqref{Q''V}.	
	\end{lem}
\begin{pf}
	From the diagram \eqref{V0V-diag} it follows that the cokernel
	$$
	C^0:=coker(V^0 \otimes \OO_X \longrightarrow Tor(S_{i_{{}_V}}))
	$$
	injects into $\Big(V/V_0 \Big) \otimes \OO_X$. But the latter is locally free and $C^0$ is a torsion sheaf. Hence $C^0=0$.
	
	On the level of global sections the diagram \eqref{V0V-diag} induces the commutative square
	\begin{equation}\label{V0-V-sq}
\xymatrix@R=12pt@C=12pt{
	0\ar[r]&V^0 \ar[r]\ar[d]& V\ar[d]\\
	0\ar[r]&H^0( Tor(S_{i_{{}_V}}))\ar[r]&H^0(S_{i_{{}_V}})	
}
\end{equation}
with exact rows. From Lemma \ref{lem:SiV-Vgen}, 2), the vertical arrow on the right is injective. Hence the arrow on the left is injective as well.

For the last assertion we observe that the exact sequence \eqref{Q''V} is part of the diagram
$$
\xymatrix@R=12pt@C=12pt{
	&&0\ar[d]&0\ar[d]&\\
		0\ar[r]&\TET \ar[r]\ar@{=}[d]&Q''_{i_{{}_V}} \ar[r]\ar[d]&Tor(S_{i_{{}_V}})\ar[d] \ar[r]&0\\
	0\ar[r]&\TET \ar[r]&Q'_{i_{{}_V}} \ar[r]\ar[d]&S_{i_{{}_V}} \ar[r]\ar[d]&0\\
	&&{\cal J}_{A_{i_{{}_V}}} (L_{i_{{}_V}})\ar@{=}[r]\ar[d]&{\cal J}_{A_{i_{{}_V}}} (L_{i_{{}_V}})\ar[d]&\\
	&&0&0		
}
$$
On the cohomology level we obtain the commutative square
$$
\xymatrix@R=12pt@C=12pt{
0\ar[r]	&H^0(Tor(S_{i_{{}_V}}))\ar[d] \ar[r]&H^1(\TET)\ar@{=}[d]\\
0\ar[r]	&H^0(S_{i_{{}_V}}) \ar[r]&H^1(\TET)
}
$$
where the exactness of the rows come from the vanishing
$$
H^0 (Q'_{i_{{}_V}})=H^0(Q''_{i_{{}_V}})=0
$$
and this because the sheaves $Q'_{i_{{}_V}}$ and $Q''_{i_{{}_V}}$ are
in the subcategory $\FF_{K_X,bK_X}$. The above diagram together with
the commutative square \eqref{V0-V-sq} imply the diagram asserted in the lemma.	
\end{pf}

Assume $V^0$, the subspace of $V$ corresponding to $H^0 (F_{L_{i_{{}_V}}})$, is nonzero. We take the extension sequence
\begin{equation}\label{V0ext}
\xymatrix@R=12pt@C=12pt{
	0\ar[r]& \TET \ar[r]& F_{V^0} \ar[r]& V^0 \otimes \OO_X \ar[r]&0
}
\end{equation}
corresponding to the inclusion $V^0 \subset H^1(\TET)$. The inclusions
$$
V^0 \subset V \subset H^1 (\TET) 
$$
give rise to the inclusion of extensions
$$
\xymatrix@R=12pt@C=12pt{
	0\ar[r]& \TET \ar[r]\ar@{=}[d]& F_{V^0} \ar[r] \ar[d]& V^0 \otimes \OO_X \ar[r]\ar@{^{(}->}[d]&0\\
	0\ar[r]& \TET \ar[r]& F_{V} \ar[r]& V \otimes \OO_X \ar[r]&0	
}
$$
 The vertical arrow in the middle is a monomorphism and can be completed to the exact sequence
$$
\xymatrix@R=12pt@C=12pt{
	0\ar[r]& F_{V^0}\ar[r]& F_{V} \ar[r] & \Big(V/V^0 \Big)\otimes \OO_X \ar[r]&0.}
$$
This implies that $F_{V^0}$ is a subbundle of $F_V$ and hence lies in $\FF_{K_X,bK_X}$.

We now make use of the exact sequence \eqref{Q''V}. It tells us
that $H^1(\TET)=Ext^1 (\OO_X,\TET)$ maps into 
$H^1 (Q''_{i_{{}_V}})=Ext^1 (\OO_X,Q''_{i_{{}_V}})$ and the cohomology classes in
$H^1(\TET)$ coming from $H^0 (Tor(S_{i_{{}_V}}))$ are precisely the extension
classes in $Ext^1 (\OO_X,\TET)$ which go to zero in \linebreak$Ext^1(\OO_X,Q''_V)$.
In particular, all extension classes in $V^0 \subset H^1 (\TET)$ are of this kind. This means that the monomorphism
$\TET \hookrightarrow Q''_V$ in \eqref{Q''V} lifts to a morphism 
$F_{V^0} \longrightarrow Q''_V$, that is, we have a commutative diagram
$$
\xymatrix@R=12pt@C=12pt{
	0\ar[r]& \TET \ar[r]\ar@{=}[d]& F_{V^0} \ar[r] \ar[d]& V^0 \otimes \OO_X \ar[r]&0\\
0\ar[r]&\TET \ar[r]	&Q''_{i_{{}_V}} &&
}
$$
This can be completed to the morphism of exact sequences
$$
\xymatrix@R=12pt@C=12pt{
	0\ar[r]& \TET \ar[r]\ar@{=}[d]& F_{V^0} \ar[r] \ar[d]& V^0 \otimes \OO_X \ar[r] \ar[d]&0\\
	0\ar[r]&\TET\ar[r]	&Q''_{i_{{}_V}}\ar[r] &Tor(S_{i_{{}_V}})\ar[r]&0
}
$$
The vertical arrow on the right is surjective, see Lemma \ref{lem:V0}, so is the arrow in the middle.
Thus we arrive to the diagram 
 \begin{equation}\label{V0ext-diag}
 	\xymatrix@R=12pt@C=12pt{
 		&&0\ar[d]&0\ar[d]&\\
 		&&F'_{V^0} \ar@{=}[r]\ar[d]&F'_{V^0}\ar[d]&\\
 		0\ar[r]& \TET \ar[r]\ar@{=}[d]& F_{V^0} \ar[r] \ar[d]& V^0 \otimes \OO_X \ar[r] \ar[d]&0\\
 		0\ar[r]&\TET\ar[r]	&Q''_{i_{{}_V}}\ar[r]\ar[d] &Tor(S_{i_{{}_V}})\ar[r]\ar[d]&0\\
 		&&0&0
 	}
 \end{equation}
where $F'_{V^0}:=ker(F_{V^0} \longrightarrow Q''_{i_{{}_V}})$.
This is the situation of rank $2$ case considered in the previous subsection. More precisely, the following holds.

\begin{lem}\label{lem:V0ext}
	In the diagram \eqref{V0ext-diag} the column in the middle is the 
	$K_X$-destabilizing sequence of $F_{V^0}$. Furthermore, the sheaf
	$Q''_{i_{{}_V}}$ is the last $K_X$-semistable factor in the HN filtration of $F_{V^0}$ with respect to $K_X$.
	
	 All
	the properties of Proposition \ref{pro:QVrk2} hold with the sheaf $Q''_{i_{{}_V}}$ and the divisor $D_{V^0}=c_1 (Tor(S_{i_{{}_V}}))$ playing respectively the roles of the sheaf $Q'_{i_{{}_V}}$ and the divisor $D_{i_{{}_V}}$ of that proposition.
	In addition, in  the extension 
	$$
	\xymatrix@R=12pt@C=12pt{
		0\ar[r]& \TET \ar[r]& F_{[\xi]} \ar[r]&  \OO_X \ar[r]& 0
	}
$$
corresponding to a line $[\xi]$ in $ V^0 \subset H^1 (\TET)$, the sheaf $ F_{[\xi]}$ in the middle is $K_X$-unstable with 
$\OO_X (-D_{V^0})$ being a $K_X$-destabilizing subsheaf of $ F_{[\xi]}$. For a general $[\xi] \in \PP(V^0)$, the sheaf $\OO_X (-D_{V^0})$ is the maximal $K_X$-destabilizing subsheaf of $ F_{[\xi]}$.	
\end{lem}
\begin{pf}
	Observe that the statements here are somewhat stronger than the ones in  Proposition \ref{pro:QVrk2}: we only need $dim(V^0) \neq 0$ and no further assumption on the ratio of the Chern numbers of $X$ is needed
	to see that $\OO_X (-D_{V^0})$ is a $K_X$-destabilizing (resp. the maximal $K_X$-destabilizing)  subsheaf of $ F_{[\xi]}$ for $[\xi]\in \PP(V^0)$ (resp. general $[\xi]\in \PP(V^0)$). The reason is that the divisor $D_{V^0} =c_1 (Tor(S_{i_{{}_V}}))$ is  a component of the divisor $c_1 (S_{i_{{}_V}})$: recall, we have
	the exact sequence
	$$
\xymatrix@R=12pt@C=12pt{
	0\ar[r]&Tor(S_{i_{{}_V}}) \ar[r]&S_{i_{{}_V}}\ar[r]& {\cal J}_{A_{i_{{}_V}}} (L_{i_{{}_V}})\ar[r]&0
}	
	$$
	which implies
	$$
	c_1 (S_{i_{{}_V}})=D_{V^0}+ L_{i_{{}_V}}
	$$
	where two terms of the sum are effective divisors. Furthermore, the upper bound
	$$
	c_1 (S_{i_{{}_V}})\cdot K_X \leq 2(3c_2(X)-K^2_X)
	$$
	in Lemma \ref{lem:F'Vrk3} implies 
	\begin{equation}\label{DV0-deg}
		D_{V^0}\cdot K_X \leq 2(3c_2(X)-K^2_X).
	\end{equation}
Using this inequality and arguing as in the proof of Lemma \ref{lem:Q'iV-stab}
give a proof of the assertion about 	$Q''_{i_{{}_V}}$ being the last $K_X$-semistable factor in the HN filtration of $F_{V^0}$ with respect to $K_X$. 
Once we know this, Proposition \ref{pro:QVrk2} applies to $Q''_{i_{{}_V}}$ and $D_{V^0}$. In particular, we obtain that $D_{V^0}$ annihilates the subspace $V^0$ of $H^1 (\TET)$ and it is minimal with respect to this property, see Lemma \ref{lem:DV}, 3). This means that
the extension 
	$$
\xymatrix@R=12pt@C=12pt{
	0\ar[r]& \TET \ar[r]& F_{[\xi]} \ar[r]&  \OO_X \ar[r]& 0
}
$$
corresponding to a line $[\xi]$ in $V^0$ admits
a monomorphism
$$
\xymatrix@R=12pt@C=12pt{
	&&0\ar[d]&0\ar[d]&\\
	&&\OO_X (-D_{V^0})\ar@{=}[r]\ar[d]&\OO_X (-D_{V^0})\ar[d]&\\
	0\ar[r]& \TET \ar[r]& F_{[\xi]} \ar[r]&  \OO_X \ar[r]& 0
}
$$
 and this shows that $\OO_X (-D_{V^0})$ is a subsheaf of $F_{[\xi]}$.
 To see that it is destabilizing we need to check the slope inequality
 $$
 -D_{V^0}\cdot K_X > -\frac{1}{3}c_1 (F_{[\xi]})\cdot K_X =-\frac{1}{3}K^2_X
 $$ 
 or, equivalently, the inequality
 $$
 D_{V^0}\cdot K_X <\frac{1}{3}K^2_X
 $$
 and this is ensured by the upper bound \eqref{DV0-deg} since
 it gives
 $$
 	D_{V^0}\cdot K_X \leq 2(3c_2(X)-K^2_X)<\frac{1}{4}K^2_X
 	$$
 	where the second inequality comes from the assumption 
 	on the ratio $\alpha_X <\frac{3}{8}$.
 	
 	We complete the monomorphism $\OO_X (-D_{V^0})\longrightarrow F_{[\xi]} $ to the exact sequence
 	$$
 \xymatrix@R=12pt@C=12pt{
 	0\ar[r]&\OO_X (-D_{V^0})\ar[r]& F_{[\xi]}\ar[r]& Q_{[\xi]}\ar[r]&0.
 }
 	$$
 	For a general $[\xi] \in \PP(V^0)$ the quotient sheaf $Q_{[\xi]}$ is torsion free, see Proposition \ref{pro:QVrk2}, 3). In addition, it differs from $Q''_{i_{{}_V}}$ on at most $0$-dimensional subscheme of $X$:
 	$$
 	\xymatrix@R=12pt@C=12pt{
 		0\ar[r]& Q_{[\xi]}\ar[r]&Q''_{i_{{}_V}}\ar[r]& S_[\xi]\ar[r]&0
 	},
 $$ 
see the middle column of the second diagram in Proposition \ref{pro:QVrk2}, 3): remember that $Q'_{i_{{}_V}}$ and $D_{i_{{}_V}}$ there must be replaced by $Q''_{i_{{}_V}}$ and $D_{V^0}$
respectively. Since  $Q''_{i_{{}_V}}$
 	is $K_X$-semistable, the exact sequence above tells us that $Q_{[\xi]}$ is $K_X$-semistable as well. This in turn implies that
 	$$
 	F_{[\xi]} \supset \OO_X (-D_{V^0}) \supset 0
 	$$
 	is the HN filtration of $F_{[\xi]}$ with respect to $K_X$. Hence
 	$\OO_X (-D_{V^0})$ is the maximal $K_X$-destabilizing subsheaf of $F_{[\xi]}$.	
\end{pf}

At this stage we have the following situation: for any subspace $V \subset H^1(\TET)$ of dimension $dim(V)\geq 2$ for which the sheaf
$Q^{l_{i_{{}_V}}}_{i_{{}_V}}$ is of rank $3$, there is a one step
filtration
$$
V^0 \subset V,
$$
where the subspace $V^0$, if nonzero, gives rise to the exact
sequence
\begin{equation}\label{Q''iVTh-seq}
\xymatrix@R=12pt@C=12pt{
	0\ar[r]& \TET \ar[r]& Q''_{i_{{}_V}} \ar[r]&  Tor(S_{i_{{}_V}}) \ar[r]& 0
}
\end{equation}
as defined in \eqref{Q''V} and this sequence is a part of the diagram \eqref{V0ext-diag} which has all the features of the rank $2$ case and thus is subject to Proposition \ref{pro:QVrk2}.

We now turn to the quotient $V/V^{0}$. It has two interpretations: one is cohomological and the other is geometric.

The cohomological interpretation comes from applying the homological functor $Hom(\OO_X,\bullet)$ to the sequence \eqref{Q''iVTh-seq} to obtain the exact sequence 
$$
\xymatrix@R=12pt@C=12pt{
	0\ar[r]& H^0 (Tor(S_{i_{{}_V}}))\ar[r]& H^1(\TET)\ar[r]& Ext^1 (\OO_X,Q''_{i_{{}_V}} ),
}
$$
where the identification $Ext^1 (\OO_X, \TET)=H^1 (\TET)$ has been used. Putting this together with the commutative square in Lemma \ref{lem:V0}
gives the inclusion
\begin{equation}\label{V-V0-Ext1}
	V/V^0 \hookrightarrow  Ext^1(\OO_X,Q''_V).
\end{equation}

The geometric interpretation comes from identifying $V/V^{0}$ with a subspace of $H^0 ({\cal J}_{A_{i_{{}_V}}}(L_{i_{{}_V}}))$:
$$
V/V^{0} \hookrightarrow  H^0 ({\cal J}_{A_{i_{{}_V}}}(L_{i_{{}_V}})),
$$
 see the diagram \eqref{V0V-diag}.
To reconcile the two inclusions above
we use the diagram 
\begin{equation}\label{Q''iV-Q'iv-diag}
\xymatrix@R=12pt@C=12pt{
	&&0\ar[d]&0\ar[d]&\\
	0\ar[r]&\TET \ar[r]\ar@{=}[d]&Q''_{i_{{}_V}} \ar[r]\ar[d]&Tor(S_{i_{{}_V}})\ar[d] \ar[r]&0\\
	0\ar[r]&\TET \ar[r]&Q'_{i_{{}_V}} \ar[r]\ar[d]&S_{i_{{}_V}} \ar[r]\ar[d]&0\\
	&&{\cal J}_{A_{i_{{}_V}}} (L_{i_{{}_V}})\ar@{=}[r]\ar[d]&{\cal J}_{A_{i_{{}_V}}} (L_{i_{{}_V}})\ar[d]&\\
	&&0&0		
}
\end{equation}
that we encountered in the proof of Lemma \ref{lem:V0}.
Namely, we apply the homological functor  $Hom(\OO_X,\bullet)$ to the middle column of the above diagram to deduce the exact sequence
\begin{equation}\label{ExtQ''-ExtQ'}
	\xymatrix@R=12pt@C=12pt{
		0\ar[r]&  H^0({\cal J}_{A_{i_{{}_V}}}(L_{i_{{}_V}}))\ar[r]&Ext^1(\OO_X,Q''_{i_{{}_V}})  \ar[r]&Ext^1(\OO_X,Q'_{i_{{}_V}}).  
	}
\end{equation}
Combining this with the inclusion in \eqref{V-V0-Ext1} we obtain
\begin{equation}\label{V-V0-Extseq}
\xymatrix@R=12pt@C=12pt{
	&&0\ar[d]&\\
	&&V/V^{0} \ar[d]&\\
	0\ar[r]&  H^0({\cal J}_{A_{i_{{}_V}}}(L_{i_{{}_V}}))\ar[r]&Ext^1(\OO_X,Q''_{i_{{}_V}})  \ar[r]&Ext^1(\OO_X,Q'_{i_{{}_V}})
}
\end{equation}
\begin{lem}\label{lem:V-V0-factors}
	The inclusion 
	$$
	V/V^0 \hookrightarrow  Ext^1(\OO_X,Q''_{i_{{}_V}})
	$$
	in \eqref{V-V0-Extseq} factors through  $H^0({\cal J}_{A_{i_{{}_V}}}(L_{i_{{}_V}}))$ 
	to form the commutative diagram
	\begin{equation}\label{V-V0-Ext-diag-lem}
\xymatrix@R=12pt@C=12pt{
	&&0\ar[d]&\\
	&&V/V^{0} \ar[d]\ar[ld]&\\
	0\ar[r]&  H^0({\cal J}_{A_{i_{{}_V}}}(L_{i_{{}_V}}))\ar[r]&Ext^1(\OO_X,Q''_{i_{{}_V}})  \ar[r]&Ext^1(\OO_X,Q'_{i_{{}_V}})
}
\end{equation}	
\end{lem}
\begin{pf}
	We apply the homological functor  $Hom(\OO_X,\bullet)$ to the diagram \eqref{Q''iV-Q'iv-diag} to obtain

	\begin{equation}\label{V-V0-twofactor}
		\xymatrix@R=12pt@C=12pt{
			&0\ar[d]&0\ar[d]&0\ar[d]&\\
			0\ar[r]&H^0(Tor(S_{i_{{}_V}})) \ar[r]\ar[d]& H^0(S_{i_{{}_V}})\ar[r]\ar[d]& H^0({\cal J}_{A_{i_{{}_V}}}(L_{i_{{}_V}}) )\ar[d]&\\
			&H^1 (\TET)\ar@{=}[r]&H^1 (\TET)\ar[r]&Ext^1(\OO_X,Q''_{i_{{}_V}})&&
		}
	\end{equation}
This and the identifications of $V$ and $V^0$ with subspaces of $H^0(S_{i_{{}_V}})$ and $H^0(Tor(S_{i_{{}_V}}))$ respectively, see the commutative square \eqref{V0-V-sq}, imply the commutative diagram
$$
\xymatrix@R=12pt@C=12pt{
	V/V^0 \ar[r] \ar@{=}[d]&H^0({\cal J}_{A_{i_{{}_V}}}(L_{i_{{}_V}}) )\ar[d]\\
	V/V^0 \ar[r]& Ext^1(\OO_X,Q''_{i_{{}_V}}) 
}
$$
where all arrows are monomorphisms.
\end{pf}

The diagram \eqref{V-V0-Ext-diag-lem} tells us:
 
$\bullet$ the quotient space $V/V^0$ parametrizes extension classes in
$Ext^1(\OO_X,Q''_{i_{{}_V}}) $ - this is the vertical arrow in the diagram,

$\bullet$   the extension classes in the previous item go to the trivial (split) extensions in
$Ext^1(\OO_X,Q'_{i_{{}_V}})$ - this is the slanted arrow of the diagram.

\vspace{0.2cm}
  The first item above gives the extension sequence
 $$
 \xymatrix@R=12pt@C=12pt{
	0\ar[r]&Q''_{i_{{}_V}} \ar[r]&G_{V/V^0} \ar[r]& \Big(V/V^0\Big)\otimes \OO_X \ar[r]&0
}
$$ 
corresponding to the inclusion $V/V^0 \hookrightarrow Ext^1 (\OO_X,Q''_{i_{{}_V}})$,
while the second item says that the diagram
$$
\xymatrix@R=12pt@C=12pt{
	&0\ar[d]&&\\
	0\ar[r]&Q''_{i_{{}_V}} \ar[r]\ar[d]&G_{V/V^0} \ar[r]& \Big(V/V^0\Big)\otimes \OO_X \ar[r]\ar@{=}[d]&0\\
0\ar[r]	&Q'_{i_{{}_V}} \ar[r]&Q'_{i_{{}_V}}\oplus\Big(V/V^0\Big)\otimes \OO_X \ar[r]&\Big(V/V^0\Big)\otimes \OO_X \ar[r] &0 
}
$$
where the bottom sequence is a split exact sequence, can be completed by a vertical arrow in the middle to a morphism of
extensions
\begin{equation}\label{Gext-to-split}
\xymatrix@R=12pt@C=12pt{
	&0\ar[d]&&\\
	0\ar[r]&Q''_{i_{{}_V}} \ar[r]\ar[d]&G_{V/V^0} \ar[r]\ar[d]& \Big(V/V^0\Big)\otimes \OO_X \ar[r]\ar@{=}[d]&0\\
	0\ar[r]	&Q'_{i_{{}_V}} \ar[r]&Q'_{i_{{}_V}}\oplus\Big(V/V^0\Big)\otimes \OO_X \ar[r]&\Big(V/V^0\Big)\otimes \OO_X \ar[r] &0 
}
\end{equation}
We also have two epimorphisms
$$
\xymatrix@R=12pt@C=12pt{
&\Big(V/V^0\Big)\otimes \OO_X \ar[d]\\	
Q'_{i_{{}_V}} \ar[r]& {\cal J}_{A_{i_{{}_V}}} (L_{i_{{}_V}})
}
$$
the horizontal one from the middle column in \eqref{Q''iV-Q'iv-diag}
and the vertical one from the diagram in \eqref{V0V-diag}. This can be completed to the commutative square 
$$
\xymatrix@R=12pt@C=12pt{
Q'_{i_{{}_V}} \oplus \Big(V/V^0\Big)\otimes \OO_X \ar[r]\ar[d]	&\Big(V/V^0\Big)\otimes \OO_X \ar[d]\\	
	Q'_{i_{{}_V}} \ar[r]& {\cal J}_{A_{i_{{}_V}}} (L_{i_{{}_V}})
}
$$
From this and the morphism of extension sequences in \eqref{Gext-to-split} we deduce commutative diagram
$$
\xymatrix@R=12pt@C=12pt{
	0\ar[r]&Q''_{i_{{}_V}} \ar[r]\ar@{=}[d]&G_{V/V^0} \ar[r]\ar[d]& \Big(V/V^0\Big)\otimes \OO_X \ar[r]\ar[d]&0\\
		0\ar[r]&Q''_{i_{{}_V}} \ar[r]&Q'_{i_{{}_V}}\ar[r]\ar[d]&{\cal J}_{A_{i_{{}_V}}} (L_{i_{{}_V}})\ar[r]\ar[d]&0\\
		&&0&0& 
}
$$

This diagram in turn can be expanded to
\begin{equation}\label{Q''Vext-diag}
	\xymatrix@R=12pt@C=12pt{
		&&0\ar[d]&0\ar[d]&\\
		&&G'_{V/V^0} \ar@{=}[r]\ar[d]&G'_{V/V^0} \ar[d]&\\
		0\ar[r]&Q''_{i_{{}_V}} \ar[r]\ar@{=}[d]&G_{V/V^0} \ar[r]\ar[d]& \Big(V/V^0\Big)\otimes \OO_X \ar[r] \ar[d]&0\\
		0\ar[r]	&Q''_{i_{{}_V}} \ar[r]&Q'_{i_{{}_V}} \ar[r]\ar[d]&{\cal J}_{A_{i_{{}_V}}}(L_{i_{{}_V}})\ar[r]\ar[d]&0\\
		&&0&0& 
	}	
\end{equation}
This diagram is similar to the one in \eqref{rkQV=3diag}, except now it involves the extension sequence for the torsion free sheaf $Q''_{i_{{}_V}}$ which differs from $\TET$ by the torsion sheaf 
$Tor(S_{i_{{}_V}})$, see the exact sequence \eqref{Q''V}. 

The point
of the above discussion is 
that there is some kind of an inductive process:

we started with the pair $(\TET, V)$, where $V$ is a subspace of $H^1(\TET)$, and produced the subspace $V^0$ of $V$ together with `bigger' sheaf $Q''_{i_{{}_V}}$ (it contains $\TET$); that sheaf controls the instability of the extension
$$
\xymatrix@R=12pt@C=12pt{
	0\ar[r]	&\TET	\ar[r]	&F_{V^0} \ar[r]	&V^0 \otimes \OO_X \ar[r]&0
}
$$
corresponding to the subspace $V^0 \subset H^1(\TET)$ as well as gives the new pair $(Q''_{i_{{}_V}}, V/V^{0})$, where 
$V/V^0 $ is a subspace of $H^1(Q''_{i_{{}_V}})$; so we can repeat our considerations with this pair. 

\vspace{0.2cm}
We summarize all the properties of the rank $3$ case in the following.

\begin{pro}\label{pro:QVrk3}
	Assume $\displaystyle{\alpha_X < \frac{3}{8} }$ and the rank of $Q^{l_{i_{{}_V}}}_{i_{{}_V}}$ in Proposition \ref{pro:QV-23} is $3$. Then either 
	$dim(V)=1$ or $dim(V) \geq 2$ and then the extension
	$$
	\xymatrix@R=12pt@C=12pt{
		0\ar[r]& \TET \ar[r]& F_V \ar[r]& V\otimes \OO_X \ar[r]& 0
	}
	$$
	corresponding to the natural inclusion $V\subset H^1 (\TET)$ fits into the diagram \eqref{rkQV=3diag}, see the middle row of that diagram. Furthermore, the following properties hold.
	
	1) The canonical divisor $K_X$ admits the decomposition
	$$
	K_X =c_1 (S_{i_{{}_V}}) + (-c_1 (Q'_{i_{{}_V}}))=L_{i_{{}_V}} +D_{V^0} +(-c_1 (Q'_{i_{{}_V}})),
	$$
	where $L_{i_{{}_V}}$ is as in Lemma \ref{lem:LV} and $D_{V^0}=c_1(Tor(S_{i_{{}_V}}))$. In particular, 
	$(L_{i_{{}_V}} +D_{V^0} )$ is an effective, nonzero divisor and it is 
	 subject to
	$$
	(L_{i_{{}_V}} +D_{V^0} ) \cdot K_X \leq 2(3c_2(X) -K^2_X).
	$$
	In addition, the divisor
	$(-c_1 (Q'_{i_{{}_V}}))$ lies in the positive cone of $X$.
	
	2) The sheaf $F_V$ is $K_X$-unstable and the column in the middle of the diagram \eqref{rkQV=3diag} is a $K_X$-destabilizing sequence of $F_V$. The quotient sheaf $Q'_{i_{{}_V}}$ is either $K_X$-semistable
	and then $Q'_{i_{{}_V}}=Q^m_V$ is the last $K_X$-semistable factor of the HN filtration of $F_V$ with respect to $K_X$, or 
	$Q'_{i_{{}_V}}$ is $K_X$-unstable and admits the destabilizing sequence
	$$
	\xymatrix@R=12pt@C=12pt{
		0\ar[r]& {\cal J}_A(-M)\ar[r]&Q'_{i_{{}_V}}\ar[r]&Q^m_V \ar[r]&0,
	}
$$
where $Q^m_V$ is the last $K_X$-semistable factor of the HN filtration of $F_V$ with respect to $K_X$.

	3) If the torsion part $Tor(S_{i_{{}_V}})$ of $S_{i_{{}_V}}$ is nonzero, then the inclusion $F'_{i_{{}_V}} \hookrightarrow V\otimes \OO_X$ in the diagram \eqref{rkQV=3diag} admits the locally free intermediate subsheaf $F_{L_{i_{{}_V}}}$, that is, one has inclusions
	$$
	F'_{i_{{}_V}} \hookrightarrow F_{L_{i_{{}_V}}} \hookrightarrow V\otimes \OO_X
	$$
	 fitting into the diagram
		$$
	\xymatrix@R=12pt@C=12pt{
		&&0\ar[d]&0\ar[d]&\\
		0\ar[r]&F'_{{i_{{}_V}}}\ar[r]\ar@{=}[d]& F_{L_{i_{{}_V}}} \ar[r] \ar[d]& Tor(S_{i_{{}_V}})\ar[r]\ar[d]&0\\
		0\ar[r]&F'_{{i_{{}_V}}}\ar[r]& V\otimes \OO_X \ar[r]\ar[d]& S_{i_{{}_V}}\ar[r]\ar[d]&0\\
		&&{\cal J}_{A_{i_{{}_V}}}(L_{i_{{}_V}})\ar@{=}[r]\ar[d]&{\cal J}_{A_{i_{{}_V}}}(L_{i_{{}_V}})\ar[d]&\\
		&&0&0&	
	}
$$
 
	4) The space $ H^0 (F_{L_{i_{{}_V}}})$ of global sections of $F_{L_{i_{{}_V}}}$ defines a proper subspace of $V$ denoted $V^0$. If nonzero, it gives rise to the extension
	$$
	\xymatrix@R=12pt@C=12pt{
		0\ar[r]& \TET \ar[r]& F_{V^0} \ar[r] & V^0 \otimes \OO_X \ar[r]&0
	}
	$$
	corresponding to the inclusion $V^0 \subset H^1(\TET)$. This extension is subject to the properties of Lemma \ref{lem:V0ext}:  the vector bundle $F_{V^0}$, the middle term in the above extension, is $K_X$-unstable whose last $K_X$-semistable factor
	in its HN filtration is the rank $2$ sheaf $Q''_{i_{{}_V}}$ fitting into the exact sequence
	$$
	\xymatrix@R=12pt@C=12pt{
		0\ar[r]& \TET \ar[r]&Q''_{i_{{}_V}} \ar[r]& Tor(S_{i_{{}_V}})\ar[r]&0;
	}
$$

\noindent
the divisor $D_{V^0}=c_1(Tor(S_{i_{{}_V}}))$ in the lemma is the one appearing in the decomposition of $K_X$ in the part 1) of the proposition.
	
	5) The sheaf ${\cal J}_{A_{i_{{}_V}}}(L_{i_{{}_V}})$, the torsion free part of $S_{i_{{}_V}}$, is generated by its global sections. More precisely, one has an epimorphism
	$$
		\xymatrix@R=12pt@C=12pt{
			\Big(V/V^0 \Big)\otimes \OO_X \ar[r]& {\cal J}_{A_{i_{{}_V}}}(L_{i_{{}_V}})\ar[r]&0
		}
	$$
which induces the inclusion
$$
	V/V^0 \subset H^0 ( {\cal J}_{A_{i_{{}_V}}}(L_{i_{{}_V}})).
$$
In particular, $dim\Big(V/V^0 \Big) \geq 2$, unless ${\cal J}_{A_{i_{{}_V}}}(L_{i_{{}_V}}) =\OO_X$. 
	
	6)	The quotient space $V/V^0$ is identified with a subspace of
	$Ext^1 (\OO_X, Q''_{i_{{}_V}} )$ and the pair $(Q''_{i_{{}_V}}, V/V^0)$ produces the diagram \eqref{Q''Vext-diag} which could be viewed as inductive `descendant' of the diagram \eqref{rkQV=3diag}.
\end{pro}

\section{$H^1(\TET)$ for surfaces with $\displaystyle{\alpha_X <\frac{3}{8}}$}
We spell out the results of the previous sections in the case $V=H^1(\TET)$.
\begin{thm}\label{th:H1Theta}
	Let $X$ be a smooth compact complex surface with the canonical divisor $K_X$ ample and the ratio $\alpha_X$ of the Chern numbers subject to $\displaystyle{\alpha_X <\frac{3}{8}}$. Assume
	$h^1 (\TET)\geq 2$. Then the extension
	sequence
	$$
\xymatrix@R=12pt@C=12pt{
	0\ar[r]& \TET \ar[r]& F_X \ar[r]& H^1(\TET)\otimes \OO_X \ar[r]& 0
}
$$
corresponding to  $id_{H^1(\TET)} \in End( H^1 (\TET))\cong Ext^1 (H^1 (\TET)\otimes \OO_X, \TET)$ fits into the diagram 
\begin{equation}\label{X-diag}
\xymatrix@R=12pt@C=12pt{
	&&0\ar[d]&0\ar[d]&\\
	&&F'_X \ar@{=}[r]\ar[d]&F'_X \ar[d]&\\
	0\ar[r]& \TET \ar[r]\ar@{=}[d]& F_X \ar[r]\ar[d]& H^1(\TET)\otimes \OO_X \ar[r]\ar[d]& 0\\
	0\ar[r]&\TET \ar[r]&Q'_X \ar[r]\ar[d]& S_X \ar[r]\ar[d]&0\\
	&&0&0&
}
\end{equation}
where the column in the middle is the $K_X$-destabilizing sequence of $F_X$ and the quotient sheaf $Q'_X$ is torsion free, $K_X$-stable and of rank $2$ or $3$. Furthermore, the following properties hold.

1) $H^1(\TET) \cong H^0(S_X)$ and the sheaf $S_X$ is globally generated.

2) The canonical divisor $K_X$ admits the decomposition
$$
K_X =c_1 (S_X) + (-c_1 (Q'_X)),
$$
where $(-c_1 (Q'_X))$ lies in the positive cone of $X$ and $c_1(S_X)$ is an effective, nonzero divisor subject to
$$
(c_1 (S_X)) \cdot K_X \leq (4-2rk(S_X))(3c_2(X) -K^2_X).
$$

3) Let $Tor(S_X)$ and $S^{tf}_X$ be respectively the torsion and torsion free parts of $S_X$; in other words $S_X$ fits into the exact sequence
$$
\xymatrix@R=12pt@C=12pt{
0\ar[r]&Tor(S_X)\ar[r]&S_X \ar[r]&S^{tf}_X \ar[r]&0.
}
$$
The subspace $H^0 (Tor(S_X)) \subset H^0 (S_X)$ is identified, via the isomorphism in 1), with the subspace denoted $H^1(\TET)^0 $ of $H^1 (\TET)$ and this gives a one step
filtration of $H^1 (\TET)$
\begin{equation}\label{1stepfilt}
H^1 (\TET) \supset H^1 (\TET)^0 \cong H^0 (Tor(S_X)).
\end{equation}
The equality $H^1 (\TET) = H^1 (\TET)^0$ holds if and only if the rank of
$Q'_X$ in \eqref{X-diag} is $2$. In this case  $S_X=Tor(S_X)$ is a torsion sheaf. The divisor $D_X=c_1 (S_X)$ is subject to
$$
D_X \cdot K_X < 4(3c_2(X) -K^2_X).
$$
In addition, $D_X$ annihilates $H^1 (\TET)$, that is, for a global section $\delta_X$ of $\OO_X (D_X)$ defining $D_X$, the homomorphism
$$
H^1 (\TET)\stackrel{\delta_X}{\longrightarrow}H^1(\Theta_X (D_X))
$$
defined by multiplication with $\delta_X$ is identically zero and one has
$$
H^1 (\TET)\cong H^0 (\TET\otimes \OO_{D_X}(D_X)).
$$

3) If the subspace $H^1(\TET)^0$ in \eqref{1stepfilt} is not equal to  $H^1 (\TET)$,  that is, we are in the case $rk(Q'_X)=3$, then the torsion free part $S^{tf}_X$ of $S_X$ is nonzero and  has the form
$$
S^{tf}_X={\cal J}_{A_X}(L_X),
$$
where ${\cal J}_{A_X}$ is the ideal sheaf of at most $0$-dimensional subscheme $A_X$ of $X$ and $L_X$ is a  divisor.  The bottom row in \eqref{X-diag} fits into the following commutative diagram
$$
\xymatrix@R=12pt@C=12pt{
	&&0\ar[d]&0\ar[d]&\\ 
	0\ar[r]	&\TET\ar[r]\ar@{=}[d]&Q''_X \ar[r]\ar[d]&Tor(S_X)\ar[r]\ar[d]&0\\
	0\ar[r]&\TET\ar[r]&Q'_X \ar[r]\ar[d]&S_X \ar[r]\ar[d]&0\\
	&&{\cal J}_{A_X}(L_X)\ar@{=}[r]\ar[d]&{\cal J}_{A_X}(L_X)\ar[d]&\\
	&&0&0&
}
$$ 
In particular, the divisor
$c_1 (S_X)$ in 2) has the form
$$
c_1(S_X)=c_1 (Tor(S_X)) +L_X.
$$
and it is subject to
$$
c_1(S_X)= (c_1(Tor(S_X))+L_X)\cdot K_X \leq 2(3c_2(X)-K^2_X).
$$  

Furthermore, the quotient space
$H^1 (\TET)/H^1 (\TET)^0$ injects into the space
$H^0 ({\cal J}_{A_X}(L_X))$ of global sections of ${\cal J}_{A_X}(L_X)$ and this sheaf is globally generated by that subspace. In particular,
$L_X =0$ if and only if $H^1(\TET)^0$ is a codimension $1$ subspace of
$H^1 (\TET)$ and in this case ${\cal J}_{A_X}(L_X)=\OO_X$. Otherwise, $L_X$ is a nonzero nef divisor with $L^2_X \geq deg(A_X)$.

4) 
If $H^1(\TET)^0 \cong H^0 (Tor(S_X)) \neq 0$ and is a proper subspace of $H^1 (\TET)$, then the extension 
$$
\xymatrix@R=12pt@C=12pt{
0\ar[r]	&\TET\ar[r]&F^0_X \ar[r]&H^1(\TET)^0 \otimes \OO_X \ar[r]&0
}
$$
corresponding to the inclusion $H^1(\TET)^0 \subset H^1 (\TET)$ fits into the diagram similar to the one in \eqref{X-diag}
\begin{equation}\label{X-diag1}
	\xymatrix@R=12pt@C=12pt{
		&&0\ar[d]&0\ar[d]&\\
		&&F^{0'}_X \ar@{=}[r]\ar[d]&F^{0'}_X \ar[d]&\\
		0\ar[r]& \TET \ar[r]\ar@{=}[d]& F^0_X \ar[r]\ar[d]& H^1(\TET)^0 \otimes \OO_X \ar[r]\ar[d]& 0\\
		0\ar[r]&\TET \ar[r]&Q''_X \ar[r]\ar[d]& Tor(S_X) \ar[r]\ar[d]&0\\
		&&0&0
	}
\end{equation}
The sheaf $Q''_X $ is the last $K_X$-semistable factor of the HN filtration of $F^0_X$ with respect to $K_X$ and the divisor
 $D^0_X:=c_1 (Tor(S_X))$ has the property of annihilating the subspace $H^1 (\TET)^0$, see part 2) of the theorem for the meaning of this. 

5) Continuing to assume $Tor(S_X)\neq 0$ and $H^1(\TET)^0 \neq H^1 (\TET)$, the quotient space $W_X:= H^1 (\TET)/ H^1 (\TET)^0$ injects into
$H^1(Q''_X)$ and the pair $(Q''_X, W_X)$ is a `descendant' of the pair 
$(\TET, H^1(\TET))$, this is to say that the pair $(Q''_X, W_X)$ in its turn defines the extension sequence
$$
\xymatrix@R=12pt@C=12pt{
	0\ar[r]&Q''_X \ar[r]&G_{W_X} \ar[r]& W_X \otimes \OO_X \ar[r]&0
}
$$
corresponding to an inclusion 
$W_X \subset Ext^1 (\OO_X,Q''_X)$. This fits into the diagram
\begin{equation}\label{X-diag2}
\xymatrix@R=12pt@C=12pt{
	&&0\ar[d]&0\ar[d]&\\
	&&G'_{W_X}\ar@{=}[r]\ar[d]&G'_{W_X}\ar[d]&\\
	0\ar[r]&Q''_X \ar[r]\ar@{=}[d]&G_{W_X} \ar[r]\ar[d]& W_X \otimes \OO_X \ar[r]\ar[d]&0\\
0\ar[r]&Q''_X \ar[r]&Q'_X \ar[r]\ar[d]& {\cal J}_{A_X}(L_X)\ar[d]\ar[r]&0\\
&&0&0&	
}
\end{equation}
analogous to \eqref{X-diag}.
\end{thm}

We now suggest the following geometric interpretation of the homological algebra considerations above.

\vspace{0.2cm}
a) The exact sequence
$$
\xymatrix@R=12pt@C=12pt{
	0\ar[r]& \TET \ar[r]& Q''_X \ar[r]& Tor(S_X)\ar[r]&0
}
$$
at the bottom of the sequence \eqref{X-diag1} suggests the infinitesimal version of
$X$ being a ramified covering of some other surface, call it $X''$, where $Q''_V$
looks like the pull back of the tangent sheaf of $X''$ and the above sequence plays the role of the differential of that ramified covering;
the formula
$$
K_X=-c_1 (Q''_X)+ c_1 (Tor(S_X))
$$
is a homological counterpart of the formula
$$
K_X =f^{\ast}(K_{X''})  +R_f
$$
of the canonical divisor
of $X$ under the ramified covering $f:X\longrightarrow X''$; 

the property saying that $(-c_1 (Q''_X))$ is in the positive cone of $X$,
see Theorem \ref{th:H1Theta}, 2), now can be interpreted as  $X''$ being of general type; furthermore, the relation of the Chern numbers of $X$
and those of $Q''_X$ given by the inequality
$$
3c_2 (X)-K^2_X\geq (3c_2 (Q''_X)-c^2_1(Q''_X)) +K_X \cdot D^0_X -2(D^0_X)^2,
$$
see Lemma \ref{lem:Q'Vstable}, 3), says that $X$ is a ramified covering of $X''$, a surface of general type whose Chern numbers give the point $(c_2(X''), K^2_{X''})=(c_2(Q''_X), c^2_1 (Q''_X))$
in the region depicted in \eqref{Pworld} which lies  closer to the BMY line than the point $(c_2(X), K^2_{X})$;

 the divisor $D^0_X=c_1 (Tor(S_X))$ could be viewed as a homological version of the ramification divisor of $f$.

  The subspace $H^1(\TET)^0 \cong H^0 (Tor(S_X))$ can now be interpreted as the part of the infinitesimal deformations of $X$ controlled by the infinitesimal deformations of the ramification divisor of the covering
$X\stackrel{f}{\longrightarrow} X''$; 
from this point of view the equality $H^1(\TET)=H^1(\TET)^0$
says that all infinitesimal deformations of $X$ 
come from the infinitesimal deformations of the ramification divisor of
the covering.

\vspace{0.2cm}
b)	The bottom row of the diagram \eqref{X-diag2} looks like
the infinitesimal version of a morphism of $X''$ in a) into a variety $X'$  of dimension $3$, where the sheaf $Q'_X$ looks like the pull back of the tangent sheaf of $X'$; the exact sequence itself plays the role of the differential of that morphism - a sort of generalized normal sequence of $X''$ in $X'$;

 the quotient space $H^1(\TET)/H^1(\TET)^0 $ together with
the inclusion 
$$
H^1(\TET)/H^1(\TET)^0  \subset H^0({\cal J}_{A_X}(L_X))
$$
 looks like
the infinitesimal deformations of the divisor $X''$ inside $X'$.

To summarize, the homological considerations suggest the following.

\vspace{0.2cm}
{\bf Meta-principle.} Any smooth compact complex surface $X$ with the canonical
class $K_X$ ample and the Chern numbers $(c_2 (X),K^2_X) $ with the ratio
$\alpha_X =\displaystyle{\frac{c_2 (X)}{K^2_X} < \frac{3}{8}}$ and the number of moduli 
$h^1 (\TET)\geq 2$ are obtained

either as a ramified covering $f:X \longrightarrow X''$ of a surface
of general type $X''$ with Chern numbers $(c_2(X''), K^2_{X''})$ closer
to BMY line than the Chern numbers of $X$ and all infinitesimal deformations of $X$ come from the infinitesimal deformations of the ramification divisor of $f$,

or $X$ admits a generically finite morphism $g:X \longrightarrow X'$ into a normal variety $X'$ of general type of dimension $3$; in this case all infinitesimal deformations of $X$ either come from the infinitesimal deformation of the divisor $g(X)$ in $X'$ or there is a proper nonzero subspace $H^1(\TET)^0$ of $H^1(\TET)$ coming from the infinitesimal deformations of the ramification divisor of $g:X \longrightarrow g(X)$, while the quotient
space
$H^1(\TET)/H^1(\TET)^0$ is a subspace of the infinitesimal deformations of
$g(X)$ in $X'$; in this latter case the smooth minimal model of $g(X)$ is a surface of general type with the Chern numbers closer to BMY line than the Chern numbers of $X$.

\section{Surfaces with the small $(3c_2 (X)-K^2_X)$}
Until now our results relied on the assumption on the ratio $\alpha_X$ of the Chern numbers of $X$. However, the geometric output of our considerations, the divisors $D_X$ and $c_1 (S_X) $ in Theorem \ref{th:H1Theta}, have their degrees bounded from above in terms
of $(3c_2 (X)-K^2_X)$, the quantity which tells us how far the point
$(c_2(X),K^2_X) $ in the region depicted in \eqref{Pworld} is from the
BMY line $K^2=3c_2$. This suggests to study surfaces of general type according to the values of $(3c_2 (X)-K^2_X)$. In this section we give several results of this nature.

In the previous section we found the decomposition of the canonical divisor of the form
\begin{equation}\label{KXdecomp}
K_X =P_X + C_X,
\end{equation}
where $P_X$ is a divisor lying in the positive cone of $X$ and $C_X$ is an effective, nonzero divisor, see Theorem \ref{th:H1Theta}, 2). 
The divisor $C_X$ and the sheaves $Tor(S_X)$ and ${\cal J}_{A_X}(L_X)$ associated to it control  the space $H^1(\TET)$ of infinitesimal deformations of $X$. So the properties of this divisor are important for the further understanding of $H^1 (\TET)$.

From Theorem \ref{th:H1Theta}, 2), we know that
$C_X$ is subject to
\begin{equation}\label{CVbound}
C_X \cdot K_X \leq (4-2rk(S_X))(3c_2 (X)-K^2_X).
\end{equation}
In this section we impose a numerical condition on $(3c_2 (X)-K^2_X)$ which allows to control the geometry of the divisor $C_X$ as well as the geometry of $X$.
The numerical condition in question is as follows:
\begin{equation}\label{qbound}
	3c_2 (X)-K^2_X <\HA \Big(K^2_X \Big)^{\HA}.
\end{equation}
As we will see shortly, this is essentially motivated by the inequality
\eqref{CVbound} and the Hodge Index inequality: the two put together imply that the intersection form restricted to the sublattice of $NS(X)$
spanned by the irreducible components of $C_X$ is negative semidefinite.

We begin by observing that the condition guaranties the assumption 
$\displaystyle{\alpha_X <\frac{3}{8}}$ on the ratio of the Chern numbers imposed in Theorem \ref{th:H1Theta}.
\begin{lem}\label{lem:qbound-ratio}
Let $X$ be a smooth compact minimal complex surface of general type of positive index $\tau_X >0$. Then the inequality in \eqref{qbound} implies $\displaystyle{\alpha_X <\frac{3}{8}}$.
	\end{lem}
\begin{pf}
	Dividing \eqref{qbound} by $K^2_X$ we obtain
	$$
	3\alpha_X -1 <\HA \Big(K^2_X \Big)^{-\HA}.
	$$
	Assume $\displaystyle{\alpha_X \geq \frac{3}{8}}$. Then the above
	implies $K^2_X \leq 15$. But the condition
	$\tau_X>0$ tells us $K^2_X >8\chi(\OO_X)$. Hence $\chi(\OO_X)=1$.
	From BMY inequality $K^2_X\leq 9\chi(\OO_X)$, we deduce
	$K^2_X=9$ and $c_2 (X)=3$ with the ratio 
	$\displaystyle{\alpha_X =\frac{1}{3}}$ and this contradicts the assumption
	$\displaystyle{\alpha_X \geq \frac{3}{8}}$. 
\end{pf}

From now on we assume:
{\small
\begin{equation}\label{Xassupm}
	\begin{gathered}
	\text{\it $X$ is a smooth compact complex surface with the canonical divisor $K_X$ ample,}
	\\
	\text{\it the index $\tau_X >0$ and its Chern numbers $(c_2 (X),K^2_X)$ are subject to the inequality \eqref{qbound}.}
\end{gathered}
\end{equation}
}

\begin{lem}\label{lem:CX-negsd}
	For a surface $X$ subject to the conditions \eqref{Xassupm}
	either $h^1 (\TET)\leq 1$ or $h^1 (\TET)\geq 2$ and then the canonical divisor $K_X$ admits
	a distinguished decomposition
	$$
	K_X =P_X +C_X
	$$
	where $P_X$ lies in the positive cone of $X$ and $C_X$ is an effective, nonzero divisor whose irreducible components span
	a negative semidefinite sublattice of $NS(X)$, the N\'eron-Severi
	group of $X$. Furthermore, the divisor $C_X$ is subject to
	$$
	C_X \cdot K_X -\HA C^2_X \leq 3c_2 (X)-K^2_X.
	$$  
	\end{lem}
\begin{pf}
	We may assume $h^1 (\TET)\geq 2$, otherwise there is nothing to prove. The conditions \eqref{Xassupm} and Lemma \ref{lem:qbound-ratio} imply that the ratio of the Chern numbers $\displaystyle{\alpha_X <\frac{3}{8}}$. So we can apply Theorem \ref{th:H1Theta} to obtain the decomposition
	$$
	K_X =P_X +C_X.
	$$
	In addition, the effective nonzero divisor $C_X$ is subject to
	$$
	C_X \cdot K_X \leq (4-2rk(S_X))(3c_2 (X)-K^2_X),
	$$
	see part 2) of Theorem \ref{th:H1Theta}. With the inequality \eqref{qbound} we can improve it to
	\begin{equation}\label{CX-bound}
		C_X \cdot K_X \leq 2(3c_2 (X)-K^2_X).
	\end{equation}
This is certainly true if $rk(S_X)=1$, so we assume $rk(S_X)=0$. In this case we only have
$$
	C_X \cdot K_X \leq 4(3c_2 (X)-K^2_X)
	$$
and we need to `upgrade' it to \eqref{CX-bound}. For this combine the above inequality with \eqref{qbound} 
	$$
	C_X \cdot K_X \leq 4(3c_2 (X)-K^2_X) <2(K^2_X)^{\HA}.
	$$	
This and the Hodge Index inequality give
	$$
	C^2_X K^2_X \leq (C_X\cdot K_X)^2 < 4K^2_X;
	$$
	hence $C^2_X\leq 3$; next recall that in $rk(S_X)=0$ case the divisor $C_X$ satisfies the inequality
	\begin{equation}\label{CX-BMYbound}
	C_X \cdot K_X -2C^2_X \leq 3c_2(X)-K^2_X,
\end{equation}
	see Lemma \ref{lem:Q'Vstable}, 4); this and the previous bound on
	 $C^2_X$ imply
	$$
	C_X \cdot K_X \leq 3c_2(X)-K^2_X + 6.
	$$ 
Since we are assuming that $(3c_2(X)-K^2_X)$ is not zero, this is because $H^1(\TET)\neq 0$, and it is divisible by $4$:
$$
3c_2(X)-K^2_X=3c_2(X)-(12\chi(\OO_X)-c_2(X))=4c_2(X)-12\chi(\OO_X),
$$
where the first equality uses the Noether's formula. Hence
\begin{equation}\label{BMYbigger4}
3c_2(X)-K^2_X \geq 4
\end{equation}
 and the preceding inequality can be rewritten as follows
$$
	C_X \cdot K_X \leq 2(3c_2(X)-K^2_X) + 2;	
$$
using the inequality \eqref{qbound} again gives
$$
C_X \cdot K_X  <(K^2_X)^{\HA} + 2;
$$
applying the Hodge Index inequality once again we obtain
$$
C^2_X K^2_X \leq (C_X\cdot K_X)^2 <((K^2_X)^{\HA} + 2)^2=K^2_X +4(K^2_X)^{\HA} +4;
$$
hence
$$
C^2_X <1+4(K^2_X)^{-\HA} +\frac{4}{K^2_X}\leq 1+\frac{4}{3} +\frac{4}{9}<3,
$$
where the second inequality uses $K^2_X \geq 9$.
Thus $C^2_X \leq 2$; returning with this to the inequality 
\eqref{CX-BMYbound} gives
$$
C_X \cdot K_X \leq (3c_2(X)-K^2_X)+4\leq 2(3c_2(X)-K^2_X),
$$
where the second inequality comes from \eqref{BMYbigger4}.
The inequality \eqref{CX-bound} is now proved. 

We go on to showing that the sublattice of the N\'eron-Severi group of $X$ generated by the irreducible components of $C_X$ is negative semidefinite.
 
For any component $C$ of $C_X$ we have 
$$
C\cdot K_X \leq C_X \cdot K_X \leq 2(3c_2(X)-K^2_X).
$$
 Combining this  with the inequality \eqref{qbound} gives
 $$
C\cdot K_X <(K^2_X)^{\HA}.
$$
This and the  Hodge Index inequality imply
 $$
 C^2K^2_X \leq  (C\cdot K_X)^2 <K^2_X.
 $$
 Hence
	$$
	C^2 \leq 0,
	$$
	for every component $C$ of $C_X$. This means that the sublattice of
	$NS(X)$ generated by the irreducible components of $C_X$ is negative semidefinite.
	
	We now turn to the inequality asserted in the lemma. This is an improvement of the bound in Theorem \ref{th:H1Theta}, 2), due to the fact that $C^2_X \leq 0$. Indeed, recall that the middle column in the diagram \eqref{X-diag} is the $K_X$-destabilizing sequence for the sheaf $F_X$. Furthermore, dualizing, the sheaf $(Q'_X)^{\ast}$ is the maximal $K_X $-destabilizing subsheaf of $F^{\ast}_X$, the dual of $F_X$. This gives the equality
	\begin{equation}\label{tauch2}
	\tau_X =2ch_2 (F_X) = 2ch_2 ((Q'_X)^{\ast}) +\sum^m_{i=2}2ch_2 (Q_i),
	\end{equation}
where the sum is over the remaining semistable factors of the HN filtration of $F^{\ast}_X$ with respect to $K_X$. Furthermore, Bogomolov-Gieseker-Miyaoka inequality for $(Q'_X)^{\ast}$ gives
$$
 2ch_2 ((Q'_X)^{\ast}) \leq \frac{1}{3}c^2_1 (Q'_X),
 $$
 while Bogomolov-Gieseker for the remaining factors reads
 $$
 2ch_2 (Q_i)\leq \frac{c^2_1 (Q_i)}{rk(Q_i)},\,\,\forall i\geq 2.
 $$
 
 	Recall that in the decomposition $K_X =P_X +C_X$, the divisor $P_X=-c_1 (Q'_X)$, while 
 	$$
 	C_X = \sum^m_{i=2}c_1 (Q_i).
 	$$
 Since all $c_1 (Q_i)$'s have positive degrees with respect to $K_X$ we deduce
 $$
 0<c_1 (Q_i)\cdot K_X \leq C_X\cdot K_X.
 $$
 From the first part of the proof we deduce
 $
 c^2_1 (Q_i) \leq 0,
 $
 for all $i\geq 2$.	Hence
 $$
 2ch_2(Q_i)\leq 0,\,\,\forall i\geq 2.
 $$
 Using these inequalities and the one for $2ch_2(Q'_X)$ in the formula \eqref{tauch2} we deduce
 $$
 \tau_X \leq \frac{1}{3}c^2_1 (Q'_X)
 $$
  	This and the relation
	$$ 
		-c_1(Q'_X)=P_X=K_X -C_X
	$$
gives the inequality	
	$$
	\tau_X \leq \frac{1}{3}(K_X -C_X)^2= \frac{1}{3} K^2_X -\frac{2}{3} K_X\cdot C_X +\frac{1}{3} C^2_X
	$$
or, equivalently,
$$
K_X \cdot C_X -\HA C^2_X \leq 3c_2 (X)-K^2_X.
$$
\end{pf}

We begin to draw the geometric consequences from the properties of the divisor $C_X$.
\begin{cor}\label{cor:f_X-fibration}
	Let $X$ be subject to \eqref{Xassupm} and assume $h^1(\TET)\geq 2$. If the filtration
	$$
	H^1 (\TET) \supset H^1 (\TET)^0
	$$ 
	in \eqref{1stepfilt} has the quotient space
	$H^1 (\TET) / H^1 (\TET)^0$ of dimension at least $2$, then
	the divisor $C_X$ in Lemma \ref{lem:CX-negsd} admits the decomposition
	\begin{equation}\label{C_X}
	C_X =D^0_X + L_X,
	\end{equation}
	where $L_X$ and $D^0_X$ are as in Theorem \ref{th:H1Theta}, 3) and 4) respectively. Furthermore, in the diagram \eqref{X-diag2} the sheaf ${\cal J}_{A_X} (L_X)=\OO_X (L_X)$, it is generated by its global sections and gives a surjective morphism
	$$
 f_X:	X\longrightarrow B
	$$
	 onto a smooth projective curve $B$ with connected fibres, that is, there is a divisor $A$ on $B$ such that 
	$\OO_X (L_X)=f^{\ast}_X(\OO_B (A))$. In addition,
	$$
	\begin{gathered}
	dim (H^1 (\TET) / H^1 (\TET)^0) \leq h^0 (\OO_X (L_X))=h^0 (\OO_B (A)),
	\\
	L_X \cdot K_X =deg(A)C \cdot K_X =2deg(A)(g_C-1) \leq 3c_2 (X)-K^2_X,
	\end{gathered}
	$$
	where $C$ is the divisor class of a fibre of $f_X$ and $g_C$ is the genus of a smooth fibre of $f_X$.
	
	In addition, if the divisor $D^0_X$ of the decomposition in \eqref{C_X} is nonzero, then it is contained in the fibres of the morphism $f_X$.	
\end{cor}
\begin{pf}
	By definition $L_X$ is a component of $C_X$ and according to Lemma \ref{lem:CX-negsd} the self-intersection $L^2_X \leq 0$. On the other hand, from Theorem \ref{th:H1Theta}, 3) we know that the sheaf
	${\cal J}_{A_X} (L_X)$ is globally generated by the space $H^1 (\TET) / H^1 (\TET)^0$. Hence
	$$
	0\geq L^2_X \geq deg(A_X).
	$$
	From this it follows
	$$
	L^2_X = deg(A_X)=0
	$$
	and ${\cal J}_{A_X} (L_X) =\OO_X (L_X)$ is globally generated.
	This implies that the line bundle defines a morphism
	$$
	\phi_{L_X}: X\longrightarrow \PP(H^0(\OO_X (L_X))^{\ast})
	$$
	whose image is a reduced irreducible curve. The Stein factorization of this morphism is the morphism $f_X$ asserted in the corollary.
	Thus the morphism $f_X: X \longrightarrow B$
	onto a smooth projective curve $B$ with connected fibres and
	the line bundle
	$
	\OO_X (L_X)=f^{\ast}_X (\OO_B (A))
	$
	is the pull back by $f_X$ of a line bundle $\OO_B (A)$ on $B$.
	Taking the direct image we obtain
	$$
	f_{X\ast} (\OO_X (L_X)) =\OO_B (A).
	$$
	Hence the identifications
	$$
	H^0(\OO_X (L_X))\cong H^0 (f_{X\ast} (\OO_X (L_X))) \cong H^0 (\OO_B (A)).
	$$
	From Theorem \ref{th:H1Theta}, 3), follows the inclusion
	$$
	H^1 (\TET)/H^1 (\TET)^0 \hookrightarrow 	H^0(\OO_X (L_X)).
	$$
	Hence the assertion
	$$
	dim(	H^1 (\TET)/H^1 (\TET)^0) \leq h^0 (\OO_X (L_X)).
	$$
	
	Let $C$ be the divisor class of a fibre of $f_X$. Then the formula
	$	\OO_X (L_X)=f^{\ast}_X (\OO_B (A))$ implies
	$$
	L_X =deg(A)C.
	$$
	Intersecting with $K_X$ gives
	\begin{equation}\label{LXgF}
	L_X \cdot K_X =deg(A)C\cdot K_X =2deg(A)(g_C -1)
	\end{equation}
	where $g_C$ is the genus of a smooth fibre of $f_X$ and the last equality uses the adjunction formula
	$$
	C\cdot K_X =2(g_C-1).
	$$
	From the decomposition \eqref{C_X} it follows
	$$
	L_X \cdot K_X \leq C_X \cdot K_X.
	$$
	This and the upper bound in Lemma \ref{lem:CX-negsd} give
	$$
	3c_2 (X)-K^2_X \geq L_X \cdot K_X =2deg(A)(g_C -1),
	$$
	where the last equality comes from \eqref{LXgF}.
	
	We now turn to the assertions about the divisor $D^0_X$ of the decomposition \eqref{C_X}. We may assume it to be nonzero and we consider
	the divisor
	$$
D_t=	tD^0_X+C,
	$$
	for $t \in \RR$. This divisor is contained in the lattice spanned
	by the irreducible components of $C_X$ and the lattice is negative semidefinite by Lemma \ref{lem:CX-negsd}. Hence
	$$
	0\geq D^2_t=(tD^0_X+C)^2=t^2 \big(D^0_X)^2 +2tD^0_X \cdot C, \,\, \forall t\in \RR.
	$$
	In particular, for all  $t>0$ we obtain
	$$
	0 \geq t \big(D^0_X)^2 +2D^0_X \cdot C.
	$$
	 Taking the limit $t\to 0^+$ gives
	$$
	D^0_X \cdot C \leq 0.
	$$
	Since $D^0_X$ is effective and $C$ is nef the above inequality
	must be equality
	$$
	D^0_X \cdot C =0
	$$
and this means that all irreducible components of $D^0_X$ are contained in the fibres of $f_X$.	 
\end{pf}

The fibration $f_X :X \longrightarrow B$ in Corollary \ref{cor:f_X-fibration} has a distinctive feature of the genus of the fibres $g_C$ being subject to the inequality
$$
2g_C -2=C\cdot K_X \leq 3c_2 (X)-K^2_X.
$$
This together with the upper bound
$$
3c_2 (X)-K^2_X <\HA (K^2_X)^{\HA}
$$
in \eqref{qbound} imply that this fibration is essentially unique. 
\begin{cor}\label{cor:fX-unique}
	Assume the setting of Corollary \ref{cor:f_X-fibration} and assume
	$$
	dim(H^1(\TET)/H^1(\TET)^0))\geq 2.
	$$
	Let $f': X\longrightarrow B'$ be a another fibration of $X$, that is, a surjective morphism onto a smooth projective curve $B'$, with connected fibres. Denote by $C'$ the divisor class of a fibre of $f'$ and assume
	$$
2g_{C'} -2=C'\cdot K_X \leq 3c_2 (X)-K^2_X.
$$
Then there is a unique isomorphism
$$
\rho:B'\longrightarrow B
$$
making the diagram
$$
\xymatrix@R=12pt@C=12pt{
	&X \ar[ld]_{f'}\ar[rd]^{f_X}&\\
	B'\ar[rr]^{\rho}&&B
}
$$
commutative, where $f_X$ is the fibration described in Corollary \ref{cor:f_X-fibration}.	
\end{cor}
\begin{pf}
The main point is to show that the divisor classes $C'$ and $C$ of the fibrations $f'$ and $f_X$ are subject to
\begin{equation}\label{F-F'}
	C'\cdot C=0.
\end{equation}
To see this consider the divisor $(C'+C)$. Intersecting with $K_X$ gives
$$
(C'+C)\cdot K_X =C'\cdot K_X + C\cdot K_X \leq 2(3c_2(X)-K^2_X).
$$
Combining with the upper bound \eqref{qbound} gives
$$
(C'+C)\cdot K_X <(K^2_X)^{\HA}.
$$
Applying the Hodge Index inequality gives
$$
(C'+C)^2 K^2 \leq \big((C'+C)\cdot K_X \big)^2 <K^2_X.
$$
From this it follows
$$
2C'\cdot C = (C'+C)^2 <1.
$$
Hence $C'\cdot C \leq 0$. Since the intersection $C'\cdot C \geq 0$, because both divisors are nef, the equality \eqref{F-F'} follows.

The equality $	C'\cdot C=0$ tells us that every fibre of the morphism
$f'$ is contained in the fibres of $f_X$ and vice-versa. In particular, a reduced irreducible fibre $C'_{b'}$ of $f'$ is contained in a unique fibre $C_b$ of $f_X$. Furthermore, since $f_X$ has at most finite number
of reducible fibres, we may assume $C_b$ is irreducible. Hence the equality
of the fibres.
$$
C'_{b'} =C_b.
$$
This implies the equality of the divisor classes
\begin{equation}\label{F'=F}
	C'=C.
\end{equation}

The above argument also defines a map
$$
\rho:B' \longrightarrow B
$$
which sends a closed point $b' \in B'$ into the closed point
$\rho (b')$ of $B$ subject to $C'_{b'} \subset C_{\rho(b')}$. Furthermore, the equality \eqref{F'=F} implies that $\rho$ has degree $1$. 

Exchanging the roles of $B$ and $B'$ gives a morphism
$$
\rho':B \longrightarrow B'.
$$
sending a closed point $b\in B$ to the closed point $\rho'(b)$ in $B'$
subject to $C_b \subset C'_{\rho'(b)}$.

The composition 
$$
g:=\rho' \circ \rho : B' \longrightarrow B',
$$
by definition, sends a closed point $b'$ to the closed point $g(b')$ of $B'$ with the property
$$
C'_{b'}\subset C'_{g(b')}.
$$
From this it follows that $g=id_{B'}$ is the identity map. Hence $\rho$ and $\rho'$ are isomorphisms  and they are inverses of each other.	
\end{pf}

From now on the genus of the fibration $f_X :X\longrightarrow B$ in Corollary \ref{cor:f_X-fibration} will be denoted $g_{{}_X}$. Set
$Sing(f_X)$ to be the subscheme of the critical {\it values} of $f_X$ and
$\text{\it \r B}=B\setminus Sing(f_X)$ its complement. Over {\it \r B} the fibration 
$$
f_X: \text{\it \r X} \longrightarrow \text{\it \r B}
$$
is smooth, here  $\text{\it \r X}=f^{-1}_X (\text{\it \r B})$.
This gives rise to a unique morphism
\begin{equation}\label{mX-morph}
	m_X : \text{\it \r B} \longrightarrow {\mathfrak{M}}_{g_{{}_X}},
\end{equation}
where ${\mathfrak{M}}_{g_{{}_X}}$ is the moduli space of curves of genus $g_{{}_X}$.
\begin{cor}\label{cor:mX-immersion}
	The morphism $m_X$ in \eqref{mX-morph} is an immersion.
\end{cor}
\begin{pf}
The differential of the morphism $f_X: X\longrightarrow B$ gives rise
to the exact sequence
$$
\xymatrix@R=12pt@C=12pt{
0\ar[r]& \Theta_{f_X} \ar[r]&\TET \ar[r]& f^{\ast}_X (\Theta_B)
},
$$
where $\Theta_{f_X}$ is the relative tangent sheaf of $f_X$. Over {\it \r B} the sequence is exact on the right as well, that is, we have the exact sequence
$$
\xymatrix@R=12pt@C=12pt{
	0\ar[r]& \Theta_{f_X}\otimes \OO_{\text{\it \r X}} \ar[r]&\TET\otimes \OO_{\text{\it \r X}} \ar[r]& f^{\ast}_X (\Theta_{\text{\it \r B}})\ar[r]& 0.
}
$$
Taking the direct image with respect to $f_{X}$ we obtain
\begin{equation}\label{fX-dirimage}
\xymatrix@R=12pt@C=12pt{
	0 \ar[r]& \Theta_{\text{\it \r B}}\ar[r]&R^1f_{X \ast} (\Theta_{f_X}\otimes \OO_{\text{\it \r X}})
	\ar[r]&R^1f_{X \ast} (\Theta_{X}\otimes \OO_{\text{\it \r X}})\ar[r]&\Theta_{\text{\it \r B}} \otimes R^1f_{X \ast} (\OO_X)\ar[r]&0.
}
\end{equation}
This is because the direct images
\begin{equation}
\text{$f_{X \ast} (\Theta_{f_X}) =0$ and $f_{X \ast} (\Theta_{X}) =0$.}	
\end{equation} 
The first vanishing follows from the fact 
$$
H^0 (\Theta_{f_X} \otimes \OO_{C_b})=0,
$$
for every fibre $C_b$ of $f_X$, and the second will be proved in a moment. Once this is done the assertion of the corollary follows, since the monomorpism
$$
\xymatrix@R=12pt@C=12pt{
	0 \ar[r]& \Theta_{\text{\it \r B}}\ar[r]&R^1f_{X \ast} (\Theta_{f_X}\otimes \OO_{\text{\it \r X}})
}
	$$
	in the exact complex \eqref{fX-dirimage} is the differential of $m_X: \text{\it \r B} \longrightarrow \mathfrak{M}_{g_{{}_X}}$ and it is  a monomorphism of vector bundles.

To complete the argument we return to the proof of the vanishing
$$
f_{X\ast}(\TET)=0.
$$
Assume that the direct image $f_{X\ast}(\TET)$ is nonzero. From the vanishing on the level of global sections
$$
H^0(f_{X\ast}(\TET))\cong H^0(\TET)=0
$$
it follows that $f_{X\ast}(\TET)$ is torsion free and
since we are on a smooth curve it is locally free of rank $h^0(\TET \otimes \OO_{C_b})$, where $C_b=f^{-1}_X (b)$ is the fibre of $f_X$ over $b\in B$. Taking $C_b$ to be a smooth fibre, consider the normal sequence of  $C_b$ in $X$
$$
\xymatrix@R=12pt@C=12pt{
	0\ar[r]& \Theta_{C_b} \ar[r]&\TET\otimes \OO_{C_b} \ar[r]& \OO_{C_b} \ar[r]&0.
}
$$ 
From this it follows 
that $H^0(\TET \otimes \OO_{C_b}) \neq 0$ if and only if the normal sequence splits
\begin{equation}\label{normCb-splits}
\TET\otimes \OO_{C_b} \cong  \Theta_{C_b} \oplus \OO_{C_b}
\end{equation}
and then the rank of $f_{X\ast}(\TET)$ is
$$
rk(f_{X\ast}(\TET))=h^0(\TET\otimes \OO_{C_b} )=h^0(\OO_{C_b}) =1.
$$
Hence $f_{X\ast}(\TET)$ is a line bundle which we denote
$\OO_B (T)$. Taking the pull back by $f_X$ gives the nonzero morphism
\begin{equation}\label{OBT-ThetaX}
f^{\ast}_X(\OO_B (T))=f^{\ast}_X (f_{X\ast} (\TET))	\longrightarrow \TET
\end{equation}
corresponding to the identity morphism
in $Hom(f_{X\ast} (\TET), f_{X\ast} (\TET))$ under the isomorphism
$$
Hom_{\OO_X}(f^{\ast}_X (f_{X\ast} (\TET)), \TET)\cong Hom_{\OO_B} (f_{X\ast} (\TET), f_{X\ast} (\TET)).
$$
The morphism in \eqref{OBT-ThetaX} is a monomorphism and we complete it to the exact sequence
$$
\xymatrix@R=12pt@C=12pt{
	0\ar[r]& f^{\ast}_X(\OO_B (T))\ar[r]& \TET \ar[r]& Q\ar[r]&0,
}
$$
where $Q=\TET /f^{\ast}_X(\OO_B (T))$ is the quotient sheaf. Factoring out the torsion part of $Q$, if nonzero, gives the exact sequence
\begin{equation}\label{Theta-fX-seq}
\xymatrix@R=12pt@C=12pt{
	0\ar[r]& \OO_X (f^{\ast}_X (T) +D)\ar[r]& \TET \ar[r]& {\cal J}_A (-K_X - f^{\ast}_X (T) -D)\ar[r]&0,
}
\end{equation}
where $D$ is an effective divisor: this is the divisorial component of the torsion part of 
$Q$, 
and ${\cal J}_A$ is the ideal sheaf of at most $0$-dimensional subscheme
$A$ of $X$.
Taking the intersection with $C$, the divisor class of a fibre of $f_X$, tells
us that the subsheaf $\OO_X (f^{\ast}_X (T) +D)$ in \eqref{Theta-fX-seq} has  the degree
$$
(f^{\ast}_X (T) +D)\cdot C=D\cdot C \geq 0
$$
on the fibres of $f_X$. We deduce 
$$
D\cdot C = 0,
$$
 since otherwise by a result of Miyaoka, see Theorem 8.4, \cite{Mi3}, $X$ is covered by a family of rational curves which is impossible for a surface of general type. The above equality implies that $D$ is contained in the fibres of $f_X$ and hence
 \begin{equation}\label{D2less0}
 (f^{\ast}_X (T) +D)^2 =D^2 \leq 0.
\end{equation}
Dualizing the exact sequence \eqref{Theta-fX-seq} gives the inclusion
$$
\xymatrix@R=12pt@C=12pt{
	0\ar[r]& \OO_X (K_X +f^{\ast}_X (T) +D)\ar[r]& \Omega_X
}
$$
Intersecting $(K_X +f^{\ast}_X (T) +D)$ with $C$ gives
$$
(K_X +f^{\ast}_X (T) +D) \cdot C=K\cdot C >0.
$$
Since $C$ is nef and hence is a limit of ample divisors, we have
$$
(K_X +f^{\ast}_X (T) +D) \cdot H >0
$$
for some ample divisor $H$. In addition, $\OO_X (K_X +f^{\ast}_X (T) +D)$ is
a 
line subsheaf of $\Omega_X$, so $(K_X +f^{\ast}_X (T) +D)$ can not lie in the positive cone of $X$. Hence
$$
(K_X +f^{\ast}_X (T) +D) ^2 \leq 0.
$$
This together with \eqref{D2less0} and the exact sequence \eqref{Theta-fX-seq} imply
$$
\begin{gathered}
\tau_X =2ch_2 (\TET)=2ch_2 (\OO_X (f^{\ast}_X (T) +D))+2ch_2({\cal J}_A (-K_X - f^{\ast}_X (T) -D))=
\\
(f^{\ast}_X (T) +D)^2 +	(K_X +f^{\ast}_X (T) +D) ^2 -2deg(A) \leq 0
\end{gathered}
$$
and this is contrary to the assumption $\tau_X >0$.
\end{pf}
 
The map $m_X$ in \eqref{mX-morph} can be extended to a unique map
$$
\overline{m}_X: B\longrightarrow \overline{\mathfrak{M}}_{g_{{}_X}}
$$
into the Deligne-Mumford stack $\overline{\mathfrak{M}}_{g_{{}_X}}$ of stable curves of genus $g_{{}_{X}}$. The pull back of the boundary divisor
$
\delta_{g_{{}_X}}:=\overline{\mathfrak{M}}_{g_{{}_X}} \setminus \mathfrak{M}_{g_{{}_X}}
$
gives us $(B, {\overline{m}}^{\ast}_X (\delta_{g_{{}_X}}))$, the curve $B$ together with the marked points ${\overline{m}}^{\ast}_X (\delta_{g_{{}_X}})$ on it. This is expected to be functorial. 

Namely,
 a smooth {\it connected} family
$$
\pi_U: {\mathfrak X}_{U} \longrightarrow U
$$
of surfaces with the Chern numbers $(c_2,K^2)$ subject to the assumptions
\eqref{Xassupm}
 and 
\begin{equation}\label{h1-cond}
dim (H^1 (\Theta_{X_u})/H^1 (\Theta_{X_u})^0)\geq 2,
\end{equation}
 for every closed point $u\in U$, produces the family of curves
 $$
 \omega_U:{\mathfrak{B}}_U \longrightarrow U
 	$$
over $U$ fitting into the commutative diagram of $U$-morphisms
$$
\xymatrix@R=12pt@C=12pt{
{\mathfrak X}_{U}\ar_{\pi_U}[dr] \ar^{f_U}[rr]& &{\mathfrak{B}}_U  \ar^{\omega_U}[dl]\\
& U&
}
$$
where over every closed point $u\in U$ the diagram becomes
$$
\xymatrix@R=12pt@C=12pt{
	\pi^{-1}_U (u)=X_{u}\ar[dr] \ar^{f_{X_u}}[rr]& &B_u=\omega^{-1}(u)  \ar[dl]\\
	& u&
}
$$
where $f_{X_u}$ is the morphism in Corollary \ref{cor:f_X-fibration}.
Thus we obtain 
$$
f_U:{\mathfrak X}_{U} \longrightarrow {\mathfrak{B}}_U
$$ 
a family of curves of genus $g=g(c_2,K^2)$ depending only on $(c_2,K^2)$
naturally attached to every connected family $\pi_U: {\mathfrak X}_{U} \longrightarrow U$ of surfaces with Chern numbers $(c_2,K^2)$ and subject to
\eqref{Xassupm} 
and \eqref{h1-cond}. That family of curves in turn gives a 
unique rational map
$$
\overline{m}_{{}_U}:{\mathfrak{B}}_U -- > \overline{\mathfrak{M}}_g
$$
and this can be viewed as a family of maps of marked curves into the Deligne-Mumford stack $\overline{\mathfrak{M}}_g$. 
Thus what comes out of this discussion is the following:
\begin{equation}\label{GW}
	\begin{gathered}
\text{\it the stack of surfaces subject to \eqref{Xassupm} and \eqref{h1-cond}}
\\
\text{\it
	defines Gromov-Witten type theory.}
\end{gathered}	
\end{equation}
From this point of view regular surfaces, that is, surfaces $X$ with
$q_X =h^1 (\OO_X)=0$, emerge as the Gromov-Witten theory of genus $g=0$:
\begin{equation}\label{GW0}
	\begin{gathered}
	\text{\it the stack of regular surfaces subject to \eqref{Xassupm} and \eqref{h1-cond}
		defines a substack}
		\\
		 \text{\it${\bf Maps}(\PP^1,\overline{\mathfrak{M}}_g)$ of maps
		from $\PP^1$ to $\overline{\mathfrak{M}}_g$. }	
	\end{gathered}
\end{equation}
This is because for $X$ regular the morphism $f_X$ in Corollary \ref{cor:f_X-fibration} maps $X$ onto $\PP^1$. 

All of the above will be elaborated in the sequel to this paper. We want to conclude our discussion here by pointing out that the images of maps in \eqref{GW0} sweep out at most a 
surface in $\overline{\mathfrak{M}}_g$.

To explain this result we return to the morphism
$$
f_X :X\longrightarrow B
$$
in Corollary \ref{cor:f_X-fibration}. For a smooth fibre $C_b =f^{-1}_X (b)$ of $f_X$ we have the pair $(X,C_b)$. This gives the closed point of the stack
${\bf Pairs}_{c_2,K^2,g}$ of
 pairs $(X,C)$ of surfaces $X$ with Chern numbers
$(c_2,K^2)$ and the fibration $f_X:X\longrightarrow B$ in Corollary \ref{cor:f_X-fibration}, and $C$ is a smooth fibre of $f_X$, a curve of genus
$g=g_X$. 
This comes with two morphisms
$$
\xymatrix@R=12pt@C=12pt{
	&{\bf Pairs}_{c_2,K^2,g}\ar_(.56)p[ld]\ar^(.56){p'}[rd]&\\
	{\mathfrak{M}}_{c_2,K^2}&&{\mathfrak{M}}_{g}
		}
	$$
	onto the stacks ${\mathfrak{M}}_{c_2,K^2}$ and ${\mathfrak{M}}_{g}$ of surfaces and curves respectively. We wish to understand the image of
	the morphism
	$$
	p':{\bf Pairs}_{c_2,K^2,g} \longrightarrow {\mathfrak{M}}_{g}.
	$$
	For this we look at the differential of $p'$ at $(X,C_b)$. Recall: 
	
	$\bullet$ 
	the space of the infinitesimal deformations of the pair $(X,C_b)$ is 
	$H^1(\Theta_X (-log(C_b)))$, where $\Theta_X (-log(C_b))$ is the sheaf of germs of holomorphic vector fields on $X$ tangent along $C_b$, 
	
	$\bullet$ 
	the space of the infinitesimal deformations of $C_b$ is $H^1 (\Theta_{C_b})$, where $\Theta_{C_b}$ is the holomorphic tangent bundle of $C_b$,
	
	$\bullet$ 
	the differential of $p'$ at $(X,C_b)$ is given by the homomorphism
	\begin{equation}\label{dp'-X-Cb}
	dp'_{(X,C_b)}:	H^1 (\Theta_X (-log(C_b)))\longrightarrow H^1(\Theta_{C_b}).
	\end{equation}
	
 We 
 claim the following.
 \begin{pro}\label{pro:imdp}
 	Assume $X$ is a regular surface subject to \eqref{Xassupm} and to the cohomological condition 
 	$$
 	dim(H^1(\TET)/H^1(\TET)^0) \geq 2.
 	$$
 	 Let $f_X:X\longrightarrow \PP^1$ be the fibration in Corollary \ref{cor:f_X-fibration}. Then for every smooth fibre $C_b$ of
 	$f_X$
 	the image of the homomorphism in \eqref{dp'-X-Cb} is at most two-dimensional.
 \end{pro}
\begin{pf}
	We begin by recording the defining exact sequence for  $\Theta_X (-log(C_b))$:
	$$
	\xymatrix@R=12pt@C=12pt{
		0\ar[r]&\Theta_X (-log(C_b))\ar[r]& \TET\ar[r]& \OO_{C_b} (C_b)\ar[r]\ar@{=}[d]&0,\\
		&&&\OO_{C_b}
	}
$$
where the epimorphism 
$$
\TET \longrightarrow \OO_{C_b} (C_b)=\OO_{C_b}
$$
comes from the composition of the restriction-to-$C_b$ morphism with the epimorphism
in 
the normal sequence of $C_b$ in $X$, the slanted arrow below:
$$
\xymatrix@R=12pt@C=12pt{
&&\TET\ar[d]\ar[dr]& &\\
0\ar[r]&\Theta_{C_b}\ar[r]& \TET \otimes \OO_{C_b} \ar[r]&\OO_{C_b} \ar[r]&0.
}
$$
This implies that the defining sequence of $\Theta_X (-log(C_b))$ is a part of the following diagram
\begin{equation}\label{log-diag-Cb}
	\xymatrix@R=12pt@C=12pt{
		&0\ar[d]&0\ar[d]&&\\
		0\ar[r]&\Theta_X (-C_b)\ar@{=}[r]\ar[d]& \TET(-C_b)\ar[d]&  &\\	
		0\ar[r]&\Theta_X (-log(C_b))\ar[r]\ar[d]& \TET\ar[r]\ar[d]& \OO_{C_b} \ar[r]\ar@{=}[d]&0\\
		0\ar[r]&\Theta_{C_b} \ar[r]\ar[d]&\TET \otimes \OO_{C_b}\ar[d] \ar[r]&\OO_{C_b} \ar[r]&0\\
		&0&0&&
	}
\end{equation}
The homomorphism in \eqref{dp'-X-Cb} is induced, on the cohomology level, by the epimorphism of the left column of the diagram. So we examine the relevant part of the diagram of cohomology groups resulting from the diagram
\eqref{log-diag-Cb}
\begin{equation}\label{log-diag-Cb-coh}
	\xymatrix@R=12pt@C=12pt{
	&	&0\ar[d]&0\ar[d]\\
	&	&H^1(\Theta_X (-C_b))\ar@{=}[r]\ar[d]& H^1(\TET(-C_b))\ar[d]  \\	
	0\ar[r]&	H^0(\OO_{C_b})\ar[r]\ar@{=}[d]&H^1(\Theta_X (-log(C_b)))\ar[r]\ar[d]& H^1(\TET)\ar[d]\\
	0\ar[r]&	H^0(\OO_{C_b})\ar[r]	&H^1 (\Theta_{C_b}) \ar[r]&H^1(\TET \otimes \OO_{C_b}) 
	}
\end{equation}
From this we deduce
\begin{equation}\label{dim-dp'}
	\begin{gathered}
	\text{\it the dimension of the image of the homomorphism in \eqref{dp'-X-Cb}}
	\\
\text{\it is equal to $\big(h^1(\Theta_X (-log(C_b))) -h^1(\Theta_X (-C_b))\big)$.}
\end{gathered}
\end{equation}
From the middle row of \eqref{log-diag-Cb-coh} it follows
$$
h^1(\Theta_X (-log(C_b))) -h^1(\Theta_X (-C_b)) \leq 1 +h^1 (\TET)-	h^1(\Theta_X (-C_b))
$$
and the assertion of the propositions amounts to showing
\begin{equation}
	\text{$h^1(\Theta_X (-C_b)) \geq h^1 (\TET)-1,$ {\it for every smooth fibre $C_b$ of $f_X$.}}
\end{equation}
The rest of the proof is devoted to showing the above statement.

We 
begin by checking that the subspace $H^1 (\TET)^0 \subset H^1 (\TET)$ lies in the kernel of the homomorphism
$$
H^1 (\TET)\longrightarrow H^1 (\TET \otimes \OO_{C_b})
$$
for every smooth fibre $C_b$ of $f_X$. For this recall the exact sequence
$$
\xymatrix@R=12pt@C=12pt{
	0\ar[r]& \TET \ar[r]& Q''_X \ar[r]& Tor (S_X)\ar[r]&0
}
$$
at 
the bottom of the diagram \eqref{X-diag1}. The torsion sheaf $ Tor (S_X)$ is supported on the divisor $D^0_X=c_1 (Tor (S_X))$. From the description of the irreducible components of $D^0_X$ in Proposition \ref{pro:nefpair}, 2), 3), we know that this is disjoint from any smooth fibre of $f_X$. So the restriction of the above exact sequence 
to $C_b$ gives the diagram
$$
\xymatrix@R=12pt@C=12pt{
	&0\ar[d]&0\ar[d]&&\\
	0\ar[r]& \TET(-C_b) \ar[d]\ar[r]& Q''_X (-C_b) \ar[r]\ar[d]& Tor (S_X)\ar[r] \ar@{=}[d]&0\\
	0\ar[r]& \TET \ar[d]\ar[r]& Q''_X \ar[d]\ar[r]& Tor (S_X)\ar[r]&0\\
	&\TET \otimes \OO_{C_b} \ar[d]\ar@{=}[r]&Q''_X \otimes \OO_{C_b}\ar[d]&&\\
	&0&0&&  
}
$$
Passing to cohomology results in the following
$$
\xymatrix@R=12pt@C=12pt{
	 0\ar[r]& H^0(Tor (S_X))\ar[r] \ar@{=}[d]&H^1( \TET(-C_b))\ar[d]\\
	0\ar[r]& H^0(Tor (S_X))\ar[r]&H^1(\TET)\ar[d]\\
	&&H^1 (\TET \otimes \OO_{C_b})  
}
$$
By definition the subspace $H^1 (\TET)^0$ is the image of $H^0(Tor (S_X))$ in the middle row of the above diagram and that diagram tells us that the
inclusion $H^0(Tor (S_X)) \hookrightarrow H^1 (\TET)$ factors through
$H^1( \TET(-C_b))$. Hence
\begin{equation}\label{H10inkerH1-Cb}
	H^1(\TET)^0 \subset ker\Big(H^1(\TET) \longrightarrow H^1 (\TET \otimes \OO_{C_b})\Big)
\end{equation} 

Next we choose a subspace $W$ in $H^1(\TET)$ complementary to $H^1(\TET)^0$.
With this choice we have
\begin{equation}
	\xymatrix@R=12pt@C=12pt{
	&&W\ar@{^{(}->}[d]\ar[rd]&&\\
	0\ar[r]&H^1(\TET)^0 \ar[r]&H^1(\TET)\ar[r]&W_X \ar[r]&0
}
\end{equation}
where $W_X =H^1(\TET)/H^1(\TET)^0$ and the slanted arrow is an isomorphism.
Our task is to show the following:
\begin{equation}\label{cl-W}
	\begin{gathered}
	\text{\it the image of $W$ under the homomorphism}
	\\
	 H^1(\TET) \longrightarrow H^1 (\TET \otimes \OO_{C_b})
	 \\
	 \text{\it
	  is at most one-dimensional.}
  \end{gathered}
\end{equation}

We 
proceed by viewing the subspace $W\subset H^1(\TET)$ as the extension sequence
$$
\xymatrix@R=12pt@C=12pt{
	0\ar[r]& \TET \ar[r]& F_W \ar[r]& W\otimes \OO_X \ar[r]&0.
}
$$
This is related to the extension
$$
\xymatrix@R=12pt@C=12pt{
	0\ar[r]& \TET \ar[r]& F_X \ar[r]& H^1(\TET)\otimes \OO_X \ar[r]&0
}
$$
in 
the middle row of \eqref{X-diag} via the diagram
$$
\xymatrix@R=12pt@C=12pt{
	&&0\ar[d]&0\ar[d]&\\
		0\ar[r]& \TET \ar[r]\ar@{=}[d]& F_W \ar[r]\ar[d]& W\otimes \OO_X\ar[d] \ar[r]&0\\
	0\ar[r]& \TET \ar[r]& F_X \ar[r]\ar[d]& H^1(\TET)\otimes \OO_X \ar[r]\ar[d]&0\\
	&&H^1(\TET)/W \otimes \OO_X \ar@{=}[r]\ar[d]&H^1(\TET)/W \otimes \OO_X \ar[d]&\\
	&&0&0& 
}
$$
We 
now recall that $F_X$ admits the $K_X$-destabilizing sequence
$$
\xymatrix@R=12pt@C=12pt{
	0\ar[r]& F'_X \ar[r]& F_X \ar[r]& Q'_X\ar[r]&0,
}
$$
see the middle column in \eqref{X-diag}. Composing the inclusion
$F_W \hookrightarrow F_X$ 
with the epimorphism in the above exact sequence gives
\begin{equation}\label{FW-Q'X}
\xymatrix@R=12pt@C=12pt{
	0\ar[r]& \TET \ar[r]\ar@{=}[d]& F_W \ar[r]\ar[d]& W\otimes \OO_X \ar[r]\ar[d]&0\\
		0\ar[r]&\TET\ar[r]&Q'_X\ar[r]&S_X\ar[r]&0
}
\end{equation}
In 
addition, the bottom sequence in the above diagram is a part of the following diagram
$$
\xymatrix@R=12pt@C=12pt{
	&&0\ar[d]&0\ar[d]&\\
	0\ar[r]& \TET \ar[r]\ar@{=}[d]& Q''_X\ar[r]\ar[d]& Tor(S_X) \ar[r]\ar[d]&0\\
	0\ar[r]&\TET\ar[r]&Q'_X\ar[r]\ar[d]&S_X \ar[r]\ar[d]&0\\
	&&\OO_X (L_X)\ar@{=}[r]\ar[d]&\OO_X (L_X)\ar[d]&\\
	&&0&0&
}
$$
see Theorem \ref{th:H1Theta}, 3). We also recall the epimorphism
$$ 
W_X \otimes \OO_X \longrightarrow \OO_X (L_X), 
$$
 see the right column in \eqref{X-diag2}. This and the choice of $W\subset H^1 (\TET)$ made so that $W \cong W_X$, ensure 	that the vertical arrow on the right in \eqref{FW-Q'X} gives the composition arrow,
 $$
 \xymatrix@R=12pt@C=12pt{
 W\otimes \OO_X \ar[r]\ar[rd]&S_X \ar[d]\\
 &\OO_X (L_X)\ar[d]\\
 &0
}
$$
the slanted arrow of the above diagram, which continue to be surjective.
This implies the following commutative diagram:
\begin{equation}\label{FW-Q'X-1}
	\xymatrix@R=12pt@C=12pt{
		&&&0\ar[d]&0\ar[d]&&\\
		&&0\ar[d]\ar[r]&F'_W\ar@{=}[r]\ar[d]&F'_W\ar[r]\ar[d]\ar[r]&0\\
		&0\ar[r]& \TET \ar[r]\ar[d]& F_W \ar[r]\ar[d]& W\otimes \OO_X \ar[r]\ar[d]&0&\\
	&	0\ar[r]&Q''_X\ar[r]\ar[d]&Q'_X\ar[r]\ar[d]&\OO_X (L_X)\ar[r]\ar[d]&0&\\
	&0\ar[r]	&Tor(S_X)\ar[d]\ar@{=}[r]&Tor(S_X)\ar[r]\ar[d]&0&\\
		&&0&0&&		
	}
\end{equation}

The epimorphism $W\otimes \OO_X \longrightarrow \OO_X (L_X)$ in the column on the  right of the above diagram gives an injection
\begin{equation}\label{WinH0LX}
W\hookrightarrow H^0 (\OO_X (L_X)).
\end{equation}

Given a smooth fibre $C_b$ of $f_X :X\longrightarrow \PP^1$ we choose   $w\in W$  
so that the corresponding global section $l_w$ of $\OO_X (L_X)$ has no zeros on $C_b$: this is possible because the linear subsystem $|W|$ of $|L_X|$ is base point free. 

Consider the extension
\begin{equation}\label{ext-w}
\xymatrix@R=12pt@C=12pt{
	0\ar[r]&\TET \ar[r]& F_{[w]} \ar[r]&\OO_X\ar[r] &0
}
\end{equation}
corresponding to the line $[w]$ in $H^1(\TET)$ spanned by $w$. This is related to the extension corresponding to $W\subset H^1(\TET)$ by the diagram
$$
\xymatrix@R=12pt@C=12pt{
	&&0\ar[d]&0\ar[d]&\\
	0\ar[r]&\TET \ar[r]\ar@{=}[d]& F_{[w]} \ar[r]\ar[d]&\OO_X\ar[r]\ar[d] &0\\
		0\ar[r]&\TET \ar[r]& F_{W} \ar[r]\ar[d]&W\otimes \OO_X\ar[r] \ar[d]&0\\
		&&W/[w] \otimes \OO_X \ar@{=}[r]\ar[d]&W/[w] \otimes \OO_X \ar[d]&\\
		&&0&0&
}
$$
The middle column is the extension sequence corresponding to the inclusion
\begin{equation}\label{W-w-ExtFw}
	W/[w] \hookrightarrow H^1 (F_{[w]}) \cong Ext^1(\OO_X, F_{[w]}).
\end{equation}
Furthermore, composing the inclusion $F_{[w]} \hookrightarrow F_W$ with the morphism
$F_W \longrightarrow Q'_X$ in the middle column of \eqref{FW-Q'X-1}
gives the diagram

\begin{equation}\label{Fwslanted}
\xymatrix@R=12pt@C=12pt{
	&&0\ar[d]&&\\ 
	&&F'_W\ar[d]\ar[rd]&&\\
0\ar[r]&F_{[w]}\ar[r]\ar[dr]&F_W \ar[r]\ar[d]&W/[w] \otimes \OO_X \ar[r]&0\\
&&Q'_X&&
}
\end{equation}
where the slanted arrows are the compositions. To understand those arrows
we 
have
$$
\xymatrix@R=12pt@C=12pt{
	&0\ar[d]&0\ar[d]&0\ar[d]&\\
	0\ar[r]&\TET \ar[r]\ar[d]& F_{[w]} \ar[r]\ar[d]&\OO_X\ar[r]\ar[d]^{l_w} &0\\
	0\ar[r]&Q''_X \ar[r]\ar[d]&Q'_X \ar[r]\ar[d]&\OO_X (L_X)\ar[r]\ar[d]&0\\
	0\ar[r]&Tor(S_X)\ar[d]\ar[r]&T_{[w]}\ar[r]\ar[d]&\OO_{D_w} (L_X)\ar[d]\ar[r]&0\\
	&0&0&0&
}
$$
where $D_w =(l_w=0)$ is the zero divisor of the global section $l_w$ of $\OO_X (L_X)$ corresponding to $w$ under the inclusion in \eqref{WinH0LX}. From this it follows that the morphism
$$
F_{[w]} \longrightarrow Q'_X
$$
fails to be an isomorphism precisely on the divisor
$$
\widetilde{D}_w=c_1(Tor(S_X)) +D_w=D^0_X +D_w.
$$	
According to our choice of $w$, the fibre $C_b$ is disjoint from $\widetilde{D}_w$. This implies that the slanted arrows in 
\eqref{Fwslanted}
are isomorphisms when restricted to $C_b$. This means that the extension sequence
$$
\xymatrix@R=12pt@C=12pt{
	0\ar[r]&F_{[w]}\ar[r]&F_W \ar[r]&W/[w] \otimes \OO_X \ar[r]&0,
}
$$
in
 \eqref{Fwslanted} splits on $C_b$. From this it follows that the inclusion
in \eqref{W-w-ExtFw} factors through the kernel of the homomorphism
$$
H^1(F_{[w]})\longrightarrow H^1(F_{[w]}\otimes \OO_{C_b}).
	$$
	coming from the restriction-to-$C_b$ morphism in the exact sequence
	$$
\xymatrix@R=12pt@C=12pt{
	0\ar[r]&F_{[w]}(-C_b)\ar[r]&F_{[w]}	 \ar[r]&F_{[w]}\otimes \OO_{C_b} \ar[r]&0.
}
$$
From this we deduce the inclusion
$$
\xymatrix@R=12pt@C=12pt{
	&0\ar[d]&\\
	&W/[w]\ar[d]\ar[dl]&\\
H^1 (F_{[w]}(-C_b)) \ar[r]&H^1(F_{[w]})&
}
$$
of $W/[w]$
 into $H^1 (F_{[w]}(-C_b))$. By definition $F_{[w]}$ is the middle term of the extension sequence \eqref{ext-w}. This implies
$$
\xymatrix@R=12pt@C=12pt{
&	&0\ar[d]&&\\
&	&W/[w]\ar[d]&&\\
0\ar[r]&H^1(\TET(-C_b))\ar[r]&	H^1 (F_{[w]}(-C_b)) \ar[r]&H^1(\OO_X (-C_b))&
}.
$$
From the vanishing of $H^1 (\OO_X)$ (the assumption of the regularity of $X$) it follows that $H^1(\OO_X (-C_b))=0$. Hence the commutative diagram
$$
\xymatrix@R=12pt@C=12pt{
		0\ar[d]&0\ar[d]\\
	W/[w]\ar@{=}[r]\ar[d]	&W/[w]\ar[d]\\
	H^1(\TET(-C_b))\ar[r]^{\cong}&	H^1 (F_{[w]}(-C_b))
}.
$$
Putting this together with the restriction-to-$C_b$ sequence gives the commutative diagram
$$
\xymatrix@R=12pt@C=12pt{
	&&0\ar[d]&\\
&&	H^1(\TET(-C_b)) \ar[r]^{\cong}\ar[d]& H^1 (F_{[w]}(-C_b))\ar[d]\\
0\ar[r]&[w]\ar[r]\ar[d]&	H^1(\TET) \ar[r]\ar[d]& H^1 (F_{[w]})\ar[d]\\
	&H^0(\OO_{C_b})\ar[r]&	H^1(\TET\otimes \OO_{C_b}) \ar[r]& H^1 (F_{[w]}\otimes \OO_{C_b})
	}
$$
From this it follows that the image of the subspace $W\subset H^1(\TET)$
under the homomorphism
$H^1 (\TET)\longrightarrow H^1 (\TET\otimes \OO_{C_b})$ is
 equal to the image of the line $[w]$ under that homomorphism. The statement \eqref{cl-W}
is now proved.
\end{pf}

\vspace{1cm}
\begin{flushright}
	Universit\'e d'Angers\\
	D\'epartement de Math\'ematiques
	\\
	2, boulevard Lavoisier\\
	49045 ANGERS Cedex 01 \\
	FRANCE\\
	{\em{E-mail addres:}} reider@univ-angers.fr
\end{flushright} 
	
\end{document}